\documentclass[12pt,reqno]{amsart}
\usepackage[toc]{appendix} 
\usepackage{amsmath, amssymb, amsfonts, tikz, graphicx, textcomp, caption, amsthm, listings, verbatim, mathrsfs, url, tipa, datetime,  stmaryrd, xypic, eucal, quiver} 
\usepackage{subfiles}
\usepackage[LGR,T1]{fontenc}

\usepackage[utf8]{inputenc}
\usepackage[colorlinks=true,hyperindex, linkcolor=magenta, pagebackref=false, citecolor=cyan,   pdfpagelabels]{hyperref}
\usepackage[inline,shortlabels]{enumitem}
\usepackage{array}
\usepackage[margin=3cm]{geometry}
\usepackage[capitalize]{cleveref}
\usepackage{mathtools}

\newcolumntype{L}{>{$}l<{$}}
\usetikzlibrary{calc, arrows}
\newcommand{\stackyq}{\text{ }/\kern-0.3cm /\text{ }}

\newcommand{\textgreek}[1]{\begingroup\fontencoding{LGR}\selectfont#1\endgroup}

\newtheorem{thmx}{Theorem}

\newtheorem{corx}[thmx]{Corollary}

\newtheorem{theorem}{Theorem}[subsection]

\makeatletter
\let\c@equation\c@theorem
\makeatother

\DeclareMathOperator{\Aut}{Aut}

\DeclareMathOperator{\End}{End}

\DeclareMathOperator{\Ind}{Ind}

\DeclareMathOperator{\Hom}{Hom}

\DeclareMathOperator{\Sym}{Sym}
\DeclareMathOperator{\Gal}{Gal}
\DeclareMathOperator{\Fun}{Fun}

\DeclareMathOperator{\Ext}{Ext}
\DeclareMathOperator{\colim}{colim}

\DeclareMathOperator{\Sp}{Sp}

\DeclareMathOperator{\Spf}{Spf}
\DeclareMathOperator{\Spa}{Spa}
\DeclareMathOperator{\Spd}{Spd}
\DeclareMathOperator{\Shv}{Shv}
\DeclareMathOperator{\St}{St}

\DeclareMathOperator{\id}{id}

\DeclareMathOperator{\Pic}{Pic}

\DeclareMathOperator{\Spec}{Spec}
\DeclareMathOperator{\Proj}{Proj}

\DeclareMathOperator{\cts}{cts}
\DeclareMathOperator{\Dist}{Dist}

\DeclareMathOperator{\BC}{BC}
\DeclareMathOperator{\LT}{LT}
\DeclareMathOperator{\RZ}{RZ}

\DeclareMathOperator{\cris}{cris}
\DeclareMathOperator{\Bun}{Bun}
\DeclareMathOperator{\Isoc}{Isoc}
\DeclareMathOperator{\mer}{mer}

\DeclareMathOperator{\VB}{VB}

\DeclareMathOperator{\val}{val}
\DeclareMathOperator{\Symm}{Symm}
\DeclareMathOperator{\Kos}{Kos}
\DeclareMathOperator{\Laz}{Laz}
\DeclareMathOperator{\Gr}{Gr}
\DeclareMathOperator{\vStk}{vStk}
\DeclareMathOperator{\un}{un}
\DeclareMathOperator{\Sht}{Sht}
\DeclareMathOperator{\ProFin}{ProFin}
\DeclareMathOperator{\Fin}{Fin}
\DeclareMathOperator{\Cond}{Cond}
\DeclareMathOperator{\Ani}{Ani}
\DeclareMathOperator{\Cat}{Cat}
\DeclareMathOperator{\CAlg}{CAlg}
\DeclareMathOperator{\Mod}{Mod}
\DeclareMathOperator{\LMod}{LMod}
\DeclareMathOperator{\Cts}{C}
\DeclareMathOperator{\Tot}{Tot}
\DeclareMathOperator{\Lie}{Lie}

\newcommand{\EE}{\mathbb{E}}
\newcommand*{\B}[1][]{\mathbb{B}_{#1}}
\newcommand*{\bA}[1][]{\mathbb{A}_{#1}}
\newcommand{\hB}{\widehat{\mathbb{B}}}
\newcommand{\op}{{\mathrm{op}}}

\newcommand*{\LKN}{\mathsf{LKN}}

\newcommand{\abs}[1]{\left\lvert #1 \right\rvert}

\newcommand{\set}[1]{\left\{ #1 \right\}}
\newcommand{\powerseries}[1]{\llbracket #1 \rrbracket}
\newcommand{\laurentseries}[1]{(\!(#1)\!)}
\newcommand{\isom}{\cong}
\newcommand{\from}{\colon}
\newcommand{\injto}{\mathrel{\hookrightarrow}}

\newcommand{\surjto}{\mathrel{\twoheadrightarrow}}

\renewcommand{\xto}[1]{\mathbin{\xrightarrow{#1}}} 
\newcommand{\xfrom}[1]{\mathbin{\xleftarrow{#1}}} 
\newcommand{\isoto}{\xto{\smash{\raisebox{-.4ex}{\ensuremath{\scriptstyle\sim}}}}}

\newcounter{Partcounter}

\usepackage{tikz-cd}\tikzset{node distance=2cm, auto}

\newcommand{\C}{\mathbb{C}}
\newcommand{\CC}{\mathcal{C}}
\newcommand{\DD}{\mathcal{D}}
\newcommand{\Z}{\mathbb{Z}}
\DeclareMathOperator{\GL}{GL}
\newcommand{\N}{\mathbb{N}}
\newcommand{\F}{\mathbb{F}}
\newcommand{\G}{\mathbb{G}}
\newcommand{\Q}{\mathbb{Q}}
\newcommand{\R}{\mathbb{R}}
\newcommand{\A}{\mathbb{A}}

\renewcommand{\SS}{\mathbb{S}}

\newcommand{\Y}{\mathcal{Y}}

\newcommand{\hS}[1]{\mathbb{S}_{K(#1),\infty}} 
\newcommand{\PP}{\mathbb{P}}

\newcommand{\eps}{\varepsilon}

\newcommand{\Einf}{\mathbb{E}_\infty}

\newcommand{\OO}{\mathcal{O}}
\newcommand{\HH}{\mathcal{H}}
\newcommand{\et}{\mathrm{\acute{e} t}}
\newcommand{\proet}{\mathrm{pro\acute{e} t}}

\newcommand{\mc}{\mathcal}

\def\alg{\mathrm{alg}}

\def\SW{\mathbb{SW}}
\def\notimes{\boxtimes}
\def\cl{\mathrm{cl}}
\def\un{\mathrm{un}}
\def\Map{\mathrm{Map}}

\definecolor{codegray}{gray}{0.9}

\def\cond{\mathrm{cond}}
\def\Prof{\mathrm{Prof}}
\def\Comod{\mathrm{Comod}}
\def\Coalg{\mathrm{coAlg}}
\def\hyp{\mathrm{hyp}}
\def\triv{\mathrm{triv}}

\def\Grp{\mathrm{Grp}}
\def\coMon{\mathrm{coMon}}

\def\mix{\mathrm{mix}}
\def\for{\mathrm{for}}
\def\lax{\mathrm{lax}}

\numberwithin{equation}{subsection}

\theoremstyle{definition}
\newtheorem{defn}[equation]{Definition}
\newtheorem{notations}[equation]{Notations}
\newtheorem{construction}[equation]{Construction}
\DeclareFontFamily{U}{min}{}
\DeclareFontShape{U}{min}{m}{n}{<-> dmjhira}{}

\theoremstyle{plain}
\newtheorem{thm}[equation]{Theorem}
\newtheorem{lemma}[equation]{Lemma}
\newtheorem{cor}[equation]{Corollary}
\newtheorem{conj}[equation]{Conjecture}
\newtheorem{proposition}[equation]{Proposition}
\newtheorem*{thm*}{Theorem}

\theoremstyle{remark}
\newtheorem{rmk}[equation]{Remark}
\newtheorem{example}[equation]{Example}
\newtheorem{ntn}[equation]{Notation}

\newtheorem{theorem*}{Theorem}
\newtheorem{conv}[equation]{Convention}

\theoremstyle{remark}

\newcommand{\E}{\mathcal{E}}

\newcommand{\Perf}{\mathrm{Perf}}

\newcommand{\Pro}{\mathrm{Pro}}
\renewcommand{\solid}{{\scalebox{0.5}{$\square$}}}
\newcommand{\LC}{\mathrm{LC}} 
\newcommand{\ad}{\mathrm{ad}} 
\newcommand{\fg}{\Gamma}

\newcommand{\sm}{\mathrm{sm}}

\newcommand{\height}{\mathrm{height}}

\newlist{thmenum}{enumerate}{1}
\setlist[thmenum]{label=(\roman*), ref=\thetheorem.(\roman*)}
\crefalias{thmenumi}{theorem}

\newlist{propenum}{enumerate}{1}
\setlist[propenum]{label=(\roman*), ref=\theproposition.(\roman*)}
\crefalias{propenumi}{proposition}

\newlist{corenum}{enumerate}{1}
\setlist[corenum]{label=(\roman*), ref=\thecor.(\roman*)}
\crefalias{corenumi}{cor}

\newlist{lemenum}{enumerate}{1}
\setlist[lemenum]{label=(\roman*), ref=\thelemma.(\roman*)}
\crefalias{lemenumi}{lemma}

\newlist{defenum}{enumerate}{1}
\setlist[defenum]{label=(\alph*), ref=\thedefn.(\alph*)}
\crefalias{defenumi}{defn}

\AddToHook{env/theorem/begin}{     \crefalias{equation}{theorem}}
\AddToHook{env/thm/begin}{     \crefalias{equation}{thm}}
\AddToHook{env/example/begin}{     \crefalias{equation}{example}}
\AddToHook{env/proposition/begin}{    \crefalias{equation}{proposition}}
\AddToHook{env/lemma/begin}{     \crefalias{equation}{lemma}}
\AddToHook{env/cor/begin}{     \crefalias{equation}{cor}}
\AddToHook{env/conj/begin}{    \crefalias{equation}{conj}}
\AddToHook{env/rmk/begin}{    \crefalias{equation}{rmk}}
\AddToHook{env/defn/begin}{    \crefalias{equation}{defn}}

\title{On the chromatic splitting conjecture in coheight 1}

\author[Barthel]{Tobias Barthel}
\address{Max Planck Institute for Mathematics, Vivatsgasse 7, 53111 Bonn, Germany}
\email{tbarthel@mpim-bonn.mpg.de}

\author[Mann]{Lucas Mann}
\address{Faculty of Mathematics and Computer Science, Einsteinstra{\ss}e 62, 48149 M\"unster, Germany}
\email{mann.lucas@uni-muenster.de }

\author[Ray]{Rin Ray}
\address{Faculty of Mathematics and Computer Science, Einsteinstra{\ss}e 62, 48149 M\"unster, Germany}
\email{cray@uni-muenster.de}

\author[Schlank]{Tomer M.~Schlank}
\address{Department of Mathematics, Harvard University, 1 Oxford Street, Cambridge, MA, USA}
\email{schlank@math.harvard.edu}

\author[Senger]{Andrew Senger}
\address{Department of Mathematics, University of Maryland, College Park, MD, USA}
\email{senger@umd.edu}

\author[Weinstein]{Jared Weinstein}
\address{Boston University Department of Mathematics and Statistics, 665 Commonwealth Avenue, Boston, MA, USA}
\email{jsweinst@bu.edu}

\author[Zhou]{Xinyu Zhou}
\address{Boston University Department of Mathematics and Statistics, 665 Commonwealth Avenue, Boston, MA, USA}
\email{xyz6@bu.edu}

\date{\today}

\usepackage[dvipsnames]{xcolor}
\usepackage[disable]{todonotes}
\usepackage{tcolorbox}
\usepackage{soul}
\newif\ifcomments
\commentstrue

\ifcomments
\DeclareRobustCommand{\xinyu}[1]{{\sethlcolor{magenta!30}\quad \hl{Xinyu: #1}}}

\DeclareRobustCommand{\jared}[1]
{{\sethlcolor{blue!30}\quad \hl{Jared: #1}}}
\DeclareRobustCommand{\tomer}[1]
{{\sethlcolor{violet!30}\quad \hl{Tomer: #1}}}
\DeclareRobustCommand{\tobi}[1]{{\sethlcolor{green!60}\quad \hl{Tobi: #1}}}
\DeclareRobustCommand{\andy}[1]{{\sethlcolor{teal!30}\quad \hl{Andy: #1}}}
\else 
\DeclareRobustCommand{\jared}[1]{}
\DeclareRobustCommand{\tomer}[1]{}
\DeclareRobustCommand{\xinyu}[1]{}
\DeclareRobustCommand{\tobi}[1]{}
\DeclareRobustCommand{\andy}[1]{}
\fi 

\begin{document}

\begin{abstract}
    This article confirms predictions made by Hopkins' {\em chromatic splitting conjecture} in stable homotopy theory.  A consequence of our main theorem is that a certain map of $K(h-1)$-local $\mathbb{E}_\infty$-rings $\mathbb{S}_{K(h-1)}^{B\mathbb{Z}}\to L_{K(h-1)}\mathbb{S}_{K(h)}$ admits a retract, so long as $p-1$ does not divide $h-1$.  If $h=2$ and $p$ is odd, we show this map is an isomorphism. 

    We use methods from the theory of $p$-divisible groups and $p$-adic geometry.  Along the way we prove an entirely new kind of ``two-tower isomorphism'', in analogy with Faltings' isomorphism between the Lubin--Tate and Drinfeld towers, which parametrizes $p$-divisible groups over characteristic $p$ perfectoid spaces with fixed generic fiber in terms of negative Banach--Colmez spaces.  Using the two-tower isomorphism, we gain control over the continuous cohomology of the Morava stabilizer group $\mathbb{G}_h$ with coefficients in $\pi_*L_{K(h-1)}E_h$. 

    To study the continuous action of $\mathbb{G}_h$ on $L_{K(h-1)} E_h$, it was necessary to develop the theory of solid Lubin--Tate spectra.  Adapting Lurie's constructions, we attach a solid Lubin--Tate spectrum to $p$-divisible groups over a broad class of solid rings in characteristic $p$.
\end{abstract}

\maketitle

\begin{center}  \includegraphics[scale=0.075]{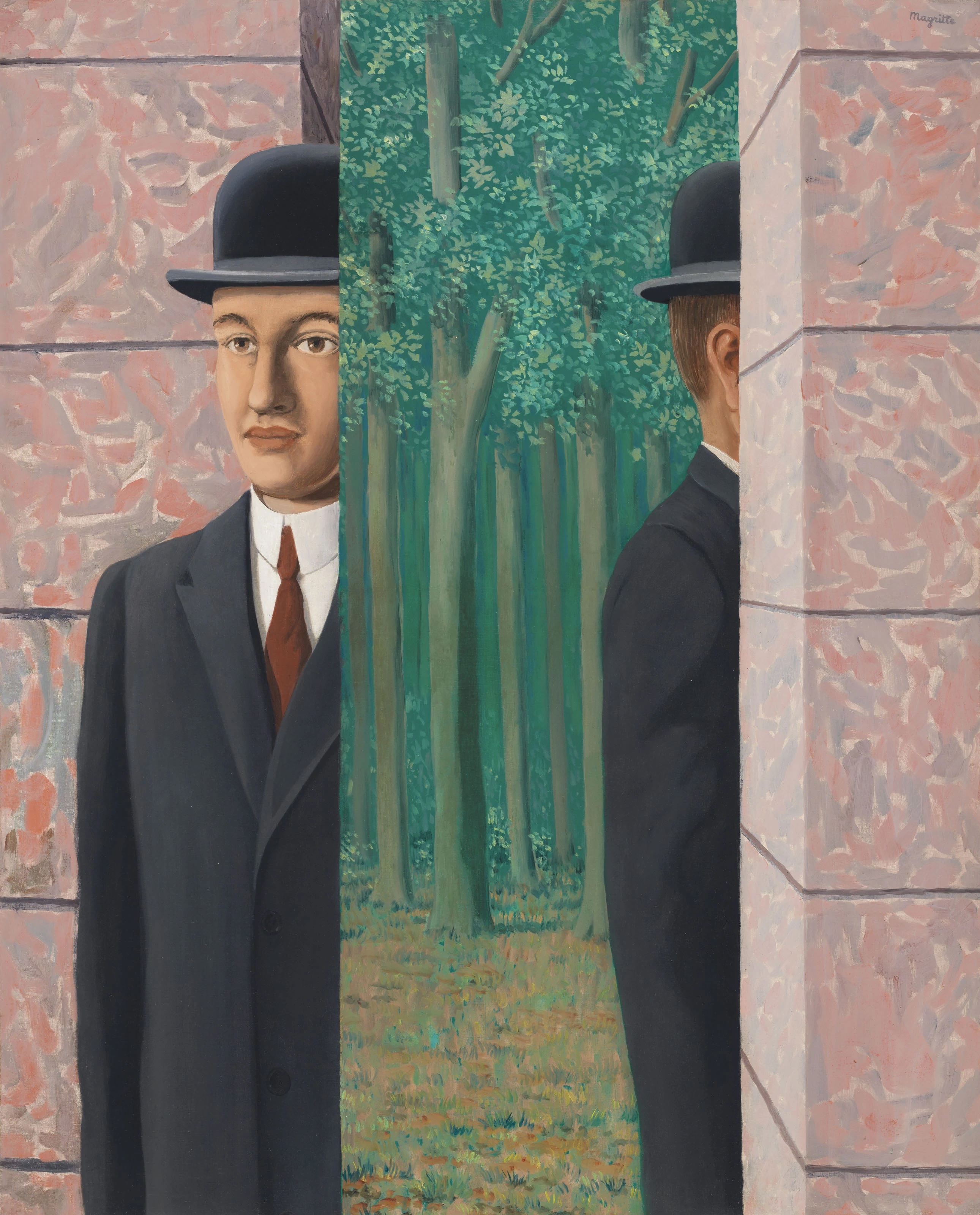} \\
Ren\'e Magritte, The Common Place, 1964
\copyright 2026 C. Herscovici / Artists Rights Society (ARS), New York

\end{center}

\newpage
\setcounter{tocdepth}{1}

\newpage
\section{Introduction and main results}

Chromatic homotopy theory studies the stable homotopy groups of the sphere spectrum $\SS$ through a divide-and-conquer approach: First decompose $\SS$ into ``irreducible'' pieces $\SS_{K(h)}$, study these individually, and then reassemble them. Here, $K(h)$ is the height $h$ Morava $K$-theory spectrum at the implicit prime $p$, forming the residue fields of $\SS$, and $\SS_{K(h)}$ denotes the corresponding completion. Writing $L_h$ for the localization at $K(0), K(1), \ldots, K(h)$ jointly, the inductive assembly process is formally captured by a homotopy pullback square, displayed on the right, which forms the heart of the subject:
    \[
        \xymatrix{\Z_{(p)} \ar[r] \ar[d] & \Z_{p} \ar[d] && L_h\SS \ar[r] \ar[d] & \SS_{K(h)} \ar[d] \\
        \Q \ar[r] & \Q_p  && L_{h-1}\SS \ar[r]^-{\iota_h} & L_{h-1}\SS_{K(h)}.}
    \]
This height $h$ chromatic fracture square can be thought of as a chromatic analogue of the familiar Hasse square reconstructing the ring of $p$-local integers from its places, as in the square on the left. Based on low height evidence, Hopkins' {\em chromatic splitting conjecture} predicts an explicit form of chromatic assembly which would have important theoretical consequences; see \cite{MR1320994} or \cite[\S 2.4]{BarthelBeaudry} for a detailed exposition. At first approximation, it comes in two different flavors:
    \begin{itemize}
        \item \emph{Qualitatively}, the \emph{weak chromatic splitting conjecture} predicts that $\iota_h$ splits.
        \item \emph{Quantitatively}, the \emph{strong chromatic splitting conjecture} then describes the fiber of $\iota_h$ in terms of localizations of familiar spectra. 
    \end{itemize}
Both conjectures have been resolved for all $h \leq 2$ through substantial calculations, while for $h \geq 3$ only very sporadic information has been known as computational techniques quickly become intractable, necessitating a different approach. 

In recent work \cite{BSSW}, it is suggested that the problem becomes more amenable via ideas from modern arithmetic geometry. Indeed, utilizing techniques from $p$-adic geometry, specifically the two-tower isomorphism of Faltings and Fargues in the form established by Scholze--Weinstein \cite{ScholzeWeinstein} as well as integral $p$-adic Hodge theory \cite{BhattMorrowScholze}, the main result of \cite{BSSW} solves the chromatic splitting conjecture rationally for all heights $h$ and primes. In other words, this determines precisely the coheight $h$ monochromatic part of $\SS_{K(h)}$, i.e., the iterated localization $L_{K(0)}\SS_{K(h)}$. Naturally, we may probe the conjecture through other iterated localizations $L_{K(h-i)}\SS_{K(h)}$ of $\SS_{K(h)}$. The subject of the present paper is, as the title suggests, the coheight 1 part $L_{K(h-1)}\SS_{K(h)}$.

The coheight 1 weak chromatic splitting conjecture predicts that the map
    \[
        L_{K(h-1)}\iota_h\colon L_{K(h-1)}\SS \to L_{K(h-1)}\SS_{K(h)}
    \]
splits, while its strong form 
identifies the cofiber as $\Sigma^{-1}L_{K(h-1)}\SS$, with suitable modifications for small primes. As a consequence of the results of this paper, we prove:
    \begin{thm*}
        The coheight 1 weak chromatic splitting conjecture is true for all $p-1\nmid h-1$.
    \end{thm*} 
In fact, as we explain in the next part, our theorems are significantly sharper and also partially address arbitrary primes. The proofs rely on novel results in $p$-adic geometry which are of independent interest.

\subsection{Main results}

In order to state our main chromatic results, it will be convenient to work with the {\em Iwasawa extension} $\hS{h}$ of the $K(h)$-local sphere, also known as the half sphere: this is a $K(h)$-local Galois extension $\hS{h}/\SS_{K(h)}$ with Galois group $\Z_p$, see \S\ref{ssec:csc} for a detailed discussion. Now assume $h \geq 2$ and consider the localization $L_{K(h-1)}\hS{h}$, so there exists a unique map of $\EE_\infty$-rings $\SS_{K(h-1)}\to L_{K(h-1)}\hS{h}$. Based on our findings, we entertain the following conjecture which is consistent with the strong form of Hopkins' chromatic splitting conjecture in coheight $1$ and in fact implies it.

\begin{conj}[The half sphere chromatic splitting conjecture in coheight 1]
\label{ConjStrongSplitting} 
Assume that $p > h \geq 2$.   The unit map $\SS_{K(h-1)} \to L_{K(h-1)}\hS{h}$ of $\EE_\infty$-ring spectra is an isomorphism. 
\end{conj}

We approach this conjecture by constructing a proxy in the form of a $\Z_p$-equivariant map of $\EE_{\infty}$-ring spectra $L_{K(h-1)}\hS{h} \to \bA[h]$, where $\bA[h]$ will be defined rigorously in \eqref{eq:introAph}. Here is the main result of the article:

\begin{thmx}\label{thmx:main}
Let $h \geq 2$ and let $p$ be an odd prime.
    \begin{enumerate}
        \item If $p-1\nmid h-1$, then the unit map $L_{K(h-1)}\SS \to \bA[h]$ is an equivalence.
        \item If $p-1\mid h-1$, then the cofiber of $L_{K(h-1)}\SS \to \bA[h]$ is in $\Pic(\Sp_{K(h-1)})$.
    \end{enumerate}
\end{thmx}

In particular, in case $p-1\nmid h-1$, passing to $\Z_p$-fixed points induces a retraction of $\SS_{K(h-1)}^{B\Z_p}\to L_{K(h-1)}\SS_{K(h)} \simeq L_{K(h-1)}\hS{h}^{h\Z_p}$. We remark that for any $p$-complete $\Einf$-ring spectrum $X$, the underlying spectrum of the $\EE_{\infty}$-ring spectrum $X^{B\Z_p}$ is $X\oplus \Sigma^{-1}X$. For $X = \SS_{K(h-1)}$, the two components of the map $\SS_{K(h-1)}^{B\Z_p}\to L_{K(h-1)}\SS_{K(h)}$ may be identified with the canonical unit map and the map $\zeta_h$ constructed by Devinatz--Hopkins \cite{DevinatzHopkins2004}, respectively. Thus:

\begin{corx}\label{corx:main}
     The map $\SS_{K(h-1)}^{B\Z_p} \to L_{K(h-1)}\SS_{K(h)}$ admits a retract when $p-1\nmid h-1$.
\end{corx}

In particular, this establishes the coheight 1 weak chromatic splitting conjecture at all heights and sufficiently generic primes as stated above. It is the first result in this direction beyond height 2. In the case $p-1\mid h-1$, the $K(h-1)$-locally invertible element we exhibit provides a potential obstruction to the chromatic splitting conjecture in its original form, compatible with the known divergence for $h=p=2$, see \cite{beaudry2022chromatic}. For $h = 2$ and odd primes, our results are sharper and recover the known cases of the strong chromatic splitting conjecture:
 
\begin{thmx}\label{thmx:height2} 
    Assume $p\neq 2$. The unit map $\SS_{K(1)} \to L_{K(1)}\SS_{K(2),\infty}$ is an equivalence. In particular, it induces a canonical isomorphism of $\EE_{\infty}$-ring spectra $\SS_{K(1)}^{B\Z_p} \simeq L_{K(1)}\SS_{K(2)}$.
\end{thmx}

\Cref{thmx:height2} was first proved by Shimomura--Yabe (\cite{ShimomuraYabe1995}, \cite{Behrens2012}) for $p>3$ and by Goerss, Henn, Karamanov, Mahowald, Rezk, and Shimomura \cite{Shimomura2000,GHMR2005,HKM2013,GHM2014} for $p=3$, using intense explicit calculations of both sides.  In contrast, our proof of \cref{thmx:height2} is conceptual and requires only comparatively little computational input.  In fact, the proof of \cref{thmx:main} is independent of any intricate group cohomology calculations; instead, both proofs make essential use of the geometry of moduli spaces of $p$-divisible groups. 

\begin{rmk} 
    When $p=2$, it is no longer true that  $L_{K(1)}\hS{2}$ is isomorphic to $\SS_{K(1)}$: Instead, there is an additional $L_{K(1)}\SS/2$-summand in degree $(-2)$, as found and described by Beaudry, Goerss and Henn in \cite[Theorem 1.1.9]{beaudry2022chromatic}.  While we have not been able to recover this result using our methods, we obtain a partial explanation in terms of the so-called Heyer--Jacquet module. 

    In fact, in the general situation where $p-1$ divides $h-1$, we can point to the Heyer--Jacquet module as an explanation for the failure of \cref{ConjStrongSplitting}.  We consider this observation to be one of the major contributions of this paper to the analysis of iterated chromatic localizations.  
\end{rmk}

The rest of the introduction lays out the main ideas in the proof of \cref{thmx:main} and \cref{thmx:height2}.

\subsection{Overview of the proof and main algebraic result}

The conventional approach to access the chromatic information contained in $\SS_{K(h)}$ is by descent from Morava $E$-theory $E_h$, where $E_h$ is an even periodic $K(h)$-local $\EE_{\infty}$-ring spectrum equipped with a continuous action by the Morava stabilizer group $\G_h$. To recall its basic structure, let $k\subset \overline{\F}_p$ be a subfield containing $\F_{p^h}$, and let $\Gamma_h$ be the $p$-divisible group of the Honda formal group of height $h$ over $k$, for which $\pi_0E_h$ is the deformation ring of $\Gamma_h$.  Let $\Gamma^{\un}_h$ be the universal deformation of $\Gamma_h$.  Then $\pi_{-2}E_h=\Lie\Gamma_h^{\un}$ is the Lie algebra; it is an invertible $\pi_0E_h$-module. Explicitly:
\begin{eqnarray*}
    \pi_0E_h &=& W(k)\powerseries{u_1,\dots,u_{h-1}} \\ 
    \pi_*E_h &=& \pi_0E_h[u^{\pm 1}],
\end{eqnarray*}
where $u$ is a generator in degree $-2$.

Devinatz and Hopkins \cite{DevinatzHopkins2004} show that the canonical unit map $\SS_{K(h)} \to E_h$ exhibits $E_h$ as a $K(h)$-local $\G_h$-Galois extension of $\SS_{K(h)}$. In particular, this begets a descent spectral sequence with excellent convergence properties, of signature 
    \begin{equation}
        \label{eq:DH}
        H_{\cts}^*(\G_h,\pi_*E_h) \Longrightarrow \pi_*\SS_{K(h)}.
    \end{equation}
This permits one to study $\pi_*\SS_{K(h)}$ via the representation $\G_h \curvearrowright \pi_*E_h$ and its continuous cohomology. 
Explicit computations of $H_{\cts}^*(\G_h,\pi_*E_h)$, however, quickly increase in complexity and are essentially out of reach for $h >2$. 

Nevertheless, \cite{BSSW} calculates the rational cohomology $\Q\otimes H_{\cts}^*(\G_h,\pi_*E_h)$ for all primes  and all $h\geq 1$.  The method of \cite{BSSW} uses 
Faltings' isomorphism \cite{Faltings2002} between the Lubin--Tate tower and the Drinfeld tower, appearing in the form of an isomorphism between v-stacks
\begin{equation}
    \label{eq:introFaltings}
[\LT_{h,0}^{\Diamond}/\G_h] 
\isom [\mathcal{H}^{h-1,\Diamond}/\GL_h(\Z_p)]. 
\end{equation}
Here $\LT_h=\Spf \pi_0E_h$ and $\LT_{h,0}^{\Diamond}$ is its diamond generic fiber, and $\mathcal{H}^{h-1}$ is Drinfeld's symmetric space of dimension $h-1$.  The next step is to take $\OO^+$-cohomology on either side of \eqref{eq:introFaltings} and then rationalize.  Using some heavy machinery from integral $p$-adic Hodge theory, one arrives at
\[ 
\Q\otimes H^*_{\cts}(\G_h,\pi_0E_h)\isom H^*_{\cts}(\GL_h(\Z_p),\Q_p). 
\] 
The right-hand side is well-known by a result of Lazard.  
From there it is straightforward to give a full computation of $\Q \otimes \pi_*\SS_{K(h)}=\pi_*L_{K(0)}\SS_{K(h)}$.

Faltings' isomorphism \eqref{eq:introFaltings} does not grant access to the intermediate localizations $L_{K(t)}\SS_{K(h)}$ for $0<t<h$, essentially for the reason that it is purely a characteristic 0 phenomenon.  To access transchromatic information contained in iterated localizations such as $L_{K(h-1)}\SS_{K(h)}$, we follow Torii \cite{ToriiComparison}, \cite{MR2794609},\cite{Torii2010}. Let us assume that $k$ contains both $\F_{p^{h-1}}$ and $\F_{p^h}$. Then $E_{h-1}$ and $E_h$ are commutative algebras over $\SS\mathbb{W}(k)$, the $\EE_\infty$-ring of spherical Witt vectors relative to the perfect field $k$.  Define
    \begin{equation}
    \label{eq:B}
        \B \coloneqq E_{h-1}\hat{\otimes}_{\SS\mathbb{W}} L_{K(h-1)}E_h, 
        \end{equation}  
where the tensor product is computed in the category of commutative $K(h-1)$-local $\SS\mathbb{W}$-algebras. By construction, $\B$ admits an action of the group
\[
    \G_{h,h-1}
    \coloneqq
    \G_{h-1}\times_{\Gal(k/\F_p)}\G_h,
\]
induced by the factorwise actions of these groups on their respective
Morava $E$-theories. 

Torii \cite{ToriiComparison} describes the $\pi_*E_h$-algebras $\pi_*L_{K(h-1)}E_h$ and $\pi_*\B$ as follows.  Let $I_{h-1}\subset \pi_0E_h$ be the ideal $(p,u_1,\dots,u_{h-2})$.  Then 
    \[ 
        \pi_*L_{K(h-1)}E_h \isom (\pi_*E_h)[u_{h-1}^{-1}]^\wedge_{I_{h-1}}, 
    \]
meaning the $I_{h-1}$-adic completion of the ring $(\pi_*E_h)[u_{h-1}^{-1}]$.  In particular $\pi_0L_{K(h-1)}E_h$ is a local   ring with residue field
    \begin{equation} 
    \label{eq:K}
        K\coloneqq (\pi_0E_h/I_{h-1})[u_{h-1}^{-1}]\isom k\laurentseries{u_{h-1}}.
    \end{equation}
The base change $\Gamma_{h,K}^{\un}=\Gamma_h^{\un}\otimes_{\pi_0E_h} K$ has identity component $(\Gamma_{h,K}^{\un})^{\circ}$ of height $h-1$. 

\begin{defn}[The Gross--Torii extension $L/K$, {\cite[\S 2.3]{MR2004428}}]
\label{def:torii_ext} Let $L/K$ be the separable extension which classifies isomorphisms $\Gamma_{h-1,L}\isom (\Gamma^{\un}_L)^\circ$.      
\end{defn}

By a result of Gross \cite[Theorem 3.5]{GrossRamification}, $L/K$ is a totally ramified Galois extension with Galois group $\overline{\G}_{h-1}\coloneqq \ker(\G_{h-1}\to\Gal(k/\F_p))$. The continuous action of $\G_h$ on $K$ extends to a continuous action of $\G_{h,h-1}$ on $L$.

We can now state Torii's description of $\pi_*\B$:  it is the unique lift of the Galois extension $L/K$ to a $\pi_*L_{K(h-1)}E_h$-algebra which is complete with respect to the maximal ideal of $\pi_0L_{K(h-1)}E_h$.
Thus $(\pi_*\B)/I_{h-1}\isom L[u^{\pm 1}]$.

The field $L$, being algebraic over $K$, carries a unique valuation extending the $u_{h-1}$-adic valuation on $K$.  Since $L/K$ is infinite, $L$ is not complete with respect to this valuation.  We define $\widehat{L}$ to be the completion of $L$.  Importantly, $\widehat{L}$ is a perfect field;  it is isomorphic to the field $k\laurentseries{t^{1/p^\infty}}$ of perfect Laurent series over $k$ (\cref{PropL}).  Thus $\widehat{L}$ is an example of a perfectoid field of characteristic $p$. The action of $\G_{h,h-1}$ on $L$ extends to a continuous action on $\widehat{L}$. 

The key point now is that the continuous cohomology of $\widehat{L}$ is accessible via $p$-adic geometry.  (Indeed, this entire project began by imagining what a Faltings-style isomorphism involving $\widehat{L}$ would have to look like.)  Before explaining our main algebraic result, let us fix some notation for certain subgroups of the Morava stabilizer group $\G_h$ which appear throughout the article.  We make essential use of the determinant homomorphism $\det\from\G_h\to\G_1$, which commutes with the maps toward $\Gal(k/\F_p)$.  We have $\G_1=\Z_p^\times \times\Gal(k/\F_p)$, and 
$\Z_p^\times \isom \mu \times \Z_p$ for $\mu$ the torsion subgroup of $\Z_p^\times$; thus $\mu$ is cyclic of order $p-1$ or $2$ as $p$ is odd or even.

We may then define closed normal subgroups of $\G_h$ as in the following diagram
\begin{equation}\label{eq:allgroups} 
\begin{split}
\xymatrix{
S\overline{\G}_h \ar[r] \ar[d] & 
\overline{\G}_h^1 \ar[r] \ar[d] & 
\overline{\G}_h \ar[d] \ar[r] &
\G_h\ar[d]^{\det} \\ 
1 \ar[r] & \mu \ar[r] & \Z_p^\times \ar[r] & \G_1 
}
\end{split}
\end{equation} 
in which all squares are cartesian.  Each of the upper groups is $\Gal(k/\F_p)$-stable, and so we may define 
    \begin{equation}\label{eq:allgroups2}
        S\G_h=S\overline{\G}_h\rtimes \Gal(k/\F_p), \quad \G_h^1=\overline{\G}_h^1\rtimes\Gal(k/\F_p).
    \end{equation}
Importantly, the Iwasawa extension $\SS_{K(h),\infty}$ is the subextension of the Galois extension $\SS_{K(h)}\to E_h$ corresponding to the closed normal subgroup $\G_h^1\subset \G_h$, with quotient group $\Z_p$.  

As the main algebraic input en route to \cref{thmx:main}, we establish:

\begin{thmx}\label{thmx:cohom} 
    Assume $h\geq 2$. As representations of $\G_{h-1} \times \Z_p^{\times}$ we have:
        \[ 
            H^i_{\cts}(S\overline{\G}_h,\widehat{L}) \isom 
                \begin{cases}
                    k, & i=0,\\
                    k(h-1), & i=h,\\
                    0. & \text{else},
                \end{cases}
        \] 
        The action of $\Z_p^\times$ on $k(h-1)$ is through the character $a\mapsto a^{h-1}$ mod $p$.   
\end{thmx}

The module $k(h-1)$ appearing in cohomological degree $h$ is what we call the Heyer--Jacquet module, for it arises in our setting by taking the ``higher Jacquet module'' of a Steinberg module, in the sense of Heyer's work \cite{Heyer2023Left}; see \cref{sec:Heyer-Jacquet}. 

From here it is straightforward to access $H^i_{\cts}(\overline{\G}_h^1,\widehat{L})$ by taking derived invariants under the action of $\mu=\overline{\G}_h^1/S\overline{\G}_h$.  The action of $\mu$ on $k(h-1)$ is nontrivial if and only if $p-1\nmid h-1$.  If we assume $p-1\nmid h-1$, the $\mu$-cohomology of $k(h-1)$ vanishes, and we find:
\begin{equation}
    \label{eq:coho_intro}
    H^*_{\cts}(\overline{\G}_h^1,\widehat{L})\isom k. 
\end{equation}

\subsection{Solid spectra and solid Lubin--Tate theories}

We now discuss how to derive the main homotopical result \cref{thmx:main} from the main algebraic result \cref{thmx:cohom}, by combining the following two steps:

{{\em Step 1:  Recover $L_{K(h-1)}\SS_{K(h)}$ as the homotopy fixed points of $\B$.}} 
Consider the diagram of $\EE_\infty$-ring spectra
\begin{equation}
\label{eq:double_descent}
        \xymatrix{L_{K(h-1)}\SS_{K(h)} \ar[r]^-{\G_h} & L_{K(h-1)}E_h \ar[r]^-{\overline{\G}_{h-1}} & \B}
\end{equation}
The first (respectively, second) arrow is derived from the Devinatz--Hopkins extension at height $h$ (respectively, $h-1$) by localizing at height $h-1$ (respectively, base changing to $E_{h-1}$ in the $K(h-1)$-local category).  Since the Devinatz--Hopkins extension is Galois, it is reasonable to think that each arrow is a Galois extension for the indicated group, which should imply that the source of each arrow appears as the homotopy fixed points of the target.

To make this precise, we work with the category of condensed spectra \cite{Clausen2025arxiv}.  There is a notion of continuous homotopy fixed points for a continuous action of a profinite group on a condensed spectrum.  In fact all of the condensed spectra we encounter will be solid spectra, meaning that their homotopy groups are solid abelian groups.  

The precise statement of Step 1 is therefore:  construct a diagram of solid $\EE_\infty$-ring spectra, whose value on a point is \eqref{eq:double_descent}, and for which the source of each arrow is the continuous homotopy fixed points of the target.  
In particular, the Iwasawa extension $L_{K(h-1)}\SS_{K(h),\infty}$ should appear by taking homotopy fixed points of (the solid enhancement of) $\B$ under the subgroup $\G_h^1\times_{\Gal(k/\F_p)} \G_{h-1}$ of $\G_{h,h-1}$.

It is insufficient to give the objects in \eqref{eq:double_descent} the topology inherent to their structure as $K(h-1)$-local spectra.  Indeed, while the rings $\pi_*L_{K(h-1)}E_h$ and $\pi_*\B$ are complete for the $I_{h-1}$-adic topology, their topologies cannot simply be the $I_{h-1}$-adic ones:  the quotients of $\pi_0$ by $I_{h-1}$ are the fields $K$ and $L$ respectively, and these carry a nontrivial $u_{h-1}$-adic topology.  Thus it really is necessary to work with condensed spectra.  

{{\em Step 2: Construct a lift $\B\to\widehat{\B}$ of $L\to\widehat{L}$.}}  More precisely, the idea is to construct a $\G_{h,h-1}$-equivariant map of solid spectra, whose source is the condensed enhancement of $\B$ from the previous paragraph, and which recovers the map of solid rings $L\to\widehat{L}$ upon passage to $\pi_0$ modulo $I_{h-1}$.  Such a lift ought to exist by functoriality of Morava E-theory;  in particular the uncondensed spectrum $\widehat{\B}$ ought to be the height $h-1$ theory associated to the connected part of $\Gamma_{h,\widehat{L}}^{\un}$, which is identified with $\Gamma_{h-1,\widehat{L}}$.  It will be vital to define the solid structure on $\widehat{\B}$ compatibly with the metric topology on $\widehat{L}$, so as to be able to apply the continuous cohomology result in \cref{thmx:cohom}.

Granting the above two constructions, we apply continuous homotopy fixed points by the group 
$\G_h^1\times_{\Gal(k/\F_p)} \G_{h-1}$ to the map $\B\to\widehat{\B}$ (or rather to its condensed enhancement).  The result will be the map $L_{K(h-1)}\SS_{K(h),\infty}\to \mathbb{A}_h$ appearing in \cref{thmx:main}.  The description of the continuous cohomology of $\overline{\G}_h^1$ on $\widehat{L}$ in 
\cref{thmx:cohom} gets used to control the behavior of $\mathbb{A}_h$ as seen in \cref{thmx:main}.

To execute both steps, we generalize Lurie's construction \cite{ECII} of Lubin--Tate theories via derived algebraic geometry.  In that construction, the input is a pair $(R,\Gamma)$ consisting of an $\F_p$-algebra $R$ and a $p$-divisible group $\Gamma$ over $R$ satisfying certain conditions rendering its deformation theory well-behaved (see \cite[Theorem 0.0.8]{ECII}), and the output is an $\EE_\infty$-ring spectrum $E(R,\Gamma)$ whose $\pi_0$ is the universal deformation ring of $\Gamma$.  In our desired generalization, $R$ is a solid ring (such as $K$, $L$, or $\widehat{L}$), and the output is a solid $\EE_\infty$-ring spectrum whose value on a profinite set $S$ is $E(R(S),\Gamma)$.  

We are not able to apply \cite[Theorem 0.0.8]{ECII} ``out of the box'' due to a technical issue:  Lurie's original construction requires the input ring to be Noetherian, whereas $R(S)$ is almost never Noetherian.  To deal with this issue, we propose an axiomatization of the input data required in Lurie's construction which is adapted to handle all cases of interest to us.  In particular we introduce (\cref{def:solidLKNpair}) a category of solid Luriean Kalypto--Noetherian ($\LKN$) pairs $(R,\Gamma)$ consisting of a solid $\F_p$-algebra $R$ and a $p$-divisible group $\Gamma$ over $R(\ast)$ satisfying a number of conditions which guarantee that each $E(R(S),\Gamma)$ exists. The resulting category $\LKN^{\solid}$  contains the fields $K$, $L$, and $\widehat{L}$, and is closed under suitable operations, such as localizations and a certain variant on the tensor product. See \cref{thm:solidLT} for a precise statement of the following theorem.

\begin{thmx}\label{thmx:solidLT}
    The Morava $E$-theory functor extends to a natural functor enriched over condensed anima
    \[
        E^{\solid}(-) \colon  \LKN^{\solid} \to \CAlg^{\solid}
    \]
    taking values in solid $\EE_{\infty}$-ring spectra, compatible with certain tensor products and localizations. Moreover, the condensed functor $E^\solid(-)$ is computed by applying the uncondensed functors of the same name levelwise. 
\end{thmx}

As a basic example, $E_h^\solid\coloneqq E^\solid(k,\Gamma_h)$ is Morava E-theory equipped with its natural topology as a $K(h)$-local spectrum (\cref{exm:usualLT}).  The next example is $E^\solid(K,\Gamma_{h,K}^{\un})$, which turns out to be isomorphic to  the localized solid E-theory $L_{K(h-1)}^{\solid} E_h^\solid$ (\cref{exm:kh-1}).  Finally, the solid spectrum $\B^\solid$ is defined as the $K(h-1)$-localization of $(L_{K(h-1)} ^{\solid} E_h ^{\solid} \otimes_{\SW(k)} ^{\solid} E_{h-1} ^{\solid})$, and this turns out to be isomorphic to $E^{\solid} (L, \Gamma^{\un} _{h,L})$ (\cref{exm:B}).

We now explain how to execute the two steps described above.  For Step 1:  we find a diagram of $K(h-1)$-local solid $\EE_{\infty}$-ring spectra
\begin{equation}
\xymatrix{L^\solid_{K(h-1)}\SS_{K(h)} \ar[r]^-{\G_h} & E^\solid(K,\Gamma_{h,K}^{\un}) \ar[r]^-{\overline{\G}_{h-1}} & E^\solid(L,\Gamma_{h,L}^{\un})}\isom \B^{\solid}
\end{equation}
where $K$ and $L$ are interpreted as solid rings, whose value on $(\ast)$ recovers \eqref{eq:double_descent}.  
Regarding the first arrow: the compatibility of $E^\solid$ with $K(h-1)$-localization shows directly that the source is the homotopy fixed points of the target.  For the second arrow, this is also true, but the argument is more subtle, relying on the equivalence between smooth and condensed fixed points for $p$-adic Lie groups acting on solid, characteristic $p$ coefficients (\cref{rslt:smooth-F-p-reps-embed-fully-faithfully}).  

For Step 2:  define $\widehat{\B}^{\solid}=E^{\solid}(\widehat{L},\Gamma_{h-1,\widehat{L}})$, using the solid structure on $\widehat{L}$ coming from its metric topology.  The map $\B^{\solid}\to \widehat{\B}^{\solid}$ is constructed in \cref{ExBtoBhat}.  

The computation of the $S\overline{\G}_h$-cohomology of $\widehat{L}$ in \cref{thmx:cohom} lifts to a description of the $\G_{h-1}\times\Z_p^{\times}$-equivariant solid $E_{h-1}$-algebra $(\widehat{\B}^{\solid})^{hS\overline{\G}_h}$:  it sits in a $\G_{h-1}\times\Z_p^\times$-equivariant fiber sequence 
\[ E_{h-1}^{\solid}\to (\widehat{\B}^{\solid})^{hS\overline{\G}_h} \to \mathcal{H}_h, \] 
where $\mathcal{H}_h$ is an invertible $E_{h-1}^\solid$-module which lifts the Heyer--Jacquet module.  After taking derived $\mu$-invariants, this module survives exactly when $p-1\vert h-1$.  Therefore:
\begin{itemize}
    \item If $p-1\nmid h-1$, then $E_{h-1}^{\solid}\to (\widehat{\B}^{\solid})^{h\overline{\G}_h^1}$ is an equivalence;  after taking $\G_{h-1}$-fixed points and evaluating on $(*)$ we find an isomorphism 
    \[ 
        \SS_{K(h-1)}\isom (\widehat{\B}^{\solid})^{h (\G_h^1\times_{\Gal(k/\F_p)} \G_{h-1})}(*). 
    \] 
    \item If $p-1\mid h-1$ and $p \neq 2$, then the cofiber of $E_{h-1}^{\solid}\to (\widehat{\B}^{\solid})^{h \overline{\G}_h^1}$ is an invertible $E_{h-1}^{\solid}$-module, and so there is a fiber sequence
    \[ 
        \SS_{K(h-1)}\to  (\widehat{\B}^{\solid})^{h (\G_h^1\times_{\Gal(k/\F_p)} \G_{h-1})}(*) \to \mathcal{C}, 
    \]
    where $\mathcal{C}$ is an invertible $\SS_{K(h-1)}$-module.
\end{itemize}

To conclude, we construct the map of $\EE_{\infty}$-ring spectra appearing in \cref{thmx:main} as the homotopy invariants of the morphism $\B \to \widehat{\B}$:
    \begin{equation}\label{eq:introAph}
        (L_{K(h-1)}\hS{h} \to \bA[h]) \coloneqq (\B^\solid  \to \widehat{\B}^\solid)^{h (\G_h^1\times_{\Gal(k/\F_p)} \G_{h-1})}(*)
    \end{equation}
equipped with its residual $\Z_p$-action. \Cref{thmx:main} is then an immediate consequence of the points above, thereby finishing the homotopical part of our exposition.

\subsection{The moduli interpretation of $\widehat{L}$, and the period domain $\Omega_{h,h-1}$}\label{ssec:intro_moduli}

We now discuss the $p$-adic geometry result powering \cref{thmx:cohom}. 
The idea is to interpret the continuous group cohomology $H^*_{\cts}(S\overline{\G}_h,\widehat{L})$ as the cohomology of the structure sheaf on an appropriate quotient by $S\overline{\G}_h$ of a geometric object associated to $\widehat{L}$.

There is a ready-made category of v-stacks which allows for such quotients. The field $\widehat{L}$ is a perfectoid field.  For a perfectoid field $M$, let $\Spd M$ be the corresponding diamond, meaning the sheaf of sets represented by $M$ on the v-site of perfectoid spaces in characteristic $p$.  Define the v-stack
 \[ \Omega_{h,h-1}\coloneqq (\Spd \widehat{L})/S\overline{\G}_h, \] 
which carries a residual action of $\G_{h,h-1}/S\overline{\G}_h\isom \G_{h-1}\times\Z_p^\times$.  (In fact this quotient is itself a diamond.) 

The theory of the Fargues--Fontaine curve enters the picture at this point.  The proof of the following theorem will be discussed in \S\ref{ssec:intro_twotowers}.  

\begin{proposition}[\Cref{PropL}]
 \label{ThmMeaningOfL} 
 As a diamond over $\Spd k$ with $\G_{h,h-1}$-action, $\Spd \widehat{L}$ assigns to a perfectoid space $S$ over $k$ the set of exact sequences
\[ 0\to\OO\to \OO(1/h) \to \OO(1/(h-1))\to 0  \] 
of vector bundles on the Fargues--Fontaine curve $X_S$ which preserve a fixed isomorphism $\det\OO(1/h)\isom \det\OO(1/(h-1))$, where the group $\G_h$ acts via automorphisms of $\OO(1/h)$ and the group $\G_{h-1}$ via automorphisms of $\OO(1/(h-1))$.  Accordingly, the $S$-points of $\Omega_{h,h-1}=(\Spd\widehat{L})/S\G_h$ classify extensions of vector bundles on $X_S$ of the form
\[ 0 \to \OO\to \E\to \OO(1/(h-1)) \to 0, \] 
where $\E$ is semistable of slope $1/h$ at every geometric point of $S$.
\end{proposition}

The diamond $\Omega_{h,h-1}$ can be interpreted as the complement of the origin in the Banach--Colmez space $\BC(-1/(h-1))$ which classifies extensions of $\OO$ by $\OO(-1/(h-1))$. 
\cref{thmx:cohom} can be reduced to the following calculation of the coherent cohomology of $\Omega_{h,h-1}$, which states that this punctured Banach--Colmez space is a homology sphere of dimension $h$.  

\begin{thmx}[\cref{PropCohoOfBC}] \label{ThmAcyclicity} There is an isomorphism of graded modules with $\G_{h-1}$-action:
\[ H^*(\Omega_{h,h-1},\OO) \isom k\oplus k(-1)[-h]. \] 
Here the $(-1)$ refers to the characteristic $p$ Tate twist (\cref{DefnTateTwist}).
\end{thmx}

To prove \cref{ThmAcyclicity}, we leverage a strange phenomenon concerning negative Banach--Colmez spaces:  when base changed to an algebraically closed perfectoid field $C$ of characteristic 0, they become isomorphic to quotients of affine spaces by finite-dimensional $\Q_p$-vector spaces.   In our situation, the relevant isomorphism reads (\cref{rslt:BaseChangeOfOmegaToC}) 
\[ \Omega_{h,h-1}\times \Spd C\isom \HH^{1,\Diamond}_{F,C}/F,\]
where $F\subset C$ is an arbitrary extension of $\Q_p$ of degree $h-1$, and 
where $\HH^1_F=\mathbb{A}^1_F\backslash F$ is 1-dimensional Drinfeld symmetric space relative to the field $F$.  

Therefore the $\OO$-cohomology of $\Omega_{h,h-1}\times\Spd C$ is isomorphic to the derived $F$-fixed points of the cohomology $R\Gamma(\HH^1_{F,C},\OO^\flat)$, where $\OO^\flat$ refers to the tilted structure sheaf on the v-site of the rigid-analytic space $\HH^1_{F,C}$. 

The $\OO^\flat$-cohomology of $\HH^1_{F,C}$ is computed in \S\ref{sec:oflat}.  It relies on the 6-functor formalism for $\OO^\flat$-cohomology developed in \cite{ALBM-6ff-on-FF-curve}.  

\begin{thmx}[\cref{thm:OFlatCohomDrinfeld}]
\label{thm:OflatCohomologyIntro}
 There are isomorphisms of solid $C^\flat$-modules with $\GL_2(F)$-action:
    \[ 
    H^i(\HH^1_{F,C},\OO^\flat)=\begin{cases}
        C^\flat, & i=0,\\
        \St(C^\flat)^*, & i=1,\\
        0, & \text{otherwise.}
    \end{cases}
    \]
    Here $\St(C^\flat)^*$ is the $C^\flat$-linear dual of the Steinberg representation $\St(C^\flat)$, defined as the quotient of the space of locally constant functions $\mathbb{P}^1(F)\to \C^\flat$ by the constant functions.  
\end{thmx}

The description of the cohomology in Theorem \ref{thm:OflatCohomologyIntro} has the same shape as the result of Schneider--Stuhler \cite{SchneiderStuhler}, which calculates the cohomology of Drinfeld symmetric space with respect to any cohomology theory (such as de Rham or $\ell$-adic \'etale) satisfying certain axioms, including homotopy-invariance by the closed disk $D$.  We are unable to apply \cite{SchneiderStuhler} directly, as in fact $H^1(D,\OO^\flat)$ is quite large.  Key to Theorem \ref{thm:OflatCohomologyIntro} is the observation that the restriction map $H^1(D,\OO^\flat)\to H^1(D',\OO^\flat)$ is zero whenever $D'\subset D$ is a sub-disk of sufficiently small radius.

Deriving \cref{ThmAcyclicity} from \cref{thm:OflatCohomologyIntro} is a matter of taking derived invariants under the unipotent subgroup $F$ of $\GL_2(F)$, and then descending the result from $C$ to the base field $k$.  For the first step, we encounter the problem of taking the derived $F$-invariants of $\St(C^\flat)^*$, or dually the derived $F$-coinvariants of $\St(C^\flat)$, which we can think of as a derived Jacquet module of the Steinberg representation.  Over the discrete base $k$, this problem was solved by Heyer \cite{Heyer2023Left}, who shows that the module of coinvariants is invertible and supported in degree $[F:\Q_p]=h-1$.  We adapt Heyer's computation to the setting of $C^\flat$-coefficients to compute $H^*(\HH^1_{F,C}/F,\OO^\flat)\isom H^*(\Omega_{h,h-1}\times\Spd C,\OO)$ in \cref{PropCohoOfHModF}; the cohomology is 1-dimensional in degrees 0 and $h$.

For the second step of descending along $\Spd C\to\Spd k$, the idea is to assume that $k=\mathbb{F}_q$ is finite and to apply derived invariants to $R\Gamma(\Omega_{h,h-1}\times \Spd C,\OO)$ under the $q$th power Frobenius map on $C^\flat$.  This operation recovers $R\Gamma(\Omega_{h,h-1},\OO)$ by a phenomenon reminiscent of Drinfeld's lemma in the context of diamonds, see \cref{LemmaDrinfeld}.

\subsection{A family of two-towers isomorphisms, after Faltings}\label{ssec:intro_twotowers}

In fact we found a family of geometric results for each height $0\leq t\leq h$, of which \cref{ThmMeaningOfL} is the case $t=h-1$.  These may be of independent interest, so we discuss them here.  It may well be possible to use these results to investigate $L_{K(t)}\SS_{K(h)}$.

Let $A_h=\pi_0E_h$, and let $\LT_h=\Spf A_h$ be the corresponding affine formal scheme, so that $\LT_h$ is the formal moduli space of deformations of $\Gamma_h$.  There is a descending chain of $\G_h$-stable closed formal subschemes $\Spf A_h/I_t$ for $0\leq t\leq h$, corresponding to an ascending chain of $\G_h$-invariant ideals $I_t\subset A_h$.  We put $I_0=0$, and for each $1\leq t\leq h$, $\Spf A_h/I_t$ classifies those deformations $\Gamma$ in characteristic $p$ for which the $t$th power of the Frobenius isogeny divides $p$ in $\End \Gamma$.  

Explicitly, one can choose formal coordinates $p=u_0,u_1,\dots,u_{h-1}\in A_h$ for which 
\[ A_h = W(k)\powerseries{u_1,\dots,u_{h-1}},\]
and $I_t=(u_0,\dots,u_{t-1})$ for $1\leq t\leq h$.  

Let $\LT_h^{\ad}=\Spa(A_h,A_h)$ denote the affinoid adic space associated to the formal scheme $\LT_h$.  The points of $\LT_h^{\ad}$ are continuous valuations $f\mapsto \abs{f}$ on $A_h$.  In light of the remarks above, there is a $\G_h$-stable height stratification
\[ \LT_h^{\ad}=\coprod_{t=0}^h \LT_{h,t}, \]
where $\LT_{h,t}$ is the locally closed subspace of $\LT_h^{\ad}$ defined by the conditions 
\[ 
    \abs{u_0}=\cdots=\abs{u_{t-1}}=0,\quad \abs{u_t}\neq 0; 
\]
here, by conventions, we set $u_h=1$. In particular $\LT_{h,0}$ is defined by $\abs{p}\neq 0$;  this is isomorphic to the $(h-1)$-dimensional unit ball over $\Spa W(k)[1/p]$. At the other extreme, $\LT_{h,h}=\Spa k$.  For intermediate values $0<t<h$, the adic space $\LT_{h,t}$ is isomorphic to the $(h-t-1)$-dimensional ball over $\Spa k\laurentseries{u_t}$.  
In the coheight 1 case, $\LT_{h,h-1}=\Spa K$ for our field $K=k\laurentseries{u_{h-1}}$.  

We now pass to the associated diamond $\LT_{h,t}^{\Diamond}$.  
For every geometric point $\Spa(M,M^+)$ of $\LT_{h,t}^{\Diamond}$ with $M$ algebraically closed, the base change $\Gamma^{\un}_{h,M}$ is isomorphic to $\Gamma_{t}\oplus  (\Q_p/\Z_p)^{h-t}$.  Therefore there exists a $\G_h$-equivariant pro-\'etale torsor
\[ \LT_{h,t}^{\infty}\to \LT_{h,t}^{\Diamond} \]
for the group $\overline{\G}_t\times \GL_{h-t}(\Z_p)$ which classifies such trivializations.  We note that the group $\G_{h,t}\times \GL_{h-t}(\Z_p)$ acts on $\LT_{h,t}^{\infty}$, where $\G_{h,t}=\G_h\times_{\Gal(k/\F_p)} \G_t$.

We prove the following alternative characterization of $\LT_{h,t}^\infty$. 

\begin{thmx}[\cref{thm:twotowersinchapter}]   
\label{thm:twotowersintro}
Let $0< t<h$.  As a diamond with $\G_{h,t}\times\GL_{h-t}(\Z_p)$-action, $\LT_{h,t}^{\infty}$ is isomorphic to the sheaf which associates to a perfectoid space $S$ the set of exact sequences of vector bundles on the Fargues--Fontaine curve $X_S$
\[ 0 \to \OO^{\oplus(h-t)}\to \OO(1/h)\to \OO(1/t) \to 0\] 
such that the induced isomorphism on determinants $\wedge^h \OO(1/h)\isom \wedge^t \OO(1/t)$ is a $\Z_p^\times$-multiple of a fixed one.  
\end{thmx}

Whereas \cite{ScholzeWeinstein} shows that infinite-level Lubin--Tate space $\LT_{h,0}^{\infty}$ in characteristic 0 classifies modifications of vector bundles on $X_S$ along the divisor $S^\sharp$, Theorem \ref{thm:twotowersintro} shows that infinite-level Lubin--Tate space in characteristic $p$ classifies extensions of vector bundles.

We summarize here the main points of the proof of \cref{thm:twotowersintro}.  Consider an affinoid perfectoid space $S=\Spa(R,R^+)$ in characteristic $p$.   Since $R^+$ is perfect, the categories of $p$-divisible groups over $R^+$ and Dieudonn\'e modules over $R^+$ are anti-equivalent \cite{LauSmoothness}. From a Dieudonn\'e module over $R^+$, we build a {\em meromorphic vector bundle} over $S$ in the sense of \cite{GleasonIvanovZillinger}:  this means a vector bundle over the Fargues--Fontaine curve relative to $S$ whose Frobenius is meromorphic along $p=0$.  The stack $\Bun^{\text{mer}}$ classifying meromorphic vector bundles admits a ``special fiber'' morphism to the moduli stack of vector bundles $\Bun$, and a ``generic fiber'' morphism to the moduli stack of isocrystals $\mathrm{Isoc}$.  

Up to v-sheafifying, a morphism $S\to \LT_{h,t}^\infty$ corresponds to a $p$-divisible group over $R$ whose special and generic fibers have been identified with $\Gamma_h$ and $\Gamma_t\oplus (\Q_p/\Z_p)^{h-t}$, respectively.  This in turn corresponds to an $S$-point of the fiber of $\Bun^{\text{mer}}$ over the point of $\Bun\times\mathrm{Isoc}$ determined by those respective $p$-divisible groups.  The main theorem of \cite{GleasonIvanovZillinger} identifies this fiber with a moduli space of extensions of vector bundles on the Fargues--Fontaine curve, which is why such extensions make an appearance in \cref{thm:twotowersintro}.

If $\LT_{h,t}^{\infty}\to \LT_{h,t}$ is an analogue of the Lubin--Tate tower, the analogue of the Drinfeld tower appears when one quotients by $\overline{\G}_h$.  This quotient $\Omega_{h,t}\coloneqq \LT_{h,t}^\infty/\overline{\G}_h$ is the diamond which classifies those extensions of vector bundles on the Fargues--Fontaine curve
\[ 
0 \to \OO^{\oplus (h-t)} \to \E\to \OO(1/t) \to 0, 
\] 
where $\E$ is semistable of slope $1/h$. Dually, $\Omega_{h,t}$ classifies extensions of $\OO^{\oplus (h-t)}$ by $\OO(-1/t)$ which are semistable of slope $-1/h$.  As such, there is an open immersion
\[ \Omega_{h,t}\hookrightarrow \mathrm{BC}(-1/t)^{\oplus(h-t)}\]
where $\mathrm{BC}(-1/t)$ is the Banach--Colmez space sending $S$ to $H^1(X_S,\OO(-1/t))$.  The image is exactly those vectors not lying in any $\Q_p$-rational hyperplane.  We think of $\Omega_{h,t}$ as a period domain in characteristic $p$, analogous to the $p$-adic period domains appearing in the theory of Rapoport--Zink spaces.  

According to \cref{thm:twotowersintro}, there is a ``two-towers'' diagram of diamonds
\[
\xymatrix{
& \LT_{h,t}^{\infty} \ar[dl]_{\overline{\G}_t\times\GL_{h-t}(\Z_p)} \ar[dr]^{\overline{\G}_h} & \\
\LT_{h,t}^{\Diamond} && \Omega_{h,t}
}
\] 
in which each morphism is a pro-\'etale torsor for the indicated group.  The diagram recovers Faltings' original picture in the case $t=0$, if one interprets $\Omega_{h,0}$ as the diamond associated to Drinfeld symmetric space $\HH^{h-1}_{\Q_p}$, base changed to $W(\overline{\F}_p)[1/p]$.  

The two-tower diagram witnesses an isomorphism of v-stacks:
\begin{equation}
    \label{eq:TwoTowerIntro}[\LT_{h,t}^{\Diamond}/\G_h] \isom [\Omega_{h,t}/(\G_t\times \GL_{h-t}(\Z_p))]. 
\end{equation} 

We offer some remarks on the prospect of applying \eqref{eq:TwoTowerIntro} to produce any new results about $L_{K(t)}\SS_{K(h)}$ for $0<t<h-1$, say by taking $\OO$-cohomology on either side.  An immediate difficulty arises:  $H^0(\LT_{h,t},\OO)$ is the ring of convergent series on a (positive-dimensional) open ball over $k\laurentseries{u_t}$, and this is much larger than $\pi_0L_{K(t)}E_h/I_t=k\powerseries{u_t,u_{t+1},\dots,u_{h-1}}[u_t^{-1}]$. In \cite{BSSW} this difficulty was skirted in the $t=0$ case by taking $\OO^+$-cohomology instead, and then inverting $p$ on the Lubin--Tate side to access the $\G_h$-cohomology of $\pi_0L_{K(0)}E_h$.
To carry out the corresponding strategy for $0<t<h-1$, note that the analogue of the $\G_h$-invariant element $p$ is the height $t$ Hasse invariant in $\pi_{2(p^t-1)}E_h/I_t$.  There is a $\G_h$-equivariant invertible $\OO^+$-module $\mathcal{L}$ on $\LT_{h,t}$ coming from the Lie algebra of the universal deformation, whose global sections are exactly $\pi_{-2} E_h/I_t$.   
So one could take cohomology in \eqref{eq:TwoTowerIntro} with coefficients in powers of $\mathcal{L}^{\Diamond}$, and then invert the Hasse invariant on the Lubin--Tate side to access the $\G_h$-cohomology of $\pi_* L_{K(t)}E_h/I_t$.  We now run into another difficulty:  it is not immediately clear how to describe the invertible $\OO^+$-module on the right side of \eqref{eq:TwoTowerIntro} corresponding to $\mathcal{L}^{\Diamond}$, let alone compute its cohomology on $\Omega_{h,t}$.

\subsection{The strong splitting conjecture, and the proof of \cref{thmx:height2}} 

Assume that $p-1\nmid h-1$.  Our \cref{ConjStrongSplitting} states that the unit map $\SS_{K(h-1)}\to L_{K(h-1)}\SS_{K(h),\infty}$ is an isomorphism.  We observed that the composition of $\G_{h-1}$-equivariant $\EE_\infty$-ring spectra 
\[ E_{h-1}^\solid\to (\B^{\solid})^{h\overline{\G}_h^1}\to (\widehat{\B}^\solid)^{h\overline{\G}_h^1}\]  is an isomorphism.  After taking $\G_{h-1}$-fixed points and performing various manipulations (see Lemma \ref{lem:TFAEsCSC}), one finds that the composition 
\[ \SS_{K(h-1)}\to ((\B^{\solid})^{h\overline{\G}_{h-1}})^{h\G_h^1}(*) \to ((\widehat{\B}^{\solid})^{h\overline{\G}_{h-1}})^{h\G_h^1}(*) \] 
is also an isomorphism.  Regarding the middle term, we have $(\B^{\solid})^{h\overline{\G}_{h-1}}\isom L_{K(h-1)}^{\solid}E_h^{\solid}$, and the $\G_h^1$-fixed points of this evaluated at $(*)$ is $L_{K(h-1)}\SS_{K(h),\infty}$.  Therefore Conjecture \ref{ConjStrongSplitting} is equivalent to the statement that the second arrow is an isomorphism.  The naturality properties of our $E$-theory functor show that $(\widehat{\B}^{\solid})^{h\overline{\G}_{h-1}}$ is the height $h-1$ Morava $E$-theory
over the ring of derived $\overline{\G}_{h-1}$-invariants in $\widehat{L}$. 
In fact this ring agrees with the ring of un-derived invariants, which is the  
 completed perfection $\widehat{K}_p=k\laurentseries{u_{h-1}^{1/p^\infty}}$ of $K$.  

We find therefore that \cref{ConjStrongSplitting} reduces to the statement that the map
\[ L_{K(h-1)}^{\solid}E_h^{\solid}\to E^{\solid}(\widehat{K}_p,\Gamma^{\un}_{h,\widehat{K}_p})
\] 
induces an isomorphism after passage to $\G_h^1$-fixed points, and this homotopical statement is in turn reducible to an algebraic statement, namely that 
    \begin{equation}\label{eq:comparisonmap}
        H^*_{\cts}(\G_h^1,\pi_*E_h\otimes_{\pi_0E_h} K)\to H^*_{\cts}(\G_h^1,\pi_*E_h\otimes_{\pi_0E_h} \widehat{K}_p)
    \end{equation}
is an isomorphism (see \cref{PropDeperfection}).

Therefore we are interested in criteria under which a topological $\F_p$-vector space $V$ with action by a $p$-adic Lie group $G$ would have $H^*_{\cts}(G,V)=0$, with the intent of applying such criteria to $G$ a convenient closed normal subgroup of $\G_h^1$ and $V=\pi_*E_h\otimes_{\pi_0E_h} \widehat{K}_p/K$.  

Here we were inspired by Tate's original technique \cite{TatePDiv} for showing acyclicity of the action of $\Z_p$ on the module $V=\widehat{F}_\infty/F$, where $F$ is a finite extension of $\Q_p$ and $F_\infty$ is a ramified $\Z_p$-extension of $F$.  Letting $g$ be a topological generator of $\Z_p$, Tate shows that $(g-1)$ is invertible on $F_\infty/F$, and bounds the operator norm of $(g-1)^{-1}$, which then extends to a continuous operator on $V$ which is a nullhomotopy for the complex computing $H^*_{\cts}(\Z_p,V)$. 

In the criteria for acyclicity we devised, $G$ is a uniform pro-$p$ group of dimension $d$, with topological generators $g_1,\dots,g_d$, which acts on an $\F_p$-vector space $V$ which is separated and complete for an ultrametric valuation $\val_V$, in such a way that each $(g_i-1)$ is a strict contraction on $V$.  

\begin{thmx}[Theorem \ref{ThmAcyclicityOne}] \label{thmx:acyclicity}
Assume that the operator $(g_1-1)$ has bounded inverse $(g_1-1)^{-1}$ on $V$.  Then $H^*_{\cts}(G,V)=0$.
\end{thmx}

The proof of \ref{thmx:acyclicity} makes essential use of Lazard's study \cite{Lazard} of continuous cohomology of uniform pro-$p$ groups, whereby this cohomology can be computed by a finite complex which lifts the Koszul complex for a polynomial algebra.  The existence of the operator $(g_1-1)^{-1}$ allows us to write down a nullhomotopy for Lazard's complex.

Using some explicit formulas of Salch for the Morava action at height 2, we were able to verify the acyclicity criteria in \cref{thmx:acyclicity} in the setting of the strong splitting conjecture (\cref{ThmDeperf2});  as a consequence we proved \cref{thmx:height2}.   

Regarding the strong splitting conjecture in heights $h\geq 3$, we find it unlikely that a similar strategy using \cref{thmx:acyclicity} will prevail. Instead it may be more fruitful to seek a reformulation of the two-tower isomorphism in a category where test objects are not required to arise from perfectoid rings, thus obviating the need for the completed spectrum $\widehat{\B}$ in the first place.

\begin{rmk} \label{rmk:astra-compute-kernel}
    Preliminary calculations with an AI system indicate that 
    \[
        H^*_{\cts}(\G_h^1,\pi_*E_h\otimes_{\pi_0E_h} K)\to H^*_{\cts}(\G_h^1,\pi_*E_h\otimes_{\pi_0E_h} \widehat{K}_p)
    \]
    is \emph{not} an isomorphism when $p \leq h-1$.
    As a consequence, we expect the statement of \Cref{ConjStrongSplitting} to be false when $p \leq h-1$.
\end{rmk}

\subsection{Outline of the paper}

The first two sections form the $p$-adic backbone for the chromatic results of this paper. The main subject of \Cref{sec:twotowers} is our intermediate height generalization of the two-tower isomorphism (\cref{thm:twotowersintro}). In \cref{sec:oflat}, we first compute the coherent cohomology of $\Omega_{h,h-1}$ (\cref{ThmAcyclicity}) and then translate it via the two-tower isomorphism to deduce \cref{thmx:cohom}. \Cref{sec:solidLT} deals with our generalization (\cref{thmx:solidLT}) of Lurie's construction of Lubin--Tate theories and uses it to construct the completion $\widehat{\B}$ of Torii's spectrum $\B$. After a short recollection on chromatic theory, \cref{sec:CSC} then assembles our results to prove \cref{thmx:main} and its consequences. The last section, \cref{sec:height2}, establishes the acyclicity criterion (\cref{thmx:acyclicity}) and applies it to reprove the strong form of the chromatic splitting conjecture for height 2 and odd primes, \cref{thmx:height2}. 
Finally, Appendix~\ref{app:actions} collects some material on actions by profinite groups which will be used crucially in the main body of the paper.

\subsection*{Notation}

For convenience, we provide here a list of frequently used symbols.

\begin{itemize}

\item $A_h$, the Lubin--Tate deformation ring in height $h$, isomorphic to $W(k)\powerseries{u_1,\dots,u_{h-1}}$, see \cref{defn:LTring}.
\item $\mathbb{A}_h$, the ring spectrum defined in \eqref{eq:introAph}.
\item $\B$, the Torii spectrum, see \eqref{eq:B}.
\item $\BC(\lambda)$, the Banach--Colmez space of slope $\lambda$, see \cref{def:BCspaces}.
\item $\Bun$, the v-stack of vector bundles on the Fargues--Fontaine curve, see \cref{def:IsocBun}.
\item $\G_h$, the Morava stabilizer group, constructed relative to the field $k$. There are several variants and subgroups used throughout the paper, including:
\begin{itemize}
    \item $S\G_h = \ker(\G_h \to \Z_p^{\times})$ and  
    $\G_h^1 =\ker(\G_h \to \Z_p^\times\to \Z_p)$, see \eqref{eq:allgroups} and \eqref{eq:allgroups2}.
    \item $\overline{\G}_h=\ker(\G_h\to \Gal(k/\F_p))$, and similarly for $S\overline{\G}_h$ and $\overline{\G}_h^1$.
\end{itemize}
\item $\Gamma_h$, a 1-dimensional formal group of height $h$ over $k$, often viewed as a $p$-divisible group.
\item $\HH_F^n$, Drinfeld symmetric space of dimension $n$ relative to a finite extension $F/\Q_p$, see \cref{defn:DrinfeldSymmSpace}.
\item $I_t$, the Landweber ideal $(p,u_1,\dots,u_{t-1})$ in $A_h$.
\item $\Isoc$, the v-stack of isocrystals, see \cref{def:IsocBun}.
\item $k$, an algebraic extension of $\F_p$, assumed to contain $\F_{p^h}$ whenever the context references a height $h$.  
\item $K$, the residue field of $\pi_0L_{K(h-1)}E_h$, isomorphic to $k\laurentseries{u_{h-1}}$, see \eqref{eq:K}.
\item $L$, the Gross--Torii extension of $K$, see \cref{def:torii_ext}.
\item $\LT_h$, the Lubin--Tate deformation space $\Spf A_h$, see \cref{DefLT}. 
\item $\Omega_{h,t}$, the characteristic $p$ period domain, see \cref{defn:Omegaht}.
\item $\RZ_h$, Rapoport--Zink space in height $h$, see \cref{DefRZ}. 
\item $\SS$, the sphere spectrum;  also its localization $\SS_{K(h)}=L_{K(h)}\SS$.
\item $\Sht_{h,\geq 1}$, the v-stack of shtukas with a paw at $p$, see \cref{DefnShtukas}.
\end{itemize}

\subsection*{Statement on the use of AI}
Most of this article was completed before the age of AI in mathematics. We have used AI in the following ways:
\begin{itemize}
    \item To write magma code to assist with the computations of \Cref{sec:height2}.
    \item To locate references for various facts in higher algebra in \Cref{sec:oflat}.
    \item To navigate some of the higher categorical technicalities in the appendix. In particular, AI suggested the use of \cite[Variant 2.2.3]{Arpon}, saving us from the finer points of the theory of lax symmetric monoidal comonads, as well as \Cref{lem:chatgpt} and its proof.
    \item To check the manuscript for mathematical and typographical errors. This resulted in many edits, most of which were minor; the most significant of these was a correction to and strengthening of the statement of \Cref{ThmAcyclicityOne}.
    \item To explore whether the comparison map \eqref{eq:comparisonmap} is an isomorphism, resulting in \Cref{rmk:astra-compute-kernel} and a comment at the end of \Cref{sec:height2}, as well as a modification of \Cref{ConjStrongSplitting} to exclude all $p \leq h$.
\end{itemize}

\subsection*{Acknowledgements}

This article arose out of interactions that occurred during a December 2024 workshop hosted by the American Institute of Mathematics, organized by three of us (TB, TS, JW) together with Nathaniel Stapleton.  We would like to thank Agn\`es Beaudry, Mike Hopkins, Jeremy Hahn, Piotr Pstragowski, Padmavathi Srinivasan, Yujie Xu, Zheng Yang, and Ningchuan Zhang for many helpful conversations.  Special thanks are due to Andrew Salch for sharing formulas on the Morava stabilizer action which appear in Chapter 6.  
XZ thanks Dylan Pentland and Niven Achenjang for the opportunity to present the results in this manuscript at the Harvard Number Theory Seminar. 

TB is supported by the European Research Council (ERC) under Horizon Europe (grant No.~101042990) and would like to thank the Max Planck Institute for its hospitality. JW is supported by NSF Grant DMS-2401472 and a grant from the Simons Foundation.

\tableofcontents

\newpage
\section{The two tower isomorphism}\label{sec:twotowers}

This chapter provides the proof of \cref{thm:twotowersintro}, which is the main arithmetic geometry result powering the homotopical results in this paper.  For the reader's convenience, we review the context and restate the theorem.

For an integer $h\geq 1$, the formal scheme $\LT_h$ classifies deformations of a 1-dimensional formal group of height $h$ over a perfect field $k$ of characteristic $p$, which for the purposes of this introduction will be $\overline{\F}_p$. The associated adic space $\LT_h^{\ad}$ 
admits a stratification by locally closed subspaces $\LT_{h,t}$ for $0\leq t\leq h$.  The space $\LT_{h,t}$ is the locus of those deformations whose generic fiber has height $t$.  

Let $\LT_{h,t}^{\diamondsuit}$ be the diamond associated to the adic space $\LT_{h,t}$.  Lying above $\LT_{h,t}^{\diamondsuit}$ is the infinite-level space $\LT_{h,t}^{\infty}$, which is the v-sheaf classifying deformations together with a trivialization of the generic fiber of the associated $p$-divisible group.  Then as a v-sheaf over $\Spd k$, the infinite-level space $\LT_{h,t}^{\infty}$ admits an action of $\G_{h,t}\times\GL_{h-t}(\Z_p)$ lying over the action on $\Spd k$ via $\Gal(k/\F_p)$.   

The open stratum $\LT_{h,0}$ is simply the characteristic-0 locus in $\LT_h^{\ad}$, and $\LT_{h,0}^\infty\to \LT_{h,0}^{\diamondsuit}$ is the classical Lubin--Tate tower.  The infinite-level space $\LT_{h,0}^{\infty}$ admits an alternative description in terms of linear algebra.  For a test object $S$ (an affinoid perfectoid space over $W(k)[1/p]$), the $S$-points of $\LT_{h,0}^{\infty}$ classify morphisms of vector bundles $\OO^{\oplus h} \to \OO(1/h)$ on the Fargues--Fontaine curve $X_{S^\flat}$ whose cokernel is a line bundle over the distinguished divisor on $X_{S^\flat}$ corresponding to the untilt $S$ of $S^\flat$ \cite[\S 6.4]{ScholzeWeinstein}.

When $t\geq 1$, the locus $\LT_{h,t}$ lies over $k$.  We also found an alternative interpretation of $\LT_{h,t}^\infty$ in terms of vector bundles on the Fargues--Fontaine curve.  In this context, 
the test objects $S$ are affinoid perfectoid spaces in characteristic $p$, and so $S=S^\flat$; the curve $X_S$ has no distinguished divisor. Once again, here is Theorem \ref{thm:twotowersintro}:

\begin{theorem}
\label{thm:twotowersinchapter}
Assume $0<t\leq h$.  As a v-sheaf over $\Spd k$ with $\G_{h,t}\times\GL_{h-t}(\Z_p)$-action, $\LT_{h,t}^{\infty}$ is isomorphic to the sheaf which associates to a perfectoid space $S$ the set of exact sequences of vector bundles on the Fargues--Fontaine curve $X_S$:
\[ 0 \to \OO^{\oplus(h-t)}\to \OO(1/h)\to \OO(1/t) \to 0\] 
such that the induced isomorphism on determinants $\wedge^h \OO(1/h)\isom \wedge^t \OO(1/t)$ is a $\Z_p^\times$-multiple of a fixed one.  
\end{theorem}

The proof of Theorem \ref{thm:twotowersinchapter} follows a chain of equivalences starting from $\LT_{h,t}$ and leading to a moduli space of extensions of vector bundles.  Along the way we encounter a sequence of closely related moduli spaces. For clarity's sake, we have elected to give each of these a distinct notation.  We offer the reader the following guide to the moduli spaces appearing in this chapter.

\begin{itemize}
    \item $\LT_h$, the deformation space of a formal group of height $h$.  Related:  the adic space $\LT_h^{\ad}$ and its strata $\LT_{h,t}$ for $0\leq t\leq h$, and in particular the characteristic $p$ locus $\LT_{h,\geq 1}$.  
    \item $\RZ_h$, the Rapoport--Zink deformation space of a formal 1-dimensional $p$-divisible group of height $h$ up to quasi-isogeny. Related:  the adic space $\RZ_h^{\ad}$ and its strata $\RZ_{h,t}$ for $0\leq t\leq h$, and in particular the characteristic $p$ locus $\RZ_{h,\geq 1}$. 
\item $\Sht_{h,\geq 1}$, the v-sheaf classifying {\em shtukas with a paw at $p$} (Definition \ref{DefnShtukas}), which are vector bundles on the curve $\Y_{S,[0,\infty)}$ together with a $\phi$-linear automorphism which is meromorphic at $p$, plus some additional data.  Related:  the strata $\Sht_{h,t}$ for $1\leq t\leq h$. 
\item Infinite-level spaces $\LT_{h,t}^\infty$, $\RZ_{h,t}^\infty$, and $\Sht_{h,t}^\infty$, which appear by ``trivializing the generic fiber''.  
\end{itemize}

$\RZ_{h}$ is isomorphic to a disjoint union of $\Z$ copies of $\LT_{h}$ (Proposition \ref{PropRZ}), and $\Sht_{h,\geq 1}$ is isomorphic to the diamond version of $\RZ_{h,\geq 1}$ (\cref{prop:RZtoSht}), in a manner compatible with passage to the infinite-level spaces.  Thereby Theorem \ref{thm:twotowersinchapter} may be reduced to the question of finding an alternative description for $\Sht_{h,t}^\infty$.

At this point our work makes contact with \cite{GleasonIvanovZillinger} and the stack $\Bun^{\mer}$ of {\em meromorphic vector bundles on the Fargues--Fontaine curve}.  The motivation for \cite{GleasonIvanovZillinger} is an ambitious program to relate two $\Q_p$-linear v-stacks which appear in distinct approaches to the categorical local Langlands conjecture:

\begin{itemize}
    \item $\Isoc$, the stack classifying isocrystals over a perfectoid $\F_p$-algebra.
    \item $\Bun$, the stack of vector bundles on the Fargues--Fontaine curve.
\end{itemize}
The stacks $\Isoc$ and $\Bun$ have the same underlying set of points (corresponding to lists of slopes), but  otherwise do not seem to be directly related. 

The stack $\Bun^{\mer}$ classifies vector bundles on the Fargues--Fontaine curve together with a meromorphic extension to the locus $p=0$ (Definition \ref{DefnBunMer}).  There is a correspondence $\Bun^{\mer}\to \Isoc\times \Bun$ carrying such a meromorphic vector bundle to its generic and special fiber, respectively (Definition \ref{DefnBunMerCorresp}, see Remark \ref{rmk:genericspecial} for an explanation for this terminology.)

We use a Beauville-Lazslo gluing argument to identify $\Sht_{h,t}^{\infty}$ with the pullback of $\Bun^{\mer}$ over a specific point $\ast\to \Isoc\times\Bun$ (Theorem \ref{ThmShtVSMer}).  On the other hand, one of the main results of \cite{GleasonIvanovZillinger} identifies the pullback of $\Bun^{\mer}$ over any point of $\Isoc$, in terms of vector bundles on the Fargues--Fontaine curve equipped with a slope filtration (Theorem \ref{ThmBunMerMain}).  Combining these two results shows that $\Sht_{h,t}^\infty$ classifies two-step filtrations on $\OO(1/h)$, where the quotients are identified with $\OO^{\oplus (h-t)}$ and $\OO(1/t)$, respectively (Theorem  \ref{ThmShtInf}).  From there, Theorem \ref{thm:twotowersinchapter} follows easily. 

The chapter concludes with a discussion of the case $t=h-1$.  Recall in that setting the Gross--Torii extension $L/K$ \ref{def:torii_ext}.  We have $\LT_{h,h-1}=\Spa K$, and the diamond $\Spd \widehat{L}$ appears as a connected component of $\LT_{h,h-1}^{\infty}$.  Theorem \ref{ThmMeaningOfL}, on the interpretation of $\Spd \widehat{L}$ in terms of extensions of vector bundles on the Fargues--Fontaines curve, appears as Proposition \ref{PropL}(2).  

\subsection{The Lubin--Tate and Rapoport--Zink deformation spaces}
Let $h\geq 1$, and let $k\subset \overline{\F}_p$ be a subfield containing $\F_{p^h}$.  
Let $\fg_h$ be the Honda formal group of height $h$ over $k$.  With respect to a parameter $X$, the $p$-series in $\fg_h$ is $[p]_{\fg_h}(X)=X^{p^h}$.  The ring of endomorphisms $\End \fg_h=\OO_{D_h}$ is the ring of integers in a division algebra $D_h$ over $\Q_p$ of invariant $1/h$, with unit group $\Aut \fg_h=\OO_{D_h}^\times=\overline{\mathbb{G}}_h$.  

\begin{defn} \label{DefLT} The Lubin--Tate deformation functor $\LT_h$ inputs a complete Noetherian local ring $R$ with residue field $k$, and outputs the set of pairs $(\Gamma,\iota)$, where $\Gamma$ is a formal group over $R$, and $\iota\from \fg_h\isom \fg\otimes_R k$ is an isomorphism.  The action of $\G_h$ on $k$ lifts to an action on $\LT_h$. 
\end{defn}

By a classical result of Lubin--Tate \cite{LubinTate}, $\LT_h$ is representable by an affine formal scheme $\Spf A_h$, where
\begin{equation}
    \label{defn:LTring}
A_h\isom W(k)\powerseries{u_1,\dots,u_{h-1}}
\end{equation} 
is the ring of power series over the ring of Witt vectors $W(k)$ in formal parameters $u_1,\dots,u_{h-1}$.  Let $\fg^{\un}_h$ denote the universal deformation of $\fg_h$ to $A_h$.   
 
It is possible to choose the parameters $u_1,\dots,u_{h-1}$ in such a way that 
\begin{equation}
    \label{EqPSeries}
    [p]_{\fg^{\un}_h}(X)=pX+_{\fg^{\un}_h} u_1X^p
+_{\fg^{\un}_h} \cdots+_{\fg^{\un}_h} u_{h-1}X^{p^{h-1}}+_{\fg^{\un}_h}  X^{p^h}, 
\end{equation}
by specializing the parameters in the logarithm of the universal $p$-typical formal group law \cite{Araki}.

Having chosen the parameters in this manner, we define a chain of ideals $I_t\subset A_h$ for $0\leq t\leq h$ by $I_0=(0)$ and for $1\leq t\leq h$, 
\[ I_t=(p,u_1,\dots,u_{t-1}). \]
For $1\leq t\leq h$, the vanishing locus $V(I_t)$ is a formal subscheme of $\LT_h$, which lives entirely in the characteristic $p$ locus.  It classifies those deformations $(\fg,\iota)$ in characteristic $p$ for which $[p]_{\fg}$ factors through the $t$th power of the relative Frobenius map $\fg\to \fg^{(p)}$.  As this condition does not depend on $\iota$, the ideals $I_t$ are stable under the $\G_h$-action.

Over such $R$ as in Definition \ref{DefLT}, there is an equivalence between formal groups of finite height and connected $p$-divisible groups \cite{TatePDiv}.  For many reasons it will be profitable to immediately translate our discussion into the language of $p$-divisible groups.  We will use the same notation $\Gamma_h$ for the $p$-divisible group of the Honda formal group of height $h$ over $k$.  The deformation problem $\RZ_h$ appearing below, which replaces the isomorphism in \cref{DefLT} with a quasi-isogeny, was introduced by Rapoport--Zink \cite{RapoportZink}.   One could define $\RZ_h$ as a functor on $W(k)$-algebras in which $p$ is nilpotent, or equivalently, as we do here, on adic $W(k)$-algebras.  (An adic ring $R$ is a topological ring which is complete for the $I$-adic topology for some ideal of definition $I$.  If $R$ is an adic $W(k)$-algebra, one may choose $I$ to contain $p$.)

The following material is discussed in a more leisurely manner in \cite[\S3.2]{HopkinsPicard}.

\begin{defn}
    \label{DefRZ} Let $\RZ_h$ be the functor which inputs an adic $W(k)$-algebra $R$ and outputs the set of isomorphism classes of pairs $(\Gamma,\alpha)$, where $\Gamma$ is a $p$-divisible group over $R$, and $\alpha\from \Gamma_h\otimes_k R/I\to \Gamma\otimes_R R/I$ is a quasi-isogeny.  Here $I$ is an unspecified ideal of definition of $R$ containing $p$.
\end{defn}

\begin{rmk} Given $\Gamma$, the set of $\alpha$ as in \cref{DefRZ} is independent of the choice of ideal of definition $I$: this is due to the rigidity property of quasi-isogenies \cite{DrinfeldCoverings}.  
\end{rmk}

The space $\RZ_h$ admits more symmetries than $\LT_h$.  Let $\G_h^0$ be the group of automorphisms of the pair $(\Gamma_h,k)$ in the isogeny category. 
Then $\G_h^0$ fits into a diagram
\[
\xymatrix{
& 0 \ar[d] & 0 \ar[d] & & \\
0 \ar[r] & \overline{\G}_h \ar[r] \ar[d] & \overline{\G}_h^0 \ar[r]^{\height}\ar[d] & \Z \ar[r] \ar[d]^{=} & 0 \\
0 \ar[r] &  \G_h \ar[r] \ar[d] & \G_h^0 \ar[r]^{\height} \ar[d] & \Z \ar[r] & 0\\
& \Gal(k/\F_p) \ar[d] \ar[r]_{=} & \Gal(k/\F_p) \ar[d] & & \\
& 0 & 0 &&
}
\]
in which rows and columns are exact.

The group $\G_h^0$ acts on $\RZ_h$.  
The height of the quasi-isogeny $\alpha$ is a locally constant $\Z$-valued function on $\RZ_h$, so there is a decomposition $\RZ_h=\coprod_{n\in\Z} \RZ_h^{(n)}$. 
The action of $\G_h^0$ on $\RZ_h$ permutes the components transitively; the subgroup $\G_h$ preserves each component.

\begin{proposition}[{\cite[Proposition 3.79]{RapoportZink}}] \label{PropRZ} $\RZ_h$ is representable by a formal scheme over $\Spf W(k)$.  There is a $\G_h$-equivariant isomorphism of formal schemes:  $\LT_h\isom \RZ_h^{(0)}$.  Reciprocally, $\RZ_h$ with its $\G_h^0$-action can be expressed in terms of $\LT_h$:  $\RZ_h=(\LT_h \times \G_h^0)/\G_h$.  
\end{proposition}

\subsection{The adic spaces $\LT_h^{\ad}$ and $\RZ_h^{\ad}$, and their height stratifications}

In order to state the main results of the chapter, it will be necessary to pass from formal schemes to their associated adic spaces, and then ultimately to diamonds.
Let $\LT_h^{\ad}=\Spa A_h$ be the adic space associated to the topological ring $A_h$. 

For a complete affinoid $W(k)$-algebra $(R,R^+)$, morphisms $\Spa(R,R^+)\to \LT_h^{\ad}$ are the same as continuous $W(k)$-algebra homomorphisms $A_h\to R^+$; such a homomorphism determines a formal group $\Gamma$ over $R^+$ by specialization of the universal deformation $\Gamma_h^{\un}$. 

The ascending chain of ideals $I_t\subset A_h$ gives rise to a descending chain of Zariski closed subsets $\LT_{h,\geq t}\subset\LT_h^{\ad}$.  Let $\LT_{h,t}=\LT_{h,\geq t}\backslash \LT_{h,\geq t+1}$, so that we have a decomposition of $\LT_h^{\ad}$ into $\G_h$-stable locally closed subsets:
\[ \LT_h^{\ad}=\coprod_{t=0}^h \LT_{h,t}. \]
Put another way, $\LT_{h,t}$ is the locus where $p=u_1=\cdots=u_{t-1}=0$, but where $u_t\neq 0$. 

In particular $\LT_{h,0}$ is defined by $p\neq 0$;  this is isomorphic to the $(h-1)$-dimensional unit ball over $\Spa W(k)[1/p]$. At the other extreme, $\LT_{h,h}=\Spa k$ is the unique non-analytic point of $\LT_h^{\ad}$.  For intermediate values $1\leq t\leq h-1$, the adic space $\LT_{h,t}$ is isomorphic to the $(h-t-1)$-dimensional ball over $\Spa k\laurentseries{u_t}$.  

We state below a moduli-theoretic interpretation of the $\LT_{h,t}$.  

\begin{proposition}  \label{PropAffField} Let $0\leq t\leq h$.  Let $(R,R^+)$ be a complete affinoid $W(k)$-algebra, and let $\Spa(R,R^+)\to \LT_h^{\ad}$ be a morphism, determining a formal group $\Gamma$ over $R^+$.  The following are equivalent:
\begin{enumerate} 
\item The morphism $\Spa(R,R^+)\to\LT_h^{\ad}$ factors through $\LT_{h,t}$.
\item The formal group $\Gamma\otimes_{R^+} R$ has height $t$.   
\end{enumerate}
\end{proposition}

\begin{proof} First we dispense with the $t=0$ case:  condition (1) is equivalent to $p$ being invertible in $R$, which is in turn equivalent to $\fg\otimes_{R^+} R$ having height 0. 

Now suppose $1\leq t\leq h$.  Let $w_1,\dots,w_{h-1}\in R^+$ be the images of $u_1,\dots,u_{h-1}$
under $A_h\to R^+$.  Condition (1) means that $p=w_1=\cdots=w_{t-1}=0$ and $w_t$ is a unit in $R$, so that in light of 
\[ [p]_{\fg}(X) =pX+_{\fg} w_1X^p+_{\fg} 
 \cdots+_{\fg} w_{h-1}X^{p^{h-1}}+_{\fg}  X^{p^h}, \]
the $p$-series of $\fg\otimes_{R^+} R$ has the form $f(X^{p^t})$, where $f(X)=(\text{unit})X+\cdots$. 
This means exactly that the formal group law $\fg$ has exact height $t$ over $R$, which is equivalent to (2).
\end{proof}

Similarly, one has a decomposition of the adic Rapoport--Zink space $\RZ_h^{\ad}$ into locally closed subpaces $\RZ_{h,t}$ which are stable for the action of $\G_h^0$.  

Let $\Spa(R,R^+)\to \RZ_{h,t}$ be a morphism, determining a $p$-divisible group $\Gamma$ over $R^+$. By \Cref{PropAffField}, the height of the identity component of $\fg_R=\fg\otimes_{R^+} R$ is constant on $\Spec R$, so that by \cite[Proposition II.4.9]{MessingCrystals} there is a connected-\'etale exact sequence of $p$-divisible groups over $R$:
\[ 0\to \fg_R^{\circ}\to \fg_R\to \fg_R^{\et}\to 0,\] 
where $\fg_R^{\et}$ is \'etale of height $h-t$, and $\fg_R^{\circ}$ has connected fibers of height $t$.  

One obtains a tower of finite \'etale covers of $\RZ_{h,t}$, by trivializing $\fg_R^{\circ}$ and $\fg_R^{\et}$.  First consider the $t=0$ case:  there $\fg_R=\fg_R^{\et}$, and one obtains a tower of finite \'etale covers $\RZ_{h,0}^H\to \RZ_{h,0}$ indexed by open subgroups $H\subset \GL_h(\Z_p)$;  over an $(R,R^+)$-point represented by a pair $(\fg,\alpha)$ as above, $\RZ_{h,0}^H$ parametrizes isomorphisms $\fg_R\isom (\Q_p/\Z_p)^{\oplus h}$ (equivalently, isomorphisms of the Tate module of $\fg_R$ with $\Z_p^{\oplus h}$), modulo the action of $H$.  In \cite{ScholzeWeinstein} it is shown that the limit of the tower is (in a certain sense known as a tilde limit) a perfectoid space $\RZ_{h,0}^\infty$;  in more recent language, $\RZ_{h,0}^\infty$ is the limit in v-sheaves of the diamonds $\RZ_{h,0}^{H,\Diamond}$.   

Now suppose $1\leq t\leq h-1$.  We impose the assumption that $k$ contains $\F_{p^{ht}}$;  let $\fg_t$ be the Honda formal group of height $t$ over $k$.  We obtain a tower of finite \'etale covers $\RZ_{h,t}^H\to \RZ_{h,t}$ indexed by open subgroups $H\subset \overline{\G}_t\times \GL_{h-t}(\Z_p)$, which parameterize pairs of isomorphisms $\fg_R^{\circ}\isom \fg_{t,R}$ and $\fg_R^{\et}\isom (\Q_p/\Z_p)^{\oplus (h-t)}$ modulo $H$.  Similarly there is a tower over $\LT_{h,t}$, which we call the {\em Lubin--Tate tower for the pair $(h,t)$}.  

To define $\RZ_{h,t}^\infty$ and $\LT_{h,t}^\infty$, we pass once again from rigid-analytic spaces to their associated diamonds, and take the limit in v-sheaves.  The cost of this passage is that, unlike in the characteristic 0 case, $X\mapsto X^\Diamond$ is far from being fully faithful.  As an example, if $X=\Spa(R,R^+)$ is a rigid space in characteristic $p$, then the ring of global sections of $X^\Diamond$ is not $R$ but rather the completion of the perfection of $R$.  Ultimately it is this issue which prevented us from proving the strong splitting conjecture in coheight 1 (\cref{ConjStrongSplitting}) for $h\geq 3$, and which required the special techniques of \S\ref{sec:height2} in the case $h=2$.

\subsection{From $p$-divisible groups to shtukas}
\label{ssec:dieudonne}

The next step is to consider the various linear-algebra structures one can associate to a $p$-divisible group $\fg$ over $R^+$, when $(R,R^+)$ is an affinoid perfectoid algebra in characteristic $p$.  Ultimately these will all come from Dieudonn\'e theory, which we now review.

Suppose $R$ is an f-semiperfect ring in characteristic $p$.  
 (See \cite[Definition/Proposition~4.1.2]{ScholzeWeinstein} for the definition of f-semiperfect ring;  the quotient of a perfect ring by a finitely generated ideal is an example.)  Then there is a universal $p$-adically complete PD thickening $\mathrm{A}_{\cris}(R)$ of $R$.  By functoriality, the Frobenius endomorphism of $R$ lifts to an endomorphism $\varphi$ of $\mathrm{A}_{\cris}(R)$.  If $R$ is perfect, then $\mathrm{A}_{\cris}(R)=W(R)$.

We adopt the following definition of Dieudonn\'e modules from the paragraph before \cite[Theorem~4.1.1]{ScholzeWeinstein}; in the case of perfect rings, the definition is the same as \cite[Definition~6.1]{LauSmoothness}.
\begin{defn}[Dieudonn\'e modules]\label{Def:Dieudonne-modules}
    Let $R$ be a $\mathrm{f}$-semiperfect ring of characteristic $p$. Let $\varphi$ be the Frobenius on $\mathrm{A}_{\cris}(R)$. 
    \textit{A Dieudonn\'e module} over $R$ is a finite projective $\mathrm{A}_{\cris}(R)$-module $M$ equipped with maps 
    \begin{equation*}
        F\colon M\otimes_{\mathrm{A}_{\cris}(R),\varphi} \mathrm{A}_{\cris}(R)\to M
    \end{equation*}
and
\[ V\colon M\to M\otimes_{\mathrm{A}_{\cris}(R),\varphi}\mathrm{A}_{\cris}(R) \]
    satisfying 
    \begin{equation*}
        FV=VF=p.
    \end{equation*}
\end{defn}

For an f-semiperfect ring $R$,  let $\fg\mapsto M(\fg)$ denote the (contravariant) Dieudonn\'e module functor, so that $M(\fg)$ is a Dieudonn\'e module over $R$.

\begin{thm} 
\label{ThmFullFaithfulness}
Let $R$ be an f-semiperfect ring.
\begin{enumerate}
    \item The functor $\Gamma\mapsto M(\Gamma)$ is fully faithful up to isogeny.
    \item If $R$ is perfect, then $\Gamma\mapsto M(\Gamma)$ is an anti-equivalence.  
\end{enumerate} 
\end{thm}

\begin{proof}
(1) is \cite[Theorem A]{ScholzeWeinstein}, and (2) is \cite[Theorem D]{LauSmoothness}.  
\end{proof}

Now let us specialize to the situation that $S=\Spa(R,R^+)$ is an affinoid perfectoid space in characteristic $p$, with $\varpi\in R$ a pseudo-uniformizer, so that both $R$ and $R^+$ are perfect rings, and $R^+/\varpi$ is f-semiperfect.  

Recall from \cite[\S12.2]{BerkeleyLectures} the analytic adic spaces
$Y_{S,I}$ derived from $S$, where $I\subset [0,\infty]$ is an interval.  In particular $Y_{S,[0,\infty]}$ is the analytic locus in $\Spa W(R^+)$, this being the complement in $\Spa W(R^+)$ of the vanishing locus $V(J)$ of the ideal of definition $J=(p,[\varpi])$ of the topological ring $W(R^+)$.  Furthermore
\begin{eqnarray*}
\Y_{S,[0,\infty)} &=& \Spa W(R^+) \setminus V([\varpi]) \\
\Y_{S,(0,\infty)} &=& \Y_{S,[0,\infty)} \setminus V(p).
\end{eqnarray*}

We write $\varphi$ for the Frobenius automorphism of $W(R^+)$, $\Y_{S,(0,\infty)}$, or $\Y_{S,[0,\infty)}$ arising from the Frobenius automorphism of $(R,R^+)$.  The categories of $\varphi$-equivariant vector bundles on these adic spaces play a prominent role in $p$-adic Hodge theory.  Importantly, $\varphi$ acts discontinuously on $\Y_{S,(0,\infty)}$, so that one may define the quotient 
\[ X_S\coloneqq \Y_{S,(0,\infty)}/\varphi^{\Z}. \] 
known as the relative Fargues--Fontaine curve. 
 
Some further notational conventions are in order.  Let $\Proj A$ denote the category of finite projective modules over a ring $A$ (disregarding any topology on $A$). Let $\VB(X)$ denote the category of vector bundles on an adic space $X$ (i.e., locally free $\OO_X$-modules).  For an endomorphism $\varphi$ of $A$, we let $\Proj A/\varphi$ denote the category of finite projective $A$-modules $M$ together with an isomorphism $\varphi^*M\stackrel{\sim}{\to} M$, and similarly for $\VB(X/\varphi)$.  

\begin{lemma}  \label{LemCrisBundles} Let $S=\Spa(R,R^+)$ be an affinoid perfectoid space in characteristic $p$.  
\begin{enumerate}
    \item  Recall that $\mathrm{B}_{\cris}^+(R^+/\varpi)=\mathrm{A}_{\cris}(R^+/\varpi)[1/p]$.  There is an equivalence of categories $\Proj B_{\cris}^+(R^+/\varpi)/\varphi \isom \VB(Y_{S,(0,\infty]}/\varphi)$.    
    \item If $S$ is a product of points \cite[Definition 2.14]{GleasonIvanovZillinger}, then the restriction map $\VB(Y_{S,(0,\infty]}/\varphi)\to \VB(Y_{S,(0,\infty)}/\varphi)$ is fully faithful.     
\end{enumerate}
\end{lemma}

\begin{proof} For the proof, we use the notation $B_{S,I}=H^0(Y_{S,I},\OO)$ for a subinterval $I\subset [0,\infty]$.  We start by noting the following inclusions of $\Q_p$-algebras:
\[ B_{S,[1/p,\infty]} \subset \mathrm{B}^+_{\cris}(R^+/\varpi) \subset B_{S,[1,\infty]}, \]
The first inclusion holds because $[\varpi^p]/p\in A_{\cris}(R^+/\varpi)$, and the second holds because for all $n\geq 1$, $\varpi^n/n!\in H^0(Y_{S,[1,\infty]},\OO^+)$, ultimately because $v_p(n!)\leq n$.  (In fact $v_p(n!)\leq n/(p-1)$, so the $[1,\infty]$ can be improved to $[1/(p-1),\infty]$, cf. \cite[p. 8]{FarguesFontaine2014} in the case that $S=\Spa(C,C^\circ)$ is a rank 1 geometric point.)  

Given an object $M$ of $\Proj {\mathrm{B}_{\cris}^+(R^+/\varpi)/\varphi}$, one can base change to $B_{S,[1,\infty]}$ to obtain a vector bundle $\E(M)$ on $Y_{S,[1,\infty]}$.  The pullback $\varphi^*\E(M)$ is a vector bundle on $Y_{S,[1/p,\infty]}$, and we obtain an isomorphism  $\varphi^*\E(M)\vert_{Y_{S,[1,\infty]}}\isom \E(M)$.  Using the iterated pullbacks $(\varphi^n)^*\E(M)$, we obtain a unique $\varphi$-equivariant extension of $\E(M)$ to $Y_{S,(0,\infty]}$.  Conversely, given a $\varphi$-equivariant vector bundle $\E$ on $Y_{S,(0,\infty]}$, one recovers $M$ by taking sections over $Y_{S,[1/p,\infty]}$ followed by base change along $B_{S,[1/p,\infty]} \subset \mathrm{B}^+_{\cris}(R^+/\varpi)$.  This is (1).

Part (2) is \cite[Corollary 7.3(1)]{GleasonIvanovZillinger}.  (The case when $S$ is a rank 1 geometric point is \cite[Proposition 7.20]{FarguesFontaine2014}.)
\end{proof}

In particular if $S$ lies over $k$, the Dieudonn\'e module of the constant $p$-divisible group $\Gamma_h\otimes_k R^+/\varpi$ determines a rank $h$ vector bundle $\OO_{X_S}(1/h)$ on $X_S$.   
 
We describe now a construction which assigns to a $p$-divisible group over $R^+$ its ``shtuka with a paw at $p$''.  The corresponding construction in characteristic 0 and its relation to $p$-divisible groups is detailed in \cite[\S12]{BerkeleyLectures}.  

\begin{defn} 
\label{defn:shtukapawatp} 
Let $S=\Spa(R,R^+)$ be an affinoid perfectoid space in characteristic $p$.  A {\em shtuka with a paw at $p$} over $S$ is a pair $(\E,\varphi_\E)$, where $\E$ is a vector bundle over $\Y_{S,[0,\infty)}$, and $\varphi_{\E}\from \varphi^*\E\vert_{\Y_{S,(0,\infty)}}\to \E\vert_{\Y_{S,(0,\infty)}}$ is an isomorphism which is meromorphic along $p=0$.

Given such a $(\E,\varphi_{\E})$, we define some associated objects:
\begin{itemize}
    \item The completion $\E^\wedge_p$ along $p=0$ is a projective $W(R)$-module which comes equipped with an isomorphism $\varphi_{\E^\wedge_p}\from \varphi^*\E^\wedge_p[1/p]\to \E^\wedge_p[1/p]$.  The pair consisting of $\E^\wedge_p[1/p]$ together with its Frobenius is an isocrystal over $R$, meaning an object of $\Proj W(R)[1/p]/\varphi$.  We call this the generic fiber of $(\E,\varphi_{\E})$ and denote it by $\eta(\E,\varphi_{\E})$. 
    \item The special fiber $\sigma(\E,\varphi_{\E})$, defined as the descent of $\E\vert_{\Y_{S,(0,\infty)}}$ via $\varphi$.  This is a vector bundle over $X_S$.
\end{itemize}
\end{defn}

\begin{rmk}
\label{rmk:genericspecial}
The names ``generic fiber'' and ``special fiber'' are a bit misleading in \cref{defn:shtukapawatp}, as $\gamma(\E,\varphi_{\E})$ is the localization at (a punctured formal neighborhood of) the locus $p=0$ while $\sigma(\E,\varphi_{\E})$ is the localization away from $p=0$. Nonetheless, as the theorem below shows, there is a functor from $p$-divisible groups $\Gamma$ over $R^+$ towards shtukas with a paw at $p$, in which the generic and special fibers of the shtuka are determined by $\Gamma\otimes_{R^+} R$ and $\Gamma\otimes_{R^+} R^+/\varpi$, respectively.
\end{rmk}

\begin{theorem}
\label{thm:RZtoShtuka}
     Let $S=\Spa(R,R^+)$ be an affinoid perfectoid space over $k$.  There is a natural transformation $(1)\to (2)$ between the following two presheaves of sets on the category of affinoid perfectoid spaces over $k$:
\begin{enumerate}
    \item The functor which assigns to $S$ the set $\RZ_h(R^+)$ classifying pairs $(\Gamma,\alpha)$, where $\Gamma$ is a $p$-divisible group over $R^+$ and $\alpha\from \Gamma_h\otimes_{k} R^+/\varpi \to \Gamma \otimes_{R^+} R^+/\varpi$ is a quasi-isogeny, where $\varpi$ is an unspecified pseudo-uniformizer.
    \item The functor which assigns to $S$ the set of isomorphism classes of quintuples $(\E,\varphi_{\E},N,\beta_1,\beta_2)$, where:
    \begin{enumerate}
        \item $(\E,\varphi_{\E})$ is a shtuka with a paw at $p$ over $S$,
        \item $N$ is a Dieudonn\'e module over $R$,
        \item $\beta_1\from N\to \E^\wedge_p$ 
        is an isomorphism of $W(R)$-modules satisfying $\beta_1\circ F=\varphi_{\E^\wedge_p}\circ \beta_1$, and 
        \item $\beta_2\from \OO_{X_S}(1/h)\to \sigma(\E,\varphi_{\E})$ is an isomorphism of vector bundles over $X_S$.
    \end{enumerate}
\end{enumerate}
If $S$ is a product of points, the map $(1)\to (2)$ on $S$ is a bijection.
\end{theorem}

\begin{proof}
We first describe the natural transformation $(1)\to (2)$.  Let $M=M(\Gamma)$ be the Dieudonn\'e module, a finite projective $W(R^+)$-module.  Associated to $M$ is a vector bundle on $\Spa W(R^+)$; let $\E$ be the restriction of this vector bundle to $\Y_{S,[0,\infty)}$.  The Frobenius operator $F$ on $M$ induces a map $\varphi_{\E}\from \varphi^*\E\to \E$;  since $F$ divides $p$ this is an isomorphism on $\E\vert_{Y_{S,(0,\infty)}}$ which is meromorphic along $p=0$, and so $(\E,\varphi_{\E})$ is a shtuka with a paw at $p$ over $S$.  Let $N=M(\Gamma\otimes_{R^+} R)$;  then by functoriality of the Dieudonn\'e module, we have an isomorphism $\beta_1\from N\to \E^\wedge_p=M\otimes_{W(R^+)}W(R)$ commuting with Frobenius.  Finally, the quasi-isogeny $\alpha$ induces an isomorphism $M(\Gamma_h)[1/p] \to M(\Gamma\otimes_{R^+} R^+/\varpi)[1/p]$ in $\Proj \mathrm{B}_{\cris}^+(R^+/\varpi)/\varphi$;  from \cref{LemCrisBundles} we have a functor $\Proj \mathrm{B}_{\cris}^+(R^+/\varpi)/\varphi\isom \VB(\Y_{S,(0,\infty]}/\varphi)\to \VB(X_S)$, from which we derive the isomorphism $\beta_2\from \OO_{X_S}(1/h)\to\sigma(\E,\varphi)$.
 
Now suppose $S=\Spa(R,R^+)$ is a product of points, and that data $(\E,\varphi_\E,N,\beta_1,\beta_2)$ constituting an object of (2) are given. Restrict this data to $S'=\Spa(R,R^\circ)$.   
By \cite[Corollary 1.9]{GleasonIvanovZillinger}, see also \cite[Theorem 8.14]{MingjiaZhang}, the shtuka $(\E,\varphi_{\E})$ arises from a Breuil-Kisin-Fargues module:  that is, there is a finite projective $W(R^\circ)$-module $M$ together with an isomorphism $\varphi_M\from \varphi^*M[1/p]\stackrel{\sim}{\to} M[1/p]$ whose restriction to $\Y_{S',[0,\infty)}$ is $(\E,\varphi_{\E})$.

The idea now is to show that $\varphi_M$ gives $M$ the structure of a Dieudonn\'e module over $R^\circ$, which is to say that both $\varphi_M$ and $p\varphi_M^{-1}$ preserve $M$.  Since the completion of $W(R^\circ)[1/[\varpi]]$ at $p$ is $W(R)$, 
Beauville--Laszlo gluing (see e.g., \cite[Lemma~5.2.9]{BerkeleyLectures}) gives an equivalence of categories
\[ \Proj W(R^\circ)[1/[\varpi]] \isom \Proj W(R) \times_{\Proj W(R)[1/p]} \Proj W(R^\circ)[1/p[\varpi]]. \] 
The isomorphism $\beta_1$ identifies $M\otimes_{W(R^\circ)} W(R)$ with $N$;  since $F$ and $V$ preserve $N$, the operators $\varphi_M$ and $p\varphi_M^{-1}$ preserve $M[1/[\varpi]]$.  

Therefore the operators $\varphi_M$ and $p\varphi_M^{-1}$, a priori defined on $M[1/p[\varpi]]$, preserve both $M[1/[\varpi]]$ and $M[1/p]$; we claim they must also preserve $M$.  It suffices to check this claim on an open covering of $\Spec W(R^\circ)$ by affines $\Spec T$ for which $M\otimes_{W(R^\circ)} T$ is free. With respect to a basis, our operators have 
matrix entries in the fiber product of $T[1/[\varpi]]$ and $T[1/p]$ over $T[1/p[\varpi]]$.  The key observation now is that $[\varpi]\in W(R^\circ)$ is a non-zero divisor modulo $p$ in $W(R^\circ)$;  the same is true for $T$, and consequently this fiber product is just $T$, proving the claim.

Therefore $M$ is a Dieudonn\'e module;  by \cref{ThmFullFaithfulness}, $M=M(\Gamma)$ for some $p$-divisible group $\Gamma$ over $R^\circ$.

Finally, passing $\beta_2$ through the equivalence of categories in \cref{LemCrisBundles} produces an isomorphism $M\otimes_{W(R^\circ)} \mathrm{B}_{\cris}^+(R^\circ/\varpi) \stackrel{\sim}{\to} M(\Gamma_h\otimes_k R^\circ/\varpi)[1/p]$ compatible with the Frobenius actions; by \cref{ThmFullFaithfulness} this arises from a quasi-isogeny 
$\alpha\from \Gamma_h\otimes_k R^\circ/\varpi \to \Gamma \otimes_{R^\circ} R^\circ/\varpi$.  The pair $(\Gamma,\alpha)$ therefore determines a point of $\RZ_h(R^\circ)$.  As noted above, the map $\LT_h(R^\circ)\to \LT_h(R^+)$ is a bijection; since $\RZ_h$ is a disjoint union of copies of $\LT_h$, the same is true for $\RZ_h$.  Therefore we have constructed the required $S$-point of (1).  
\end{proof}

\begin{rmk} All the shtukas we need in this article have one paw (or leg) which lies in characteristic $p$.  See \cite[Proposition 23.3.1]{BerkeleyLectures} for a definition of shtukas with a paw in characteristic 0 and their relation to the space $\RZ_{h,0}$ \cite[Theorem 24.2.5]{BerkeleyLectures}. 
\end{rmk}

\subsection{The v-sheaf of shtukas, and the v-stacks $\Isoc$ and $\Bun$}

Here we define v-sheaves in parallel to the Rapoport--Zink spaces $\RZ_{h,\geq 1}$, $\RZ_{h,t}$, and $\RZ_{h,t}^\infty$.  

\begin{defn}
\label{DefnShtukas} Let $h\geq 1$, and let $k$ be a field which is algebraic over $\mathbb{F}_{p^h}$.  
\begin{enumerate}
    \item 
Let $\Sht_{h,\geq 1}$ denote the v-sheafification of the presheaf on $\Perf_k$ which assigns to an affinoid perfectoid space $S=\Spa(R,R^+)$ the set of isomorphism classes of data $(\E,\varphi,N,\beta_1,\beta_2)$ as in \cref{thm:RZtoShtuka}.
\item Let $1\leq t\leq h$, and write $\Gamma_{h,t}=\Gamma_t\oplus (\Q_p/\Z_p)^{h-t}$.  Let $\Sht_{h,t}\subset \Sht_{h,\geq 1}$ denote the subsheaf imposing the condition that the Dieudonn\'e module $N$ over $W(R)$ is isomorphic to $M(\Gamma_{h,t})$ at every geometric point of $\Spec R$.
\item Keep $t$ as in (2), and assume that $k$ contains $\mathbb{F}_{p^{ht}}$.  Let $\Sht_{h,t}^\infty$ be the v-sheafification of the presheaf on $\Perf_k$ which assigns to $S=\Spa(R,R^+)$ the set of isomorphism classes of data $(\E,\varphi,\beta_1,\beta_2)$, where $\E,\varphi,\beta_2$ are as in \cref{thm:RZtoShtuka}, but now
\[ \beta_1\from M(\Gamma_{h,t}\otimes_k R)\isom \E_p^\wedge \] is an isomorphism of $W(R)$-modules satisfying $\beta_1\circ F=\varphi_{\E_p^\wedge}\circ \beta_1$.  
\end{enumerate}
\end{defn}

The group $\overline{\G}_h^0$ acts on all three spaces defined above, by operating on $\beta_2$.  The group $\overline{\G}_t\times \GL_{h-t}(\Z_p)$, being the group of automorphisms of $\Gamma_{h,t}$, acts on $\Sht_{h,t}^\infty$ by operating on $\beta_1$.  Thus $\Sht_{h,t}^\infty$ admits an action of $\overline{\G}_h^0\times\overline{\G}_t\times \GL_{h-t}(\Z_p)$.  This action even extends to a larger group $\G_h^0\times_{\Gal(k/\F_p)} \G_t\times\GL_{h-t}(\Z_p)$.

Recall the functor $X\mapsto X^{\diamondsuit}$ from adic spaces over $\Z_p$ to v-sheaves over $\Spd \Z_p$. If $X$ lives in characteristic $p$, then $X^\diamondsuit$ is the restriction 
of the sheaf represented by $X$ from characteristic $p$ adic spaces to characteristic $p$ perfectoid spaces.

\begin{proposition} \label{prop:RZtoSht} There is a $\G_h^0$-equivariant isomorphism $\RZ_{h,\geq 1}^{\Diamond}\to \Sht_{h,\geq 1}$ which carries $\RZ_{h,t}$ isomorphically onto $\Sht_{h,t}$ for each $1\leq t\leq h$.  
\end{proposition}

\begin{proof} The v-sheafification of the natural transformation described in \cref{thm:RZtoShtuka} gives a $\G_h^0$-equivariant morphism $\RZ_{h,\geq 1}^{\Diamond}\to \Sht_{h,\geq 1}$.  Since products of points form a basis for the v-topology (\cite[Remark~2.15]{GleasonIvanovZillinger}), part (2) of the theorem shows this is an isomorphism.  With respect to this isomorphism, $N$ is the Dieudonn\'e module of $\Gamma$, and so the respective conditions for membership in $\RZ_{h,t}$ and $\Sht_{h,t}$ agree.    
\end{proof}

\begin{lemma} \label{lem:caraianischolze}
The map $\Sht_{h,t}^\infty\to \Sht_{h,t}$ sending $(\E,\varphi_{\E},\beta_1,\beta_2)$ to $(\E,\varphi_{\E},M(\Gamma_{h,t}),\beta_1,\beta_2)$ is a pro-finite \'etale torsor for the group $\overline{\G}_t\times\GL_{h-t}(\Z_p)$.  Pulling back the map from $\LT_{h,t}^{\Diamond}$ (respectively $\RZ_{h,t}^{\Diamond}$) through this torsor produces $\LT_{h,t}^{\infty}$ (respectively, $\RZ_{h,t}^{\infty}$).  
\end{lemma}

\begin{proof} Given data $(\E,\varphi_{\E},N,\beta_1,\beta_2)$ defining a point of $\Sht_{h,t}$ over $S=\Spa(R,R^+)$, the fiber of $\Sht_{h,t}^\infty$ classifies all trivializations of $N$.  These form a pro-\'etale $\Aut N$-torsor by \cite[Proposition 2.2.6]{CaraianiScholze}.  When the data arose from a $p$-divisible group $\Gamma$ over $R^+$, such trivializations correspond exactly to trivializations of $\Gamma\otimes_{R^+} R$.
\end{proof}

As constructed in the previous subsection, a shtuka $(\E,\varphi_{\E})$ has a generic point $\gamma(\E,\varphi_\E)$ and a special point $\sigma(\E,\varphi_\E)$.  We review here the v-stacks $\Isoc$ and $\Bun$, which are the targets of these constructions.

\begin{defn}[The stacks $\Isoc$ and $\Bun$]
\label{def:IsocBun}
Let $S = \Spa(R, R^+)$ be an affinoid perfectoid space over $k$.
\begin{enumerate}
    \item An \emph{analytic isocrystal} over $S$ is a pair $(M, \varphi_M)$ consisting of a finite projective $W(R)[\frac1p]$-module $M$ and an isomorphism $\varphi_M\colon \varphi^*M \isoto M$. We denote by $\Isoc'(S)$ the groupoid of analytic isocrystals over $S$ and by $\Isoc$ the v-sheafification of the functor $S\mapsto \Isoc'(S)$.
\item Define a presheaf of groupoids $\Bun$ on $\Perf_k$, which sends $S$ to the groupoid of vector bundles on $X_S$.  This is in fact already a v-stack \cite[I.4.1]{FarguesScholze}.
\end{enumerate}
\end{defn}

The underlying topological spaces of $\Isoc$ and $\Bun$ are understood through Newton polygons.

\begin{defn}  We denote by $\mc N$ the set of {\em Newton polygons}, which we index by non-increasing tuples $b=(\lambda_1,\dots,\lambda_n)$ of rational numbers $\lambda_i=d_i/h_i$, with $h_i\geq 1$ and $\gcd(d_i,h_i)=1$.  To such a $b$ there is an associated convex polygon in the plane, obtained by starting at $(0,0)$ and then drawing a line segment of slope $\lambda_i$ and horizontal extent $h_i$ for $i=1,\dots,n$.  The $\lambda_i$ are the {\em slopes} of $b$, if they are all equal we call $b$ {\em isoclinic}.  We equip $\mc N$ with a partial order:  $b\leq b'$ if the polygons have the same endpoint, and if the $b$-polygon is always on or below the $b'$-polygon.  Thus isoclinic polygons are minimal in $\mc N$.
\end{defn}

\begin{thm} \label{rslt:topology-on-Bun-and-Isoc}
There are canonical bijections
\begin{align*}
   \abs{\Isoc} \isom \abs{\Bun} \isom \mc N
\end{align*}
and we have the following explicit description of the topologies:
\begin{thmenum}
    \item Given $b \in \abs{\Isoc}$, the closure of $\{ b \}$ in $\abs{\Isoc}$ is the set of all $b' \in \mc N$ with $b' \le b$.
    \item Given $b \in \abs{\Bun}$, the closure of $\{ b \}$ in $\abs{\Bun}$ is the set of all $b' \in \mc N$ with $b' \ge b$.
\end{thmenum}
\end{thm}
\begin{proof}
The claim about $\abs{\Isoc}$ is \cite[Propositions 6.9 and 6.11]{GleasonIvanovZillinger}.  The claim about the underlying set $\abs{\Bun}$ is \cite[I.4.1]{FarguesScholze}, and the claim about its topology is \cite{hansen2017degenerating}.
\end{proof}

With the above description of $\abs{\Bun}$ and $\abs{\Isoc}$ it is now natural to consider the following locally closed substacks:

\begin{defn}
For a Newton polygon $b \in \mc N$ we denote by
\begin{align*}
    \Bun^b \subseteq \Bun, \qquad \Isoc^b \subseteq \Isoc
\end{align*}
the full substacks given by the locally closed point $b \in \abs{\Bun} = \abs{\Isoc}$. More explicitly, for a small v-stack $S$, $\Bun^b(S)$ resp. $\Isoc^b(S)$ is the stack of those $\mc E \in \Bun(S)$ resp. $M \in \Isoc(S)$ whose Newton polygon $\abs S \to \mc N$ is the constant function $b$. 
\end{defn}

\begin{proposition} \label{rslt:strata-of-Bun-and-Isoc}
\begin{thmenum}
\item For $b \in \mc N$ we have isomorphisms
\begin{eqnarray*}
    \Bun^b &\isom& */\widetilde G_b,\\
    \Isoc^b &\isom&  */G_b,
\end{eqnarray*}
where $G_b$ is the (locally profinite) group of isomorphisms of the isocrystal associated to $b$ and $\widetilde G_b$ is a semidirect product of $G_b$ with a certain unipotent group diamond $\widetilde G^{>0}_b$.
\item For $b$ isoclinic, we have $\widetilde G_b=G_b$, and so $\Bun^b\isom \Isoc^b$.
\item Consider the category of $\Q$-graded vector bundles in $\Bun$ with Newton polygon $b$ whose $\lambda$-graded piece is isoclinic of slope $\lambda$.  This category is equivalent to $\Isoc^b$. 
\end{thmenum}
\end{proposition}
\begin{proof}
Regarding (1):  For $\Bun$ see \cite[Proposition~III.5.3]{FarguesScholze}, for $\Isoc$ see \cite[Proposition~6.11]
{GleasonIvanovZillinger}.  (2) is \cite[Theorem III.4.5]{FarguesScholze}, and (3) is \cite[Proposition III.4.7]{FarguesScholze}.  
\end{proof}

\begin{example}
Suppose that $b \in \mc N$ is isoclinic. Then $\Bun^b\subset \Bun$ is an open immersion, while $\Isoc^b\subset \Isoc$ is a closed immersion. We have $\widetilde G_b = G_b$ by \cref{rslt:strata-of-Bun-and-Isoc}(ii), so in this case $\Bun^b \isom \Isoc^b \isom [*/G_b]$. 

An important special case is $b = (1/h)$.  The Newton polygon is isoclinic, meaning that $\Bun^b\subset \Bun$ is open.  We have $\Bun^b=[\ast/G_b]$, where $G_b=\overline{\G}_h^0=D_h^\times$ is the group of units in the division algebra over $\Q_p$ of invariant $1/h$.  For $S\in \Perf_k$, we write $\OO_{X_S}(1/h)$ for the vector bundle on $X_S$ corresponding to $S\to \ast\to [\ast/G_b]\isom \Bun^b$ (in agreement with the convention of the last subsection) and we write $\OO(1/h)$ to refer to the morphism $\ast\to \Bun^b$ itself.

Another special case is $b=(1/t,0,\dots,0)$ (with $(h-t)$ zeros).   We have $\Isoc^b=[\ast/G_b]$, where $G_b=\overline{\G}_t^0\times \GL_{h-t}(\Q_p)$.  For $S=\Spa(R,R^+)$ affinoid in $\Perf_k$, we write $\mathbf{D}_R(1/t,0,\dots,0)$ for the isocrystal over $R$ corresponding to $S\to\ast\to [\ast/G_b]\isom \Isoc^b$; note this is the same as $M(\Gamma_{h,t,R})[1/p]$, where once again $\fg_{h,t}=\fg_t\oplus (\Q_p/\Z_p)^{\oplus (h-t)}$. We write $\mathbf{D}(1/t,0,\dots,0)$ to refer to the morphism $\ast\to \Isoc^b$ taking this value at every $S$.
\end{example}

\subsection{Meromorphic vector bundles}

In \cite{GleasonIvanovZillinger} the authors construct a stack $\Bun^{\mer}$ of ``meromorphic vector bundles'' on the Fargues--Fontaine curve, which turn out to be closely related to our moduli spaces of shtukas.   After establishing this relationship, we can leverage the main result of \cite{GleasonIvanovZillinger} to produce a two-tower isomorphism for $\Sht_{h,t}$.

We summarize the construction of $\Bun^{\mer}$ from a slightly different viewpoint than in \cite{GleasonIvanovZillinger}.  Let $S=\Spa(R,R^+)$ be an affinoid perfectoid space in characteristic $p$ with pseudouniformizer $\varpi \in R^+$.  The adic space $\mc Y_{S,[0,\infty)}$ admits an affinoid covering by 
\[ \mc Y_{S,[0,p^n]} = \Spa B_{S,[0,p^n]},\]
as $n$ ranges through $\Z$.  
Explicitly, $B_{S,[0,p^n]}$ has the following presentation as a subring of $W(R)$:
\begin{equation}
    \label{EqPresentationOfB}
 B_{S,[0,p^n]} = W(R^+)\left[\frac{p}{\varpi^{1/p^n}}\right]^{\wedge}_{(p)}\left[\frac{1}{\varpi}\right] = 
\set{ \sum_{i\geq 0} [\rho_i]p^i\in W(R) \biggm\vert\rho_i\varpi^{i/p^n}\to 0} 
\end{equation}
\cite[proof of Proposition 11.2.1]{BerkeleyLectures}.  The automorphism $\varphi$ of $\mc Y_{S,[0,\infty)}$ restricts to an isomorphism $\mc Y_{S,[0,p^n]}\to \mc Y_{S,[0,p^{n+1}]}$.  

We now introduce the notion of a meromorphic vector bundle, which is meant to be an intermediate notion between vector bundles on $\mc Y_{S,[0,\infty)}$ and on $\mc Y_{S,(0,\infty)}$.

\begin{defn} \label{def:meromorphic-vector-bundle-over-aff-perfectoid}
Given an affinoid perfectoid space $S = \Spa(R, R^+)$ in characteristic $p$, let
$\VB(\mc Y_{S,\langle 0,\infty)})$ denote the colimit as $n\to -\infty$ of the categories 
\begin{align*}
    \Proj {B_{S,[0,p^n]}[\tfrac1p]} \times_{\VB(\mc Y_{S,(0,p^n]})}  \VB(\mc Y_{S,(0,\infty)})
\end{align*}
Then $\varphi$ is an automorphism of $\VB(\mc Y_{S,\langle 0,\infty)})$;  we let
\begin{align*}
    \VB(\mc Y_{S,\langle 0,\infty)}/\varphi^{\Z})
\end{align*}
denote the category of objects $\E$ of $\VB(\mc Y_{S,\langle 0,\infty)})$ equipped with an isomorphism $\varphi^* \mathcal E \isoto \mathcal E$.  We call this the category of {\em meromorphic vector bundles on the Fargues--Fontaine curve over $S$}, thinking of the formal symbol $\mc Y_{S,\langle 0, \infty)}$ as ``$\mc Y_{S,(0,\infty)}$ together with an infinitesimal punctured neighborhood at $0$.'' 
\end{defn}

\begin{rmk}
Using the theory of analytic stacks by Clausen--Scholze, one can make literal sense of the symbol $\mc Y_{S,\langle 0,\infty)}$.  Namely, given $S = \Spa(R, R^+)$ as above, we can define $\mc Y_{S,\langle 0, \infty)} \coloneqq  \mc Y_{S,[0,\infty)} \times_{\Spec \Z_p} \Spec \Q_p$. Then there are natural inclusions
\begin{align*}
    \mc Y_{S,(0,\infty)} \subseteq \mc Y_{S,\langle 0, \infty)} \subseteq \mc Y_{S,[0,\infty)}
\end{align*}
and a meromorphic vector bundle on the Fargues--Fontaine curve over $S$ is exactly a vector bundle on $\mc Y_{S,\langle 0,\infty)}/\varphi^{\Z}$.
\end{rmk}

Since meromorphic vector bundles extend to a punctured infinitesimal neighbourhood of $0$, one can use Beauville-Laszlo gluing in order to glue them to a vector bundle on all of $\mc Y_{[0,\infty)}$:

\begin{lemma} \label{rslt:BL-glueing-on-Y-at-0}
Let $S = \Spa(R, R^+)$ be an affinoid perfectoid space of characteristic $p$. 
\begin{enumerate}
\item     Beauville-Laszlo gluing at the point $p = 0$ induces a natural $\varphi$-equivariant equivalence of categories
\begin{align*}
    \VB(\mc Y_{S,[0,\infty)}) \isom \VB(\mc Y_{S,\langle 0,\infty)}) \times_{\Proj {(W(R)[\frac1p])}} \Proj {W(R)}.
\end{align*}
\item  The restriction map $\VB(\mathcal{Y}_{S,[0,\infty)})\to \VB(\mathcal{Y}_{S,\langle 0,\infty)})$ induces a $\varphi$-equivariant functor between $\Q_p$-linear categories
\[ \VB(\mathcal{Y}_{S,[0,\infty)})[1/p] \to  \VB(\mathcal{Y}_{S,\langle 0,\infty)}). \]
If $S$ is a product of points in the sense of \cite[Definition 2.14]{GleasonIvanovZillinger}, then this functor is an equivalence.
\end{enumerate}
\end{lemma}
\begin{proof}
For (1):  By the definition of $\VB(\mc Y_{S,\langle 0,\infty)})$ (and a similar formula for $\VB(\mc Y_{S,[0,\infty)})$) the claim immediately reduces to showing that
\begin{align*}
    \Proj {B_{S,[0,1]}} = \Proj (B_{S,[0,1]}[\tfrac1p]) \times_{\Proj (W(R)[\frac1p])} \Proj {W(R)}.
\end{align*}
We observe that $W(R)$ is the $p$-adic completion of $B_{S,[0,1]}$. (Explanation:  By the presentation given in \eqref{EqPresentationOfB}, the subring $B_{S,[0,1]}$ of $W(R)$ has $B_{S,[0,1]}/p\isom R^+[p/\varpi][1/\varpi]=R$;  therefore its $p$-completion must be $W(R)$.) Thus the claim follows immediately from Beauville-Laszlo gluing.

For (2): Under the hypothesis that $S=\Spa(R,R^+)$ is a product of points, we claim that the functor 
\begin{equation}
    \label{eq:projW}
(\Proj W(R))[\tfrac1p]\to \Proj (W(R)[\tfrac1p]) 
\end{equation}
is an equivalence.  That it is fully faithful is true without any hypotheses on $R$:  this is the statement that $\Hom(M,N)[1/p]\to \Hom(M[1/p],N[1/p])$ is an isomorphism for finite projective $W(R)$-modules $M$ and $N$, and it reduces immediately to the case that $M\isom W(R)^{\oplus n}$ is free, in which case both sides identify with $N[1/p]^{\oplus n}$.  For essential surjectivity:  Suppose $M$ is a finite projective $W(R)[1/p]$-module.  The result \cite[Proposition 2.19]{GleasonIvanovZillinger} shows that the perfect scheme $\Spec R$ is a {\em comb}, meaning that every connected component of $\Spec R$ is the spectrum of a valuation ring with algebraically closed fraction field.  Then by \cite[Theorem 6.1]{Ivanov-arc-descent}, there are finitely many rings $R_i$, and for each $i$ a finite free module $M_i$ over $W(R_i)[1/p]$, such that $R\isom \prod_i R_i$ and $M\isom \prod_i M_i$.  Since each $M_i$ obviously descends to a finite projective $W(R_i)$-module, we find that $M$ descends to a finite projective $W(R)$-module. 

Inverting $p$ in the equivalence in (1), we then have equivalences:
\begin{eqnarray*}
    \VB(\mathcal{Y}_{S,[0,\infty)})[\tfrac1p]
    &\isom& 
    \VB(\mathcal{Y}_{S,\langle 0,\infty)})
    \times_{\Proj (W(R)[\tfrac1p])} (\Proj W(R))[\tfrac1p] \\
    &\isom& 
    \VB(\mathcal{Y}_{S,\langle 0,\infty)}).
\end{eqnarray*}
\end{proof}

We do not know if the assignment $S \mapsto \VB(\mc Y_{S,\langle 0,\infty)})$ satisfies v-descent, so we need to sheafify it to make the following definition.

\begin{defn}
\label{DefnBunMer}
We denote by $\Bun^{\mer}$ the v-sheafification of the prestack $S \mapsto \VB(\mc Y_{\langle 0,\infty)}/\varphi^\Z)$, and call this the stack of {\em meromorphic vector bundles}.
\end{defn}

We now look at two localizations on $\mc Y_{S,\langle 0,\infty)}$. Firstly we can localize to $\mc Y_{S,(0,\infty)}$, i.e., throw away any infinitesimal information around $0$. Secondly we can localize to the punctured infinitesimal neighbourhood of $0$, which is given by the ring $W(R)[\frac1p]$ (see \cref{rslt:BL-glueing-on-Y-at-0}).  These localizations will result in functors out of $\Bun^{\mer}$ and into $\Bun$ and $\Isoc$, respectively.

\begin{defn}\label{DefnBunMerCorresp}
We denote by
\begin{align*}
    \Isoc \xfrom{\gamma} \Bun^{\mer} \xto{\sigma} \Bun
\end{align*}
the maps of v-stacks induced by pullbacks along the two localizations on $\mc Y_{S,\langle0,\infty)}$ discussed above. More precisely, given an affinoid perfectoid $S = \Spa(R, R^+)$ and an object $\mathcal E \in \VB(\mc Y_{S,\langle 0,\infty)})$ then by definition $\mathcal E$ is formally glued from a finite projective $B_{S,[0,p^n]}[\frac1p]$-module $\mathcal E_0$ for some $n\in\Z$ and a vector bundle $\mathcal E_1$ on $\mc Y_{S,(0,\infty)}$. We let
\begin{align*}
    \gamma(\mc E) = \mathcal E_0 \otimes_{B_{S,[0,p^n]}[\frac1p]} W(R)[\tfrac1p], \qquad \sigma(\mc E) \coloneqq  \mc E_1,
\end{align*}
together with their respective $\varphi$-equivariant structures. After passing to v-sheafification we obtain the desired morphisms of small v-stacks. We call $\gamma$ the \emph{generic polygon map} and $\sigma$ the \emph{special polygon map}.
\end{defn}

Let us now check that the above definitions are compatible with the ones in \cite{GleasonIvanovZillinger}, so that we can apply their results.

\begin{lemma}
Our above constructions produce the same results as the ones in \cite{GleasonIvanovZillinger}. More precisely:
\begin{lemenum}
    \item The v-stack $\Bun^{\mer}$ is isomorphic to the one from \cite[Definition~4.19]{GleasonIvanovZillinger}.

    \item The v-stack $\Isoc$ is isomorphic to  the v-stack $\mathfrak B^\Diamond$ from \cite[Definition~4.12]{GleasonIvanovZillinger}.

    \item With respect to the isomorphisms above, the special and generic polygon maps $\gamma$ and $\sigma$ coincide with the ones from \cite[Definition~4.20]{GleasonIvanovZillinger}.
\end{lemenum}
\end{lemma}
\begin{proof}
For (i): Let us temporarily denote by $\Bun^{\mer'}$ the small v-stack from \cite{GleasonIvanovZillinger}.  By \cite[Definition 4.18 and Proposition 4.21(2)]{GleasonIvanovZillinger}, $\Bun^{\mer'}$ is the v-sheafification of the functor sending $S$ to the category of objects $\E$ in $\VB(\mathcal{Y}_{[0,\infty)})[1/p]$ equipped with an isomorphism $\varphi^*\E\stackrel{\sim}{\to}\E$.  The $\varphi$-equivariant functor $\VB(\mathcal{Y}_{[0,\infty)})[1/p]\to \VB(\mathcal{Y}_{\langle 0,\infty)})$ from \cref{rslt:BL-glueing-on-Y-at-0}(2) therefore induces a functor $\Bun^{\mer}\to\Bun^{\mer'}$ which is an equivalence on products of points.  These form a basis for the v-topology (see \cite[Remark~2.15]{GleasonIvanovZillinger}), so the functor is an equivalence.  

For (ii):  By \cite[Definition 4.4 and Proposition 4.14(2)]{GleasonIvanovZillinger}, $\mathfrak{B}^\Diamond$ is the v-sheafification of the functor sending $S=\Spa(R,R^+)$ to the category of objects $M$ in $(\Proj W(R))[\tfrac1p]$ equipped with an isomorphism $\varphi^*M\stackrel{\sim}{\to} M$.  In light of the $\varphi$-equivariant functor \eqref{eq:projW}, there is a functor $\mathfrak{B}^{\Diamond}\to \Isoc$.  By the proof of \cref{rslt:BL-glueing-on-Y-at-0}(2), the functor \eqref{eq:projW} is an equivalence when $S$ is a product of points, which implies that $\mathfrak{B}^{\Diamond}\to \Isoc$ is an equivalence.  

Part (iii) follows easily by going through the above identifications and the definition of $\gamma$ and $\sigma$ on both sides.
\end{proof}

We can now state a relationship between our space of shtukas and $\Bun^{\mer}$.  

\begin{thm} \label{ThmShtVSMer} 
There is a cartesian diagram of small v-stacks
\begin{equation*}\begin{tikzcd}
    \Sht_{h,t}^{\infty} \arrow[r] \arrow[d] & \Bun^{\mer} \arrow[d,"{(\gamma,\sigma)}"]\\
    * \arrow[r] & \Isoc \times \Bun
\end{tikzcd}\end{equation*}
where the bottom horizontal map is the product of $\mathbf{D}(1/t,0,\dots,0)$ and $\OO(1/h)$.  
\end{thm}
\begin{proof}
The upper horizontal arrow in the diagram is restriction of shtukas from $\Y_{S,[0,\infty)}$ to $\Y_{S,\langle 0,\infty)}$; the definition of $\Sht_{h,t}^\infty$ makes it fairly clear that the diagram commutes.

To see that the diagram is cartesian:  Let $S=\Spa(R,R^+)$ be an affinoid perfectoid space over $k$. 
Then up to sheafification, the $S$-points of the fiber product parameterize data $(\E^{\mer},\varphi_{\E^{\mer}},\alpha_1,\alpha_2)$, where $\E^{\mer}$ is an object of $\VB(\Y_{S,\langle 0,\infty)})$, $\varphi_{\E^{\mer}}\from \varphi^*\E^{\mer}\to\E^{\mer}$ is an isomorphism, $\alpha_1\from M(\Gamma_{h,t})\otimes_{W(k)} W(R)[1/p] \to \gamma(\E^{\mer},\varphi_{\E^{\mer}})$ 
is an isomorphism of isocrystals over $R$, and $\alpha_2\from \OO_{X_S}(1/h)\to \sigma(\E^{\mer},\varphi_{\E^{\mer}})$ 
is an isomorphism of vector bundles over $X_S$.  By \cref{rslt:BL-glueing-on-Y-at-0}(1), the meromorphic vector bundle $\E^{\mer}$ and the Dieudonn\'e module $M(\Gamma_{h,t})$ glue together to give an object $\E\in \VB(\Y_{S,[0,\infty)})$. 
The isomorphism $\varphi_{\E^{\mer}}$ (respectively, $p\varphi_{\E^{\mer}}^{-1}$)
glues together with Frobenius (respectively, Verschiebung)
to give an operator $\varphi_{\E}\from \varphi^*\E\to \E$ (respectively, an operator $V\from \E\to \varphi^*\E$) for which $\varphi_{\E}V=V\varphi_{\E}=p$; in particular $\varphi_{\E}$ is meromorphic at $p=0$.  Thus $(\E,\varphi_{\E})$ defines a shtuka with a paw at $p$, whose $p$-completion and special fiber are trivialized by $\alpha_1$ and $\alpha_2$, respectively.  This defines an $S$-point of $\Sht_{h,t}^\infty$. 
\end{proof}

\cref{ThmShtVSMer} shows that the action of $\G_h^0\times_{\Gal(k/\F_p)} \G_t \times\GL_{h-t}(\Z_p)$ on $\Sht_{h,t}^\infty$ extends to the larger group $\G_h^0\times_{\Gal(k/\F_p)} \G_t^0 \times\GL_{h-t}(\Q_p)$, since the lower horizontal map in the diagram is a torsor onto its image relative to this group.

\subsection{The two-tower isomorphism}

In \cite[Definition~V.3.2]{FarguesScholze} Fargues--Scholze introduce the v-stack $\mc M$ of vector bundles $\mc E$ on the Fargues--Fontaine curve together with an increasing filtration $(\mc E_\lambda)_{\lambda \in \Q}$ such that for each $\lambda\in\Q$ the quotient $\mc E_\lambda/\mc E_{<\lambda}$ has a single slope $\lambda$.  The morphism $\mc M\to \Bun$ forgetting the filtration serves as a smooth uniformization of $\Bun$ by a v-stack representable in diamonds, showing that $\Bun_G$ is an Artin v-stack.  The v-stack $\mc M$ decomposes as a disjoint union of v-stacks $\mc M_b$ over Newton polygons $b\in\mathcal{N}$, where $\mc M_b$ classifies those $\Q$-filtered vector bundles whose associated graded has Newton polygon $b$.  

There is an analogue of $\mc M$ for meromorphic vector bundles:

\begin{defn}
We denote by $\mc M^{\mer}$ the stack of \emph{meromorphic vector bundles with slope filtration}: It is the stack of $\Q$-filtered objects $(\mc E_\lambda)_{\lambda\in\Q}$ such that the quotients $\mc E_\lambda/\mc E_{<\lambda}$ have a generic Newton polygon with single slope $\lambda$ (cf. \cite[Definition~5.15]{GleasonIvanovZillinger}).  
\end{defn}

In the following we summarize some of the main results in \cite{GleasonIvanovZillinger}, relating the stacks $\mc M^{\mer}$ and $\mc M$.
\begin{thm}\label{ThmBunMerMain}
\begin{thmenum}
    \item The special polygon map $\sigma\colon \Bun^{\mer} \to \Bun$ induces an isomorphism
    \begin{align*}
        \mc M^{\mer} \isoto \mc M.
    \end{align*}
    In other words, every vector bundle on the Fargues--Fontaine curve equipped with a slope filtration extends uniquely to a meromorphic vector bundle equipped with a slope filtration. 
    \item Let $b \in \mc N$ be a Newton polygon. There is a cartesian square:
    \begin{equation*}\begin{tikzcd}
        \mc M_b \arrow[r,hook] \arrow[d] & \Bun^{\mer} \arrow[d,"\gamma"] \\
        \Isoc^b \arrow[r,hook] & \Isoc
    \end{tikzcd}\end{equation*}
    Here, the top horizontal arrow composes the restriction of the inverse of the isomorphism in (1) to $\mc M_b$ with the map $\mc M^{\mer}\to \Bun^{\mer}$ forgetting the filtration.  
    The left vertical arrow inputs a $\Q$-filtered vector bundle, passes to its associated $\Q$-graded vector bundle, and passes it through the equivalence to $\Isoc^b$ in \cref{rslt:strata-of-Bun-and-Isoc}. 
\end{thmenum}
\end{thm}
\begin{proof}
Part (i) is proved in \cite[Proposition~5.16]{GleasonIvanovZillinger} and part (ii) is proved in \cite[Theorem~6.13]{GleasonIvanovZillinger}.
\end{proof}

We derive from this our key result.

\begin{thm} \label{ThmShtInf} Let $1\leq t\leq h$.
\begin{enumerate}
    \item As a v-sheaf with $\overline{\G}_h^0\times \overline{\G}_t^0 \times \GL_{h-t}(\Q_p)$-action, the infinite-level shtuka space $\Sht_{h,t}^\infty$ is isomorphic to the v-sheaf over $k$ whose $S$-points classify exact sequences of vector bundles on $X_S$:
\[ 0 \to \OO^{\oplus(h-t)}\to \OO(1/h)\to \OO(1/t) \to 0.\] 
\item The infinite-level Lubin--Tate space $\LT_{h,t}^\infty$ classifies the same exact sequences, such that the induced isomorphism on determinants $\wedge^h \OO(1/h)\isom \wedge^t \OO(1/t)$ is a $\Z_p^\times$-multiple of a fixed one.  
\end{enumerate}
\end{thm}

\begin{proof}  For (1):  Let $b=(1/t,0,\dots,0)$.  Combining \cref{ThmShtVSMer} and \cref{ThmBunMerMain} gives a cartesian diagram
\[ 
\xymatrix{
\Sht_{h,t}^\infty \ar[d] \ar[r] & \mathcal{M}_b \ar[d] \\
\ast \ar[r] & \Isoc_b \times \Bun 
}
\]
That is: $\Sht_{h,t}^\infty$ classifies objects $\E$ of $\Bun$ equipped with (1) an isomorphism $\E\isom \OO(1/h)$ together with (2) an increasing $\Q$-filtration of $\E$, whose graded pieces are identified with $\OO^{\oplus h-t}$ and $\OO(1/t)$.  This is exactly the data given by an exact sequence as in the theorem.

For (2):  Given the description of $\Sht_{h,t}^\infty$ in (1), there is a map from $\Sht_{h,t}^\infty$ 
towards the space of isomorphisms $\wedge^h \OO(1/h)\isom \wedge^t \OO(1/t)$, 
which is a $\Q_p^\times$-torsor;  
this map is compatible with the determinant homomorphism $\delta\from \overline{\G}_h^0\times \overline{\G}_t^0\times\GL_{h-t}(\Q_p)\to \Q_p^\times$.  
Quotienting by $\overline{\G}_t\times\GL_{h-t}(\Z_p)$, we arrive at a map from $\Sht_{h,t}$ towards the $\Z$-torsor $T$ of isomorphisms $\wedge^h \OO(1/h)\isom \wedge^t \OO(1/t)$ considered up to $\Z_p^\times$, 
which is compatible with the height map $\overline{\G}_h^0\to \overline{\G}_h^0/\overline{\G}_h\isom \Z$.  
Since $\Sht_{h,t}\isom (\LT_{h,t}^{\Diamond}\times \overline{\G}_h^0)/\overline{\G}_h$, 
it must be that $\LT_{h,t}^\Diamond$ is one of the fibers of $\Sht_{h,t}\to T$.  Therefore its preimage $\LT_{h,t}^\infty$ in $\Sht_{h,t}^\infty$ is as described.
\end{proof}

In our version of the two-towers isomorphism, the analogue of Drinfeld symmetric space is obtained by quotienting $\Sht_{h,t}^\infty$ by $\overline{\G}_h^0$ to obtain a ``characteristic $p$ period domain''.  First we introduce some notations.

\begin{defn} [Banach--Colmez spaces] \label{def:BCspaces} Let $\lambda\in\Q$.  Assume that $k$ contains $\F_{p^d}$, where $d$ is the denominator of $\lambda$.  Define a sheaf of $\Q_p$-vector spaces $\BC(\lambda)$ on $\Perf_k$ as follows.
\begin{enumerate}
    \item If $\lambda\geq 0$, $\BC(\lambda)(S)=H^0(X_S,\OO_{X_S}(\lambda))$. 
    \item If $\lambda<0$, $\BC(\lambda)(S)=H^1(X_S,\OO_{X_S}(\lambda))$.
\end{enumerate}
Note that $\BC(\lambda)$ admits endomorphisms by the division algebra over $\Q_p$ of invariant $\lambda$.  
\end{defn}

\begin{defn} [The characteristic $p$ period domains $\Omega_{h,t}$]  
\label{defn:Omegaht}
Considering that the direct sum $\BC(-1/t)^{\oplus (h-t)}$ classifies extensions
\[ 0 \to \OO^{h-t} \to \E\to \OO(1/t)\to 0, \] 
we let 
\[ \Omega_{h,t}\subset \BC(-1/t)^{\oplus (h-t)} \]
be the open subsheaf classifying those extensions for which $\E$ is isoclinic of slope $1/h$, i.e., it is the pullback of $\Bun^{1/h}$ under $\BC(-1/t)^{\oplus (h-t)}\to \Bun$.  
\end{defn}

We are ready to present our main two-tower isomorphism statement.

\begin{proposition}[The two-tower isomorphism in characteristic $p$] Recall the notation $\overline{\G}_h=\Aut \Gamma_h$ and $D_h=\End \fg_h[1/p]$ for each $h\geq 1$.  
\begin{enumerate}
    \item 
There is a diagram of v-sheaves over $k$
\[
\xymatrix{
& \Sht_{h,t}^{\infty} \ar[dl]_{\overline{\G}_t  \times\GL_{h-t}(\Z_p)} \ar[dr]^{\overline{\G}_h^0} & \\
\Sht_{h,t} && \Omega_{h,t}
}
\] 
in which each morphism is a torsor for the indicated group.

\item  There is a similar diagram of v-sheaves over $k$:
\[
\xymatrix{
& \LT_{h,t}^{\infty}\ar[dl]_{\overline{\G}_t\times\GL_{h-t}(\Z_p)} \ar[dr]^{\overline{\G}_h} & \\
\LT_{h,t}^{\Diamond} && \Omega_{h,t}
}
\] 
\end{enumerate}
\end{proposition}

\begin{proof}  The left morphism in (1) is a torsor as indicated by \cref{lem:caraianischolze}.  The right morphism comes from the description of $\Sht_{h,t}^\infty$ given in \cref{ThmShtInf}:  $\Sht_{h,t}^\infty$ classifies the same data as in $\Omega_{h,t}$ together with a trivialization $\OO(1/h)\isom \E$.  Such trivializations form a torsor under $\Aut \OO(1/h)\isom \overline{\G}_h^0$.

(2) is obtained from (1) by restricting over the inclusion $\LT_{h,t}^\Diamond\to \Sht_{h,t}$.  
\end{proof}

\subsection{The case $t=h-1$}

We provide some more details in the coheight 1 case. Let us review the basic setup:  $k$ is a perfect field  containing $\F_{p^{h(h-1)}}$, and we have fixed formal groups $\Gamma_h$ and $\Gamma_{h-1}$ of dimension 1 and heights $h$ and $h-1$ over $k$.  Then the adic space $\LT_{h,h-1}$ is one point $\Spa K$, where $K\isom k\laurentseries{u_{h-1}}$.  Recall that $\Gamma_h^{\un}$ is the universal deformation of $\Gamma_h$ to $\LT_h$.  The identity component $(\Gamma_{h,K}^{\un})^{\circ}$ of $\Gamma_{h,K}^{\un}$ has height $h-1$, so it becomes isomorphic to $\Gamma_{h-1}$ over a separable extension of $K$.  As in \cref{def:torii_ext}, $L/K$ be the separable extension which classifies isomorphisms $\Gamma_{h-1,L}\isom (\Gamma^{\un}_L)^\circ$.      

Let $\widehat{L}$ be the $u_{h-1}$-adic completion of $L$.  Then $\widehat{L}$ is a topological field admitting a continuous action of $\overline{\G}_h\times \overline{\G}_{h-1}$.  To state our results, it will be convenient to fix isomorphisms of Dieudonn\'e modules $\det M(\Gamma_h)\isom \det M(\Gamma_{h-1})\isom M(\Gamma_1)$. (A quick calculation shows that for the Honda formal group $\Gamma_n$ of height $n$, the determinant $\det M(\Gamma_n)$ is isomorphic to the Dieudonn\'e module of the multiplicative group over $\F_{p^2}$.  Since $h(h-1)$ is even, $\F_{p^2}$ is contained in $k$, so such an isomorphism exists, and the set of all such isomorphisms is a $\Z_p^\times$-torsor.)
Such an isomorphism induces an isomorphism of vector bundles $\det \OO(1/h)\isom \det \OO(1/(h-1))$, which we call $\iota$.

\begin{proposition} \label{PropL} Regarding the completion of the Gross--Torii extension $L/K$, we have:
\begin{enumerate}
    \item The completion $\widehat{L}$ is a perfect (and thus perfectoid) field. 
      
    \item $\Spd \widehat{L}$ classifies exact sequences 
    \[ 0 \to \OO\to \OO(1/h)\to \OO(1/(h-1))\to 0\] 
    which induce the fixed isomorphism $\iota\from \det\OO(1/h)\isom \det \OO(1/(h-1))$.
    \item As a v-sheaf with $\overline{\G}_h\times \overline{\G}_{h-1}$-action, $\Spd \widehat{L}$ classifies surjective maps $\OO(1/h)\to \OO(1/(h-1))$. 
    \item $\widehat{L}$ is isomorphic to the completed perfection of a Laurent series field:
\[ \widehat{L}\isom k\laurentseries{t^{1/p^\infty}}. \]
    \item Recall that $S\overline{\G}_h$ is the kernel of the determinant map $\overline{\G}_h\to \overline{\G}_1=\Z_p^\times$.  Then the quotient $(\Spd \widehat{L})/S\overline{\G}_h$ is isomorphic to the punctured negative BC space $\Omega_{h,h-1}=\BC(-1/(h-1))^\ast$.  With respect to this isomorphism, the $\overline{\G}_{h-1}\times\Z_p^\times$-action on $\Omega_{h,h-1}$ is described as follows:
    \begin{enumerate}
        \item The action of $\overline{\G}_{h-1}$ is a twist of the natural action of this group (arising from its action on $\OO(1/(h-1))$) by the character $\det_{h-1}\from \overline{\G}_{h-1}\to\Z_p^\times$.
        \item An element $a\in \Z_p^\times$ (representing the determinant of an element of $\overline{\G}_h$) acts through multiplication by the scalar $a^{-1}$.  
    \end{enumerate} 
\end{enumerate}
\end{proposition}
\begin{proof} (1): $L/K$ is a totally ramified $p$-adic Lie extension, so that by \cite[Theorem 1.2]{Wintenberger1979}
it is strictly APF.  As pointed out in \cite[Remark~2.6]{Wintenberger1979}, this implies that $\widehat{L}$ is perfect.

(2) The map of $v$-sheaves $\LT_{h,h-1}^{\infty}\to \Spd \widehat{L}$ which forgets the trivialization of the \'etale part of $\Gamma_{h,K}^{\un}$ is a pro-finite $\Z_p^\times$-torsor.  By \cref{ThmShtInf}(2), this quotient classifies exact sequences of 
vector bundles as described in the proposition, where the induced isomorphism of determinants is fixed.

(3):  Extensions as in (2) are the same as surjective maps $\OO(1/h)\to \OO(1/(h-1))$, since the kernel $\mathcal{L}$ is necessarily a line bundle which comes equipped with an isomorphism to $\det \OO(1/h)\otimes (\det \OO(1/(h-1)))^{-1}\isom \OO$ via $\iota$. 

(4):  By (3), $\Spd \widehat{L}$ classifies surjections $\OO(1/h)\to \OO(1/(h-1))$.  We claim that any map $\OO(1/h)\to \OO(1/(h-1))$ which is nonzero at every point of a test object $S$ are automatically surjective.  For this we can argue on a geometric point of $S$.  The image of such a map is a quotient of $\OO(1/h)$ and a subbundle of $\OO(1/(h-1))$. If $d$ and $r\geq 1$ are relatively prime integers and $d/r$ appears as a slope of this image, then $d/r$ lies in $[1/h,1/(h-1)]$.  In particular $d\geq 1$.  But also the rank of this image is $\leq h-1$, implying $r\leq h-1$ and therefore $d/r\geq 1/(h-1)$; the only possibility is $d/r=1/(h-1)$.  Therefore the map is surjective.

Nonzero maps $\OO(1/h)\to \OO(1/(h-1))$ are the same as nonzero global sections of $\OO(-1/h)\otimes \OO(1/(h-1))\isom \OO(1/h(h-1))$; they are parametrized by the complement of the origin in $\BC(1/h(h-1))$.  (See \cite[Example 2.2.3]{BirkbeckExtensions} for this manipulation.)  By \cref{ThmFullFaithfulness} combined with \cref{LemCrisBundles}, $\BC(1/h(h-1))$ is isomorphic to the v-sheafification of the presheaf sending $S=\Spa(R,R^+)$ to $\Hom(\Q_p/\Z_p,\Gamma_{h(h-1)})[1/p]$ in the category of $p$-divisible groups over $R^+/\varpi$.  The latter is the universal cover $\tilde{\Gamma}_{h(h-1)}=\lim_p \Gamma_{h(h-1)}$ evaluated on $R^+/\varpi$.  This universal cover is represented by $\Spf k\powerseries{t^{1/p^\infty}}$ by \cite[Proposition 3.13]{ScholzeWeinstein}, and (by the crystalline property of such universal covers) it takes the same value on $R^+$.  After v-sheafifying, we find an isomorphism $\BC(1/h(h-1))\isom \Spd k\powerseries{t^{1/p^\infty}}$.  Removing the origin, we find that $\Spd \widehat{L}$ is isomorphic to $\Spd k\laurentseries{t^{1/p^\infty}}$.  

(5):  In light of (3), the quotient $(\Spd \widehat{L})/S\overline{\G}_h$ classifies surjective maps $\E\to \OO(1/(h-1))$, where $\E$ is an unspecified rank $h$ vector bundle which is isoclinic of slope $1/h$, together with an isomorphism $\det\E \isom \OO(1)$.  The kernel of $\E\to \OO(1/(h-1))$, being a line bundle, is identified with $\det\E\otimes\det \OO(1/(h-1))^{-1}\isom \OO$.  We find that $(\Spd \widehat{L})/S\overline{\G}_h$ classifies extensions
\[ 0 \to \OO\to \E\to \OO(1/(h-1))\to 0 \] 
for $\E$ isoclinic of slope $1/h$, as indeed such an exact sequence forces a trivialization of $\det \E$.  Since such an extension is necessarily everywhere nonsplit, $(\Spd \widehat{L})/S\overline{\G}_h$ maps toward 
$\BC(-1/(h-1))^\ast$.  Conversely, a section of $\BC(-1/(h-1))^*$ corresponds to an extension as above for which $\E$ which is necessarily isoclinic of slope $1/h$, as there are no Newton polygons between $(1/(h-1),0)$ and $(1/h)$.  

The descriptions of the group actions in (a) and (b) are  required to preserve the isomorphism $\det \E \otimes \det \OO(1/(h-1))^{-1}\isom \OO$.  (As a sanity check, note that for $a\in\Z_p^\times$, the pair $(a,a)\in \overline{\G}_{h-1}\times \overline{\G}_h$ acts trivially on $\Spd \widehat{L}$, and indeed $a\det_{h-1}(a)\det_h(a)^{-1}=1$.)  
\end{proof}

\newpage
\section{The $\OO^\flat$-cohomology of characteristic $p$ period domains}\label{sec:oflat}

Let $L/K$ be the Torii extension (\cref{def:torii_ext}), so that $K\isom k\laurentseries{u_{h-1}}$ is a field of Laurent series over a perfect field $k/\F_p$, and $L/K$ is a Galois extension with group $\overline{\G}_{h-1}$.   The field $L$ admits an action of $\G_{h,h-1}$.  Let $\widehat{L}$ denote the $u_{h-1}$-adic completion of $L$.  In this section we prove \cref{thmx:cohom}, on the structure of the continuous cohomology $H^*_{\cts}(S\overline{\G}_h,\widehat{L})$, where $S\overline{\G}_h$ is the kernel of the determinant map $\overline{\G}_h\to \overline{\G}_1$.   

\cref{ThmMeaningOfL} provides an isomorphism of v-sheaves on $\Perf_k$:
\[ (\Spd \widehat{L})/S\overline{\G}_h \isom \BC(-1/(h-1))^\ast. \] 
We adopt the notation $\OO^\flat$ for the structure sheaf on $\Perf_k$, which sends $\Spa(R,R^+)$ to $R$.  Taking $\OO^\flat$-cohomology on either side of the above isomorphism, and noting that $\OO^\flat$ is acyclic on affinoid perfectoids such as $\Spa \widehat{L}$, we obtain an isomorphism of graded rings:
\[ H^*_{\cts}(S\overline{\G}_h,\widehat{L}) \isom H^*(\BC(-1/(h-1))^*,\OO^\flat). \] 

Let us outline our strategy for computing the right side, i.e., $H^*(\BC(-1/t)^*,\OO^\flat)$ for $t\geq 1$. In \cref{sec:GeometryOfOmega} we study the geometry of $\BC(-1/t)^*$ and show that after base-change to an auxiliary algebraically closed perfectoid field $C/\Q_p$, it becomes isomorphic to 1-dimensional Drinfeld symmetric space modulo the action of a finite-dimensional $\Q_p$-vector space. 
In Section~\ref{sec:6FunctorsforOflat}, we recall the 6-functor formalism for $\OO^\flat$-cohomology from \cite{ALBM-6ff-on-FF-curve} and compute the dualizing complex for smooth morphisms. This part does not depend on the other parts of the paper and may be of interest in its own right. 
In preparation for the computation of the $\OO^\flat$-cohomology of the Drinfeld symmetric space, we study the $\OO^\flat$-cohomology of closed affinoid disks over $\Spa C$ in \cref{sec:O-flat-cohom-of-disk} using a very general approach that allows similar computations for every cohomology theory on v-stacks. In \cref{sec:cohomology-of-Drinfeld-space}, we compute the $\OO^\flat$-cohomology and $\OO^\flat$-homology of the Drinfeld symmetric space. As $\GL_2(F)$-representations, the homology turns out to consist of the trivial representation and the Steinberg representation. 

In \cref{sec:Heyer-Jacquet}, we develop a general comparison between mod-$p$ representations of a $p$-adic group with discrete coefficients and those with condensed coefficients. This comparison is then used to compute the $\OO^\flat$-cohomology of $\BC(-1/t)^*$ over $C$.
Finally, in \cref{sec:DrinfeldsLemma} we prove a version of Drinfeld's lemma for $\OO^\flat$-cohomology, which allows us to deduce the cohomology of $\BC(-1/t)^*$ over $k$ from that of $\BC(-1/t)^* \times \Spa C$ computed in the previous section.

Recall that $\BC(-1/(h-1))^*=\Omega_{h,h-1}$ is an example of a characteristic $p$ period domain as in \cref{defn:Omegaht}.  It ought to be possible to extend the results of this section to compute the $\OO^\flat$-cohomology of $\Omega_{h,t}$ for general $1\leq t\leq h-1$;  this may be relevant to the splitting conjecture in general coheight. 

\subsection{The geometry of $\BC(-1/t)^*$} \label{sec:GeometryOfOmega}

Fix a (discrete) perfect field $k/\F_p$ containing $\F_{p^t}$, and use it to define the v-sheaf $\BC(-1/t)^*$ on $\Perf_k$.  Let $C/\Q_p$ be an algebraically closed perfectoid field $C$ whose residue field contains $k$, so that $\Spd C$ is fibered over $\Spd k$.  We study the geometry of the v-sheaf $\BC(-1/t)^*$ and relate it to Drinfeld symmetric spaces. Let us first review what these spaces are. 

\begin{defn}
\label{defn:DrinfeldSymmSpace} 
    Let $F$ be a finite extension of $\Q_p$ contained in $C$. The \emph{Drinfeld symmetric space} of dimension $n$ associated to $F$ is the rigid-analytic space over $F$ defined by 
    \begin{equation*}
        \HH^n_{F} \coloneqq  \PP^n_F \setminus \bigcup_H H,
    \end{equation*}
    where $H$ ranges over all $F$-rational hyperplanes in the projective space $\PP^n_F$. $\HH^n_{F}$ is a smooth Stein space of dimension $n$ over $F$. Let  $\HH^n_{F,C}$ denote the base change of $\HH^n_F$ to $C$.
\end{defn}

The main result of the present subsection is the following description of $\BC(-1/t)^*$ after base-change to $\Spd C$:

\begin{proposition} \label{rslt:BaseChangeOfOmegaToC} Let $F/\Q_p$ be a field extension of degree $t$ contained in $C$.  
\begin{enumerate}
    \item There is an isomorphism of $\Q_p$-linear diamonds over $\Spd C$:
    \[  \BC(-1/t)\times_{\Spd k}\Spd C \simeq \A^{1,\Diamond}_C/F, \]
    where $F$ acts by translation on $\A^{1,\Diamond}_C$. 
    \item There is a $\Q_p^\times$-equivariant isomorphism of diamonds over $\Spd C$:
    \[ \BC(-1/t)^\ast \times_{\Spd k}\Spd C \simeq \HH^{1,\Diamond}_{F,C}/F, \] 
    where $F$ acts by translation on $\HH^{1,\Diamond}_{F,C}\subset \A^{1,\Diamond}_C$.  
\end{enumerate}
\end{proposition}

\begin{proof} Part (1) is well-known; see the discussion after \cite[Definition II.0.6]{FarguesScholze} for the case $t = 1$.  Also, we remark that we ought to write $\A^{1,\Diamond}_{F,C}/\underline{F}^{\Diamond}$ for this quotient.  

For $S=\Spa(R,R^+)$ a strictly totally disconnected affinoid perfectoid space over $\Spa(C,\OO_C)$, let $X_{S^\flat,F}$ be the Fargues--Fontaine curve relative to the base $S^\flat$ and the local field $F$.    The affinoid perfectoid $S=\Spa(R,R^+)$, being an untilt of its tilt $S^\flat$, defines a Cartier divisor $i:S\hookrightarrow X_{S^\flat,F}$. 
    
Choose a generator of the Tate module of the Lubin--Tate formal $\OO_F$-module, with coefficients in $C$.  This choice induces a short exact sequence of sheaves on $X_{S^\flat,F}$:
    \begin{equation*}
        0\to \OO_{X_{S^\flat,F}}(-1)\to \OO_{X_{S^\flat,F}}\to i_{*}\OO_S\to 0. 
    \end{equation*}
    Let $\pi\from X_{S^\flat,F}\to X_{S^\flat}$ be the natural degree $t$ \'etale map corresponding to the degree $t$ extension $F/\Q_p$.  Pushing forward the modification by $\pi_*$, we obtain a short exact sequence of sheaves on $X_{S^\flat}$:
    \begin{equation*}
        0\to \OO_{X_{S^\flat}}(-1/t)\to \OO_{X_{S^\flat}} \otimes_{\Q_p} F\to (\pi\circ i)_* \OO_S\to 0.
    \end{equation*}
    (See \cite[\S 8.2.2]{FarguesFontaineCourbes} for a discussion of this push-forward operation on vector bundles.)
   Passing to the long exact sequence in cohomology, we obtain a short exact sequence of $\Q_p$-vector spaces:
    \begin{equation*}
        0 \to \underline{F}(S) \to \mathbb{A}^{1,\diamond}_C(S) \to H^1(X_{S^\flat},\OO_{X_{S^\flat}}(-1/t))\to 0,
    \end{equation*}
    which is to say, an isomorphism
    \begin{equation*}
        H^1(X_{S^\flat},\OO(-1/t))\isom \mathbb{A}^{1,\Diamond}_C(S)/\underline{F}(S).
    \end{equation*}
   The diamond $\BC(-1/t)\times_k\Spd C$ is now seen to be the sheaf of $\Q_p$-vector spaces whose value on a strictly totally disconnected affinoid perfectoid space $S$ over $C$ is $\mathbb{A}^{1,\diamond}_C(S)/\underline{F}(S)$.  This means exactly that 
    \begin{equation*}
        \BC(-1/t)\times_{\Spd k}\Spd C\simeq \A^{1,\diamond}_C/F,
    \end{equation*}
which is (1).

The isomorphism (2) is obtained from (1) by removing the zero section on either side.
\end{proof}

In passing, we remark that the isomorphisms in \cref{rslt:BaseChangeOfOmegaToC} are defined over the completion of the maximal abelian extension of $F$.

\subsection{A 6-functor formalism for $\OO^\flat$-cohomology}\label{sec:6FunctorsforOflat}

Let $\vStk$ denote the category of small v-stacks as defined in \cite[\S 12]{ScholzeDiamond}. One of the main results of \cite{ALBM-6ff-on-FF-curve} is the construction of a 6-functor formalism for $\OO^\flat$-cohomology on small v-stacks:

\begin{thm}\label{thm:6FunctorsOflat}
There is a certain class of morphisms $E$ in $\vStk$ and a 6-functor formalism $X \mapsto \mathcal D_{\hat\solid}(\OO^\flat_X)$ on the geometric setup $(\vStk, E)$ (i.e., with !-functors defined for all maps in $E$) such that the following is true:
\begin{thmenum}
    \item $E$ enjoys the following properties:
    \begin{enumerate}[(a)]
        \item $E$ is stable under composition and base-change. Moreover, if $f\colon Y \to X$ and $g\colon Z \to Y$ are maps of small v-stacks such that $f, fg \in E$, then $g \in E$.

        \item $E$ contains every morphism between (the diamonds of) rigid varieties over some fixed base field and every map between the classifying stacks of $p$-adic Lie groups.
    \end{enumerate}

    \item \label{rslt:6FunctorsOnAffdPerfd} If $X = \Spa(A, A^+)$ is an affinoid perfectoid space in characteristic $p$ which admits a map of finite transcendence degree to a totally disconnected space (e.g., the spectrum of perfectoid field) then
    \begin{align*}
        \mathcal D_{\hat\solid}(\OO^\flat_X) = \mathcal D_{\hat\solid}(A, A^+),
    \end{align*}
    where the right-hand side denotes the (slightly modified) category of $(A, A^+)_\solid$-modules.\footnote{Roughly, this is the derived category of topological $A$-modules with a certain completeness condition with respect to $A^+$. If $A^+$ is up to integral closure of topological finite type over $\F_p$ (e.g., $A = \C_p$), then the $\hat\solid$-modification does nothing and agrees with usual solid modules.}

    \item \label{rslt:6FunctorsComputeCohomology} Given a small v-stack $X$, let $\OO^\flat \in \mathcal D_{\hat\solid}(\OO^\flat_X)$ denote the tensor unit. Then
    \begin{align*}
        \Hom(\OO^\flat, \OO^\flat) = R\Gamma(X, \OO^\flat),
    \end{align*}
    where the cohomology on the right denotes v-cohomology.
    \item \label{rslt:6FunctorsSmoothProperMaps} We have the following results regarding smooth and proper maps:
    \begin{enumerate}[(a)]
        \item Let $f\colon Y \to X$ be a map of analytic adic spaces over $\Q_p$; we abuse notation and denote by $f$ also the diamond version. If $f$ is proper then $f$ is $\mathcal D_{\hat\solid}$-proper, i.e., $f_! = f_*$. If $f$ is smooth of pure dimension $d$ then $f$ is $\mathcal D_{\hat\solid}$-smooth and naturally $f^! = f^*[2d](d)$, where ``$(d)$'' denotes the Tate twist.

        \item Let $G$ be a $p$-adic Lie group and consider the map $f\colon */G \to *$ of v-stacks. Then $f$ is $\mathcal D_{\hat\solid}$-$!$-able, i.e., it belongs to the class $E$. The map $f$ is $\mathcal D_{\hat\solid}$-smooth and $f^! = f^* \otimes \omega_G[\dim G]$, where $\omega_G$ is a character of $G$. If $G$ is compact with $\mathrm{cd}_p(G)<\infty$, then $f$ is $\mathcal D_{\hat\solid}$-proper.
    \end{enumerate}
\end{thmenum}
\end{thm}
\begin{proof}
We let $E$ be the class of !-able maps from \cite[Definition~4.3.6]{ALBM-6ff-on-FF-curve}. By construction this class satisfies (i)(a) and by \cite[Lemma~4.3.5(ii)]{ALBM-6ff-on-FF-curve}) and \cite[Proposition~4.3.8]{ALBM-6ff-on-FF-curve} it satisfies (i)(b). In \cite[Theorem~3.2.8]{ALBM-6ff-on-FF-curve} a 6-functor formalism $X \mapsto \mathcal D_{\hat\solid}(\OO^{\flat,+,a}_X)$ is constructed. We obtain the desired one by restricting to the full subcategory of $\OO_X$-modules (i.e., inverting a pseudouniformizer). Claim (ii) then follows from \cite[Theorem~3.1.2]{ALBM-6ff-on-FF-curve}. Claim (iii) is obvious on affinoid perfectoid spaces as in (ii) (then both sides compute $A$) and follows in general because $\mathcal D_{\hat\solid}$ satisfies v-descent by \cite[Theorem~3.1.2]{ALBM-6ff-on-FF-curve}. The claim about properness in (iv)(a) is true by construction. The claim about smooth maps in (iv)(a) reduces to the 6-functor formalism $X \mapsto \mathcal D_{\hat\solid}(\OO^{\flat,+,a}_X)$ because morphisms of 6-functor formalisms preserve smoothness and commute with taking the dualizing sheaf for smooth maps (this follows from the interpretation of smoothness in terms of the category of kernels, see \cite[Definition~4.5.1]{HeyerMann6FF}). Then \cite[Lemma~3.3.2]{ALBM-6ff-on-FF-curve} reduces the question to $\OO^{\flat,+,a}/\pi$ for some pseudouniformizer $\pi$, where we can deduce it from \cite[Theorem~3.10.17]{MannThesis} and \cite[Theorem~3.10.20]{MannThesis}. The claims in (iv)(b) can similarly be reduced modulo $\pi$ and then follow from \cite[Corollary~3.12]{HansenMannJL} and \cite[Theorems~3.16,~3.18]{HansenMannJL}.
\end{proof}

\begin{rmk}
    Let $G$ be a compact $p$-adic Lie group, and $f:*/G\to *$ be the associated map of v-stacks. 
    Then the following statements are equivalent:
    \begin{enumerate}
        \item $f$ is $\mathcal D_{\hat\solid}$-proper;
        \item $\mathrm{cd}_p(G)<\infty$; and
        \item $G$ is $p$-torsion free.
    \end{enumerate}
    The equivalence between (2) and (3) is a result by Serre \cite{SerreCohDim}, and \cref{rslt:6FunctorsSmoothProperMaps} shows that (2) implies (1). To show that (1) implies (2), it is enough to show that the pullback of $f:*/G\to *$ to a geometric point $\Spd C\to *$ is not $\mathcal D_{\hat\solid}$-proper if $G$ has nontrivial $p$-torsion. Assume $H\simeq\Z/p$ is a closed subgroup of $G$. Then there are induced maps between small v-stacks 
    \[\begin{tikzcd}
	{\Spd C/H} & {\Spd C/G} & {\Spd C}.
	\arrow["i", from=1-1, to=1-2]
	\arrow["g"', curve={height=12pt}, from=1-1, to=1-3]
	\arrow["f", from=1-2, to=1-3]
    \end{tikzcd}\]
    Here $i$ is $\mathcal D_{\hat\solid}$-proper: pulling $i$ back along the natural map $\Spd C\to \Spd C/G$ yields $\underline{G/H}\to \Spd C$, which is $\mathcal D_{\hat\solid}$-proper since $G/H$ is a profinite set. Let $M\coloneqq i_* \mathbf{1}\simeq i_! \mathbf{1}$, where $\mathbf{1}$ is the tensor unit in $\mathcal D_{\hat\solid}(\Spd C/H)$ and corresponds to the trivial representation of $H$ over $C$.
    It is then enough to show $f_*M\neq f_!M$, or $g_* \mathbf{1}\neq g_! \mathbf{1}$. A simple calculation in group cohomology shows $H^1(\Spd C, g_* \mathbf{1})\neq H^1(\Spd C, g_! \mathbf{1})$.
\end{rmk}

The next result discussed a canonical topology, more precisely condensed structure, that can be put on the cohomologies in the above 6-functor formalism. This will become important in the Drinfeld lemma below.

\begin{proposition} \label{rmk:solid-v-cohomology}
For every small v-stack $X$ there is a canonical functor
\begin{align*}
    R\Gamma_\solid(X, -)\colon \mathcal D_{\hat\solid}(\OO^\flat_X) \to \mathcal D_\solid(\F_p)
\end{align*}
with the following properties:
\begin{propenum}
    \item Given $M \in \mathcal D_{\hat\solid}(\OO^\flat_X)$, the underlying $\F_p$-vector space of $R\Gamma_\solid(X, M)$ is given by
    \begin{align*}
        R\Gamma_\solid(X, M)(*) = \Hom(\OO^\flat, M).
    \end{align*}
    In particular, for $M = \OO^\flat$ we get $R\Gamma_\solid(X, \OO^\flat)(*) = R\Gamma(X, \OO^\flat)$.

    \item \label{rmk:solidCohomologyoverC} Suppose there is a map $f\colon X \to S$, where $S = \Spa(R, R^+)$ is an affinoid perfectoid space as in \cref{rslt:6FunctorsOnAffdPerfd}. Then there is a canonical isomorphism $R\Gamma_\solid(X, -) = f_*$, where on the right-hand side we implicitly apply the forgetful functor $\mathcal D_{\hat\solid}(\OO^\flat_S) = \mathcal D_{\hat\solid}(R, R^+) \to \mathcal D_\solid(\F_p)$.

    \item \label{rslt:compute-solid-v-cohomology-via-complex} Suppose that $X$ is qcqs and pick a v-hypercover $g_\bullet\colon Y_\bullet \to X$ by affinoid perfectoid spaces $Y_n = \Spa(B_n, B_n^+)$ that satisfy the condition in \cref{rslt:6FunctorsOnAffdPerfd}; then
    \begin{align*}
        R\Gamma_{\solid}(X, \OO^\flat) = \varprojlim_{n\in\Delta} B_n = [\dots \to 0 \to B_0 \to B_1 \to B_2 \to \dots],
    \end{align*}
    where on the right-hand side we have the usual resolution coming out of the hypercover, viewed as a complex in $\mathcal D_\solid(\F_p)$ for the natural condensed structures on $B_n$.
\end{propenum}
\end{proposition}
\begin{proof}
Note that for every affinoid perfectoid space $Y = \Spa(B, B^+)$ there is a canonical pullback functor $\mathcal D_\solid(\F_p) \to \mathcal D_{\hat\solid}(B, B^+)$ and this pullback functor commutes with pullbacks along maps of affinoid perfectoid spaces (cf.\ \cite[Lemma~2.7]{AM-Descent-on-vStacks}). Passing through descent, we see that $\mathcal D_{\hat\solid}$ together with its pullback and tensor product can be seen as a functor
\begin{align*}
    \mathcal D_{\hat\solid}\colon \vStk^\op \to \CAlg(\mathrm{Pr}^L)_{\mathcal D_\solid(\F_p)/},
\end{align*}
i.e., there is a natural symmetric monoidal map $\alpha^*_X\colon \mathcal D_\solid(\F_p) \to \mathcal D_{\hat\solid}(\OO^\flat_X)$ for every small v-stack $X$. The functor $R\Gamma_\solid(X, -)$ is defined to be the right adjoint of $\alpha_X^*$.

To prove (i), we note that for any $M \in \mathcal D_{\hat\solid}(\OO^\flat_X)$ we have
\begin{align*}
    R\Gamma_{\solid}(X, M)(*) = \Hom(\F_p, R\Gamma_{\solid}(X, M)) = \Hom(\alpha_X \F_p, M) = \Hom(\OO^\flat, M),
\end{align*}
as desired. For (ii) we observe that both sides of the claimed identity are right adjoint to $\alpha_X$. Finally, for (iii) we note that by descent and \cref{rslt:6FunctorsOnAffdPerfd} we have
\begin{align*}
    \mathcal D_{\hat\solid}(\OO^\flat_X) = \varprojlim_{n\in\Delta} \mathcal D_{\hat\solid}(B_n, B_n^+),
\end{align*}
and so the right adjoint of $\alpha_X = (\alpha_{Y_n})_n$ is given by the limit over $\Delta$ of the right adjoints of the $\alpha_{Y_n}$'s (see \cite[Lemma~D.4.7(ii)]{HeyerMann6FF}).
\end{proof}

\begin{rmk} \label{rmk:solid-v-cohomology-over-k}
The results in \cref{rmk:solid-v-cohomology} hold verbatim if one replaces $\F_p$ by an arbitrary discrete perfect field $k$ of characteristic $p$ and works with small v-stacks over $\Spd k$.
\end{rmk}

\begin{notations}
    Let $X$ be a small v-stack. For any $M\in \mathcal D_{\hat\solid}(\OO^\flat_X)$, we write $H^i_{\solid}(X,M)$ for $H^i(R\Gamma_{\solid}(X,M))\in\Mod_{\solid}(\F_p)$. Note that $H^i_{\solid}(X,M)(\ast)=H^i(X,M)$ since $\F_p$ is projective in $\Mod_{\solid}(\F_p)$.
\end{notations}

An important tool in the computation of various cohomology theories is that of the excision sequence: Given a space $X$ with some closed subspace $i\colon Z \injto X$ and open complement $j\colon U \injto X$ and a 6-functor formalism $\mathcal D$, one often has a natural fiber sequence $j_! j^* M \to M \to i_* i^* M$ in $\mathcal D(X)$ for any $M \in \mathcal D(X)$. This principle works for example for sheaves on topological spaces and for étale sheaves on schemes and it is a powerful tool to compute cohomology in these settings. In the coherent setting, excision sequences do not hold naively, but one can often find appropriate versions that do hold if one is willing to reinterpret the meaning of $i_* i^* M$ slightly. The following result demonstrates this in the setting of $\OO^\flat$-sheaves:

\begin{proposition} \label{rslt:excision-for-O-flat}
Let $X$ be a small v-stack, $Z \subseteq \abs X$ a closed subset and $j\colon U \injto X$ its open complement. Assume that:
\begin{itemize}
    \item[($*$)] $Z = \bigcap_{i\in I} \abs{U_i}$ for a cofiltered system of open substacks $j_i\colon U_i \subseteq X$ such that for each $i \in I$ there is some $i' \ge i$ and a closed subset $Z_i \subseteq \abs X$ with factorization $\abs{U_{i'}} \subseteq Z_i \subseteq \abs{U_i}$.
\end{itemize}
Then for every $M \in \mathcal D_{\hat\solid}(\OO^\flat_X)$ there is a natural fiber sequence
\begin{equation*}
    j_! j^* M \to M \to \varinjlim_{i \in J} \: (j_i)_* j_i^* M
\end{equation*}
in $\mathcal D_{\hat\solid}(\OO^\flat_X)$.
\end{proposition}
\begin{proof}
One can argue exactly as in \cite[Corollary 6.1.4]{ALBM-6ff-on-FF-curve} (which is stated for a slightly different 6-functor formalism) and the comments after it.
\end{proof}

\begin{rmk}
The third term in the fiber sequence of \cref{rslt:excision-for-O-flat} should be seen as an overconvergent approximation to $i_* i^* M$, where $i$ is the closed immersion of $Z$. It is necessary already for $M = \OO^\flat_X$: In this case the cofiber of $j_! \OO^\flat_U \to \OO^\flat_X$ is not complete, while $i_* \OO^\flat_Z$ is. If one looks at $\OO^{\flat,+,a}/\pi$-modules instead, then the excision sequence works with $i_* i^*$ (see \cite[Lemma~6.2.1]{ZavyalovAbstractPoincare} for a version of this).
\end{rmk}

We establish the following lemma on base-change properties of cohomology of certain qcqs v-stacks for later use.
\begin{lemma}\label{lem:base-change-identity}
    Let $k$ be a perfect discrete field containing $\F_p$, and let $k'$ be a perfect discrete extension of $k$. Then:
    \begin{enumerate}
        \item For any affinoid perfectoid $S=\Spa(R,R^+)$ over $k$, $S\times_{k}\Spd k'$ is again affinoid perfectoid and is isomorphic to $\Spd(R\widehat{\otimes}_{k}k',R^+\widehat{\otimes}_{k} k')$, where $\widehat{\otimes}$ denotes the classical completed tensor product. 
        \item Let $X$ be a qcqs\footnote{To avoid ambiguity, this means that $X$ itself is qcqs, not that the map $X\to\Spd k$ is qcqs.} small v-stack over $k$ such that $H^i_{\solid}(X,\OO^\flat)$ is a perfect solid $k$-module for all $i$. Then there are isomorphisms $H^i_{\solid}(X\times_{\Spd k}\Spd k',\OO^\flat)\simeq H^i_{\solid}(X,\OO^\flat)\otimes^{\solid}_k k'$ of solid $k'$-modules for all $i$.
    \end{enumerate}
\end{lemma}
\begin{proof}
    To show (1), we only need to show both $S\times_{\Spd k}\Spd k'$ and $\Spa(R\widehat{\otimes}_{k} k',R^+\widehat{\otimes}_{k} k')$ define the same functor of points on the category $\Perf_{k'}$ of perfectoid spaces over $k'$, which is straightforward to check.

    For (2), we first set up some notations. If $(R,R^+)$ is a Huber pair over $k$, we write $(R\otimes_kk',R^+\otimes_kk')$ for its (ordinary) base-change to $k'$, and write $(R\widehat{\otimes}_kk',R^+\widehat{\otimes}_kk')$ for the completion of the base-change.
    Take a hypercover $f_\bullet\colon Y_\bullet \to X$ by affinoid perfectoid spaces $Y_n = \Spa(A_n, A_n^+)$ that satisfy the condition in \cref{rslt:6FunctorsOnAffdPerfd}. By \cref{rslt:compute-solid-v-cohomology-via-complex} the cohomology $R\Gamma_{\solid}(X,\OO^\flat)$ is given as 
    \begin{equation*}
        R\Gamma_{\solid}(X, \OO^\flat) = \varprojlim_{n\in\Delta} A_n
    \end{equation*}
     Write $Y_{n}\widehat{\otimes}_k k'\coloneqq \Spa(A_n\widehat{\otimes}_kk',A_n^+\widehat{\otimes}_kk')$. By (1), this is isomorphic to $Y_n\times_{ k}\Spd k'$. The natural map $Y_{\bullet}\widehat{\otimes}_k k'\to X\times_{k}\Spd k'$ is again a hypercover.
    
    Let $Z^i$ (resp. $B^i$) be the condensed $k$-module of $i$-cocycles (resp. $i$-coboundaries) in the complex associated to the cosimplicial solid $k$-module $R\Gamma(Y_\bullet,\OO^\flat)$; so $B^i\subset Z^i\subset A_i$. Let $I$ be an ideal of definition of $A^+_i$. The sequence $(I^mA^+_i\otimes_k^{\solid}k')_{m\in\N}$ of solid $k'$-modules defines an exhausted, separated filtration of $A^+_i\otimes_k^{\solid}k'$. We note there is a natural isomorphism 
    \begin{equation*}
        A\widehat{\otimes}_kk'\simeq \lim_m (A_i\otimes_k^{\solid}k')/(I^mA^+_i\otimes_k^{\solid}k'),
    \end{equation*}
    where the left side is viewed as a solid $k'$-module. Using the same filtration, we can complete $Z^{i}_{k'}\coloneqq Z^i\otimes_k^{\solid}k'$ (resp. $B^{i}_{k'}\coloneqq B^i\otimes_k^{\solid}k'$) to $\widehat{Z}^i_{k'}$ (resp. $\widehat{B}^i_{k'}$), i.e., we define 
    \begin{equation*}
        \widehat{Z}^i_{k'}\coloneqq \lim_m Z^i_{k'}/(Z^i_{k'}\cap I^mA^+_i\otimes_k^{\solid}k')
    \end{equation*}
    and similarly for $\widehat{B}^i_{k'}$. So in fact, $\widehat{Z}^i_{k'}$ (resp. $\widehat{B}^i_{k'}$) is the condensed $k'$-module of $i$-cocycles (resp. $i$-coboundaries) in the complex associated to $R\Gamma(Y_\bullet\otimes_kk',\OO^\flat)$, and we have $\widehat{H}^i_{k'}\coloneqq  H^i(X\times_k\Spd k',\OO^\flat)\simeq \widehat{Z}^i_{k}/\widehat{B}^i_{k'}$.
    
    Denote $H^i\coloneqq  H^i(X,\OO^\flat)$. By assumption, $H^i$ is a perfect solid $k$-module. So in fact, $H^i$ is a finite dimensional $k$-vector space viewed as a solid module, and there is a splitting $Z^i\simeq B^i\oplus H^i$. Since $H^i\otimes_k^{\solid} k'$ is a finite dimensional $k'$-vector space, $(Z^i_{k'}\cap I^mA^+_i\otimes_k^{\solid}k')\cap (H^i\otimes_k^{\solid} k')=0$ for $m\gg 0$. 
    This implies the completion of $Z^i_{k'}$ is $\widehat{Z}^i_{k'}=\widehat{B}^i_{k'}\oplus  (H^i\otimes_k^{\solid} k')$, giving that $\widehat{H}^i_{k'}\simeq  H^i\otimes_k^{\solid} k'$.
\end{proof}

\subsection{The \texorpdfstring{$\OO^\flat$}{O-flat}-cohomology of an affinoid disc} \label{sec:O-flat-cohom-of-disk}

We fix an algebraically closed non-archimedean extension $C$ of $\Q_p$ and consider the 1-dimensional rigid-analytic disk over $\Spa C$:
\begin{align*}
    D \coloneqq  \Spa(C\langle T\rangle, \OO_C\langle T \rangle)
\end{align*}
Our goal is to compute $H^*(D,\OO^{\flat,+,a})$.  Let $\tilde{D}$ be the 1-dimensional perfectoid affinoid disc:
\[ \tilde{D} = \Spa(C\langle T^{1/p^\infty}\rangle, \OO_C\langle T^{1/p^\infty} \rangle)  \] 
Then $\tilde{D}\to D$ is a v-cover admitting an action of 
\[ G=\Z_p(1).\]
The idea is to compute $H^*(D,\OO^{\flat,+,a})$ using v-descent along $\tilde{D}\to D$. However, we have to overcome the following obstacle: The cover $\tilde D \to D$ is not a $G$-torsor, this is only true away from $0 \in D$. Indeed, the fiber of $\tilde D$ over $0$ is the single point $0 \in \tilde D$, with $G$ acting trivially. We remedy the situation in \cref{rslt:disk-as-pushout-from-perfectoid-disk} below. Before proving this result, let us state a few basic observations about small v-stacks. The first result allows us to reduce our computation to ordinary v-stacks instead of $\infty$-stacks:

\begin{lemma} \label{rslt:infty-pushout-along-mono-is-truncated}
Let $Y \leftarrow X \to Z$ be maps of small v-stacks and assume that $X \injto Y$ is a monomorphism. Then the pushout $W \coloneqq  Y \sqcup_X Z$, formed in the $\infty$-category of $\infty$-v-stacks, is an ordinary small v-stack.
\end{lemma}
\begin{proof}
We need to show that $W$ is 1-truncated stack, i.e., for every perfectoid space $S$ the $\infty$-groupoid $W(S)$ is an ordinary groupoid. To see this, we observe that $W$ can be computed in two steps: First compute the pushout in the category of $\infty$-prestacks, i.e., separatedly on $S$-valued points for all affinoid perfectoid $S$, and then take the sheafification. Note that sheafification is computed by repeatedly applying the $+$-construction, which itself is computed by taking filtered colimits of descent data along coverings -- these operations preserve truncatedness, so it is enough to show that the prestack pushout is 1-truncated. But this follows immediately from the fact that the map $X(S) \injto Y(S)$ is a monomorphism for every perfectoid space $S$, because this implies that $Y(S) = X(S) \sqcup (Y(S) \setminus X(S))$ and hence $W(S) = Z(S) \sqcup (Y(S) \setminus X(S))$, which is indeed 1-truncated.
\end{proof}

For the following result recall that a perfectoid space is called \emph{w-contractible} if it is strictly totally disconnected and $\pi_0(S)$ is extremally disconnected. These spaces form a basis of the v-site and every pro-étale cover of a w-contractible perfectoid space splits.

\begin{lemma} \label{rslt:mono-after-sheafification-checked-with-+-construction}
Let $X$ be a small v-stack and $\alpha\colon Y \to X$ a map from a prestack $Y$. Suppose that $\alpha$ satisfies the following condition:
\begin{itemize}
    \item[($*$)] Let $S$ be a w-contractible perfectoid space and pick objects $y, y' \in Y(S)$ and an isomorphism $\varphi\colon \alpha(y) \isoto \alpha(y')$ in $X(S)$. Then there is a cover $S = U_1 \sqcup \dots \sqcup U_n$ by disjoint open subsets $U_i \subseteq S$ such that $\varphi|_{U_i}$ admits a unique preimage $y|_{U_i} \isoto y'|_{U_i}$ for all $i$.
\end{itemize}
Then $\alpha$ becomes a monomorphism after v-sheafification.
\end{lemma}
\begin{proof}
Since v-sheafification is left exact and hence preserves monomorphisms, it is enough to show that $\alpha$ becomes a monomorphism after pro-étale sheafification. In fact, it is even enough to check that $\alpha$ becomes a monomorphism after applying the $+$-construction (replacing $Y(S)$ with the colimit of descent data along pro-étale covers of $Y$). Moreover, it is enough to check this isomorphism on w-contractible spaces. Altogether we have reduced the claim to showing that $Y^+(S) \to X(S)$ is a monomorphism for every w-contractible space $S$. Since every pro-étale cover of $S$ splits, one checks easily that $Y^+(S)$ is the colimit of $Y(U_1) \times \dots \times Y(U_n)$ for all decompositions $S = U_1 \sqcup \dots \sqcup U_n$ as in ($*$). It is then not hard to see that the condition in ($*$) is exactly saying that $Y^+(S) \to X(S)$ is a monomorphism.
\end{proof}

With the above preparations at hand, we can now come to the promised pushout description of the disk in terms of the perfectoid disk:

\begin{proposition} \label{rslt:disk-as-pushout-from-perfectoid-disk}
Let $\tilde D \to D$ be the v-cover of the disk by the perfectoid disk over $C$ and consider the action of $G = \Z_p(1)$ on $\tilde D$ over $D$. Then the square
\begin{equation*}\begin{tikzcd}
    */G \arrow[r] \arrow[d,"0"] & * \arrow[d]\\
    \tilde D/G \arrow[r] & D
\end{tikzcd}\end{equation*}
is a pushout square in the $\infty$-category of $\infty$-v-stacks on perfectoid spaces. Here $* = \Spa(C)$ and the quotient $\tilde D/G$ is to be understood in the stacky sense.
\end{proposition}
\begin{proof}
Let us denote $Y \coloneqq  (\tilde D/G) \sqcup_{(*/G)} *$ the pushout in the $\infty$-category of $\infty$-stacks (i.e., sheaves of anima), so that there is a natural map $\alpha\colon Y \to D$. Since the map $*/G \injto \tilde D/G$ is a monomorphism, it follows from \cref{rslt:infty-pushout-along-mono-is-truncated} that $Y$ is an ordinary small v-stack. Moreover, $\alpha$ is surjective because the composition $\tilde D/G \to Y \to D$ is surjective. It thus remains to show that $\alpha$ is a monomorphism. To see this, let us denote by $Y'$ the prestack pushout of $(\tilde D/G) \leftarrow */G \to *$. We get an induced map $\alpha'\colon Y' \to D$ which becomes $\alpha$ after sheafification. It is thus enough to show that $\alpha'$ satisfies the assumption ($*$) in \cref{rslt:mono-after-sheafification-checked-with-+-construction}.

Fix a w-contractible space $S$. Note that $(\tilde D/G)(S)$ is the groupoid of $G$-torsors $\pi\colon P \to S$ together with a $G$-equivariant map $f\colon P \to \tilde D$. By the proof of \cref{rslt:infty-pushout-along-mono-is-truncated} we see that $Y'(S)$ admits a very similar description, only the morphisms are different: Given two $G$-torsors $P, P'$ as above, if the maps $f\colon P \to \tilde D$ and $f'\colon P' \to \tilde D$ factor over $0$ then $\Hom(P, P') = *$ in $Y'(S)$, otherwise $\Hom(P, P')$ is as in $(\tilde D/G)(S)$. The map $\alpha'$ sends $(P, f) \in Y'(S)$ to the map $S \to D$ which sends $x \in S$ to the image of $f(p) \in \tilde D$ along the projection $\tilde D \to D$, for any $p \in P_x$ (the image does not depend on the choice of $p \in P_x$, because $f(p)$ runs through a single $G$-orbit of $\tilde D$). To verify that ($*$) in \cref{rslt:mono-after-sheafification-checked-with-+-construction} is satisfied, we fix $P$ and $P'$ such that $\alpha'(P) = \alpha'(P')$. Thus for every $x \in S$, $f$ and $f'$ map the fibers $P_x$ and $P'_x$ to the same $G$-orbit in $\tilde D$. 

Let $V \subseteq S$ be the open subset of those $x \in S$ such that the fibers $P_x$ and $P'_x$ do not get mapped to $0 \in \tilde D$. Then $V$ is closed under specializations, because $f$ and $f'$ preserve generalizations (like any continuous map) and $0 \in \tilde D$ does not have any generalizations. Hence $V$ corresponds to the open subset $\pi_0 V \subseteq \pi_0 S$. Let $U_1 \subseteq S$ be the closure of $V$ in $S$. Then $U_1$ corresponds to the closure of $\pi_0 V$ in $\pi_0 S$, which is open because $\pi_0 S$ is extremally disconnected. Consequently, denoting $U_2 = S \setminus U_1$ we have constructed a decomposition $S = U_1 \sqcup U_2$ into disjoint open subsets.

It remains to show that the decomposition $S = U_1 \sqcup U_2$ satisfies the condition in ($*$) for our given $P$ and $P'$. By construction, the restriction of $f$ and $f'$ to $U_2$ factors over $0 \in \tilde D$, hence $\Hom(P|_{U_2}, P'|_{U_2}) = *$ by the above description of $Y'$. Hence it only remains to see that $\Hom(P|_{U_1}, P'|_{U_1}) = *$. To simplify notation we now replace $U_1$ by $S$, i.e., we assume that $V$ is dense in $S$. Since $f$ and $f'$ do not factor over $0$, the maps $P \to P'$ in $Y'(S)$ are the same as the maps $P \to P'$ in $(\tilde D/G)(S)$. In summary, we have reduced the statement to the following key claim:
\begin{itemize}
    \item[($**$)] Fix $(P, f), (P', f') \in (\tilde D/G)(S)$, i.e., $P$ and $P'$ are $G$-torsors over $S$ and $f\colon P \to \tilde D$ and $f'\colon P \to \tilde D$ are $G$-equivariant maps. Assume that the induced maps $S \to D$ (given by sending $x \in S$ to the image point of $f(P_x)$ resp.\ $f'(P_x)$ under $\tilde D \to D$) agree and that there is a dense open subset $V \subseteq S$ such that $f(P_V)$ does not contain $0 \in \tilde D$. Then there is a unique isomorphism $P \isoto P'$ in $(\tilde D/G)(S)$.
\end{itemize}
To prove this claim, note that since every pro-étale cover of $S$ splits, we have $P \isom P' \isom S \times G$ as $G$-torsors over $S$ (disregarding the maps $f$ and $f'$). Thus, maps of $G$-torsors $P \to P'$ coincide with continuous maps $\varphi\colon \abs S \to G$. We need to see that there is exactly one such $\varphi$ that is compatible with $f$ and $f'$. This is clear on $V \subseteq S$ because then $f$ and $f'$ factor through $\tilde D \setminus \{ 0 \}$, where the $G$-action is free. Let us denote $\varphi_V\colon \abs V \to G$ the resulting map. Any extension of $\varphi_V$ to $\abs S$ (as a continuous map) produces a map $\varphi$ compatible with $f$ and $f'$ (i.e., with $f = \varphi f'$), as the locus where $f$ and $\varphi f'$ agree is closed, $G$-equivariant and contains $P_V$. We thus need to see that $\varphi_V$ admits exactly one extension to $\abs S$. By the density of $V$ in $S$, there is at most one extension. On the other hand, continuous maps $\abs S \to G$ coincide with continuous maps $\pi_0 S \to G$ (and similarly for $V$). Now $\varphi_V$ extends to a continuous map $\overline{\varphi_V}\colon \beta(\pi_0 V) \to G$, where $\beta(\pi_0 V)$ denotes the Stone-\v{C}ech compactification of the set $\pi_0 V$ (equipped with the discrete topology). The inclusion $\pi_0 V \injto \pi_0 S$ induces a cover $\beta(\pi_0 V) \surjto \pi_0 S$, which splits because $\pi_0 S$ is extremally disconnected. We thus obtain the desired map $\varphi\colon \pi_0 S \to \beta(\pi_0 V) \to G$. Note that $\varphi$ indeed extends $\varphi_V$: This reduces to the observation that the restriction of $\overline{\varphi_V}$ to $\beta(\pi_0 V) \times_{\pi_0 S} \pi_0 V$ coincides with the composition of the projection to $\pi_0 V$ and $\varphi_V$ -- but both maps $\beta(\pi_0 V) \times_{\pi_0 S} \pi_0 V \to G$ are compatible with $\varphi_V$ after restriction to $\pi_0 V$, and there can be at most one such map because $\pi_0 V$ is dense in the source.
\end{proof}

The pushout formula in \cref{rslt:disk-as-pushout-from-perfectoid-disk} provides us with an easy way to access the cohomology of the disk:

\begin{cor}
\label{cor:fibersequenceD}
In the setup of \cref{rslt:disk-as-pushout-from-perfectoid-disk}, there is a canonical fiber sequence
\begin{align*}
    R\Gamma(D, \OO^{\flat,+,a}) \to \OO^{\flat,+,a}(\tilde D)^{hG} \oplus \OO^a_{C^\flat} \to (\OO^a_{C^\flat})^{hG},
\end{align*}
where in the third term $G$ acts trivially on $\OO^a_{C^\flat}$.
\end{cor}
\begin{proof}
The assignment $X \mapsto R\Gamma(X, \OO^{\flat,+,a})$ defines a $\mathcal D_\solid^a(\OO_{C^\flat})$-valued $\infty$-sheaf on the category of small v-stacks over $\Spa C$. In particular this assignment sends colimits of stacks to limits of modules, hence by \cref{rslt:disk-as-pushout-from-perfectoid-disk} we have a pullback square
\begin{equation*}\begin{tikzcd}
    R\Gamma(D, \OO^{\flat,+,a}) \arrow[r] \arrow[d] & R\Gamma(*, \OO^{\flat,+,a}) \arrow[d]\\
    R\Gamma(\tilde D/G, \OO^{\flat,+,a}) \arrow[r] & R\Gamma(*/G, \OO^{\flat,+,a})
\end{tikzcd}\end{equation*}
in $\mathcal D_\solid^a(\OO_{C^\flat})$. The claimed fiber sequence follows easily from this pullback square.
\end{proof}

\cref{cor:fibersequenceD} reduces the study of $R\Gamma(D,\OO^{\flat,+,a})$ to the study of the continuous cohomology $H^*(G,\OO^{+,a}(\tilde{D}^\flat))$.  As an initial reduction, the inclusion $\OO_{C^\flat}\to \OO^+(\tilde{D}^\flat)$ has a section, namely evaluation at 0, so let us write:
\[ \OO^+(\tilde{D}^\flat)=\OO_C\oplus \OO^+(\tilde{D}^\flat)^0, \] 
where $\OO^+(\tilde{D}^\flat)^0$ consists of those fractional power series with no constant term.

We will need to know how this cohomology behaves under the restriction maps between discs of different radii.  Therefore let us fix some notations: for $r=\abs{\rho}>0$ in the value group of $C$, let $D_r$ be the disc of radius $r$:
\[ D_r = \Spa C\langle {T/\rho} \rangle. \]
When $0<r'\leq r$, the inclusion $D_{r'}\to D_r$ lifts to a $G$-equivariant morphism $\tilde{D}_{r'}\to \tilde{D}_r$. 

\begin{lemma} \label{LemCohoD} Concerning the continuous cohomology $H^*_{\cts}(G,\OO^{+}(\tilde{D}_r^\flat))$, we have:
\begin{enumerate}
    \item $H^0_{\cts}(G,\OO^{+}(\tilde{D}_r^\flat))=\OO_C$.
    \item Let $0<r'\ll r$.  Then the map
    \[ H^1_{\cts}(G,\OO^+(\tilde{D}_r^\flat)^0) \to H^1_{\cts}(G,\OO^+(\tilde{D}_{r'}^\flat)^0) \] is zero.
    \item $H^i_{\cts}(G,\OO^+(\tilde{D}_r^\flat))=0$ for $i\neq 0,1$.
\end{enumerate}
\end{lemma}

\begin{proof} It suffices to handle the case $r=1$, in which case we write $D=D_1$.  Let us fix an isomorphism $G=\Z_p(1)\isom \Z_p$, which is to say a topological generator 
\[ e=(1,\zeta_p,\zeta_{p^2},\dots)\in \Z_p(1).\] 
Then $\OO^{+}(\tilde{D}^\flat)^{hG}$ is computed by the complex
\[ \OO^{+}(\tilde{D}^\flat) \stackrel{e-1}{\to} \OO^{+}(\tilde{D}^\flat), \]
and in particular the cohomology is 0 for $i\neq 0,1$. 

We can be completely explicit about the action of $e$ on $\OO^+(\tilde{D}^\flat)$.  We have
\[ \OO^+(\tilde{D}^\flat) \isom \OO_{C^\flat}\langle U^{1/p^\infty} \rangle, \]
where $U=(T,T^{1/p},T^{1/p^2},\dots)$. 
The sequence defining $e$ determines an element  $\varepsilon\in C^\flat=\varprojlim_{x\mapsto x^p} C$.  This element is close to 1;  let 
\[ t = \eps-1,\]
so that $t$ is a pseudo-uniformizer of $C$.  The action of our generator $e$ on $\OO^+(\tilde{D}^\flat)$ is determined by what it does to $U$, which is:
\[ e(U)=\eps U = (1+t)U.\]
Therefore $e-1$ has this effect on $\OO^+(\tilde{D}^\flat)$:
\[ U^\alpha\mapsto [(1+t)^\alpha-1]U^\alpha, \; \alpha\in \Z[1/p]_{\geq 0}. \] 
Evidently the kernel of $e-1$ on $\OO^+(\tilde{D}^\flat)$ is $\OO_{C^\flat}$, which is part (1).
For part (2):  a class $[f(U)]$ in $H^1_{\cts}(G,\OO^+(\tilde{D}^\flat)^0)$ is represented by a series $f(U)=\sum_\alpha c_\alpha U^\alpha$, where $\alpha$ ranges through $\Z[1/p]_{>0}$, and $\abs{c_\alpha}\to 0$ as $\alpha\to\infty$, in the sense that for each $\varepsilon > 0$, there are at most finitely many $\alpha$ such that $\abs{c_\alpha} \geq \varepsilon$.  We claim that the power series
\[ g(U)=\sum_\alpha \frac{c_\alpha}{(1+t)^{\alpha}-1} U^\alpha \] 
converges on $\tilde{D}_{\abs{t}}$.  This would prove the desired claim, since the image of $[f(U)]$ under $H^1_{\cts}(G,\OO^+(\tilde{D}^\flat)^0)\to H^1_{\cts}(G,\OO^+(\tilde{D}_{\abs{t}}^\flat)^0)$ equals $[(e-1)g(U)]=0$.  
A computation shows that if $\alpha=bp^k$ where $p\nmid b$, then 
\[ \abs{(1+t)^\alpha-1}=\abs{t}^{p^k},\] 
so that 
\[ \abs{\frac{c_\alpha}{(1+t)^{\alpha}-1}}\abs{t}^\alpha = \abs{c_\alpha} \abs{t}^{\alpha-p^k} =\abs{c_\alpha} \abs{t}^{(b-1)p^k} \leq \abs{c_\alpha}\to 0, \] 
proving the claim.  
\end{proof}

\begin{lemma}\label{LemRGammaOFlat} Regarding the pro-\'etale cohomology $H^*(D,\OO^{\flat,+,a})$ of affinoid discs, we have:
\begin{enumerate}
    \item $H^0(D,\OO^{\flat,+,a})=\OO_{C^\flat}^{a}$.
    \item For a disc $D'\subset D$ of sufficiently small radius, the restriction map $H^1(D,\OO^{\flat,+,a})\to H^1(D', \OO^{\flat,+,a})$ is 0.
    \item $H^i(D,\OO^{\flat,+,a})=0$ for $i\neq 0,1$.
\end{enumerate}
\end{lemma}

\begin{proof} The claims about $H^i(D,\OO^{\flat,+,a})$ for $i\neq 1$ follow easily by combining \cref{cor:fibersequenceD} and \cref{LemCohoD}.  Whereas in the $i=1$ case, \cref{cor:fibersequenceD} shows that 
\[ H^1(D,\OO^{\flat,+,a})\isom H^1_{\cts}(G,\OO^{+,a}(\tilde{D}^\flat)^0), \] 
which has the claimed restriction behavior by \cref{LemCohoD}(2). 
\end{proof}

\subsection{Computation of $R\Gamma(\HH_{F,C}^{1,\Diamond},\OO^\flat)$}\label{sec:cohomology-of-Drinfeld-space}

Here we compute the cohomology $R\Gamma(\HH^{1,\Diamond}_{F,C},\OO^\flat)$.  The plan is to first compute the compactly supported cohomology $R\Gamma_c(\HH^{1,\Diamond}_{F,C},\OO^{\flat})$ and then to apply Poincar\'e duality.

\begin{thm} \label{thm:OflatCompactlyDrinfeld} 
Let $F$ be a finite extension of $\Q_p$, and let $C$ be an algebraically closed perfectoid field containing $F$.  There are isomorphisms of $C^\flat$-modules with $\GL_2(F)$-action:
\begin{equation*}
    H^i_c(\HH^1_{F,C},\OO^{\flat})\simeq  
\begin{cases}
   \St(C^{\flat}), & i = 1\\
   C^\flat, & i = 2 \\
   0, & \text{otherwise.}
\end{cases}
\end{equation*}
Here $\St(C^{\flat})\coloneqq \LC(\PP^1(F),C^\flat)/C^\flat $ is the Steinberg representation of $\GL_2(F)$ (cf. \cite[\S 4.1.1]{CDNIntegral}), and $\LC(-,-)$ denotes the space of locally constant functions.
\end{thm}

\begin{proof}
We apply the excision sequence in \cref{rslt:excision-for-O-flat} to the situation of the rigid-analytic space $X=\PP^1_C$, its open subspace $U=\HH^{1,\Diamond}_{F,C}$, and the complement $Z=\PP^1(F)$. We choose a cofiltered system $\{U_i\}_{i\in \Z_{>0}}$ such that 
\begin{itemize}
    \item $U_i$ is a descending family of open affinoids in $X=\PP^1_C$, with $\bigcap_i U_i=\PP^1(F)$,
    \item each $U_i$ is a disjoint union of finitely many closed discs, and 
    \item each inclusion $U_{i+1}\to U_i$ restricts on each connected component to an inclusion of discs $D'\to D$ such that the radius of $D'$ is \textit{sufficiently small} compared to that of $D$ (This is made precise in the proof of \cref{LemRGammaOFlat}). 
\end{itemize}
To see such a family $U_i$ exists, one can simply put a closed disc $D(a)$ of small radius $r$ at each $a\in \PP^1(F)\subset \PP^1_C$. For two rational points $a,b\in\PP^1(F)$, the discs $D(a)$ and $D(b)$ are either the same or disjoint. Furthermore, one only needs finitely many such discs to cover $\PP^1(F)$ since $D(a)\cap \PP^1(F)$ is open in the compact space $\PP^1(F)$. One may then take $U_i$ to be the union of these finitely many discs.

Applying \cref{rslt:excision-for-O-flat}, we find a fiber sequence in $\mathcal{D}_{\solid}(C^\flat)$:
    \begin{equation*}
        R\Gamma_c(\HH^{1,\Diamond}_{F,C},\OO^\flat)\to R\Gamma(\PP^1_C,\OO^\flat)\to \colim_i R\Gamma(U_i,\OO^\flat),
    \end{equation*}
By \cref{LemRGammaOFlat} we have:
\begin{enumerate}
    \item $H^0(U_i,\OO^{\flat,+,a})\isom C_{\cts}(\pi_0(U_i),\OO_{C^\flat}^a)$.
    \item The image of $H^1(U_i,\OO^{\flat,+,a}) \to H^1(U_{i+1},\OO^{\flat,+,a})$ is 0.
    \item $H^i(U_i,\OO^{\flat,+,a})=0$ for $i\neq 0,1$.
\end{enumerate}
Therefore we have an isomorphism of $\OO_{C^\flat}^a$-modules:
\begin{equation*}
    \colim_i H^j(U_i,\OO^{\flat,+,a}) =\begin{cases}
        \LC(\PP^1(F), \OO_{C^\flat}^a) & j=0,\\ 
        0 &\text{otherwise}.
    \end{cases} 
\end{equation*}
After base changing from $\OO_{C^\flat}^a$ to $C^\flat$, we obtain an isomorphism of $C^\flat$-modules:
\begin{equation*}
    \colim_i H^j(U_i,\OO^{\flat}) = \begin{cases}
        \LC(\PP^1(F), C^\flat) & j=0,\\ 
        0 &\text{otherwise}.
    \end{cases}
\end{equation*}
Subsituting this calculation into our exact triangle above, we obtain a long exact sequence of $C^\flat$-vector spaces with $\GL_2(F)$-action (here we suppress the $\OO^\flat$-coefficients)
\[ 0 \to H^0_c(\HH^{1,\Diamond}_{F,C})\to H^0(\PP^1_C)\to 
\LC(\PP^1(F),C^\flat) \to H^1_c(\HH^{1,\Diamond}_{F,C}) \to 0, \] 
as well as isomorphisms $H^i_c(\HH^{1,\Diamond}_{F,C})\isom H^i(\PP^1_C)$ for $i\geq 2$. 

Note $R\Gamma(\PP^1_C, \OO^{\flat}) \simeq C^{\flat} \oplus C^{\flat}[-2](-1)$ as a consequence of the Primitive Comparison Theorem (\cite[Corollary~3.9.24]{MannThesis}) and the fact that for $f\colon\PP^1_C\to\Spa C$, the functor $f_*=f_!$ commutes with limits and colimits.  For any $i$, the map $H^0(\PP^1_C)\to H^0(U_i)$ is the restriction of the constant global functions $C^\flat$ on $\PP^1_C$ to $U_i$, and the map is injective. Thus, the map $H^0(\PP^1_C)\to \colim_i H^0(U_i)$ is injective. The isomorphisms claimed in the theorem follow from these observations and from $H^2(\PP^1_C)\isom C^\flat$.  
\end{proof}

Consider a cohomologically smooth map $f:X\to Y$ of small v-stacks. The functor $f^*$ has a left adjoint $f_{\natural}$. For any $M\in \mathcal{D}_{\hat\solid}(\OO^\flat_X)$, we define the \textit{relative homology} of $M$ to be $f_{\natural}M\in \mathcal{D}_{\hat\solid}(\OO^\flat_Y)$.
Applying this to $f\colon\HH^1_{F,C}\to \Spd C$ and $M=\OO^\flat$, we get by adjunction that 
\begin{equation}
     f_*\OO^\flat\simeq \underline{\Hom}(f_{\natural}\OO^\flat,\OO^\flat).
\end{equation}
We write $H_i(\HH^1_{F,C},\OO^\flat)\coloneqq H^{-i}(f_{\natural}\OO^\flat)$.
\begin{cor}\label{cor:OFlat-Homology-Drinfeld}
    There are isomorphisms of solid $C^\flat$-modules with $\GL_2(F)$-action:
    \[ 
    H_i(\HH^1_{F,C},\OO^\flat)=\begin{cases}
        C^\flat, & i=0,\\
        \St(C^\flat), & i=1,\\
        0, & \text{otherwise.}
    \end{cases}
    \]
\end{cor}
\begin{proof}
    By the adjunction between $f_\natural$ and $f^*$ and the Poincar\'e duality for smooth rigid-analytic spaces (Theorem~\ref{thm:6FunctorsOflat}(iv)), we have 
    \begin{equation*}
        f_{\natural}\OO^\flat\simeq (f_!\OO^\flat)[2].
    \end{equation*}
    Here we omit the Tate twist since it will not be used in subsequent discussions. The desired result then follows from Theorem~\ref{thm:OflatCompactlyDrinfeld}.
\end{proof}
\begin{rmk}
    The homology group $H_i(\HH^1_{F,C},\OO^\flat)$ carries actions of three groups: $\GL_2(F)$, $\Gal(C/F)$, and $\varphi_q^{\Z}$. Here $\varphi_q$ is the $q$-Frobenius on $\Spd C$ and will be considered in \cref{sec:DrinfeldsLemma}. The $\Gal(C/F)$-action comes from the fact that $\HH^1_{F,C}$ is the base change of $\HH^1_{F}$ from $\Q_p$ to $C$, and the Tate twist mentioned in the proofs of \cref{thm:OflatCompactlyDrinfeld} and \cref{cor:OFlat-Homology-Drinfeld} refers to the $\Gal(C/F)$-action via cyclotomic characters. The $\varphi_q^{\Z}$-action commutes with the $\Gal(C/F)$-action, and in fact, we will not use the $\Gal(C/\Q_p)$-action in any subsequent discussion. There is a characteristic-$p$ Tate twist showing up in \cref{DefnTateTwist} and results thereafter, which should not be confused with the usual Tate twist. In fact, the characteristic-$p$ Tate twist comes from the $\GL_2(F)$-action, which also commutes with the $\Gal(C/F)$-action.
\end{rmk}
\begin{defn}\label{def:distributions}
    The space $\Dist(\PP^1(F),C^\flat)$ of distributions is defined as the solid $C^\flat$-linear dual of $\LC(\PP^1(F),C^\flat)$, the space of locally constant functions. Let $\Dist(\PP^1(F),C^\flat)^0$ denote the subspace of distributions which assign a value of $0$ to the constant function $1$. In other words, $\Dist(\PP^1(F),C^\flat)^0$ is defined by the short exact sequence
    \begin{equation*}
        0\to \Dist(\PP^1(F),C^\flat)^0\to \Dist(\PP^1(F),C^\flat)\xrightarrow{\varphi\mapsto \int_{\PP^1(F)}\varphi} C^\flat\to 0.
    \end{equation*}
    Note that $\Dist(\PP^1(F),C^\flat)^0$ is the solid $C^\flat$-linear dual of $\St(C^\flat)=\LC(\PP^1(F),C^\flat)/C^\flat$. We thus also write
    \begin{equation*}
        \St(C^\flat)^*=\Dist(\PP^1(F),C^\flat)^0.
    \end{equation*}
\end{defn}

\begin{thm}\label{thm:OFlatCohomDrinfeld}
 There are isomorphisms of solid $C^\flat$-modules with $\GL_2(F)$-action:
    \[ 
    H^i(\HH^1_{F,C},\OO^\flat)=\begin{cases}
        C^\flat, & i=0,\\
        \St(C^\flat)^*, & i=1,\\
        0, & \text{otherwise.}
    \end{cases}
    \]
\end{thm}
\begin{proof}
    In this proof, we will deviate from our earlier conventions and write $R\underline{\Hom}(-,-)$ for the derived internal Hom (in any appropriate category). So $\underline{\Hom}(-,-)$ will denote the underived internal Hom.
    Let $f\colon\HH^1_{F,C}\to \Spd C$ be the structure map. 
    By the adjunction between $f_\natural$ and $f^*$, we get
\begin{equation*}
    f_*\OO^\flat\simeq R\underline{\Hom}(f_{\natural}\OO^\flat,\OO^\flat).
\end{equation*}
    As shown in the proof of Theorem~\ref{thm:OflatCompactlyDrinfeld}, $R\Gamma_c(\HH^1_{F,C},\OO^\flat)$ is quasi-isomorphic to 
    \[
    R\Gamma(\PP^1_{C},\OO^\flat)[-1]\to \LC(\PP^1(F),C^\flat)[-1].
    \]
    Combined with the fact $f_{\natural}\OO^\flat\simeq (f_!\OO^\flat)[2]$, this implies $R\underline{\Hom}(f_{\natural}\OO^\flat,\OO^\flat)$ is isomorphic to 
    \begin{equation*}
        R\Hom(\LC(\PP^1(F),C^\flat),C^\flat)[-1]\to R\Hom(R\Gamma(\PP^1_{C},\OO^\flat),C^\flat)[-1]
    \end{equation*}
    For the second term, we have 
    \[
        R\Hom(R\Gamma(\PP^1_{C},\OO^\flat),C^\flat)[-1]\simeq R\Hom(C^\flat[0]\oplus C^\flat[-2],C^\flat)[-1]\simeq C^\flat[-1]\oplus C^\flat[1].
    \]
    By \cite[Theorem~5.4]{ClausenScholzeCondensed} (after base-change from $\Z$ to $C^\flat$), $\LC(\PP^1(F),C^\flat)$ is a free solid $C^\flat$-module. Hence, 
    \[
R\Hom(\LC(\PP^1(F),C^\flat),C^\flat)=\Hom(\LC(\PP^1(F),C^\flat),C^\flat)=\Dist(\PP^1(F),C^\flat).
    \]
    The desired statement for $H^i(\HH^1_{F,C},\OO^\flat)$ then follows.
\end{proof}

\subsection{The Heyer--Jacquet module}\label{sec:Heyer-Jacquet}

In this subsection we explain how to compute the cohomology of the Steinberg representation $\St(C^\flat) = \LC(\PP^1(F), C^\flat)/C^\flat$ (or rather its dual). Since this is essentially the smooth Steinberg representation (as we use locally constant functions $\PP^1(F) \to C^\flat$) we would expect the cohomology to be the one in the setting of smooth representations, which was previously computed by Heyer in \cite[\S4.3]{Heyer-GeometricLemma}. However, some care is required as the relevant 6-functor formalism for us is $\mathcal D_{\hat\solid}(\OO_X^\flat)$ and this 6-functor formalism computes continuous group cohomology, taking the topology of $C^\flat$ into account. In the following we show how smooth and continuous group cohomology compare.

Given a ring $\Lambda$, recall the 6-functor formalism $\mathcal D(-,\Lambda)$ for sheaves of $\Lambda$-modules on condensed anima from \cite[Definition~3.5.17]{HeyerMann6FF}. Here the category $\Cond(\Ani)$ of condensed anima is the category of sheaves of anima (i.e., $\infty$-stacks) on the site of profinite sets, where we implicitly impose a fixed cardinality bound on the profinite sets. Note that there is a functor $\ProFin \to \vStk$ sending a profinite set $S$ to the small v-stack $S \times \Spa C^\flat = \Spa(\Cts(S, C^\flat))$. This functor preserves covers, hence we can extend it by colimits in order to obtain a functor
\begin{align*}
    \Cond(\Ani) \to \vStk, \qquad X \mapsto X \times \Spa C^\flat.
\end{align*}
Since this functor commutes with colimits by construction, it will send the classifying stack $*/G$ of a locally profinite group $G$ to the classifying stack of $G$ viewed as a v-stack. Precomposing $\mathcal D_{\hat\solid}$ with the above functor, we obtain the 6-functor formalism $\mathcal D_{\hat\solid}(-)$ on $\Cond(\Ani)$ that sends a profinite set $S$ to
\begin{align*}
    \mathcal D_{\hat\solid}(S) = \mathcal D_{\hat\solid}(\OO_{S \times \Spa C^\flat}) = \mathcal D_{\hat\solid}(\Cts(S, C^\flat)).
\end{align*}
The following result allows us to compare the 6-functor formalism from \cite{HeyerMann6FF} (for $\Lambda = C^\flat$, viewed as a discrete ring) with the 6-functor formalism $\mathcal D_{\hat\solid}$.

\begin{lemma} \label{rslt:trafo-from-smooth-reps-to-D-solid-hat}
There is a natural transformation
\begin{align*}
    \alpha\colon \mathcal D(-, C^\flat) \to \mathcal D_{\hat\solid}(-)
\end{align*}
of functors $\Cond(\Ani)^\op \to \CAlg(\Cat)$. In other words, for every condensed anima $X$ there is a functor $\alpha_X\colon \mathcal D(X, C^\flat) \to \mathcal D_{\hat\solid}(\OO_{X \times \Spa C^\flat})$ compatible with pullback and tensor products.
\end{lemma}
\begin{proof}
Both $\mathcal D(-, C^\flat)$ and $\mathcal D_{\hat\solid}(-)$ satisfy descent and are thus the right Kan extensions of their restriction to the category $\ProFin^\op$ of profinite sets. It is therefore enough to construct $\alpha$ on these restrictions. For a profinite set $S$ we compute
\begin{align*}
    \mathcal D(S, C^\flat) &= \mathcal D(\LC(S, C^\flat)) = \Mod_{\LC(S, C^\flat)}(\mathcal D(C^\flat)),\\
    \mathcal D_{\hat\solid}(S) &= \mathcal D_{\hat\solid}(\Cts(S, C^\flat)) = \Mod_{\Cts(S, C^\flat)}(\mathcal D_{\hat\solid}(C^\flat))
\end{align*}
(for the last identity see \cite[Proposition~2.21(ii)]{ALBM-6ff-on-FF-curve} or note that in the case $C = \C_p$ we have $\mathcal D_{\hat\solid}(C^\flat) = \mathcal D_\solid(C^\flat)$ by \cite[Lemma~2.9]{ALBM-6ff-on-FF-curve}). Here $\mathcal D(C^\flat)$ denotes the derived category of the category of (classical) $C^\flat$-vector spaces. Note that $\mathcal D_{\hat\solid}(C^\flat)$ is a $C^\flat$-linear category, hence there is a uniquely determined $C^\flat$-linear symmetric monoidal embedding $\mathcal D(C^\flat) \injto \mathcal D_{\hat\solid}(C^\flat)$. Moreover, for every profinite set $S$ there is a natural map $\LC(S, C^\flat) \to \Cts(S, C^\flat)$ of rings in $\mathcal D_{\hat\solid}(C^\flat)$, compatible with maps between profinite sets (as these rings live in the heart of the t-structure, there are no non-trivial homotopies to verify here). We can thus define $\alpha_S$ as the composition
\begin{align*}
    \alpha_S\colon \Mod_{\LC(S, C^\flat)}(\mathcal D(C^\flat)) \to \Mod_{\LC(S, C^\flat)}(\mathcal D_{\hat\solid}(C^\flat)) \to \Mod_{\Cts(S, C^\flat)}(\mathcal D_{\hat\solid}(C^\flat)),
\end{align*}
and this clearly defines the desired natural transformation.
\end{proof}

\begin{rmk}
We warn the reader that the natural transformation $\alpha$ from \cref{rslt:trafo-from-smooth-reps-to-D-solid-hat} is \emph{not} compatible with lower-$!$ functors in general. In fact, this already fails in the case of the map $f\colon S \to *$ for a profinite set $S$: We have $\alpha f_! (\LC(S, C^\flat)) = \LC(S, C^\flat)$ and $f_! \alpha(\LC(S, C^\flat)) = f_! \Cts(S, C^\flat) = \Cts(S, C^\flat)$, and these two rings do not agree unless $S$ is finite. We expect that there is a lax compatibility of $\alpha$ with !-functors, in the sense of a natural transformation $\alpha f_! \to f_! \alpha$ for all $!$-able maps $f$. Making this precise requires more $(\infty,2)$-categorical machinery than we intend to use in this paper, so we content ourselves with a more ad-hoc approach in the following.
\end{rmk}

Fix a $p$-adic Lie group $G$ and let $V$ be a (derived) smooth $G$-representation on $C^\flat$-vector spaces, i.e., an object in $\mathcal D(*/G, C^\flat)$ (cf. \cite[Proposition~5.3.10]{HeyerMann6FF}). Via $\alpha$ we can view $V$ also as a continuous $G$-representation with respect to the induced topology on $V$ from $C^\flat$. Denoting $f\colon */G \to *$, we now get two pushforward functors:
\begin{itemize}
    \item $f_*\colon \mathcal D(*/G, C^\flat) \to \mathcal D(C^\flat)$ computes the smooth $G$-cohomology, i.e., $f_* V$ is given by a complex with terms $\LC(G^n, V)$.
    
    \item $f_*\colon \mathcal D_{\hat\solid}(*/G) \to \mathcal D_{\hat\solid}(C^\flat)$ computes continuous group cohomology, i.e., $f_*(\alpha V)$ is given by a complex with terms $\Cts(G^n, \alpha(V))$, where $\alpha(V) = V \otimes_{C^{\flat,\mathrm{disc}}} C^\flat$ is the condensed $C^\flat$-module obtained from $V$ by writing $V$ as a colimit of finitely generated $C^\flat$-modules and putting the canonical topology on each.
\end{itemize}
The following result shows that the above two $G$-cohomologies agree when $G$ is compact. Moreover, instead of $f_*$ one can also consider the \emph{left} adjoint of $f^*$, usually denoted $f_\natural$. This functor computes $G$-\emph{homology} and has the advantage of being compatible with $\alpha$ without the compactness assumption on $G$.

\begin{proposition}\label{prop:group-homology-compatible}
Let $\varphi\colon G \to H$ be a map of $p$-adic Lie groups and $f\colon */G \to */H$ the associated map of condensed anima.
\begin{propenum}
    \item \label{rslt:smooth-group-homology-compatible-with-D-solid} Assume that $\varphi$ is either (1) surjective or (2) an open immersion. Then $f_\natural$ exists in both of the above sheaf theories and the natural map
    \begin{align*}
        f_\natural \alpha \isoto \alpha f_\natural
    \end{align*}
    is an isomorphism of functors $\mathcal D(*/G, C^\flat) \to \mathcal D_{\hat\solid}(*/H)$.

    \item \label{rslt:smooth-group-cohomology-compatible-with-D-solid} Assume that $\varphi$ is surjective and its kernel is compact with finite $p$-cohomological dimension. Then the natural map
    \begin{align*}
        \alpha f_* \isoto f_* \alpha
    \end{align*}
    is an isomorphism of functors $\mathcal D(*/G, C^\flat) \to \mathcal D_{\hat\solid}(*/H)$.
\end{propenum}
\end{proposition}
\begin{proof}
We first prove (i), so assume that $f$ is as in the claim. Note that $f$ is cohomologically smooth for both 6-functor formalisms. Namely, this can be checked after base-change along $* \to */H$. In the second case of $\varphi$ we then end up with the map $H/G \to *$ and $H/G$ is a discrete set, which is even cohomologically étale in both 6-functor formalisms. In the first case of $\varphi$, we end up with the map $*/(\ker\varphi) \to *$. Then cohomological smoothness for $\mathcal D(-, C^\flat)$ is proved in \cite[Example~5.3.22]{HeyerMann6FF} and for $\mathcal D_{\hat\solid}(-)$ in \cite[Theorems~3.16,~3.18]{HansenMannJL} together with \cite[Lemma~3.3.2]{ALBM-6ff-on-FF-curve}. The cohomological smoothness implies that $f^*$ preserves limits and hence $f_\natural$ exists. From the compatibility of $\alpha$ with pullbacks we obtain the natural transformation $f_\natural \alpha \to \alpha f_\natural$. It remains to show that it is an isomorphism.

It is enough to show the claimed isomorphism after pullback along the map $* \to */H$ (since this pullback is conservative). Note that $f_\natural$ commutes with this pullback in both 6-functor formalisms: This follows formally from the fact that $f$ is cohomologically smooth, since after passing to right adjoints the claimed base-change is stated in \cite[Lemma~4.5.13(i)]{HeyerMann6FF}. We can thus replace $f$ by its base-change along $* \to */H$. In the case that $\varphi$ is an open immersion, we now deal with the map $H/G \to *$ and $H/G$ is a discrete set -- the claim is clear. We thus restrict to the case that $\varphi$ is surjective and hence are allowed to assume $H = *$ (by replacing $G$ with $\ker\varphi$), i.e., $f$ is now the map $*/G \to *$. Since both sides of the claimed isomorphism of functors commute with colimits, it is enough to verify the claim on compact generators. We are thus left to show that the natural map
\begin{align*}
    f_\natural \alpha (V) \to \alpha f_\natural (V)
\end{align*}
is an isomorphism for all $V$ of the form $V = g_! 1$ for some $p$-torsion free open compact subgroup $K \subseteq G$ with induced map $g\colon */K \to */G$ (here $1$ denotes the trivial representation on $C^\flat$). Note that $g$ is étale and hence $g_! = g_\natural$. By the case of open immersions for $\varphi$ we can commute $g_\natural$ past $\alpha$. Thus, replacing $G$ by $K$ we are reduced to the following key claim:
\begin{itemize}
    \item[($*$)] Suppose $G$ is $p$-torsion free and compact, and let $f\colon */G \to *$. Then the natural map $f_\natural 1 \to \alpha f_\natural 1$ is an isomorphism in $\mathcal D_{\hat\solid}(C^\flat)$.
\end{itemize}
Here is the intuition behind this claim: Since $f$ is cohomologically smooth and proper for both 6-functor formalisms, $f_\natural$ essentially computes group cohomology (up to shift and twist). Thus the claimed isomorphism states that the smooth group cohomology of the trivial $G$-representation (given by a complex with terms $\LC(G^n, C^\flat)$) coincides with the continuous group cohomology (given by a complex with terms $\Cts(G^n, C^\flat)$).

To prove the claim, we make further reductions. First, note that there is an integral version of both 6-functor formalisms, with coefficients in $\OO_{C^\flat}^a$ and $\OO^{\flat,+,a}_X$ respectively (the latter is the source of all 6-functor formalisms in \cite{ALBM-6ff-on-FF-curve}). Since all functors in the claimed isomorphism commute with colimits and hence with the base-change $- \otimes_{\OO_C^{\flat a}} C^\flat$, it is enough to prove the isomorphism in ($*$) with integral coefficients. It clearly holds modulo a pseudouniformizer $\pi$, as in this case $\alpha$ is indeed a morphism of 6-functor formalisms (this is easy to see on the full subcategory of profinite sets, where lower-$!$ is right adjoint to pullback, and then follows on all condensed anima by the uniqueness of the extension results in \cite[Theorem~3.4.11]{HeyerMann6FF}). It is thus enough to show that both $f_\natural 1$ and $\alpha f_\natural 1$ are $\pi$-adically complete. We use cohomological properness and smoothness of $f$ (in both 6-functor formalisms; for the smooth one we crucially use that $C^\flat$ has characteristic $p$ in order to apply \cite[Example~5.3.22]{HeyerMann6FF} to $\OO_C^\flat$) to compute
\begin{align*}
    f_\natural 1 = f_! f^! 1 = f_* \omega_f,
\end{align*}
where $\omega_f\coloneqq  f^! 1$ denotes the dualizing complex -- this is a character shifted by the dimension of $G$ (see \cite[Example~5.3.22]{HeyerMann6FF} for $\mathcal D(-, \OO_C^{\flat,a})$ and \cite[Theorem~3.18]{HansenMannJL} together with \cite[Lemma~3.3.2]{ALBM-6ff-on-FF-curve} for $\mathcal D_{\hat\solid}$) Since $f_*$ preserves limits and $\omega_f$ is invertible, clearly $f_* \omega_f$ is $\pi$-complete in both 6-functor formalisms, hence the same is true for $f_\natural 1$. It remains to show that $\alpha f_\natural 1$ is $\pi$-complete. But note that $f_\natural 1$ is perfect (as follows from the above computation and the fact that $\omega_f$ is admissible (cf. \cite[Definition 5.3.12, Proposition 5.3.14]{HeyerMann6FF}), hence a finite colimit of copies of $\OO_C^{\flat a}$, and the same is then true for $\alpha f_\natural 1$. But such objects are obviously $\pi$-complete. This finishes the proof of (i).

Part (ii) can be proved very similarly. The natural map $\alpha f_*\to f_*\alpha$ is given by the compatibility of $\alpha$ with pullbacks and the adjunction between $f^*$ and $f_*$. We note that $f$ is cohomologically proper for both 6-functor formalisms. Again, this can be checked after the pullback along $*\to */H$, which results in $*/\ker\varphi\to *$. So we may replace $f$ by $f:*/G\to *$ with $G$ profinite with $\mathrm{cd}_p(G)<\infty$. For $\mathcal D_{\hat\solid}(-)$, the cohomological properness follows from \cref{rslt:6FunctorsSmoothProperMaps}; and for $\mathcal D(-,C^\flat)$, it follows from \cite[Proposition 5.3.2]{HeyerMann6FF}. 

Now since $f_*=f_!$ in both 6-functor formalisms, by the same idea in the proof of (i), it is enough to show 
\begin{itemize}
    \item[($**$)] Suppose $G$ is $p$-torsion free and compact, and let $f\colon */G \to *$. Then the natural map $ \alpha f_! 1\to f_! 1 $ is an isomorphism in $\mathcal D_{\hat\solid}(C^\flat)$.
\end{itemize}
This claim then follows from (i) since $f_\natural(-)=f_!(-\otimes\omega_f)$.
\end{proof}

\begin{rmk}
In the setting of \cref{rslt:smooth-group-homology-compatible-with-D-solid} it is also true that $f_! \alpha \isom \alpha f_!$, because $f_!$ and $f_\natural$ differ by a twist that is the same in both 6-functor formalisms. However, it is difficult to construct an a priori map between $f_! \alpha$ and $\alpha f_!$ which is why we stated the result with $f_\natural$ instead.
\end{rmk}

Using \cref{rslt:smooth-group-homology-compatible-with-D-solid} we can now compute the cohomology with compact support of the smooth Steinberg representation in $\mathcal D_{\hat\solid}(*/F)$ by using a version of Heyer's computation in the setting of smooth representations.  

\begin{defn}[The characteristic $p$ Tate twist] 
\label{DefnTateTwist}
Define a character $\varepsilon\from \Q_p^\times\to \F_p^{\times}$ as the composition of the map $\Q_p^\times\to \Z_p^\times$ (identity on $\Z_p^\times$ and trivial on $p\in\Q_p^\times$) with the reduction map $\Z_p^\times\to \F_p^\times$.  For a finite extension $F/\Q_p$, define $\varepsilon_F\coloneqq \varepsilon\circ N_{F/\Q_p}$, where $N_{F/\Q_p}$ is the norm map.  For an integer $n$, the characteristic $p$ Tate twist $\F_p(n)$ is the representation of $F^\times$ on $\F_p$ whose character is $\varepsilon_F^n$.  For any $\F_p$-vector space $V$, let $V(n)\coloneqq V\otimes_{\F_p}\F_p(n)$.      
\end{defn}

\begin{proposition}
\label{PropCohoOfHModF} 
Let $F$ be a finite extension of $\Q_p$ of degree $d$. Then there are isomorphisms of $C^\flat$-modules with $(F^\times)^2$-action
\begin{align*}
    H^i(\HH^1_{F,C}/F, \OO^\flat) \isom \begin{cases}
        C^\flat, & i = 0,\\
        \chi^{-1}, & i = d+1,\\
        0,& \text{else.}
    \end{cases}
\end{align*}
where $\chi$ is the character $\chi\colon (F^\times)^2 \to C^{\flat,\times}$, $(t_1, t_2) \mapsto \varepsilon_F(t_1 t_2^{-1})$. 
\end{proposition}
\begin{proof}
By the results in the previous subsection, we essentially need to compute the $F$-cohomology of the dual Steinberg representation and the trivial representation. Even though the claim is about cohomology, it is more convenient to do the computation on homology (e.g., because in \cref{rslt:smooth-group-cohomology-compatible-with-D-solid} the group is required to be compact, but $F$ is not). Denote $f\colon \HH^1_{F,C}/B(F) \to */(F^\times)^2$ for $* = \Spa C^\flat$, where $B(F) \subset \GL_2(F)$ is the subgroup of upper triangular matrices. Let $f_\natural$ denote the left adjoint of $f^*$ in the $\mathcal D_{\hat\solid}$-formalism ($f_\natural$ exists by cohomological smoothness of $f$). We need to compute $f_* \OO^\flat$ and by adjunctions we have
\begin{align*}
    f_* \OO^\flat = (f_\natural \OO^\flat)^\vee.
\end{align*}
In order to compute $f_\natural \OO^\flat$, we factor $f$ into the maps
\begin{align*}
    g\colon \HH^1_{F,C}/B(F) \to */B(F), \qquad h\colon */B(F) \to */(F^\times)^2.
\end{align*}
Then $f_\natural = h_\natural g_\natural$. Moreover, $g_\natural$ was computed in \cref{thm:OflatCompactlyDrinfeld} and is given by $\St(C^\flat)$ (with the natural action by $B(F)$) in degree $-1$ and $C^\flat$ (with the trivial $B(F)$-action) in degree $0$. We are thus reduced to proving the following two statements:
\begin{enumerate}[(a)]
    \item $h_\natural C^\flat = C^\flat$.
    \item $h_\natural \St(C^\flat) = \chi[d]$.
\end{enumerate}
Note that both $C^\flat$ and $\St(C^\flat)$ lie in the image of the comparison functor $\alpha$ considered in \cref{prop:group-homology-compatible}. So by \cref{rslt:smooth-group-homology-compatible-with-D-solid}, it is enough to prove the analogous claims in the setting of smooth representations.  Therefore from now on all six operations are to be understood in $\mathcal D(-, C^\flat)$. Claim (a) is well-known, see e.g., \cite[Example~3.4.24]{Heyer2023Left}. Claim (b) is proved in \cite[\S4.3]{Heyer-GeometricLemma} in much more generality using the Geometric Lemma, but in our case there is a quick direct computation as follows. First note that there is an isomorphism
\begin{align*}
    \St(C^\flat) = \LC(\PP^1(F), C^\flat)/C^\flat \isoto \LC_c(F, C^\flat), \qquad \alpha\mapsto \alpha - \alpha(\infty),
\end{align*}
where $\LC_c(F, C^\flat)$ denotes the space of compactly supported locally constant maps $F \to C^\flat$. The above isomorphism is compatible with the $B(F)$-action, and in fact shows that $\St(C^\flat) = s_! C^\flat$, where $s\colon */(F^\times)^2 \to */B(F)$ is the natural section of $h$. Thus
\begin{align*}
    h_! \St(C^\flat) = h_! s_! C^\flat = C^\flat.
\end{align*}
It remains to translate between $h_!$ and $h_\natural$, which is achieved by a certain twist. Namely, note that $h$ is cohomologically smooth (see \cite[Example~5.3.22]{HeyerMann6FF}) and hence $h^! = h^* \otimes \omega_h$ for $\omega_h = h^! C^\flat$. By passing to left adjoints we deduce $h_\natural = h_!(- \otimes \omega_h)$. To compute $\omega_h$, we first note that $\omega_h = \omega_{B(F)} \otimes h^*(\omega_{(F^\times)^2}^{-1})$ and then observe that \cite[Example~5.3.22]{HeyerMann6FF} and \cite[\S2.3.2]{Heyer-GeometricLemma} provide an explicit description of the dualizing complex $\omega_G$ for any $p$-adic Lie group $G$ over $F$: It is given as $\omega_G = \mathfrak d_G[\dim_{\Q_p} G]$, where $\mathfrak d_G$ is the character $\varepsilon_F \circ \det_F(\mathfrak g)\colon G \to C^{\flat,\times}$ with $\det_F(\mathfrak g)$ being the $F$-linear determinant of the adjoint action of $G$ on its Lie algebra $\mathfrak g$. An easy explicit computation shows that $\omega_{(F^\times)^2} = C^\flat[2d]$ and $\omega_{B(F)} = h^*(\chi[3d])$ with $\chi$ as in the claim. Altogether we deduce
\begin{align*}
    &h_\natural \St(C^\flat) = h_! (\St(C^\flat) \otimes \omega_h) = h_!(\St(C^\flat) \otimes h^*(\chi[3d]) \otimes C^\flat[-2d]) \\&\qquad= h_!(\St(C^\flat)) \otimes \chi[d] = \chi[d],
\end{align*}
as desired.
\end{proof}

\begin{rmk}
In the classical setting of a smooth representation of $B(F)$ on a vector space $V$ in characteristic 0, the Jacquet module of $V$ is the space of coinvariants under the unipotent subgroup $F$;  the Jacquet module of the Steinberg representation is 1-dimensional.  In the characteristic $p$ setting, we learned from \cite{Heyer2023Left} that the correct analogue of the Jacquet module is the functor $h_\natural$ of derived coinvariants of $F$;  we also learned the calculation of $h_\natural \St(k)$ ($k$ any discrete field of characteristic $p$) in the case $F=\Q_p$ from \cite[Corollary 5.3.4]{Heyer2023Left}.  Owing to these developments, we use the name {\em Heyer--Jacquet module} to refer to $h_{\natural} \St(k)$ and related objects.
\end{rmk}

\subsection{Application of Drinfeld's Lemma, proof of \cref{thmx:cohom}}\label{sec:DrinfeldsLemma}
In this section, we discuss how to descend our result on $\OO^\flat$-cohomology from $\BC(-1/t)^*\times_k \Spd C$ to $\BC(-1/t)^*$.  
The idea is modeled on the diamond version of Drinfeld's lemma \cite[\S 16.3]{BerkeleyLectures}, which states that for $X=\Spd \Q_p$, the category of finite \'etale covers of $X$ is equivalent to the category of finite \'etale covers of $X\times_{\Spd \F_p} \Spd C/\varphi$. We give in \cref{LemmaDrinfeld} a version of Drinfeld's lemma concerning $\OO^\flat$-cohomology and in \cref{cor:solidDrinfeldLemma} a solid enhancement. But first, we need some discussions on quasi-Stein spaces, which will be used in the proof of \cref{LemmaDrinfeld}.

 Classically, a very useful property of quasi-Stein spaces is they have vanishing higher coherent cohomology in the analytic topology (\cite[Theorem~2.6.5]{KedlayaLiuImperfect}). However, the vanishing property does not hold for $\OO^\flat$-cohomology as we have seen from the cohomology of Drinfeld symmetric spaces in \cref{sec:cohomology-of-Drinfeld-space}. We must therefore impose additional conditions to ensure such vanishing. It would be interesting to investigate which classes of quasi-Stein spaces satisfy the vanishing property for the $\OO^\flat$-cohomology. But here, we only consider a rather ad-hoc class of quasi-Stein spaces needed in subsequent discussions, namely, \textit{quasi-Stein perfectoid spaces}.
\begin{defn}\label{def:quasi-Stein}
    A \textit{quasi-Stein} perfectoid space is a perfectoid space $X$ over $\F_p$ satisfying the following properties:
    \begin{enumerate}
        \item $X$ admits a map of finite transcendence degree to a totally disconnected space;
        \item $X$ can be written as the union of an ascending sequence $U_0 \subseteq U_1 \subset \cdots$ of affinoid subspaces such that for each $i \geq 0$, the map $\OO^\flat(U_{i+1}) \to \OO^\flat(U_i)$ has dense image.
    \end{enumerate}
\end{defn}

\begin{lemma}\label{thm:Stein-vanishing}
    Let $X$ be a quasi-Stein perfectoid space. Then $H^i_{\solid}(X,\OO^\flat_X)=0$ for all $i>0$. In particular, the non-solid cohomology group $H^i(X,\OO^\flat)$ vanishes for $i>0$. 
\end{lemma}
\begin{proof}
    Let $\{U_i\}_{i\geq 0}$ be a family of affinoid analytically open subspaces of $X$ as in Definition~\ref{def:quasi-Stein}.
    View $M_i\coloneqq \OO^b_{U_i}$ as an object in $\mathcal D_{\hat\solid}(\OO^\flat_X)$. Then $\OO^{\flat}_X\simeq R\lim_i M_i$, and we have 
    \begin{equation*}
        R\Gamma_{\solid}(X,\OO^\flat_X)= R\Gamma_{\solid}(X,R\lim_i M_i)\simeq R\lim_i R\Gamma_{\solid}(X,M_i)=R\lim_i R\Gamma_{\solid}(U_i,\OO^\flat_{U_i})
    \end{equation*}
    But each $U_i=\Spa(A_i,A_i^+)$ is an affinoid perfectoid space admitting a map of finite transcendence degree to a totally disconnected space. Therefore, we have 
    \begin{equation*}
        R\Gamma_{\solid}(X,\OO^\flat_X)\simeq R\lim_i A_i=R\lim_i \OO^\flat_X(U_i),
    \end{equation*}
    where each $A_i=\OO^\flat_X(U_i)$ is considered as an object in $\mathcal D_\solid(\F_p)$. By the assumptions in \cref{def:quasi-Stein}(2) and the topological Mittag-Leffler lemma (\cite[Lemma A.37]{Bosco}), $R^j\lim_i \OO^\flat_X(U_i)=0$ for any $j\geq 1$. The result for $H^i(X,\OO^\flat)$ follows since $H^i_{\solid}(X,\OO^\flat)(\ast)=H^i(X,\OO^\flat)$.
\end{proof}

\begin{lemma}[Drinfeld's lemma for $\OO^\flat$-cohomology] 
\label{LemmaDrinfeld}
Let $k/\F_p$ be a finite field of cardinality $q$.  Let $E=k\laurentseries{t^{1/p^\infty}}$, and let $X=\Spa E$. Let $C$ be an algebraically closed perfectoid field whose residue field contains $k$.  Finally, let $\varphi_q$ be the $q$th power Frobenius map on $\Spd C$.  Then the natural map
\[ R\Gamma(X,\OO^\flat)\to R\Gamma(X\times_{\Spd k} (\Spd C)/\varphi_q,\OO^\flat) \] 
is an isomorphism.    
\end{lemma}

\begin{proof} 
Without loss of generality, we may assume $C$ has characteristic $p$. We have $X\times_{\Spd k} \Spd C\isom \widetilde{B}_C^\ast$, the 1-dimensional punctured perfectoid open unit ball over $C$.  

Since $X$ and $\widetilde{B}_C^\ast$ are quasi-Stein perfectoid spaces, they have no higher cohomology with coefficients in $\OO^\flat$ by \cref{thm:Stein-vanishing}, and the map in the proposition is the same as the map
\[ E\to H^0(\widetilde{B}_C^\ast, \OO^\flat)^{h\varphi_q}. \] 
$H^0(\widetilde{B}_C^\ast,\OO^\flat)$ is the ring of power series $\sum_\alpha \gamma_\alpha t^\alpha$, where $\alpha$ runs over $\Z[1/p]$ and the series is convergent on $\widetilde{B}_C^\ast$.  The kernel of $\varphi_q-1$ is $E$, and the cokernel is 0:  for each $\alpha$, let $\delta_\alpha\in C$ be a preimage of $\gamma_\alpha$ under $\varphi_q-1$ of least absolute value.  Examining Newton polygons, we find that $\abs{\delta_\alpha}=\abs{\gamma_\alpha}$ or $\abs{\delta_\alpha}=\abs{\gamma_\alpha}^{1/q}$ as $\abs{\gamma_\alpha}\leq 1$ or $\abs{\gamma_\alpha}>1$, respectively.  Either way, $\abs{\delta_\alpha}\leq \abs{\gamma_\alpha}^{1/q}$, which shows that $\sum_\alpha \delta_\alpha t^\alpha$ converges on $\widetilde{B}_C^\ast$.
\end{proof}
In the lemma above, we did not consider the condensed structures on the two cohomology groups. Namely, the map $R\Gamma(X,\OO^\flat)\to R\Gamma(X\times_{\Spd k} (\Spd C)/\varphi_q,\OO^\flat)$ is really an isomorphism between discrete groups. Our next step is to enhance this isomorphism to an isomorphism of condensed groups. But before that, it is necessary to collect some results on Fr\'{e}chet spaces. Variants of these results appear in the literature, but none seems to provide exactly the formulation we need.

In the following, let $K$ be a complete non-archimedean field, i.e., a topological field that is complete with respect to a nontrivial non-archimedean norm. In particular, $K$ can be a local field or a perfectoid field.
\begin{defn}
    \begin{enumerate}
        \item A $K$-Banach space is a $K$-vector space complete with respect to a norm (cf. \cite[Definition~2.1.1]{LocallyConvexSpacesBook}).
        \item A $K$-Fr\'echet space is a complete metrizable locally convex $K$-vector space (cf. \cite[Definition~3.5.1]{LocallyConvexSpacesBook}).
    \end{enumerate}
\end{defn}

\begin{lemma}
    A topological $K$-vector space $V$ is Fr\'echet if and only if it is a countable cofiltered limit of $K$-Banach spaces. 
\end{lemma}
\begin{proof}
    This is a consequence of \cite[Corollary~3.5.7]{LocallyConvexSpacesBook} and \cite[Ch.II \S 5.4 Corollary~1]{TopologicalVectorSpaces-Schaefer-Wolff}.
\end{proof}
\begin{example}
    \begin{enumerate}
        \item A $K$-Banach space is Fr\'echet. In particular, a perfectoid $K$-algebra is Fr\'echet.
        \item Let $X$ be a quasi-Stein perfectoid space over a perfectoid field $K$ of characteristic $p$. Then $H^0(X,\OO_X^\flat)$ is a $K$-Fr\'echet space.
    \end{enumerate}
\end{example}

For any $K$-Fr\'echet space $V$, write $\underline{V}$ for its condensification as defined in \cite[Example 1.3]{ClausenScholzeCondensed}. Any $K$-Fr\'echet space is compactly generated since it is metrizable. Thus, the condensification functor $V\mapsto\underline{V}$ from the category of $K$-Fr\'echet spaces to condensed $K$-modules is fully faithful by \cite[Proposition 1.7]{ClausenScholzeCondensed}. The condensed $K$-module $\underline{V}$ is in fact solid by \cite[Proposition A.31]{Bosco}.
\begin{lemma}\label{lem:condensification-surjective-Frechet}
    Let $V\to W$ be a surjection of $K$-Fr\'echet spaces. Then its condensification $\underline{V}\to\underline{W}$ is a surjection of solid $K$-modules. 
\end{lemma}
\begin{proof}
    This follows from the same proof of \cite[Lemma~A.33]{Bosco}. Note that Bosco requires the field $K$ (or $F$, in Bosco's notations) to be discretely valued and to have finite residue field. However, these additional assumptions are not used in the proof of \cite[Lemma~A.33]{Bosco} and only matter in some other results in \cite{Bosco}.
\end{proof}

Our digression into general Fr\'echet algebras is now over, and we are ready to prove a solid version of \cref{LemmaDrinfeld}. In the notation of Lemma~\ref{LemmaDrinfeld}, we can describe explicitly the solid cohomology groups of $X=\Spa E$ and $X\times_{\Spd k} (\Spd C)/\varphi_q$ as defined in \cref{rmk:solid-v-cohomology} and \cref{rmk:solid-v-cohomology-over-k}. The equivalence $\mathcal{D}_{\hat\solid}(\OO^\flat_X)\simeq \mathcal{D}_{\solid}(E)$ endows $R\Gamma(X,\OO^\flat)$ with a solid structure, which is the group $R\Gamma_{\solid}(X,\OO^\flat)$. Similarly, $R\Gamma(X\times_{\Spd k} (\Spd C)/\varphi_q,\OO^\flat)$ has a solid structure $R\Gamma_{\solid}(X\times_{\Spd k} (\Spd C)/\varphi_q,\OO^\flat)$ via the identification
\begin{equation*}
    R\Gamma(X\times_{\Spd k} (\Spd C)/\varphi_q,\OO^\flat)\simeq f_\ast\OO^\flat
\end{equation*}
where $f:X\times_{\Spd k} (\Spd C)/\varphi_q\to X$ is the projection map. 

As argued in \cref{LemmaDrinfeld}, $X\times\Spd C$  is quasi-Stein and $R\coloneqq H^0(X\times\Spd C,\OO^\flat)$ is isomorphic to $ R\Gamma(X\times\Spd C,\OO^\flat)$ and is a Fr\'echet $C$-algebra. Via the projection $X\times\Spd C\to X$, $R$ is also a Fr\'echet $E$-algebra. We have the following solid enhancement of Lemma~\ref{LemmaDrinfeld}.

\begin{cor}[Solid Drinfeld's lemma]\label{cor:solidDrinfeldLemma}
    The natural map 
    \begin{equation*}
         R\Gamma_{\solid}(X,\OO^\flat)\to R\Gamma_{\solid}(X\times_{\Spd k} (\Spd C)/\varphi_q,\OO^\flat)
    \end{equation*}
    is an isomorphism in $\mathcal D_\solid(\F_p)$.
\end{cor}
\begin{proof}
    The solid cohomology group $R\Gamma_{\solid}(X\times (\Spd C)/\varphi_q,\OO^\flat) $ is computed by the condensification of the complex 
    \begin{equation*}
        R\xrightarrow{\varphi_q-\mathrm{id}}R,
    \end{equation*}
    which will be viewed as a complex of $E$-modules via the projection $X\times\Spd C\to X$.
    As shown in the proof of Lemma~\ref{LemmaDrinfeld}, the map $\varphi_q-\mathrm{id}$ is surjective. Since $R$ is a Fr\'echet $E$-algebra, the condensification of the map $\varphi_q-\mathrm{id}$ is again surjective by \cref{lem:condensification-surjective-Frechet}. Since the condensification functor preserves kernels, we then get the desired isomorphism by using Lemma~\ref{LemmaDrinfeld}.
\end{proof}

\begin{proposition}  
\label{PropCohoOfBC}
Let $k/\F_p$ be an algebraic extension of $\F_{p^t}$, and let $\BC(-1/t)$ be the negative BC space of slope $-1/t$ over $k$. For each $i\in\Z$ there is an isomorphism of solid $k$-modules with $\overline{\G}_t^0$-action:
\[ H^i_{\solid}(\BC(-1/t)^*,\OO^\flat) \isom 
\begin{cases}
k, & i=0,\\
k(-1), & i=t+1,\\
0, & \text{else}.
\end{cases}
\] 
Here $k(-1)$ is the representation of $\overline{\G}_t^0$ on $k$ with character $\varepsilon^{-1}\circ \det$, where $\det\from\overline{\G}_t^0\to\overline{\G}_1^0= \Q_p^\times$ is the determinant map, and $\varepsilon$ is defined in \cref{DefnTateTwist}. The same result holds for $H^i(\BC(-1/t)^*,\OO^\flat)$.
\end{proposition}

 \begin{proof} For $k=\F_{p^{h(h-1)t}}$, by \cref{PropL} ~(5), $\BC(-1/t)^*\simeq \Spd\widehat{L}/SG_h$, where $\Spd\widehat{L}$ is affinoid perfectoid and $SG_h$ is profinite. Thus, $\BC(-1/t)^*$ is qcqs. Since $\Spd \F_{p^{h(h-1)t}}\to \Spd \F_{p^t}$ is finite \'etale, $\BC(-1/t)^*$ is also qcqs when defined over $k=\F_{p^t}$. By \cref{lem:base-change-identity} (2), it suffices to show that $H^i_{\solid}(\BC(-1/t)^*,\OO^\flat)$ is a finite-dimensional $k$-vector space viewed as a solid $k$-module for $k=\F_q$, where $q=p^t$. Any element of $\overline{\G}_t^0=D_t^\times$ lies in a subfield $F\subset D_t$ of degree $t$ over $\Q_p$, with $\det$ agreeing with then norm map relative to $F/\Q$.  Therefore it suffices to establish the claimed isomorphism as representations of $F^\times$, where $F$ is an arbitrary maximal subfield of $D_t$.   

Applying \cref{cor:solidDrinfeldLemma} to $E=\widehat{L}$ yields a natural isomorphism
\[ R\Gamma_{\solid}(\Spd \widehat{L},\OO^\flat)\isom R\Gamma_{\solid}(\Spd \widehat{L}\times_{\Spd k} (\Spd C)/\varphi_q,\OO^\flat). \]
Passing to derived $SG_h$-invariants (note the $SG_h$-action commutes with $\varphi_q$) and applying \cref{PropL}(5) followed by \cref{rslt:BaseChangeOfOmegaToC}, we obtain isomorphisms:
\begin{eqnarray*}
R\Gamma_{\solid}(\BC(-1/t)^\ast ,\OO^\flat)&\isom& R\Gamma_{\solid}(\BC(-1/t)^\ast\times_{\Spd k} \Spd C)^{h\varphi_q}  \\
&\isom& R\Gamma_{\solid}(\mathcal{H}^1_{F,C}/F,\OO^\flat)^{h\varphi_q}  
\end{eqnarray*}
Here, as discussed in \cref{rmk:solidCohomologyoverC}, we ignore the difference between $R\Gamma_{\solid}(\mathcal{H}^1_{F,C}/F,\OO^\flat)$ and $R\Gamma(\mathcal{H}^1_{F,C}/F,\OO^\flat)$ since it is defined over the perfectoid field $C$.
The cohomology of $R\Gamma(\mathcal{H}^1_{F,C}/F,\OO^\flat)$ is computed in  \cref{PropCohoOfHModF}.  We conclude by observing that $(C^\flat)^{h\varphi_q}=k$ as solid $k$-modules, because $C^\flat$ is algebraically closed. The result for $H^i(\BC(-1/t)^*,\OO^\flat)$ follows since $H^i_{\solid}(\BC(-1/t)^*,\OO^\flat)(*)=H^i(\BC(-1/t)^*,\OO^\flat)$ and since $k$ and $k(-1)$ are discrete as solid $k$-modules.
\end{proof}

We are now ready to present the main results of the section.
\begin{thm}\label{ThmCohomologyOfLBody} 
    Assume $h\geq 2$.  Let $k/\F_p$ be an algebraic extension field containing $\F_{p^{(h-1)h}}$.   Let $L/K$ be the Torii extension, let $\G_{h-1}$ and $\G_h$ be the Morava stabilizer groups (relative to $k$), let $\overline{\G}_h$ denote the kernel of $\G_h\to\Gal(k/\F_p)$, and let $S\overline{\G}_h$ denote the kernel of the determinant map $\overline{\G}_h\to \Z_p^\times$.
    
    Then there is a fiber sequence of derived $k$-modules with $\G_{h-1}\times \Z_p^\times$-action:
        \begin{equation}\label{EqFiberSeqk}
            k \to \widehat{L}^{hS\G_h} \to k[-h] 
        \end{equation}
    Regarding the action of $\overline{\G}_{h-1}\times \Z_p^\times$:  it is trivial on $k$, and a pair $(g,a)$ acts on the $k[h]$ term by $(\det g)^{-h}a^{h-1}$.  
\end{thm}
\begin{proof} 
    Recall \cref{PropL}(5):  Under the isomorphism $(\Spd \widehat{L})/S\overline{\G}_h\isom \BC(-1/(h-1))^*$, 
   a pair $(g,a)$ acts via  $g(\det g)a^{-1}$ under the natural action of $\overline{\G}_{h-1}$ on $\BC(-1/(h-1))$.  Taking $\OO^\flat$-cohomology and applying \cref{PropCohoOfBC}, we arrive at the claim.  
\end{proof}

\begin{rmk}\label{rem:HJMquestions}
    We leave for the future some questions concerning the Heyer--Jacquet module, including:        
        \begin{enumerate}
            \item By \cref{PropCohoOfBC}, the cohomology $R\Gamma(\BC(-1/t)^*,\OO^\flat)$ is a $\G_t$-equivariant extension of $k[-t-1](-1)$ by $k$, whose extension class determines a canonical element of $H^{t+2}_{\cts}(\G_t,k(1))$. It would be interesting to know a concrete description of this ``Heyer--Jacquet class'', or at least tell when it is nontrivial. 
            \item The techniques of this chapter ought to extend to the calculation of the mixed-characteristic module $R\Gamma(\BC(-1/t)^*,W(\OO^\flat))$, where $W(\OO^\flat)$ means Witt vectors of the structure sheaf. 
        \end{enumerate}
\end{rmk}

\newpage
\section{Solid Lubin--Tate theories}\label{sec:solidLT}

In Lurie's framework \cite{ECII}, the basic data to which one associates a Lubin--Tate theory is that of a pair $(R, \Gamma)$ of a commutative $\F_p$-algebra $R$ and a $p$-divisible group $\Gamma$ over $R$.
In the case that $R$ is Noetherian and $F$-finite and $\Gamma$ is nonstationary over $R$ (see \Cref{def:nonstationary}), Lurie has constructed a well-behaved Lubin--Tate theory $E(R,\Gamma)$.
\footnote{In \cite[Section 5]{ECII}, Lurie also treats the case of a perfect $\F_p$-algebra $R$ and a connected $p$-divisible group $\Gamma$ over $R$---in this case, one does not need to assume that $R$ is Noetherian. However, this comes at the cost of having to additionally $K(h)$-localize at the end of the construction. This section will focus on extending the Noetherian constructions of \cite[Section 6]{ECII}.} \footnote{In fact, Lurie only refers to these as Lubin--Tate theories in the aforementioned case where $R$ is perfect and $\Gamma$ is connected. However, since ``orientation classifier of the universal spectral deformation of $\Gamma$'' is a mouthful, we allow ourselves to use ``Lubin--Tate theory'' in the more general situation.}

In this section, we will identify a category $\LKN$ of \emph{Luriean Kalypto--Noetherian pairs} whose Lubin--Tate theories are well-behaved. In this category, we relax the assumption that $R$ be Noetherian: instead, we ask that there exist a flat and relatively perfect homomorphism $R_0 \to R$ with $R_0$ Noetherian and $F$-finite, and that $\Gamma$ be base changed from a nonstationary $p$-divisible group $\Gamma_0$ over $R_0$.
We extend Lurie's construction of Lubin--Tate theories to this context, proving the following theorem:

\begin{thm}\label{thm:Discrete-LKN-functor}
    There is a Lubin--Tate theory functor
    $E(-) \colon  \LKN \to \CAlg$
    satisfying the following properties:
    \begin{itemize}
        \item There is a natural isomorphism $\pi_0 E(R,\Gamma) \cong R^{\mathrm{cl}}_{(R,\Gamma)},$ where $R^{\mathrm{cl}}_{(R,\Gamma)}$ is the classical deformation ring of $(R,\Gamma)$.
        \item There are natural isomorphisms $\pi_{2n} E(R,\Gamma) \cong \omega^n_{\Gamma^{\mathrm{un}}},$ where $\omega_{\Gamma^{\mathrm{un}}}$ is the invertible module of invariant differentials of the universal $p$-divisible group $\Gamma^{\mathrm{un}}$ living over $R^{\mathrm{cl}}_{(R,\Gamma)}$.
    \end{itemize}
\end{thm}

The following constructions applied to Lubin--Tate theories will play an important role for us:
\begin{itemize}
    \item Given a Lubin--Tate theory $E(R,\Gamma)$, we may consider its $K(h)$-localization $L_{K(h)} E(R,\Gamma)$.
    \item Given a $K(h)$-local Lubin--Tate theory $E(R_1, \Gamma_1)$, we may consider the $K(h)$-localized tensor product $L_{K(h)} (E(R_1, \Gamma_1) \otimes_{\SW(k)} E_h)$.
\end{itemize}

We will define a full subcategory $\LKN_h \subseteq \LKN$ of $\LKN$ pairs of height $h$ in \Cref{defn:LKN-ht-h}. If $(R,\Gamma)$ lies in $\LKN_h$, then in particular the identity component of $\Gamma$ is of exact height $h$ and $E(R,\Gamma)$ is $K(h)$-local.

\begin{thm}
    For each integer $h \geq 1$, there is a functor
    \[L^{\alg}_{K(h)} \colon  \LKN \to \LKN_h\]
    and natural isomorphisms
    \[E(L_{K(h)}^{\alg} (R,\Gamma)) \simeq L_{K(h)} E(R,\Gamma).\]
\end{thm}

\begin{thm}
    Let $k$ denote a Noetherian weakly \'etale $\F_p$-algebra and let $\LKN_{h,k}$ denote the variant of $\LKN_h$ where $R$ is a $k$-algebra, and let $\LKN_{k,h} ^{p} \subseteq \LKN_{h,k}$ denote the full subcategory of $(R,\Gamma)$ such that $R$ is perfect and $\Gamma$ is base changed from $k$.
    Then there is a functor
    \[ - \notimes_k - \colon  \LKN_{h,k} \times \LKN^{p}_{h,k} \to \LKN_{h,k}\]
    such that for $(R_1, \Gamma_1)\in \LKN_h$ and $(R_2, \Gamma_2) \in \LKN_h ^{p}$, there is a natural isomorphism
    \[E((R_1,\Gamma_1) \notimes_k (R_2, \Gamma_2)) \simeq L_{K(h)} (E(R_1, \Gamma_1) \otimes_{\SW (k)} E(R_2, \Gamma_2)).\]
\end{thm}

In the context of height $h$ formal groups over perfect $\F_p$-algebas, Clausen has enriched the Lubin--Tate functor to take values in the category of solid $\Einf$-rings \cite[Section 12]{Clausen2025arxiv}.
In this section, we use the formalism of $\LKN$ pairs to extend this construction to the context of $p$-divisible groups over solid rings.

We will define a category $\LKN^{\solid}$ of solid $\LKN$-pairs $(R, \Gamma)$---these consist of a solid ring $R$ and a $p$-divisible group $\Gamma$ over its underlying ring $R(\ast)$, subject to certain conditions. These imply, for example, that the pairs $(R(S), \Gamma)$ are $\LKN$.
It will almost never be true that $R(S)$ is Noetherian when $S$ is infinite, and this is our motivation for working in the greater generality of $\LKN$ pairs.
We again define a full subcategory $\LKN ^{\solid} _h \subseteq \LKN^{\solid}$ of solid $\LKN$ pairs of height $h$, and given a Noetherian weakly \'etale $\F_p$-algebra $k$ define a category $\LKN ^{\solid} _{h,k}$.

Then we are able to generalize the above results to the solid context.

\begin{thm} \label{thm:solidLT}
    The categories $\LKN^\solid$, $\LKN ^{\solid} _{h}$, $\LKN ^{\solid} _{h,k}$ are naturally enriched over condensed anima, and there are natural functors enriched over condensed anima
    \begin{align*}
        &E^{\solid} (-) \colon  \LKN^{\solid} \to \CAlg^{\solid} \\
        &L_{K(h)} ^{\solid, \alg} \colon  \LKN^{\solid} \to \LKN_h ^{\solid} \\
        &\notimes_{k^\delta} ^{\solid} \colon  \LKN_{h,k} ^{\solid} \times \LKN_{h,k}^{p} \to \LKN_{h,k} ^\solid,
    \end{align*}
    where we view $\LKN_{h,k} ^p$ as enriched over discrete condensed anima, and natural isomorphisms of solid $\Einf$-rings
    \begin{align*}
        &E^{\solid}(L_{K(h)}^{\solid, \alg} (R,\Gamma)) \simeq L_{K(h)} ^\solid E^{\solid}(R,\Gamma) \\
        &E^\solid ((R_1,\Gamma_1) \notimes_{k^\delta}^\solid (R_2, \Gamma_2)) \simeq L_{K(h)} ^\solid (E^\solid(R_1, \Gamma_1) \otimes_{\SW (k)}^\solid E^\solid(R_2^\delta,\Gamma_2)).
    \end{align*}
    Moreover, the condensed functors $E^\solid(-)$, $L_{K(h)}^{\solid, \alg}$, and $\notimes_{k^\delta} ^{\solid}$ are computed by applying the corresponding uncondensed functors of the same name for each profinite set. For example $E^\solid (R,\Gamma)$ is the hypersheaf
    \[S \mapsto E(R(S),\Gamma)\]
    for $S \in \Prof$.
\end{thm}

Finally, we will discuss the following examples, which are relevant for the rest of this paper, in detail in \Cref{subsec:exapp}.

\begin{example}
    Let $k$ denote an algebraic extension of $\F_{p^h}$. Then the pair $(k ^\delta, \Gamma_h)$ is an example of a solid $\LKN$ pair of height $h$, and $E_h^{\solid} = E^{\solid} (k ^\delta,\Gamma_h)$ agrees with Clausen's solid Lubin--Tate theory.
    Then $L_{K(h-1)} ^{\solid, \alg} (k,\Gamma_h) = (K,\Gamma_{h,K} ^{\mathrm{un}})$, where $K \cong k((u_{h-1}))$, so that by \Cref{thm:solidLT} there is an isomorphism
    \[L_{K(h-1)} ^\solid E^\solid _h \simeq E^\solid (K, \Gamma_{h,K}^{\un}).\]

    Furthermore, one may identify
    $(K, \Gamma_{h,K}^{\mathrm{un}}) \notimes_{k ^\delta} ^\solid (k, \Gamma_{h-1}) \cong (L,\Gamma_{h,L}^{\mathrm{un}}),$
    where $L$ is the Gross--Torii extension of $K$,
    so that
    \[\B = L_{K(h-1)} ((L_{K(h-1)} E_h) \otimes_{\mathbb{SW}(k)} E(k, \Gamma_{h-1}) \simeq E(L, \Gamma_{h,L}^{\mathrm{un}})(\ast).\]

    Using this, we are able to construct a $\G_{h,h-1}$-equivariant map of solid $\Einf$-rings
    $\mathbb{B} \to \widehat{\mathbb{B}}$
    using the map of solid $\LKN$ pairs
    \[(L, \Gamma_{h,L}^{\mathrm{un}}) \to (\widehat{L}, \Gamma_{h,\widehat{L}}^{\mathrm{un}})\]
    and the fact that, since $\widehat{L}$ is perfect by \Cref{PropL}(1), there is a natural map
    $E^{\solid} (\widehat{L}, \Gamma_{h,\widehat{L}}^{\mathrm{un}}) \to E^{\solid} (\widehat{L}, (\Gamma_{h,L}^{\mathrm{un}})^{\circ}).$
    The construction of this map is main output of this section that will be used in the rest of the paper.
\end{example}

\subsection{Background on Lubin--Tate theories} \label{sec:LT-back}

We begin by recalling Lurie's approach to the construction of Lubin--Tate theories. This is based on an extension of the definition of formal and $p$-divisible group to the context of $\Einf$-rings.

\begin{defn}
    Recall that an adic ring is a topological ring $R$ for which there exists a finitely generated ideal of definition $I$, and that a morphism of adic rings is a continuous ring homomorphism.
    An adic $\Einf$-ring is an $\Einf$-ring $R$ with the structure of adic ring on $\pi_0 (R)$, and a morphism of adic $\Einf$-rings is a morphism of $\Einf$-rings that restricts to a continuous map on $\pi_0$.
    More formally, we define $\CAlg_{\ad}$ as the pullback
    \begin{center}
        \begin{tikzcd}
            \CAlg_{\ad} \ar[r, "\pi_0"] \ar[d] & \CAlg_{\ad} ^{\heartsuit} \ar[d] \\
            \CAlg \ar[r, "\pi_0"] & \CAlg ^{\heartsuit}.
        \end{tikzcd}
    \end{center}

    We say that an adic $\Einf$-ring is \emph{complete} if it is $I$-complete for one (hence any) choice of ideal of definition $I \subseteq \pi_0 (R)$.
    We let $\CAlg^{\mathrm{cmpl}}_{\ad} \subseteq \CAlg_{\ad}$ denote the full subcategory of complete adic $\Einf$-rings.
\end{defn}

\begin{defn}[{\cite[Definition 1.6.1]{ECII}}]
    Let $R$ denote a connective $\Einf$-ring. A formal group (of dimension one) over $R$ is a functor
    \[\widehat{\Gamma} \colon  \CAlg^{\mathrm{cn}}_{R} \to D(\Z)_{\geq 0}\]
    with the property that the composite
    \[\CAlg^{\mathrm{cn}}_{R} \xrightarrow{\widehat{\Gamma}} D(\Z)_{\geq 0} \xrightarrow{\Omega^\infty} \Ani\]
    is of the form $\Spf(A)$ for $A$ an adic $\Einf$-ring satisfying the following property:
    \begin{itemize}
        \item There exists an invertible $\pi_0 (R)$-module $M$ such that
        \[
            \pi_* (A) \cong \prod_{n \geq 0} (\Sym^n_{\pi_0 (R)} (M) \otimes_{\pi_0 (R)} \pi_* (R))
        \]
        as graded $\pi_* (R)$-algebras, and $M$ generates the ideal of definition in $\pi_0 (A)$.
    \end{itemize}
\end{defn}

\begin{example}
  Recall that an $\Einf$-ring $R$ is said to be weakly periodic if $\pi_2 (R)$ is an invertible $\pi_0 (R)$-module and the map $\pi_n (R) \otimes_{\pi_0 (R)} \pi_{2} (R) \to \pi_{n+2} (R)$ is an isomorphism for all $n$.
  We say that an $\Einf$-ring $R$ is complex periodic if it is weakly periodic and complex orientable.

  If $R$ is a complex periodic $\Einf$-ring, then $R^{\mathbb{CP}^\infty}$ is a complete adic $\Einf$-ring with ideal of definition equal to the kernel of the map $\pi_0 (R^{\mathbb{CP}^\infty}) \to \pi_0 (R)$.
  The product on $\mathbb{CP}^\infty$ endows $\Spf (R^{\mathbb{CP}^\infty})$ with the structure of a formal group over $R$.
  This formal group is called the \emph{Quillen formal group} over $R$, and we denote it by $\Gamma^Q _R$.
\end{example}

\begin{rmk}
    Given a general $\Einf$-ring $R$, we define a formal group or $p$-divisible group over $R$ to be the respective object over its connective cover $\tau_{\geq 0} R$.
\end{rmk}

\begin{defn}[{\cite[Definition 2.0.2]{ECII}}]
    Let $R$ denote a connective $\Einf$-ring. A $p$-divisible group over $R$ is a functor 
    $\Gamma\colon \CAlg^{\mathrm{cn}}_R \to D(\Z)_{\geq 0}$ with the following properties:
    \begin{enumerate}
        \item For every $A \in \CAlg^{\mathrm{cn}}_R$, the $\Z$-module $\Gamma(A)$ is $p$-nilpotent: we have $\Gamma(A)[1/p] \simeq 0$.
        \item For every finite abelian $p$-group $M$, the functor
        \[A \mapsto \Omega^\infty \Hom(M,\Gamma(A))\]
        is representable by a finite flat $R$-algebra.
        \item The map $p \colon  \Gamma \to \Gamma$ is locally surjective with respect to the finite flat topology.
    \end{enumerate}
\end{defn}

\begin{defn}
    We let $\mathcal{M}_{\mathrm{BT}}^{\Einf}$ denote the moduli stack of $p$-divisible groups:
    given a connective $\Einf$-ring $R$, then $\mathcal{M} ^{\Einf}_{\mathrm{BT}} (R)$ is the $\infty$-groupoid of $p$-divisible groups over $R$.\footnote{Here the subscript ``BT'' refers to ``Barsotti--Tate'', which is another name for $p$-divisible groups.}
    We let $\mathcal{M}_{\mathrm{BT}}$ denote its underlying classical stack.

    Similarly, we let $\mathcal{M}_{\mathrm{fg}}^{\Einf}$ denote the moduli stack of formal groups over $\Einf$-rings, and let $\mathcal{M}_{\mathrm{fg}}$ denote its underlying classical stack.
\end{defn}

\begin{rmk}
    Given a $p$-divisible group $\Gamma$ over a $p$-complete base $\Einf$-ring $R$, Lurie has constructed a formal group $\Gamma^{\circ}$ over $R$ called the identity component of $\Gamma$ \cite[Theorem 2.0.8]{ECII}.
\end{rmk}

Our next goal is to state Lurie's results about the deformation theory of $p$-divisible groups over $\Einf$-rings. First, we need some preliminaries on the cotangent complex.

\begin{ntn}
    Given a map of $\Einf$-rings $R \to S$, we let $L^{\Einf}_{S/R}$ denote the $\Einf$ relative cotangent complex, defined in terms of square-zero extensions of $\Einf$-rings.

    We also make use of the relative topological cotangent complex $L_{X/Y}^{\Einf}$ of a natural transformation $X \to Y$ between functors
    $X, Y \colon  \CAlg^{\mathrm{cn}} \to \Ani$,
    which is uniquely determined if it exists (see \cite[Definition 17.2.4.2]{SAG}). 

    Given a map of commutative rings $R \to S$, we let $L^{\alg}_{S/R}$ denote the algebraic relative cotangent complex, defined in terms of square-zero extensions of animated commutative rings.
\end{ntn}

Given a map of commutative rings $R \to S$, there is a natural homomorphism $L^{\Einf}_{S/R} \to L^{\alg}_{S/R}$ which in general is not an isomorphism.
However, the following result is reassuring.

\begin{proposition} \label{prop:cot-compare}
\begin{enumerate}
    \item Let $R$ denote a commutative ring. Then the absolute $\Einf$ cotangent complex $L^{\Einf}_R$ is almost perfect as an $R$-module if and only if the absolute algebraic cotangent complex $L^{\alg}_R$ is almost perfect as an $R$-module.
    \item Let $R \to S$ denote a map of commutative rings. Then $L^{\Einf} _{S/R} \simeq 0$ if and only if $L^{\alg} _{S/R} \simeq 0$.
\end{enumerate}
\end{proposition}

\begin{proof}
    Part (1) is part of \cite[Proposition 3.3.4]{ECII}.

    To prove part (2), let $S^+$ denote the $\mathbb{E}_1$-ring of \cite[Corollary 25.3.3.3]{SAG}.  By \cite[Remark 25.3.3.7]{SAG}, there is an equivalence $L^{\alg}_{S/R} \simeq S \otimes_{S^+} L^{\Einf}_{S/R}$. Now, $L^{\Einf} _{S/R}$ is connective and the map $S^+ \to S$ has $2$-connective fiber by the proof of \cite[Proposition 25.3.5.1]{SAG}. It follows that $L^{\Einf}_{S/R}$ vanishes if and only if $L^{\alg} _{S/R}$ does.
\end{proof}

Next, we discuss the cotangent complex of $\mathcal{M}_{\mathrm{BT}}$.

\begin{proposition}[{\cite[Propositions 3.2.2 (3) \& 3.2.5]{ECII}}]
    The moduli stack of $p$-divisible groups $\mathcal{M}_{\mathrm{BT}} ^{\Einf}$ admits a $(-1)$-connective cotangent complex $L^{\Einf}_{\mathcal{M}_{\mathrm{BT}^{\Einf}}}$.
    If $R$ is a connective $\Einf$-ring for which $p$ is nilpotent in $\pi_0 (R)$ and $\eta \colon  \Spec (R) \to \mathcal{M}_{\mathrm{BT}^{\Einf}}$ corresponds to a $p$-divisible group over $R$, then $\eta^* L^{\Einf}_{\mathcal{M}_{\mathrm{BT}^\Einf}}$ is furthermore connective and almost perfect as an $R$-module.
\end{proposition}

\begin{defn}[{\cite[Definition 3.0.8 \& Remark 3.4.4]{ECII}}] \label{def:nonstationary}
    A $p$-divisible group $\Gamma$ over a commutative ring $R$ in which $p$ is nilpotent is said to be \emph{nonstationary} if the relative $\Einf$ cotangent complex $L^{\Einf} _{\Spec(R)/\mathcal{M}_{\mathrm{BT}^{\Einf}}}$ is $1$-connective, i.e., if the canonical map
    $L^{\Einf} _{\mathcal{M}_{\mathrm{BT}^{\Einf}}}\vert_{\Spec R} \to L^{\Einf}_{\Spec R}$
    is surjective on $\pi_0$.
\end{defn}

\begin{thm}[{\cite[Theorem 3.4.1]{ECII}}]
    Let $R$ be a commutative ring and let $\Gamma$ be a $p$-divisible group over $R$.
    Assume that $p$ is nilpotent in $R$, that $\Gamma$ is nonstationary, and that the cotangent complex $L^{\Einf} _{R}$ is almost perfect as an $R$-module.
    Then:
    \begin{enumerate}
        \item The $p$-divisible group $\Gamma$ admits a universal deformation to a complete adic $\Einf$-ring $R^{\mathrm{un}}_{(R,\Gamma)}$, equipped with a map 
         $R^{\mathrm{un}} _{(R,\Gamma)} \to R$. 
        \item The spectral deformation ring $R^{\mathrm{un}} _{(R,\Gamma)}$ is connective.
        \item The canonical map $R^{\mathrm{cl}}_{(R,\Gamma)} \coloneqq \pi_0 (R^{\mathrm{un}}_{(R,\Gamma)}) \to R$ is a surjection and exhibits $R^{\mathrm{cl}}_{(R,\Gamma)}$ as the classical universal deformation ring of $\Gamma$ in the realm of static rings.
        \item The canonical map $R^{\mathrm{un}} _{(R,\Gamma)} \to R$ exhibits $R$ as an almost perfect $R^{\mathrm{un}} _{(R,\Gamma)}$-module. In particular, the kernel $I$ of $R^{\mathrm{cl}}_{(R,\Gamma)} \to R$ is finitely generated.
        \item The kernel $I$ is an ideal of definition for $R^{\mathrm{cl}}_{(R,\Gamma)}$. In particular, $R^{\mathrm{un}} _{\Gamma}$ is $I$-complete.
        \item If $R$ is Noetherian, then $R^{\mathrm{un}} _{(R,\Gamma)}$ is Noetherian.
    \end{enumerate}
\end{thm}

\begin{rmk} \label{rmk:noe-f-finite}
    Recall that an $\F_p$-algebra $R$ is said to be \emph{$F$-finite} if the Frobenius morphism $F\colon R \to R$ exhibits the target as a finitely generated module over the source.

    Lurie has proven that a Noetherian $\F_p$-algebra $R$ is $F$-finite if and only if the absolute cotangent complex $L^{\alg}_{R}$ (or, equivalently, $L^{\Einf}_R$) is almost perfect as an $R$-module \cite[Theorem 3.5.1]{ECII}.
\end{rmk}

Finally, we review the notion of an \emph{orientation} of a formal group and how it may be used to define Lubin--Tate theories.

\begin{construction}[{\cite[Construction 4.2.9]{ECII}}]
    Let $R$ denote an $\Einf$-ring, and let $\widehat{\Gamma}$ denote a formal group over $R$.
    For $R$-algebras $A$, there is a natural map
    \[\mathcal{L} \colon  \Omega \widehat{\Gamma}(\tau_{\geq 0} A) \to \mathrm{Map}_{\Mod_{R}} (\omega_{\widehat{\Gamma}} [1], A)\]
    that we call the \emph{linearization map}.

    Indeed, when $R$ is connective there is a natural fiber sequence
    \[\omega_{\widehat{\Gamma}} [1] \to R \otimes_{\mathcal{O}_{\widehat{\Gamma}}} R \to R\]
    by \cite[Proposition 4.2.8]{ECII}.
    Thus we obtain
    \begin{align*}
        \Omega \widehat{\Gamma}(\tau_{\geq 0} A) &\simeq \mathrm{Map}_{\CAlg_{\tau_{\geq 0} R}} (\tau_{\geq 0} R \otimes_{\mathcal{O}_{\widehat{\Gamma}}} \tau_{\geq 0} R,\tau_{\geq 0} A) \\
        &\to \mathrm{Map}_{\Mod_{\tau_{\geq 0}R}} (\tau_{\geq 0} \omega_{\widehat{\Gamma}}[1], \tau_{\geq 0} A)
        \simeq \mathrm{Map}_{\Mod_{R}} (\omega_{\widehat{\Gamma}} [1], A)
    \end{align*}
\end{construction}

\begin{defn}[{\cite[Definition 4.3.1, Construction 4.3.7, and Definition 4.3.9]{ECII}}]
    Given a formal group $\widehat{\Gamma}$ over an $\Einf$-ring $R$, a \emph{pre-orientation} of $\widehat{\Gamma}$ is a (pointed) map $e \colon  S^2 \to \widehat{\Gamma}(\tau_{\geq 0} R).$
    Viewing $e$ as an element of $\Omega^2 \widehat{\Gamma}(\tau_{\geq 0} R)$ and applying the linearization map twice, we obtain a map
    $\beta_e \colon  \omega_{\widehat{\Gamma}}[2] \to R$, called the \emph{Bott map} of $e$.
    We say that the pre-orientation $e$ is an \emph{orientation} if the map $\beta_e$ is an isomorphism.
\end{defn}

\begin{rmk}
    If $\widehat{\Gamma}$ is an oriented formal group over an $\Einf$-ring $R$, then the Bott map is an equivalence $R[-2] \simeq \omega_{\widehat{\Gamma}}$.
    Letting $\widehat{\Gamma}^{\mathrm{cl}}$ denote the underlying formal group of $\widehat{\Gamma}$ over $\pi_0 (R)$, we find that the induced map
    $\pi_n (R) \otimes_{\pi_0 (R)} \omega_{\widehat{\Gamma}^{\mathrm{cl}}} \to \pi_{n+2}(R)$ is an isomorphism for all $n \in \Z$.
    In particular, there are canonical isomorphisms $\pi_{2n} (R) \cong \omega_{\widehat{\Gamma}^{\mathrm{cl}}} ^{\otimes n}$ for $n \in \Z$.
\end{rmk}

\begin{construction}[{\cite[Proposition 4.3.13 and Definition 4.3.14]{ECII}}]
    Let $\widehat{\Gamma}$ denote a formal group over an $\Einf$-ring $R$.
    Then there is an $R$-algebra $\mathfrak{O}_{(R,\Gamma)}$ that represents the functor taking an $R$-algebra $A$ to the anima of orientations of the base change of $\widehat{\Gamma}$ to $A$. We call $\mathfrak{O}_{(R,\widehat{\Gamma})}$ the \emph{orientation classifier} of $\widehat{\Gamma}$ over $R$.

    Indeed, the $R$-algebra $\mathfrak{P}_{(R,\widehat{\Gamma})} \coloneqq R \otimes_{R \otimes_{O_{\widehat{\Gamma}}} R} R$ classifies preorientations of $\widehat{\Gamma}$. In particular, the tautological preorientation $e$ determines a Bott map
    \[\beta_e \colon  \omega_{\widehat{\Gamma}_{\mathfrak{P}_{(R,\widehat{\Gamma})}}}[2] \to \mathfrak{P}_{(R,\widehat{\Gamma})}.\]
    Since $\omega_{\widehat{\Gamma}_{\mathfrak{P}_{(R,\widehat{\Gamma})}}}[2]$ is a line bundle, there is a universal $\mathfrak{P}_{(R,\widehat{\Gamma})}$-algebra $\mathfrak{O}_{(R,\widehat{\Gamma})}$ where the base change of $\beta_e$ becomes an isomorphism
    by \cite[Proposition 4.3.17]{ECII}.
    By definition, $\mathfrak{O}_{(R,\widehat{\Gamma})}$ classifies orientations of $\widehat{\Gamma}$, as desired.
\end{construction}

\begin{rmk}
    By \cite[Proposition 4.3.21]{ECII}, a preorientation of a formal group $\widehat{\Gamma}$ over a complex periodic $\Einf$-ring $R$ is equivalent to a map $\Gamma^Q _R \to \widehat{\Gamma}$ of formal groups over $R$.

    By \cite[Proposition 4.3.23]{ECII}, a preorientation of a formal group $\widehat{\Gamma}$ over an $\Einf$-ring is an orientation if and only if $R$ is complex periodic and the corresponding map $\Gamma^Q _R \to \widehat{\Gamma}$ is an isomorphism.
\end{rmk}

We need one more condition before we may define the Lubin--Tate theory construction.

\begin{defn}[{\cite[Definition 6.4.1]{ECII}}]
    Let $\widehat{\Gamma}$ denote a formal group over an $\Einf$-ring $R$.
    We say that $\widehat{\Gamma}$ is \emph{balanced} if the map $R \to \mathfrak{O}_{(R,\widehat{\Gamma})}$ to the orientation classifier ${\mathfrak{O}_{(R,\widehat{\Gamma})}}$ of $\widehat{\Gamma}$ satisfies the following two conditions:
    \begin{enumerate}
        \item The map $\pi_0 R \to \pi_0 {\mathfrak{O}_{(R,\widehat{\Gamma})}}$ is an isomorphism.
        \item The odd homotopy groups of ${\mathfrak{O}_{(R,\widehat{\Gamma})}}$ vanish.
    \end{enumerate}
\end{defn}

\begin{defn} \label{def:LT}
    Let $R$ denote an $\F_p$-algebra with almost perfect cotangent complex $L^{\Einf} _R$ and let $\Gamma$ denote a nonstationary $p$-divisible group over $R$. If $\Gamma^{\mathrm{un},\circ}$ over $R^{\mathrm{un}}_{(R,\Gamma)}$ is balanced, then we define the \emph{Lubin--Tate theory} of the pair $(R,\Gamma)$ to be:
    \[E(R,\Gamma) \coloneqq \mathfrak{O}_{(R^{\mathrm{un}} _{(R,\Gamma)},\Gamma^{\mathrm{un}, \circ})}.\]
    In this case, there are canonical isomorphisms
    \[\pi_n E(R,\Gamma) \cong \begin{cases} \omega_{\Gamma^{\mathrm{un}, \mathrm{cl}}}^{\otimes n/2} & \text{ if $n$ even}, \\
    0 &\text{ if $n$ odd}. \end{cases}\]
    In particular, $\pi_0 E(R,\Gamma) \cong R^{\mathrm{cl}}_{(R,\Gamma)}$.
\end{defn}

In the case that $R$ is Noetherian, Lurie has proven that $\Gamma^{\mathrm{un},\circ}$ is always balanced. 

\begin{thm}[{\cite[Theorem 6.4.6]{ECII}}] \label{thm:noe-balanced}
    Let $R$ denote a Noetherian $F$-finite $\F_p$-algebra, and let $\Gamma$ denote a nonstationary $p$-divisible group over $R$. Then the formal group $\Gamma^{\mathrm{un}, \circ}$ over $R^{\mathrm{un}}_{(R,\Gamma)}$ is balanced.
\end{thm}

As a consequence, one may assign Lubin--Tate theories $E(R,\Gamma)$ to pairs $(R, \Gamma)$ where $R$ is a Noetherian and $F$-finite $\F_p$-algebra and $\Gamma$ is a nonstationary $p$-divisible group over $R$.

We conclude this section by discussing the universal property of Lubin--Tate theory. 

\begin{rmk}[{\cite[Remark 6.0.7]{ECII}}] \label{rmk:univ}
  In the context of \Cref{def:LT}, let $I = \ker(R^{\cl}_{(R,\Gamma)} \to R)$.
    Then the Lubin--Tate theory $E(R,\Gamma)$ is a complete adic $\Einf$-ring when endowed with the $I$-adic topology, and has a universal property as a complete adic $\Einf$-ring.
    Namely, given a complete adic $\Einf$-ring $A$, the anima of maps
    $\Map_{\CAlg_\ad} (E(R,\Gamma), A)$
    is identified with the classifying anima of triples $(\widetilde{\Gamma}, \alpha, e)$, where $\widetilde{\Gamma}$ is a $p$-divisible group over $A$, $\alpha$ is an equivalence class of $\Gamma$-taggings on $\widetilde{\Gamma}$, and $e$ is an orientation of $\widetilde{\Gamma}^\circ$.

     A \emph{$\Gamma$-tagging} on $\widetilde{\Gamma}$ is a pair $(I, \alpha)$ of a finitely generated ideal of definition $I \subseteq \pi_0 (A)$ and a map of pairs $\alpha \colon (R,\Gamma) \to (\pi_0 (A) / I, \widetilde{\Gamma}_{\pi_0 (A) / I})$.
    Two taggings $(I_i, \alpha_i)$ are equivalent if there exists a finitely generated ideal of definition $I_3$ containing $I_1$ and $I_2$ such that the composites
    \[(R,\Gamma) \xrightarrow{\alpha_i} (\pi_0 (A) / I_i, \widetilde{\Gamma}_{\pi_0 (A) / I_i}) \to (\pi_0 (A) / I_3, \widetilde{\Gamma}_{\pi_0 (A) / I_3})\]
    are equal for $i = 1,2$.
\end{rmk}

\subsection{Luriean pairs}
We define \emph{Luriean pairs} as an axiomatization of the input data required for Lurie's construction of Lubin--Tate theories, and prove a permanence property of Luriean pairs that will be key in our construction of non-Noetherian Lubin--Tate theories.

\begin{defn}
    We say that a pair $(R, \Gamma)$ of an $\F_p$-algebra $R$ and a $p$-divisible group $\Gamma$ over $R$ is \emph{Luriean} if the following conditions are satisfied:
    \begin{enumerate}
        \item The cotangent complex $L_{R}$ is almost perfect as an $R$-module.
        \item The $p$-divisible group $\Gamma$ is nonstationary.
        \item The identity component $\Gamma^{\mathrm{un}, \circ}$ of the universal deformation $\Gamma^{\mathrm{un}}$ over $R^{\mathrm{un}}_{\Gamma_0}$ is balanced.
    \end{enumerate}
\end{defn}

It follows from the discussion at the end of the previous section that we may assign a well-behaved Lubin--Tate theory $E(R,\Gamma)$ to a Luriean pair $(R,\Gamma)$.

\begin{rmk} \label{rmk:Noe-good}
    By \Cref{rmk:noe-f-finite} and \Cref{thm:noe-balanced}, a pair $(R,\Gamma)$ with $R$ Noetherian is Luriean if and only if $R$ is $F$-finite and $\Gamma$ is nonstationary.
\end{rmk}

The main goal of this section is to prove the following theorem.

\begin{thm} \label{thm:base-change-Luriean1}
    Let $(R, \Gamma)$ denote a Luriean pair with $R$ Noetherian, and let
    $R \to S$ denote a flat map with $L^\alg _{S/R} \simeq 0$.
    Then $(S, \Gamma_S)$ is also a Luriean pair.
\end{thm}

The key technical input to this theorem is the following proposition.

\begin{proposition} \label{prop:un-flat1}
    Suppose that $(R, \Gamma)$ is a Luriean pair with $R$ Noetherian.
    Let $R \to S$ denote a flat map for which $L_{S/R} ^{\alg} \simeq 0$.
    Then the induced map $R_{(R,\Gamma)} ^{\mathrm{un}} \to R_{(S,\Gamma_S)} ^{\mathrm{un}}$ is flat.
\end{proposition}

Given this proposition, we may prove \Cref{thm:base-change-Luriean1}.

\begin{proof}[Proof of \Cref{thm:base-change-Luriean1}]
    Since $L^{\alg}_{S/R} \simeq 0$, we also have $L^{\Einf} _{S/R} \simeq 0$ by \Cref{prop:cot-compare}.
    It follows that $L^{\Einf}_{S} \simeq L^{\Einf}_{R} \otimes_{R} S$, which is almost perfect as an $S$-module because it is the base change of an almost perfect $R$-module.
    Similarly, we have $L^{\Einf}_{\Spec(S) / \mathcal{M}_{\mathrm{BT}}} \simeq L^{\Einf}_{\Spec(R) / \mathcal{M}_{\mathrm{BT}}} \otimes_{R} S$, so that $L^{\Einf}_{\Spec(S) / \mathcal{M}_{\mathrm{BT}}}$ is $1$-connective, i.e., $\Gamma_S$ is nonstationary.

    What remains is to show that the identity component $(\Gamma_S)^{\mathrm{un}, \circ}$ is balanced.
    However, $(\Gamma_S)^{\mathrm{un}}$ is base changed from $\Gamma^{\mathrm{un}}$ along the map $R_{(R,\Gamma)} ^{\mathrm{un}} \to R_{(S,\Gamma_S)} ^{\mathrm{un}}$.
    Since this map is flat by \Cref{prop:un-flat1} and $\Gamma^{\mathrm{un, \circ}}$ is balanced by assumption, this implies that $(\Gamma_S)^{\mathrm{un},\circ}$ is balanced by \cite[Remark 6.4.2]{ECII}.
\end{proof}

In the remainder of this section, our goal will be to prove \Cref{prop:un-flat1}.
This will be accomplished by combining the following two lemmas.

\begin{lemma}\label{lem:base-change-identity-un}
    Suppose that $(R,\Gamma)$ is a pair consisting of an $\F_p$-algebra $R$ with the property that $L^{\Einf}_R$ is almost perfect and a nonstationary $p$-divisible group $\Gamma$ over $R$.
    Suppose that $R \to S$ is a flat map such that $L^{\alg} _{S/R} \simeq 0$.
    Then $L_S ^{\Einf}$ is almost perfect as an $S$-module and $\Gamma_S$ is nonstationary, so that $R^{\mathrm{un}}_{(S,\Gamma_S)}$ exists and the map
    \[R \otimes_{R^{\mathrm{un}}_{(R,\Gamma)}} R^{\mathrm{un}}_{(S,\Gamma_S)} \to S\]
    is an isomorphism.
\end{lemma}

\begin{lemma}[Lurie] \label{lem:Neo-check-flat}
    Let $R \to S$ be a morphism of connective $\Einf$-rings, and assume that $R$ is Noetherian. Assume that $R$ is complete with respect to an ideal $I \subseteq \pi_0 R$ and that $S$ is complete with respect to the ideal $I \pi_0 (S) \subseteq \pi_0 (S)$.
    Then the following are equivalent:
    \begin{enumerate}
        \item The morphism $R \to S$ is flat.
        \item The induced morphism $\pi_0 (R)/ I \to (\pi_0 (R)) / I \otimes_R^{\mathbb{L}} S$ is flat.
    \end{enumerate}
\end{lemma}

\begin{proof}[Proof of \Cref{lem:Neo-check-flat}]
    This is \cite[Lemma 6.1.9]{ECII}, with the assumption that $S$ be Noetherian removed. This assumption was never used in Lurie's proof of \cite[Lemma 6.1.9]{ECII}.
\end{proof}

\begin{rmk}
    It follows from \Cref{lem:base-change-identity-un} that the square
    \begin{center}
        \begin{tikzcd}
            R_{(R,\Gamma)} ^{\mathrm{un}} \ar[r] \ar[d] & R_{(S,\Gamma_S)} ^{\mathrm{un}} \ar[d] \\
            R \ar[r] & S
        \end{tikzcd}
    \end{center}
    is coCartesian and that $R_{(S,\Gamma_S)} ^{\mathrm{un}}$ is $I$-complete, where $I$ is the kernel of $R^{\mathrm{cl}}_{(R,\Gamma)} \to R$.
    It therefore follows from \cite[Proposition 5.2.9]{ECII} that $R_{(R,\Gamma)} ^{\mathrm{un}} \to R_{(S,\Gamma_S)} ^{\mathrm{un}}$ is an $R_{(R,\Gamma)} ^{\mathrm{un}}$-thickening of $R \to S$ in the sense of \cite[Definition 5.2.1]{ECII}.\footnote{Lurie's assumption that $F\colon A_0 \to A_0$ is flat and that $A_0 \to B_0$ is relatively perfect may be replaced by the assumption that $L^{\Einf}_{B_0/A_0} \simeq 0$, which is all that he makes use of.}
\end{rmk}

All that remains now is the proof of \Cref{lem:base-change-identity-un}. Before we prove it, we need two preliminary lemmas.

\begin{lemma} \label{lem:cotangent-vanish}
    Let $\Gamma$ denote a nonstationary $p$-divisible group over an $\F_p$-algebra $R$.
    Suppose that $R \to S$ is a flat morphism, that $L^{\alg} _{R}$ and $L^{\alg} _{S}$ are almost perfect, and that $\Gamma_S$ is nonstationary.
    Then we have
    \[S \otimes_{R^{\mathrm{un}}_{(S,\Gamma_S)}} L^{\Einf}_{R^{\mathrm{un}}_{(S,\Gamma_S)}/R^{\mathrm{un}}_{(R,\Gamma)}} \simeq S \otimes_{R \otimes_{R^{\mathrm{un}}_{(R,\Gamma)}} R^{\mathrm{un}}_{(S,\Gamma_S)}} L^{\Einf}_{R \otimes_{R^{\mathrm{un}}_{(R,\Gamma)}} R^{\mathrm{un}}_{(S,\Gamma_S)}/ R} \simeq 0.\]
\end{lemma}

\begin{proof}
    It follows from the definition of the spectral deformation ring that
    \[L^{\Einf}_{\Spf(R^{\mathrm{un}}_{(R,\Gamma)})/\mathcal{M}_{\mathrm{BT}}} \simeq L^{\Einf}_{\Spf(R^{\mathrm{un}}_{(S,\Gamma_S)}) / \mathcal{M}_{\mathrm{BT}}} \simeq 0,\]
    so that we have $L^{\Einf}_{\Spf(R^{\mathrm{un}}_{(S,\Gamma_S)})/\Spf(R^{\mathrm{un}}_{(R,\Gamma)})} \simeq 0$. (See \cite[Proof of Lemma 6.1.7]{ECII}.)
    Letting $I \subseteq R^{\mathrm{cl}}_{(S,\Gamma_S)}$ denote the kernel of $R^{\mathrm{cl}}_{(S,\Gamma_S)} \to S$, this implies that the $I$-completion of $L^{\Einf}_{R^{\mathrm{un}}_{(S,\Gamma_S)}/R^{\mathrm{un}}_{(R,\Gamma)}}$ is trivial by \cite[Example 17.1.2.9]{SAG}, hence $S \otimes_{R^{\mathrm{un}}_{(S,\Gamma_S)}} L^{\Einf}_{R^{\mathrm{un}}_{(S,\Gamma_S)}/R^{\mathrm{un}}_{(R,\Gamma)}} \simeq 0$, as desired.
    
    Finally, the equivalence
    $S \otimes_{R^{\mathrm{un}}_{(S,\Gamma_S)}} L^{\Einf}_{R^{\mathrm{un}}_{(S,\Gamma_S)}/R^{\mathrm{un}}_{(R,\Gamma)}} \simeq S \otimes_{R \otimes_{R^{\mathrm{un}}_{(R,\Gamma)}} R^{\mathrm{un}}_{(S,\Gamma_S)}} L^{\Einf}_{R \otimes_{R^{\mathrm{un}}_{(R,\Gamma)}} R^{\mathrm{un}}_{(S,\Gamma_S)}/ R}$
    is a consequence of the base change equivalence $L^{\Einf}_{R \otimes_{R^{\mathrm{un}}_{(R,\Gamma)}} R^{\mathrm{un}}_{(S,\Gamma_S)}/ R} \simeq R \otimes_{R^{\mathrm{un}}_{(R,\Gamma)}} L^{\Einf}_{R^{\mathrm{un}}_{(S,\Gamma_S)}/R^{\mathrm{un}}_{(R,\Gamma)}}$.
\end{proof}

\begin{lemma} \label{lem:ideal-generation}
    Let $R$ denote an $\F_p$-algebra for which $L^{\Einf} _R$ is almost perfect, and let $\Gamma$ denote a nonstationary $p$-divisible group over $R$.
    Suppose that $R \to S$ is a flat map with the property that $L_{S/R}^{\alg} \simeq 0$, and let $I = \ker(R^{\mathrm{cl}}_{(R,\Gamma)} \to R)$ and $J = \ker(R^{\mathrm{cl}}_{(S,\Gamma_S)} \to S)$. The image of $I$ in $R^{\mathrm{cl}}_{(S,\Gamma_S)}$ generates $J$.
\end{lemma}

This lemma, in turn, depends on two smaller lemmas.

\begin{lemma} \label{lem:ideal-complete}
    Let $R$ denote a commutative ring that is derived $I$-complete with respect to a finitely generated ideal $I \subset R$. 
    Then any finitely generated ideal $J \subset R$ is derived $I$-complete.
\end{lemma}

\begin{proof}
    By the short exact sequence $0 \to J \to R \to R/J \to 0$, it is equivalent to prove that $R/J$ is derived $I$-complete. Fix generators $x_1, \dots, x_k \in J$ and let $\Z[x_1, \dots, x_k] \to R$ denote the map sending each $x_i$ to the corresponding element of $J$.
    Then $R \otimes_{\Z[x_1, \dots, x_k]}^{\mathbb{L}} \mathbb{Z}$ is a perfect $R$-module, hence derived $I$-complete, with $\pi_0 (R \otimes_{\Z[x_1, \dots, x_k]}^{\mathbb{L}} \mathbb{Z}) \cong R/J$. It therefore follows from \cite[Theorem 7.3.4.1]{SAG} that $R/J$ is derived $I$-complete.
\end{proof}

\begin{lemma} \label{lem:indecomp-cotangent}
    Let $R$ denote an $\F_p$-algebra with $L^{\Einf}_R$ almost perfect, and let $\Gamma$ denote a nonstationary $p$-divisible group over $R$.
    Let $I$ denote the kernel of $R^{\mathrm{cl}} _{(R,\Gamma)} \to R$.
    Then there is a canonical identification
    \[\pi_1 L^{\Einf} _{\Spec(R) / \mathcal{M}_{\mathrm{BT}}} \cong I / I^2.\]
\end{lemma}

\begin{proof}
    Since $L^{\Einf}_{\Spf (R^{\mathrm{un}}_{(R,\Gamma)})/\mathcal{M}_{\mathrm{BT}}} \simeq 0$, we have $L^{\Einf} _{\Spec(R) / \mathcal{M}_{\mathrm{BT}}} \simeq L^{\Einf} _{\Spec(R) / \Spf (R^{\mathrm{un}}_{(R,\Gamma)})}$.
    Since any $R$-module is automatically $I$-complete as a $R^{\mathrm{un}}_{(R,\Gamma)}$-module, it follows from \cite[Example 17.1.2.9]{SAG} that $L^{\Einf} _{\Spec(R) / \Spf (R^{\mathrm{un}}_{(R,\Gamma)})} \simeq L^{\Einf} _{\Spec(R) / \Spec (R^{\mathrm{un}}_{(R,\Gamma)})}$.
    Finally it follows from \cite[Theorem 7.4.3.1]{HA} that there is a natural isomorphism
    \[
        \pi_1 L^{\Einf} _{\Spec(R) / \Spec (R^{\mathrm{un}}_{(R,\Gamma)})} \simeq I/I^2. \qedhere
    \]
\end{proof}

\begin{proof}[Proof of \Cref{lem:ideal-generation}]
    Let $\widetilde{I} \subseteq J$ denote the ideal generated by the image of $I$ in $R^{\mathrm{cl}}_{(S,\Gamma_S)}$.
    It follows from \Cref{lem:ideal-complete} that $\widetilde{I}$ and $J$ are derived $J$-complete.
    Now, by \Cref{lem:indecomp-cotangent} and the fact that $L_{S/R} \simeq 0$, we have that
    \[(I / I^2) \otimes_{R} S \to J/J^2\]
    is an isomorphism.
    In particular, $\widetilde{I} / \widetilde{I} J \to J/ J^2$ is a surjection.
    Since $\widetilde{I}$ and $J$ are derived $J$-complete, it follows from derived Nakayama's lemma \cite[Tag 09B9]{stacks-project} that $\widetilde{I} \to J$ is a surjection, as desired.
\end{proof}

Using \Cref{lem:cotangent-vanish} and \Cref{lem:ideal-generation}, we prove \Cref{lem:base-change-identity-un}.

\begin{proof}[Proof of \Cref{lem:base-change-identity-un}]
    Since $L_{S/R}^{\alg} \simeq 0$, we have $L_R ^{\Einf} \otimes_R S \simeq L_S ^{\Einf}$, whence it follows that $L^{\Einf} _S$ is almost perfect and $\Gamma_S$ is nonstationary.

    To prove that
    \[R \otimes_{R^{\mathrm{un}}_{(R,\Gamma)}} R^{\mathrm{un}}_{(S,\Gamma_S)} \to S\]
    is an isomorphism, it suffices by \cite[Corollary 7.4.3.4]{HA} to show that it is an isomorphism on $\pi_0$ and has a trivial relative cotangent complex.
    On $\pi_0$, this is a consequence of \Cref{lem:ideal-generation}.
    On the other hand, to prove that
    $L^{\Einf} _{S/R \otimes_{R_{(R,\Gamma)} ^{\mathrm{un}}} R_{(S,\Gamma_S)} ^{\mathrm{un}}} \simeq 0,$
    we make use of the cofiber sequence
    \[S \otimes_{R \otimes_{R^{\mathrm{un}}_{(R,\Gamma)}} R^{\mathrm{un}}_{(S,\Gamma_S)}} L^{\Einf} _{R \otimes_{R^{\mathrm{un}}_{(R,\Gamma)}} R^{\mathrm{un}}_{(S,\Gamma_S)}/ R}  \to L^{\Einf} _{S / R} \to L^{\Einf} _{S/R \otimes_{R^{\mathrm{un}}_{(R,\Gamma)}} R^{\mathrm{un}}_{(S,\Gamma_S)}}.\]
    The middle term vanishes by assumption, and the first term vanishes by \Cref{lem:cotangent-vanish}, hence the final term vanishes and we are done.
\end{proof}

Finally, we use \Cref{lem:base-change-identity-un} and \Cref{lem:Neo-check-flat} to prove \Cref{prop:un-flat1}.

\begin{proof}[Proof of \Cref{prop:un-flat1}]
    By \Cref{lem:Neo-check-flat}, it suffices to show that the induced map
    \[R \to R \otimes_{R_{(R,\Gamma)} ^{\mathrm{un}}} R_{(S,\Gamma_S)} ^{\mathrm{un}}\]
    is flat.
    To prove this, it will suffice to prove that the map
    \[R \otimes_{R_{(R,\Gamma)} ^{\mathrm{un}}} R_{(S,\Gamma_S)} ^{\mathrm{un}} \to S\]
    is an isomorphism, since $R \to S$ is flat by assumption.
    This is the statement of \Cref{lem:base-change-identity-un}.
\end{proof}

We record the following consequence of the above results, since we will find it useful in the sequel:

\begin{lemma} \label{lem:flat-ideal-gen}
    Suppose that $(R, \Gamma)$ is a Noetherian Luriean pair.
    Let $R \to S$ denote a flat map of rings for which $L^{\alg}_{S/R} \simeq 0$. Then the induced map
    $E(R, \Gamma) \to E(S, \Gamma_S)$ is flat.
    Moreover, the induced map on $R^{\mathrm{cl}}_{(R,\Gamma)} \to R^{\mathrm{cl}}_{(S,\Gamma_S)}$ sends the kernel of $R^{\mathrm{cl}}_{(R,\Gamma)} \to R$ to generators for the kernel of $R^{\mathrm{cl}}_{(S,\Gamma_S)} \to S$.
\end{lemma}

\begin{proof}
    Since the pullback of $\Gamma^{\mathrm{un},\circ}$ along $R^{\mathrm{un}}_{(R,\Gamma)} \to R^{\mathrm{un}} _{(S,\Gamma_S)}$ is canonically isomorphic to $(\Gamma_S)^{\mathrm{un}, \circ}$, the square
    \begin{center}
        \begin{tikzcd}
            R^{\mathrm{un}}_{(R,\Gamma)} \ar[r] \ar[d] & R^{\mathrm{un}} _{(S,\Gamma_S)} \ar[d] \\
            E(R,\Gamma) \ar[r] & E(S, \Gamma_S)
        \end{tikzcd}
    \end{center}
    is a pushout of $\Einf$-rings. The flatness of the bottom is therefore a consequence of the flatness of the top arrow, proven in \Cref{prop:un-flat1}.
    The second statement is a consequence of \Cref{lem:ideal-generation}.
\end{proof}

\subsection{$\LKN$ pairs}
In this section, we introduce \emph{Luriean Kalypto--Noetheriean pairs}, which are a subclass of Luriean pairs that has better closure properties. We will frequently make use of the abbreviation \emph{$\LKN$ pairs}.

\begin{defn}[Relatively Perfect] \label{def:rel_perfect}
A map $f\colon R \to S$ of $\mathbb{F}_p$-algebras is called \emph{relatively perfect} if the diagram
\begin{equation}
\begin{tikzcd}
R \arrow[r, "f"] \arrow[d, "F_R"'] & S \arrow[d, "F_S"] \\
R \arrow[r, "f"] & S
\end{tikzcd}
\end{equation}
is a pushout in the category of commutative rings.
In other words, if the relative Frobenius
\[F_{S/R}\colon R^{F} \otimes_R \!S\to \!S^{F}\] is an isomorphism. 
\end{defn}

\begin{defn}[Kalypto--Noetherian] \label{def:kalypto}
We say that an $\mathbb{F}_p$-algebra $R$ is \emph{Kalypto--Noetherian} if there exists an $F$-finite Noetherian $\mathbb{F}_p$-algebra $R_0$ and a flat, relatively perfect map $R_0 \to R$.\footnote{We take the term Kalypto from the ancient Greek \textgreek{καλύπτω}: \emph{to cover, hide or conceal}. For example, in book 8, line 503 of the Odyssey, one finds the Achaeans \textgreek{κεκαλυμμένοι ἵππῳ}: \emph{concealed in the [Trojan] horse}. Like the Achaeans in the Trojan horse, Kalypto--Noetherian rings covertly smuggle in the properties of Noetherian $F$-finite $\F_p$-algebras that we will require.}
\end{defn}

We emphasize that the ring $R_0$ is not part of the data---it is merely required to exist. We call $R_0 \to R$ a witness to Kalypto--Noetherianity of $R$.

We now recall some basic properties of relatively perfect maps and use them to prove some closure properties of the category of Kalypto--Noetherian rings.

\begin{lemma} \label{lem:rel_perfect_props}
\begin{enumerate}
    \item If $f\colon R \to S$ is weakly \'etale, for example if it is a filtered colimit of \'etale maps, then it is relatively perfect.
    \item If $R_0 \to R_1$ and $R_1\to R_2$ are both relatively perfect, then so is their composite $R_0 \to R_2$.
    \item Let $k$ be a perfect $\F_p$-algebra with maps $k\to R_0$ and $k\to S_0$ and let $R_0 \to R_1$ and $S_0 \to S_1$ be relatively perfect, then  the map \[R_0 \otimes_k S_0 \to R_1 \otimes_k S_1\] is relatively perfect.
    \item If $R_1 \to S$ is relatively perfect and $R_1 \to R_2$ is any ring homomorphism, then the base change $R_2 \to S \otimes_{R_1}R_2$ is relatively perfect.
    \end{enumerate}
\end{lemma}

\begin{proof}
\begin{enumerate}
\item This is the statement of \cite[Tag 0F6W]{stacks-project}.
\item This is a consequence of the fact that the pasting of two pushout diagrams is a pushout diagram.
\item This follows from the fact that a tensor product of pushout diagrams is a pushout diagram. 
\item Again, this follows from the fact that a tensor product of pushout diagrams is a pushout diagram.
\qedhere
\end{enumerate}
\end{proof}

\begin{lemma}\label{lem:KN_prop}
     If $R$ is Kalypto--Noetherian and $R \to S$ is flat and relatively perfect, then $S$ is Kalypto--Noetherian.
\end{lemma}

\begin{proof}
    There exists $R_0 \to R$ flat and relatively perfect with $R_0$ $F$-finite Noetherian. If $R \to S$ is flat and relatively perfect, the composition $R_0 \to S$ is also flat and relatively perfect (by \cref{lem:rel_perfect_props}).
\end{proof}

Finally, we have the following key lemma, which we will use to show that Kalypto--Noetherianity is inherited by deformation rings.

\begin{lemma} \label{lem:check-flat-rel-perf}
    Let $R \to S$ denote a map of $\F_p$-algebras where $R$ is Noetherian and $F$-finite, and let $I \subset R$ denote a finitely generated ideal. Suppose that both $R$ and $S$ are complete with respect to $I,$
    and that the induced map $R / I \to S \otimes_R ^{\mathbb{L}} R/ I$ is flat and relatively perfect.
    Then the original map $R \to S$ is flat and relatively perfect.
\end{lemma}

\begin{proof}
    Flatness of $R \to S$ is a consequence of \Cref{lem:Neo-check-flat}.
    It remains to check relative perfectness.
    We wish to prove that the relative Frobenius
    \[F_{S/R}\colon R^F \otimes_{R} S \to S^F\]
    is an isomorphism.
    As modules over $R^F$, I claim that both sides are $I^F$-complete.
    This is clear for the right hand side. For the left-hand side, I first claim that it is $F(I)$-complete.
    Indeed, this is a consequence of the fact that $S$ is $I$-complete as an $R$-module and $R^F$ is almost perfect as an $R$-module (since $R$ is Noetherian and $F$-finite).
    Since $I$ is finitely generated, this implies that it is also $I^F$-complete, as desired.

    It therefore suffices to check that $R^F / I^F \otimes_{R/I} S/IS \cong R^F / I^F \otimes _R S \to S^F/I^FS^F$ is an isomorphism.
    This is a consequence of the assumption that $R / I \to S \otimes_R ^{\mathbb{L}} R/ I \cong S/IS$ is relatively perfect.
\end{proof}

\begin{defn}
    We say that a pair $(R,\Gamma)$ of an $\F_p$-algebra $R$ and a $p$-divisible group $\Gamma$ over $R$ is \emph{Luriean Kalypto--Noetherian} if there exists a map of pairs $(R_0, \Gamma_0) \to (R,\Gamma)$ for which $R_0$ is Noetherian and $F$-finite, $\Gamma_0$ is nonstationary, and $R_0 \to R$ is flat and relatively perfect.
\end{defn}

In particular, $R$ is a Kalypto--Noetherian ring.
Moreover, the following is an immediate consequence of \Cref{thm:base-change-Luriean1}.
\begin{thm} \label{thm:LKN-is-Luriean}
    Any Luriean Kalypto--Noetherian pair $(R,\Gamma)$ is a Luriean pair.
\end{thm}

We will abbreviate Luriean Kalypto--Noetherian pairs to $\LKN$ pairs in the remainder of the paper.

\begin{defn}
  We let $\LKN$ denote the category of $\LKN$-pairs: a morphism $(R_1, \Gamma_1) \to (R_2, \Gamma)$ consists of a map $R_1 \to R_2$ and an isomorphism $(\Gamma_1)_{R_2} \cong \Gamma_2$.
\end{defn}

\begin{proof}[Proof of \Cref{thm:Discrete-LKN-functor}]
  The oriented deformation ring construction recollected in \Cref{sec:LT-back} defines a functor
  \[E(-)\colon \LKN \to \CAlg.\]
  Since any $\LKN$ pair is Luriean by \Cref{thm:LKN-is-Luriean}, it has the desired properties.
\end{proof}

We conclude this subsection with several useful results about $\LKN$ pairs. We do not know how to prove the analogous results in the generality of Luriean pairs.

\begin{proposition} \label{prop:LKN-pair-deformed}
    Let $(R,\Gamma)$ denote an $\LKN$ pair. Then the pair $(R_{(R,\Gamma)}^{\mathrm{cl}}/p, \Gamma^{\mathrm{un}}_{R_{(R,\Gamma)}^{\mathrm{cl}}/p})$ is again an $\LKN$ pair, and the induced map $E(R_{(R,\Gamma)}^{\mathrm{cl}}/p, \Gamma^{\mathrm{un}}_{R_{(R,\Gamma)}^{\mathrm{cl}}/p}) \to E(R,\Gamma)$ is an equivalence of $\Einf$-rings.
\end{proposition}

\begin{proof}
    Let $(R_0, \Gamma_0) \to (R,\Gamma)$ be a witness to $(R,\Gamma)$ being an $\LKN$ pair.
    Then the map $R^{\mathrm{cl}} _{(R_0,\Gamma_0)} \to R^{\mathrm{cl}}_{(R,\Gamma)}$ is flat by \Cref{prop:un-flat1}, hence the induced map $R_{(R_0,\Gamma_0)}^{\mathrm{cl}}/p \to R_{(R,\Gamma)}^{\mathrm{cl}}/p$ is flat.
    Moreover, $R_{(R_0,\Gamma_0)}^{\mathrm{cl}}/p$ is Noetherian and $F$-finite and $\Gamma^{\mathrm{un}}_{R_{(R,\Gamma)}^{\mathrm{cl}}/p}$ is nonstationary by \cite[Proposition 6.1.5]{ECII}.
    It therefore suffices to show that this map is relatively perfect.
    This follows from combining \Cref{lem:check-flat-rel-perf} with \Cref{lem:flat-ideal-gen}, the latter of which shows that the following square is cocartesian:
    \begin{center}
        \begin{tikzcd}
            R_{(R_0,\Gamma_0)}^{\mathrm{cl}}/p \ar[r] \ar[d] & R_{(R,\Gamma)}^{\mathrm{cl}}/p \ar[d] \\
            R_0 \ar[r] & R. 
        \end{tikzcd}
    \end{center}
    The final statement is a consequence of \cite[Example 6.1.4]{ECII}, since the $F$-finite and Noetherian hypotheses are not necessary, just existence of the deformation ring.
\end{proof}

\begin{lemma} \label{lem:vi-regular}
    Given an $\LKN$ pair $(R, \Gamma)$, 
     let $v_0 = p,v_1,v_2,\cdots $ be a sequence of elements $v_i \in \omega_{R^{\mathrm{cl}}_{(R,\Gamma)}}^{p^i-1}$ such that 
     $v_i$ is the coefficient of $X^{p^i}$ in the $p$-series $[p]_{\Gamma^{\un, \circ}}(X)$ mod $v_0,\dots,v_{i-1}$. 
     Then the sequence $v_0 = p,v_1,v_2,\cdots $
    is regular in the graded ring $\bigoplus_{k \in \mathbb{Z}} \omega_{(R^{\mathrm{cl}}_{R,\Gamma)}}^k$.
    In other words, the $\Einf$-ring $E(R,\Gamma)$ is Landweber flat.
\end{lemma}

\begin{proof}
    Since $(R, \Gamma)$ is $\LKN$, there exists a map of pairs
    $(R_0, \Gamma_0) \to (R, \Gamma)$ with $(R_0, \Gamma_0)$ a Noetherian Luriean pair and $R_0 \to R$ flat and relatively perfect.
    By \Cref{prop:un-flat1}, the induced map $R^{\mathrm{un}} _{(R_0,\Gamma_0)} \to R^{\mathrm{un}} _{(R,\Gamma)}$, and hence the map $R^{\cl} _{(R_0,\Gamma_0)} \to R^{\cl} _{(R,\Gamma)}$, is flat, so it suffices to consider the case where $R$ is Noetherian and $F$-finite.
    
    Next, by replacing $R$ by $R^{\cl} _{(R,\Gamma)} / p$ and $\Gamma$ by $\Gamma^{\mathrm{un}}_{R^{\cl} _{(R,\Gamma)} / p}$, the classical deformation ring remains unchanged by \cite[Proposition 6.1.5]{ECII}. We may therefore assume that $R = R^{\cl} _{(R,\Gamma)} / p$.

    It suffices to show that the sequence $p, v_1, v_2, \dots$ is regular after localizing at every prime ideal $\mathfrak{p} \subseteq R^{\cl} _{(R,\Gamma)}$. If $p \notin \mathfrak{p}$, then $(R^{\cl} _{(R,\Gamma)})_{\mathfrak{p}}$ is a $\mathbb{Q}$-algebra, so our sequence is obviously regular.
    Otherwise, there is a map
    \[R^{\cl} _{(R,\Gamma)} \to R^{\cl} _{(R,\Gamma)} / p = R \to \kappa(\mathfrak{p}) \to \overline{\kappa(\mathfrak{p})},\]
    where $\kappa(\mathfrak{p})$ is the fraction field of $R^{\cl} _{(R,\Gamma)}$ at $\mathfrak{p}$.
    Moreover, we have $\Omega_{\overline{\kappa(\mathfrak{p})}/(R)} = 0 $ so by \cite[Theorem 6.1.2]{ECII} we obtain  a flat map $R^{\mathrm{un}} _{(R,\Gamma)} \to R^{\mathrm{un}} _{(\overline{\kappa(\mathfrak{p})},\Gamma_{\overline{\kappa(\mathfrak{p})}})}$.
    This factors through a faithfully flat map $R^{\mathrm{un}} _{(R,\Gamma), \mathfrak{p}} \to R^{\mathrm{un}} _{(\overline{\kappa(\mathfrak{p})},\Gamma_{\overline{\kappa(\mathfrak{p})}})}.$

    It therefore suffices to deal with the case where $R$ is an algebraically closed field. In this case, $\Gamma \cong \Gamma^{\circ} \oplus \Gamma ^{et}$, and the induced map $R^{\mathrm{un}}_{(R,\Gamma ^{\circ})} \to R^{\mathrm{un}}_{(R,\Gamma)}$ is flat by \cite[Theorem 6.2.4]{ECII}. This reduces us to the case where $\Gamma$ is connected, where it is well-known.
\end{proof}

\begin{rmk}
    The proof of \Cref{lem:vi-regular} follows closely the proof of \cite[Theorem 6.4.7]{ECII}.
\end{rmk}

\begin{proposition} \label{prop:tot-preserved}
    Let $(R, \Gamma)$ denote a Noetherian Luriean pair, and let $S^{-1} \to S^\bullet$ denote an augmented cosimplicial diagram of flat and relatively perfect $R$-algebras that identifies $S^{-1}$ with the totalization of $S^\bullet$ as derived commutative rings.
    By abuse of notation, let $\Gamma$ also denote its base change along the maps $R \to S^n$, $-1 \leq n$.
    Then the map
    \[E(S^{-1},\Gamma) \to \Tot E(S^{\bullet},\Gamma)\]
    is an equivalence.
    In fact, the maps
    \[\pi_{2n} E(S^{-1},\Gamma) \to \Tot \pi_{2n} E(S^{\bullet},\Gamma)\]
    are equivalences, where the totalization is taken in $D(\Z)$.
\end{proposition}

\begin{proof}
    It suffices to prove the statement on the level of $\pi_{2n}$, since that implies the statement on the level of spectra.
    First, we reduce to the case where $n=0$. The maps
    \[\pi_{0} E(S^{k}, \Gamma) \otimes_{\pi_0 E(R,\Gamma)} \pi_{2n} E(R,\Gamma) \to \pi_{2n} E(S^k, \Gamma)\]
    are equivalences by \Cref{lem:flat-ideal-gen}. Since the $\pi_0 E(R,\Gamma)$-module $\pi_{2n} E(R,\Gamma)$ is invertible, tensoring by it preserves totalizations, so it suffices to prove the case $n=0$.

    When $n=0$, let $I = \ker(\pi_0 E(R,\Gamma) \to R)$. It follows from \Cref{lem:flat-ideal-gen} that $E(S^n, \Gamma)$ is $I$-complete and that applying $-\otimes_{\pi_0 E(R,\Gamma)}^\mathbb{L} R$ to the augmented cosimplicial diagram we are studying produces
    \[S^{-1} \to S^{\bullet},\]
    which induces an equivalence
    \[S^{-1} \to \Tot S^{\bullet}\]
    by assumption.
    Commuting the Tot past the tensor using flatness of $\pi_0 E(S^n, \Gamma)$ over $\pi_0 E(R, \Gamma)$ (by \Cref{lem:flat-ideal-gen}), we conclude that
    \[\pi_{0} E(S^{-1},\Gamma) \to \Tot \pi_{0} E(S^{\bullet},\Gamma)\]
    becomes an equivalence after applying $-\otimes_{\pi_0 E(R,\Gamma)}^{\mathbb{L}} R$.
    Since both sides are $I$-complete, we are done.
\end{proof}

\subsection{Pairs of height $h$}
We define $\LKN$ pairs of height $h$, and describe functors on the level of $\LKN$ pairs that correspond to $K(h)$-localization and the $K(h)$-localized tensor product on the level of their Lubin--Tate theories.

\begin{defn} \label{defn:LKN-ht-h}
  We say that an $\LKN$ pair $(R,\Gamma)$ is of height $h$ if there exists a Noetherian $\LKN$ pair $(R_0, \Gamma_0)$ and a map $(R_0, \Gamma_0) \to (R,\Gamma)$ such that
  \begin{itemize}
    \item $R_0 \to R$ is flat and relatively perfect,
    \item the formal group $\Gamma_0 ^\circ$ is of exact height $h$, and
    \item the kernel of the map $R^{\cl}_{(R_0,\Gamma_0)} \to R_0$ is equal to the $h$th Landweber ideal $J_h$ of $\Gamma_0^{\un, \circ}$.
  \end{itemize}
  It follows that $\Gamma$ is also of exact height $h$; moreover, by \Cref{prop:un-flat1}, the kernel of $R^{\cl} _{(R,\Gamma)} \to R$ is also equal to $J_h$. By \cite[Proposition II.4.9]{MessingCrystals}, this implies that $\Gamma^\circ$ may be regarded as a connected $p$-divisible group, and $\Gamma$ lies in a short exact sequence
    \[0 \to \Gamma^\circ \to \Gamma \to \Gamma^{\et} \to 0,\]
    where $\Gamma^{\et}$ is an \'etale $p$-divisible group.

    We let $\LKN_h$ denote the full subcategory of $\LKN$ consisting of $\LKN$ pairs of height $h$.
\end{defn}

\begin{proposition} \label{prop:ht-n-kn-local}
  Let $(R, \Gamma)$ denote an $\LKN$ pair of height $h$. Then $E(R,\Gamma)$ is $K(h)$-local. Moreover, the ideal $J_h \subseteq \pi_0 E(R,\Gamma)$ is an ideal of definition.
\end{proposition}
\begin{proof}
    The first statement is an immediate consequence of \cite[Corollary 4.5.3]{ECII}.
    The second statement is straightforward from the definitions.
\end{proof}

\subsubsection{$\LKN$ pairs and $K(h)$-localization}

\begin{defn}
    Given an $\LKN$ pair $(R, \Gamma)$, we define a new pair as follows:
    \[L_{K(h)}^{\alg} (R,\Gamma) \coloneqq (R^{\mathrm{cl}}_{(R,\Gamma)}/ J_h[v_h^{-1}], \Gamma^{\un}_{R^{\mathrm{cl}}_{(R,\Gamma)}/ J_h[v_h^{-1}]}),\]
    where $J_h$ is the $h$th Landweber ideal.
    We will sometimes write
    \[L_{K(h)}^{\alg} (R,\Gamma) = (R^{\cl,h} _{(R,\Gamma)}, \Gamma^{un,h})\]
    for the sake of brevity of notation.
\end{defn}

\begin{proposition} \label{prop:is-ht-n}
    Given an $\LKN$ pair $(R,\Gamma)$, the pair $L_{K(h)}^\alg (R,\Gamma)$ is an $\LKN$ pair of height $h$.
\end{proposition}
\begin{proof}
    We first show that $L_{K(h)}^\alg (R,\Gamma)$ is $\LKN$.
    Let $(R_0, \Gamma_0) \to (R,\Gamma)$ denote a witness for $(R, \Gamma)$ being $\LKN$.
    Then the proof of \Cref{prop:LKN-pair-deformed} implies that $R^{\mathrm{cl}}_{(R_,\Gamma_0)}/p \to R^{\mathrm{cl}} _{(R,\Gamma)}/p$ is flat and relatively perfect.
    Modding out by $v_1, \dots, v_{h-1}$ and inverting $v_h$, we obtain a new flat and relatively perfect map
    $R^{\cl,h} _{(R_0,\Gamma_0)} \to R^{\cl,h} _{(R,\Gamma)}$ by \Cref{lem:rel_perfect_props}(4).
    Moreover, it is clear that $R^{\cl,h} _{(R_0,\Gamma_0)}$ remains Noetherian and $F$-finite, and that the base change $\Gamma^{un,h}_0$ remains nonstationary.
    It therefore follows that $L_{K(h)}^{\alg} (R_0, \Gamma_0) \to L_{K(h)}^{\alg} (R,\Gamma)$ is a witness to the fact that $L_{K(h)}^{\alg} (R,\Gamma)$ is $\LKN$.

    It remains to show that the kernel of $R^{\cl}_{L_{K(h)}^{\alg} (R_0, \Gamma_0)} \to R^{\cl,h}_{(R_0, \Gamma_0)}$ is equal to $J_h$.
    First, we note that it follows from \cite[Example 6.1.3 \& Proposition 6.1.5]{ECII} that the map
    \[R^{\cl}_{(R^{\cl} _{(R_0, \Gamma_0)}/J_h, \Gamma^{\mathrm{un}}_{R^{\cl} _{(R_0, \Gamma_0)}/J_h})} \to R^{\cl} _{(R_0, \Gamma_0)}/J_h\]
    identifies with the map 
    \[(R^{\cl} _{(R_0, \Gamma_0)})^{\wedge}_{J_h} \to R^{\cl} _{(R_0, \Gamma_0)}/J_h\]
    whose kernel is equal to the ideal $J_h$.
    To conclude, we note that since $R^{\cl} _{(R_0, \Gamma_0)}/J_h \to R^{\cl,h} _{(R_0, \Gamma_0)}$ is \'etale, it follows from \Cref{lem:flat-ideal-gen} that the kernel of $R^{\cl}_{L_{K(n)}^{\alg} (R_0, \Gamma_0)} \to R^{\cl,h} _{(R_0, \Gamma_0)}$ is also equal to $J_h$, as desired.
\end{proof}

\begin{lemma} \label{lem:alg-kn-tot}
    Let $(R,\Gamma)$ denote a Noetherian Luriean pair, and let $S^{-1} \to S^\bullet$ denote an augmented cosimplicial flat and relatively perfect $R$-algebra with the property that $S^{-1} \to \Tot S^{\bullet}$ is an equivalence when computed in derived rings.
    Then the map
    \[R^{\mathrm{cl},h}_{(S^{-1},\Gamma)} \to \Tot R^{\mathrm{cl},h}_{(S^{\bullet},\Gamma)}\]
    is an equivalence, where the totalization is computed in derived rings.
\end{lemma}

\begin{proof}
    This follows from \Cref{prop:tot-preserved} by modding out by $J_h$ and inverting $v_h$, using flatness of $R^{\mathrm{cl}}_{(S^{\bullet},\Gamma)}$ over $R^{\mathrm{cl}}_{(R,\Gamma)}$ and the fact that $J_h$ is locally generated by a regular sequence in $R^{\mathrm{cl}}_{(R,\Gamma)}$ by \Cref{lem:vi-regular}.
\end{proof}

\begin{thm} \label{thm:Kh-alg-top}
    Given an $\LKN$ pair $(R,\Gamma)$, there is a canonical equivalence of $\Einf$-rings
    \[L_{K(h)} E(R,\Gamma) \simeq E(L_{K(h)}^\alg (R,\Gamma)).\]
\end{thm}

\begin{proof}
    By definition of $L_{K(h)} ^\alg$, there is a natural map of $\LKN$ pairs
    \[(R_{(R,\Gamma)}^{\mathrm{cl}}/p, \Gamma^{\un}_{R_{(R,\Gamma)}^{\mathrm{cl}}/p}) \to L_{K(h)}^\alg (R,\Gamma)\]
    which gives rise to a map of Lubin--Tate theories
    \[E(R,\Gamma) \simeq E(R_{(R,\Gamma)}^{\mathrm{cl}}/p, \Gamma^{\un}_{R_{(R,\Gamma)}^{\mathrm{cl}}/p}) \to E(L_{K(h)}^\alg (R,\Gamma)).\]
    By \Cref{prop:ht-n-kn-local} and \Cref{prop:is-ht-n}, this factors as
    \[E(R,\Gamma) \to L_{K(h)} E(R,\Gamma) \to E(L_{K(h)}^\alg (R,\Gamma)).\]
    To show that the second map is an equivalence, it suffices to prove that it induces an isomorphism on $\pi_*$ after cofibering out by $(p, v_1, \dots, v_{h-1})$ because both source and target are $K(h)$-local.
    It follows from \cite[Proposition 7.10]{HoveyStrickland} and \Cref{lem:vi-regular} that applying $\pi_0 (-/(p,v_1, \dots, v_{h-1}))$ to the above sequence of maps produces
    \[R_{(R,\Gamma)}^{\mathrm{cl}} / J_h \to R_{(R,\Gamma)} ^{\mathrm{cl}}/J_h[v_h^{-1}] \to R^{\mathrm{cl}}_{(R,\Gamma)} /J_h [v_h^{-1}],\]
    and we note that the second map is indeed an isomorphism.
\end{proof}

\subsubsection{$\LKN$-pairs and the $K(h)$-local tensor product}

\begin{defn}
    Let $k$ denote a Noetherian weakly \'etale $\F_p$-algebra $k$, and let $\LKN_{h,k}$ denote the category of height $h$ $\LKN$ pairs $(R,\Gamma)$ where $R$ is equipped with the structure of a $k$-algebra.

    We further define $\LKN_{h,k} ^p$ denote the category of pairs $(R,\Gamma)$ where $R$ is a perfect $k$-algebra and $\Gamma$ is base-changed from a height $h$ connected $p$-divisible group over $k$.
    In particular, any such pair is $\LKN$ of height $h$, so that $\LKN_{h,k}^p$ is a full subcategory of $\LKN_{h,k}$.
\end{defn}

\begin{construction} \label{const:notimes}
    We define a functor
    \[\notimes_k\colon \LKN_{h,k} \times \LKN_{h,k} ^p \to \LKN_{h,k}\]
    as follows:
    we define $(R_1, \Gamma_1) \notimes_k (R_2, \Gamma_2)$ to corepresent the following functor on $\LKN_{k,h}$: a morphism $(R_1, \Gamma_1) \notimes_k (R_2, \Gamma_2) \to (R_3, \Gamma_3)$ corresponds to morphisms of pairs $(R_1, \Gamma_1) \to (R_3,\Gamma_3)$ and $(R_2, \Gamma_2) \to (R_3, \Gamma_3 ^\circ)$.

    Such a corepresenting pair may be constructed as follows: first, consider $(R_1 \otimes _k R_2, \Gamma_1)$, then let $R'$ denote the pro-finite \'etale cover of $R_1 \otimes_k R_2$ classifying isomorphisms $\Gamma_1 ^\circ \cong \Gamma_2.$
    Then $(R_1, \Gamma_1) \notimes_k (R_2, \Gamma_2) = (R', \Gamma_1)$, once we see that the latter is an $\LKN$ pair of height $n$.
    But this is an immediate consequence of the fact that the map $R_1 \to R_1 \otimes_k R_2 \to R'$ is a the composite of two flat and relatively perfect maps, and hence flat and relatively perfect itself.
\end{construction}

\begin{thm} \label{thm:alg-tensor}
    Given $(R_i, \Gamma_i)$ and $k$ as in Construction \ref{const:notimes},
    there is a canonical isomorphism
    \[L_{K(h)} (E(R_1, \Gamma_1) \otimes_{\SW(k)} E(R_2, \Gamma_2)) \cong E((R_1, \Gamma_1) \notimes_k (R_2,\Gamma_2))\]
    of $\Einf$-rings.
\end{thm}

\begin{proof}
    To prove the desired equivalence of $\Einf$-rings, we prove that both sides corepresent the same functor in $\CAlg_{K(h), \SW(k),\ad}^{\mathrm{cmpl}}$, the category of $K(h)$-local complete adic $\Einf$-algebras over $\SW(k)$.
    Let $A$ denote a $K(h)$-local complete adic $\Einf$-ring.
    Then $L_{K(h)} (E(R_1, \Gamma_1) \otimes_{\SW(k)} E(R_2, \Gamma_2))$ is the coproduct of $E(R_1, \Gamma_1)$ and $E(R_2, \Gamma_2)$ in the category $\CAlg_{K(h), \SW(k),\ad}^{\mathrm{cmpl}}$, so that
    \[\Map_{\CAlg_{K(h),\SW(k),\ad}^{\mathrm{cmpl}}} (L_{K(h)} (E(R_1, \Gamma_1) \otimes_{\SW(k)} E(R_2, \Gamma_2)),A)\]
    identifies with the classifying anima of tuples $(\widetilde{\Gamma}_1, \widetilde{\Gamma}_2, \alpha_1, \alpha_2, e_1, e_2)$,
    where $\widetilde{\Gamma}_i$ are $p$-divisible groups on $A$, $e_i$ are orientations of $\widetilde{\Gamma}_i$, and $\alpha_i$ are equivalence classes of $k$-linear $\Gamma_i$-taggings of $\widetilde{\Gamma}_i$.
    Since $\Gamma_2$ is connected, the existence of $\alpha_2$ implies that $\widetilde{\Gamma}_2$ is formally connected, so that it is determined by its corresponding formal group, which in turn is identified with $\Gamma_1 ^{\circ}$ by $e_1$ and $e_2$.
    
    Putting this all together, we can rewrite this mapping anima as the classifying anima of tuples $(\widetilde{\Gamma}_1, \alpha_1, \alpha_2, e_1)$, where $\alpha_2$ is now an equivalence class of $k$-linear $\Gamma_2^{\circ}$-taggings of $\Gamma_1^{\circ}$.
    This, in turn, may be rewritten as a tuple $(\widetilde{\Gamma}, \alpha, e)$, where $\widetilde{\Gamma}$ is a $p$-divisible group over $A$, $\alpha$ is an equivalence class of $k$-linear $\Gamma'$-taggings of $\widetilde{\Gamma}$, and $e$ is an orientation of $\widetilde{\Gamma}$.
    By definition, this is the same as the mapping anima
    \[\Map_{\CAlg_{K(h),\SW(k),\ad}^{\mathrm{cmpl}}} (L_{K(h)} (E((R_1, \Gamma_1) \notimes_k (R_2,\Gamma_2)),A)\]
    as desired.
\end{proof}

\subsection{Solid Lubin--Tate theories}
In this section, we define the solid version of the Lubin--Tate theory functor, and describe its interaction with $K(h)$-localization and the $K(h)$-local solid tensor product.

\subsubsection{Solid Lubin--Tate theories}

We begin with the construction of the solid Lubin--Tate theory functor, taking special care to show that this is a functor of categories enriched over the category of condensed anima $\Ani^{\cond}$.
The first step is the construction of a crude condensed Lubin--Tate theory functor.

\begin{construction}
    So far, we have constructed a functor
    \[E(-)\colon \LKN \to \CAlg.\]
    To upgrade this to a condensed functor, we may proceed as follows.
    
    Let $\mathcal{P}(\LKN)$ denote the full subcategory of $\Fun(\LKN^{op}, \Ani)$ generated by representable functors under small colimits---this is a presentable category---and let $\mathcal{P}(\LKN)^{\cond} = \mathcal{P}(\LKN) \otimes \Ani^{\cond}$.
    
    Now, since $\CAlg$ is presentable, we may factor $E(-)$ as
    \[\LKN \to \mathcal{P}(\LKN) \to \CAlg,\]
    with the latter map being a left adjoint between presentable categories.
    We may then tensor this with $\Ani^{\mathrm{cond}}$ to obtain a map
    \[E^{\mathrm{cond}} (-)\colon \mathcal{P}(\LKN)^{\mathrm{cond}} \to \CAlg^{\mathrm{cond}}\]
    of $\Ani^{\mathrm{cond}}$-modules in $\mathrm{Pr}^L$.
    In particular, $E^{\mathrm{cond}} (-)$ is a functor enriched over condensed anima.
    Writing $E^{\mathrm{cond}}(-)$ as a functor
    \[E^{\mathrm{cond}}\colon \mathrm{Shv}^{\mathrm{hyp}, \mathrm{acc}}(\mathrm{Prof}, \mathcal{P}(\LKN)) \to \mathrm{Shv}^{\mathrm{hyp}, \mathrm{acc}}(\mathrm{Prof}, \CAlg),\]
    it is computed by applying
    $\mathcal{P}(\LKN) \to \CAlg$
    levelwise and hypersheafifying.
\end{construction}

While this construction produces a well-defined functor enriched over condensed anima, it is too general to obtain good control over.
We will restrict to a subcategory of $\mathcal{P}(\LKN)^{\mathrm{cond}}$ on which we may understand $E^{\mathrm{cond}} (-)$ better.

\begin{defn}
\label{def:solidLKNpair}
    A solid $\LKN$-pair is a pair $(R, \Gamma)$ of a solid $\F_p$-algebra $R$ and a $p$-divisible group $\Gamma$ over $R(\ast)$ that satisfies the following condition:
    there is a map of pairs $(R_0, \Gamma_0) \to (R(\ast), \Gamma)$ such that
    \begin{itemize}
        \item the $\F_p$-algebra $R_0$ is Noetherian and $F$-finite.
        \item for every nonempty profinite set $S$, the composite map $R_0 \to R(\ast) \to R(S)$ is flat and relatively perfect.
        \item the $p$-divisible group $\Gamma_0$ is nonstationary.
    \end{itemize}

    As a consequence, each of the pairs $(R(S), \Gamma)$ obtained by base change is an $\LKN$ pair.

    We say that a solid $\LKN$ pair is of height $h$ if there is a map $(R_0, \Gamma_0) \to (R(\ast), \Gamma)$ as above for which $\Gamma_0 ^{\circ}$ is of exact height $h$ and the kernel of $R^{\cl}_{(R_0, \Gamma_0)} \to R_0$ is equal to $J_h$.
    In particular, this implies that each of the $\LKN$ pairs $(R(S), \Gamma)$ is an $\LKN$ pair of height $h$.
\end{defn}

\begin{example}
    Given any $\LKN$ pair $(R,\Gamma)$, the pair $(R^\delta, \Gamma)$ is a solid $\LKN$ pair. This is because any surjection of profinite sets $T \to S$ yields an ind-\'etale map $R^{\delta} (S) \to R^{\delta} (T)$. If $(R,\Gamma)$ is $\LKN$ of height $h$, then $(R^\delta, \Gamma)$ is solid $\LKN$ of height $h$.
\end{example}

\begin{construction}
Let $\LKN^\solid$ denote the category of solid $\LKN$-pairs, where a morphism $(R_1,\Gamma_1) \to (R_2,\Gamma_2)$ consists of a map of solid rings $f\colon R_1 \to R_2$ and isomorphisms $f^* \Gamma_1 \cong \Gamma_2$ over $R_2(\ast)$, and let $\LKN^{\solid} _h$ denote the full subcategory of height $h$ solid $\LKN$ pairs.
Given a solid $\LKN$-pair $(R,\Gamma)$, the assignment $S \mapsto (R(S),\Gamma)$ determines a hypersheaf $\mathrm{Prof}^{\mathrm{op}} \to \LKN$ that may be composed with the Yoneda embedding to yield a fully faithful functor
\[\LKN^\solid \to \mathcal{P}(\LKN)^{\mathrm{cond}}.\]

Viewing $\LKN^\solid$ as a full subcategory of the $\Ani^{\mathrm{cond}}$-enriched category $\mathcal{P}(\LKN)^{\mathrm{cond}}$, we may endow $\LKN^{\solid}$ with the structure of an $\Ani^{\mathrm{cond}}$-enriched category.
\end{construction}

\begin{rmk}
By presentability, $\mathcal{P}(\LKN)^{\mathrm{cond}}$ is cotensored over $\Ani^{\mathrm{cond}}$ and in particular over $\mathrm{Prof}$.
Viewing $\mathcal{P}(\LKN)^{\mathrm{cond}}$ as $\mathrm{Shv}^{\mathrm{hyp}, \mathrm{acc}}(\mathrm{Prof}, \mathcal{P}(\LKN))$, the $\mathrm{Prof}$-cotensoring can be described as follows: given $S \in \mathrm{Prof}$, the cotensoring of $A$ by $S$ is the following hypersheaf: $C(S,A)(T) = A(S\times T)$.
We note that $\LKN^{\solid}$ is closed under the $\mathrm{Prof}$-cotensoring in $\mathcal{P}(\LKN)^{\mathrm{cond}}$: given a solid $\LKN$-pair $(R,\Gamma)$, we have $C(S,(R,\Gamma)) = (C(S,R), \Gamma)$, where $C(S,R)(T) = R(S \times T)$ and $\Gamma$ abusively denotes the base change of $\Gamma$ to $R(S)$. This is indeed a solid $\LKN$-pair, with witness $(R_0, \Gamma_0) \to (R(\ast), \Gamma) \to (R(S),\Gamma)$.
\end{rmk}

\begin{proposition} \label{prop:solid-LT}
    The $\Ani^{\cond}$-enriched functor
    \[E^{\solid} (-) \colon \LKN^{\solid} \to \mathcal{P}(\LKN)^{\mathrm{cond}} \xrightarrow{E^{\mathrm{cond}}} \CAlg^{\solid}\]
    is computed explicitly by the assignment
    \[(S \mapsto (R(S),\Gamma)) \mapsto (S\mapsto E(R(S),\Gamma))\]
    and its solid homotopy groups are computed by
    \[(S \mapsto (R(S),\Gamma)) \mapsto (S\mapsto \pi_n E(R(S),\Gamma)).\]
\end{proposition}

\begin{proof}
    It follows from the definitions that the functors $E^{\solid} (-)$ and $\pi_n E^{\solid} (-)$ are computed via the hypersheafification of the formulas given in the statement of the lemma.
    It therefore suffices to prove that these presheaves are already sheaves for the profinite topology.
    To start, we choose a Noetherian Luriean pair $(R_0, \Gamma_0)$ with a map $(R_0, \Gamma_0) \to (R(\ast), \Gamma)$ such that the composition $R_0 \to R(S)$ is flat and relatively perfect for all profinite sets $S$.

    Now, let $T_\bullet \to S$ denote a hypercover of profinite sets.
    Since $\Z[S]^\solid$ is projective as a solid abelian group by \cite[Lecture V]{ClausenScholzeCondensed}, there is an equivalence
    $R(S) \xrightarrow{\sim} \Tot R(T_\bullet)$
    where the totalization is computed in derived commutative rings.
    Applying \Cref{prop:tot-preserved} to the pair $(R_0, \Gamma_0)$ and the augmented cosimplicial diagram $R(S) \to R(T_\bullet),$
    we conclude that
    \[E(R(S), \Gamma) \to \Tot E(R(T_\bullet), \Gamma))\]
    is an equivalence, as desired. The corresponding statement for the condensed homotopy groups follows from the second part of \Cref{prop:tot-preserved}.

    It remains to prove that $E^{\solid} (R,\Gamma)$ is solid.
    By definition, we must prove that its condensed homotopy groups $\pi_{n} E^{\solid}(R,\Gamma)$ are solid abelian groups.
    When $n$ is odd, these are zero.
    When $n=2m$, we saw above that they are given by 
    \[S \mapsto \pi_{2m} E(R(S),\Gamma).\]
    Choosing $(R_0, \Gamma_0) \to (R(\ast), \Gamma)$ as above and letting $I  = \ker (\pi_0 E(R_0, \Gamma_0) \to R_0)$, we see that these groups are all $I$-complete.
    Thus there is an equivalence of condensed abelian groups $\pi_{2m} E(R, \Gamma) \simeq \varprojlim_k \pi_{2m} E^{\solid}(R, \Gamma) \otimes_{\pi_0 E(R_0, \Gamma_0)^{\delta}} (\pi_0 E(R_0, \Gamma_0)/I^k)^\delta$.
    Here, this tensor product is computed pointwise: this produces a condensed abelian group because $\pi_{2m} E(R,\Gamma)$ is pointwise flat over $\pi_0 E(R_0, \Gamma_0)$.
    It therefore suffices to prove that $\pi_{2m} E(R, \Gamma) \otimes_{\pi_0 E(R_0, \Gamma_0)^{\delta}} R_0 ^{\delta}$ and $\pi_{2m} E(R, \Gamma) \otimes_{\pi_0 E(R_0, \Gamma_0)^{\delta}} (I^k/I^{k+1})^{\delta}$ are solid abelian groups.
    In this first case, this follows from the isomorphism
    \[\pi_{2m} E(R, \Gamma) \otimes_{\pi_0 E(R_0, \Gamma_0)^{\delta}} R_0 ^{\delta} \cong R \otimes_{R_0^{\delta}} (\omega_{\Gamma_0}^m)^{\delta},\]
    the fact that $\omega_{\Gamma_0}^m$ is an invertible module over $R_0$, and the fact that $R$ is solid.
    In the second case, by choosing a presentation of $I^k/I^{k+1}$ as an $R_0$-module, we may write $\pi_{2m} E(R, \Gamma) \otimes_{\pi_0 E(R_0, \Gamma_0)^{\delta}} (I^k/I^{k+1})^{\delta}$ as a colimit of direct sums of objects of the form $\pi_{2m} E(R, \Gamma) \otimes_{\pi_0 E(R_0, \Gamma_0)^{\delta}} R_0 ^{\delta}$ and hence conclude that it too is solid.
\end{proof}

\begin{rmk} \label{rmk:prof-cotensor-preserved}
    Since $E^{\solid} (-)$ is an $\Ani^{\mathrm{cond}}$-enriched functor, there is a natural transformation $E^{\solid}(C(S,(R,\Gamma))) \to C(S,E^{\solid}(R,\Gamma))$ for $(R,\Gamma) \in \LKN^{\solid}$ and $S \in \mathrm{Prof}$.
    Using \Cref{prop:solid-LT}, one verifies that this is a natural equivalence.
\end{rmk}

\subsubsection{Localization of solid Lubin--Tate theories}

We begin by recalling the interaction of accessible localizations with condensed and solid spectra from \cite[Section 13]{Clausen2025arxiv}. (He works in the light setting, but the proofs are the same.)

\begin{defn}
    Given an accessible localization functor $L\colon \mathrm{Sp} \to \mathrm{Sp}$, we obtain an accessible localization
    $L^{\mathrm{cond}}\colon \mathrm{Sp}^{\mathrm{cond}} \to \mathrm{Sp}^{\mathrm{cond}}$
    given by killing the cocomplete full subcategory of $\mathrm{Sp}^{\mathrm{cond}}$ generated by objects of the form $X \otimes A$, where $X \in \mathrm{Sp}^{\mathrm{cond}}$ and $A \in \mathrm{Sp}$ with $L(A) \simeq 0$.
    As noted in \cite[p. 114]{Clausen2025arxiv}, an object $X \in \mathrm{Sp}^{\mathrm{cond}}$ is $L^{\mathrm{cond}}$-local if and only if $X(S)$ is $L$-local for all profinite sets $S$, and the functor $L^{\mathrm{cond}}$ is computed by applying $L$ levelwise and hypersheafifying.
    By \cite[Lemma 13.5]{Clausen2025arxiv}, this restricts to a functor
    $L^{\solid}\colon \mathrm{Sp}^{\solid} \to \mathrm{Sp}^{\solid}.$
\end{defn}

\begin{lemma} \label{lem:ht-h-discrete}
  Given an $\LKN$ pair $(R,\Gamma)$ of height $h$, there is an equivalence
  \[E^{\solid} (R^\delta,\Gamma) \simeq L_{K(h)}^{\solid} E(R,\Gamma)^{\delta}.\]
\end{lemma}

\begin{proof}
  By the universal property of $(-)^{\delta}$, we have a map of solid $\Einf$-rings
  $E(R,\Gamma)^{\delta} \to E^{\solid} (R^{\delta}, \Gamma)$
  that factors through a map
  $L^{\solid} _{K(h)} E(R,\Gamma)^{\delta} \to E^{\solid} (R^{\delta}, \Gamma).$
  
  To show that this map is an equivalence, it suffices to do so after cofibering out by $(p, v_1, \dots, v_{h-1})$.
  Upon doing so, we find that
  \[\pi_{2k} (L^{\solid} _{K(h)} E(R,\Gamma)^{\delta}/(p,v_1, \dots, v_{h-1})) \cong (\pi_{2k} (E(R,\Gamma)/(p,v_1, \dots, v_{h-1})))^{\delta} \cong (\omega^{k} _{\Gamma})^{\delta}\]
  and that
  \[\pi_{2k} (E^{\solid} (R^{\delta}, \Gamma)/(p,v_1, \dots, v_{h-1})) \cong (\omega^{k} _{\Gamma})^{\delta},\]
  by computing levelwise and using the fact that $(R^{\delta} (S),\Gamma)$ is $\LKN$ of height $h$ for all profinite sets $S$.
  From this, we deduce that the map
  $L^{\solid} _{K(h)} E(R,\Gamma)^{\delta} \to E^{\solid} (R^{\delta}, \Gamma)$
  is an equivalence of solid $\Einf$-rings, as desired.
\end{proof}

\begin{proposition} \label{prop:L-Kh-solid}
    Let $(R,\Gamma)$ denote a solid $\LKN$ pair.
    Then the assignment
    \[S \mapsto L_{K(h)}^{\alg} (R(S),\Gamma)\]
    defines a solid $\LKN$ pair of height $h$, which we denote by $L_{K(h)}^{\solid,\alg} (R,\Gamma)$.
    The functor
    $L_{K(h)}^{\solid, \alg}\colon \LKN^{\solid} \to \LKN^{\solid}_{h}$
    is an $\Ani^{\mathrm{cond}}$-enriched functor that preserves the $\mathrm{Prof}$-cotensor.
\end{proposition}

\begin{proof}
    First, we check that the assignment
    \[S \mapsto L_{K(h)}^{\alg} (R(S),\Gamma)\]
    is a hypersheaf.
    This is the same as in the proof of \Cref{prop:solid-LT}, but with \Cref{lem:alg-kn-tot} replacing \Cref{prop:tot-preserved}.

    To check that the resulting condensed ring is solid, we apply \Cref{prop:solid-LT}, in which we showed that $(S \mapsto \pi_0 E(R(S), \Gamma))$ is solid.
    By definition, the underlying ring of $L_{K(h)}^{\alg} (R(S),\Gamma)$ is obtained from $\pi_0 E(R(S), \Gamma)$ by quotienting by $J_h$ and inverting $v_h$. Since the subcategory of solid abelian groups is closed under colimits, this shows that the underlying condensed ring of $L_{K(h)}^{\alg} (R,\Gamma)$ is solid.
    
    Next, to show that $L_{K(h)}^{\alg} (R,\Gamma)$ is a solid $\LKN$ pair, we choose a witness $(R_0, \Gamma_0) \to (R(\ast),\Gamma)$ to $(R,\Gamma)$ being solid $\LKN$. Then it follows from the proof of \Cref{prop:is-ht-n} that the maps
    \[L_{K(h)}^{\alg} (R_0,\Gamma_0) \to L_{K(h)}^{\alg} (R(\ast),\Gamma) \to L_{K(h)}^{\alg} (R(S),\Gamma)\]
    are flat and relatively perfect on the level of $\F_p$-algebras, and $L_{K(h)}^{\alg} (R_0,\Gamma_0)$ is Noetherian.
    It follows that $L_{K(h)}^{\alg} (R,\Gamma)$ is a solid $\LKN$ pair, as desired.
    Similarly, it follows as in the proof of \Cref{prop:is-ht-n} that $L_{K(h)}^{\alg} (R,\Gamma)$ is a solid $\LKN$ pair of height $h$.

    Finally, note that it follows from the above that there is a commutative diagram
    \begin{center}
        \begin{tikzcd}
            \LKN^{\solid} \ar[r,"L_{K(h)} ^{\solid,\alg}"] \ar[d,hook] & \LKN_h ^{\solid} \ar[d,hook] \\
            \mathcal{P}(\LKN)^{\cond} \ar[r,"L_{K(h)}^{\alg,\cond}"] & \mathcal{P}(\LKN_h)^{\cond},
        \end{tikzcd}
    \end{center}
    where the bottom functor is obtained from $L_{K(h)}^{\alg}$ by applying the construction $\P(-)^{\cond}$. (Commutativity follows from that fact that the presheaf $S\mapsto L_{K(h)}^{\alg} (R(S),\Gamma)$ is already a hypersheaf.)
    The bottom functor is $\Ani^\cond$-enriched by definition, so the diagram endows $L_{K(h)} ^{\solid,\alg}$ with the strucutre of an $\Ani^{\cond}$-enriched functor, since the vertical functors are fully faithful as $\Ani^{\cond}$-enriched functors by the definition of the enriched structures on $\LKN^{\solid}$ and $\LKN_h^{\solid}$.
\end{proof}

\begin{proposition} \label{prop:Kh-alg-top-solid}
    There is a canonical isomorphism of functors enriched over condensed anima
    \[L_{K(h)} ^{\solid} E^{\solid}(R,\Gamma) \simeq E^\solid (L_{K(h)}^{\solid,\alg} (R,\Gamma)).\]
\end{proposition}

\begin{proof}
    \Cref{thm:Kh-alg-top} determines a commutative diagram of categories
    \begin{center}
        \begin{tikzcd}
            \LKN \ar[r, "L_{K(h)}^{\alg}"] \ar[d, "E(-)"] & \LKN_h  \ar[d, "E(-)"] \\
            \CAlg \ar[r, "L_{K(h)}"] & \CAlg_{K(h)}.
        \end{tikzcd}
    \end{center}
    By left Kan extension and tensoring with $\Ani^{\cond}$, this determines a commutative diagram of $\Ani^{\cond}$-modules in $\mathrm{Pr}^L$
    \begin{center}
        \begin{tikzcd}
            \mathcal{P}(\LKN)^{\cond} \ar[r, "L_{K(h)}^{\alg, \cond}"] \ar[d, "E^\cond(-)"] & \mathcal{P}(\LKN_h)^{\cond} \ar[d, "E^\cond(-)"] \\
            \CAlg^{\cond} \ar[r, "L_{K(h)}^\cond"] & \CAlg_{K(h)} ^{\cond}.
        \end{tikzcd}
    \end{center}
    On the other hand, it follows from the definition of $L_{K(h)}^\solid$ that there is a commutative diagram of categories enriched over $\Ani^{\cond}$
    \begin{center}
        \begin{tikzcd}
            \LKN^\solid \ar[r, "L_{K(h)}^{\solid, \alg}"] \ar[d, hook] & \LKN_h ^\solid \ar[d, hook] \\
            \mathcal{P}(\LKN)^{\cond} \ar[r, "L_{K(h)}^{\alg, \cond}"] & \mathcal{P}(\LKN_h)^{\cond}.
        \end{tikzcd}
    \end{center}
    Pasting the previous two diagram, we obtain the desired equivalence of functors enriched over condensed anima.
\end{proof}

\subsubsection{Solid $\LKN$-pairs and the $K(h)$-local tensor product}

\begin{defn} \label{defn:notimes-solid}
    Given a Noetherian weakly \'etale $\F_p$-algebra $k$, we define $\LKN^{\solid} _{k^\delta,h}$ to be the category of solid $\LKN$ pairs $(R,\Gamma)$ of height $h$ where $R$ is a $k^{\delta}$-algebra.
    We define a functor
    \[\notimes_{k^{\delta}} ^{\solid} \colon \LKN_{k^\delta,h} ^{\solid} \times \LKN^{p}_{k,h} \to \LKN_{k^\delta,h} ^{\solid}\]
    as follows.
    Given a solid $\LKN$ pair $(R_1, \Gamma_1)$ of height $h$ over $k^\delta$ and a (discrete) perfect $\LKN$ pair $(R_2, \Gamma_2)$ of height $h$ over $k$ with $\Gamma_2$ connected and base changed from $k$, we define
    \[(R_1, \Gamma_1) \notimes_{k^\delta} ^{\solid} (R_2, \Gamma_2) \coloneqq (R',\Gamma')\]
    via the formula
    \[(R'(S), \Gamma') = (R_1 (S), \Gamma_1) \notimes_k (R_2, \Gamma_2).\]
    We will in the following proposition that this does indeed define a solid $\LKN$ pair of height $h$.
\end{defn}

\begin{proposition}
    Given $(R_i,\Gamma_i)$, $k$, and $k \to R_i$ as in \Cref{defn:notimes-solid}, the object
    \[(R',\Gamma') = (R_1, \Gamma_1) \notimes_{k^\delta} ^{\solid} (R_2, \Gamma_2)\]
    is a solid $\LKN$ pair of height $h$, and, if $\LKN^{p}_{k,h}$ is regarded as a discrete $\Ani^{\cond}$-enriched category, the functor
    \[\notimes_{k^{\delta}} ^{\solid} \colon \LKN_{k^\delta,h} ^{\solid} \times \LKN^{p}_{k,h} \to \LKN_{k^\delta,h} ^{\solid}\]
    is enriched over $\Ani^{\cond}$.
\end{proposition}

\begin{proof}
    We first show that $(R',\Gamma')$ is a solid $\LKN$ pair of height $h$.
  First, we show that $R'$ is a solid $\F_p$-algebra.
  Examining the definitions, we see that the following square in $\Fun(\mathrm{Prof}^{op}, \CAlg^{\mathrm{an}})$
  \begin{center}
      \begin{tikzcd}
          \underline{R_1(\ast)} \ar[r] \ar[d] & \ar[d] \underline{R'(\ast)} \\
          R_1 \ar[r] & R'
      \end{tikzcd}
  \end{center}
  is a pushout,
  where, given a commutative ring $R$, $\underline{R}$ is the constant presheaf at $R$.
  Hypersheafifying, we obtain an equivalence of condensed animated $\F_p$-algebras
  \[R_1 \otimes_{(R_1(\ast))^{\delta}} (R'(\ast))^{\delta} \xrightarrow{\sim} R'.\]
  Writing down a presentation of $R'(\ast)$ as a $R_1(\ast)$-module and using the fact that solid abelian groups are closed under colimits, we find that
  \[R' \cong R_1 \otimes_{(R_1(\ast))^{\delta}} (R'(\ast))^{\delta}\]
  is a solid $\F_p$-algebra, as desired.
  
  Now, let $(R_{10}, \Gamma_{10}) \to (R_1(\ast), \Gamma_1)$ denote a witness to $(R_1, \Gamma_1)$ being solid $\LKN$ of height $h$.
  Then it follows from \Cref{const:notimes} that the composite
  \[(R_{10}, \Gamma_{10}) \to (R_1(\ast), \Gamma_1) \to (R'(\ast), \Gamma'))\]
  is a witness to $(R_1,\Gamma_1)$ being solid $\LKN$ of height $h$.

  To construct the $\Ani^\cond$-enrichment of $\notimes^{\solid} _{k^{\delta}}$, we factor it as
  \[\LKN_{k^\delta,h} ^{\solid} \times \LKN^{p}_{k,h} \hookrightarrow \mathcal{P} (\LKN_{k^\delta,h})^{\cond} \times \mathcal{P}(\LKN^{p}_{k,h})^{\cond} \to \mathcal{P} (\LKN_{k^\delta,h})^{\cond},\]
  where the first functor is the product of
  \[\LKN_{k^\delta,h} ^{\solid} \hookrightarrow \mathcal{P} (\LKN_{k^\delta,h})^{\cond}\]
  and
  \[\LKN^{p}_{k,h} \hookrightarrow \mathcal{P}(\LKN^{p}_{k,h}) \xhookrightarrow{(-)^\delta} \mathcal{P}(\LKN^{p}_{k,h})^{\cond}\]
  and the second map is obtained from $\notimes_k$ by applying $\mathcal{P}$ and tensoring with $\Ani^{\cond}$.
\end{proof}

\begin{proposition} \label{prop:solid-tensor-thm}
    Given $(R_i,\Gamma_i)$, $k$, and $k \to R_i$ as in \Cref{defn:notimes-solid}, there is a canonical equivalence
    \[L_{K(h)}^{\solid} (E^{\solid} (R_1,\Gamma_1) \otimes_{\SW(k)^{\delta}} ^{\solid} E^{\solid}(R_2^{\delta}, \Gamma_2)) \simeq  E^{\solid}((R_1, \Gamma_1) \notimes_{k^\delta} ^{\solid} (R_2, \Gamma_2))\]
    that arises from an equivalence of $\Ani^\cond$-enriched functors
    \[\LKN_{k^\delta,h} ^{\solid} \times \LKN^{p}_{k,h} \to \CAlg_{K(h)} ^{\solid}.\]
\end{proposition}

\begin{proof}
    It follows from \Cref{lem:ht-h-discrete} that there is a canonical equivalence
  \[L_{K(h)}^{\solid} (E^{\solid} (R_1,\Gamma_1) \otimes_{\SW(k)^{\delta}} ^{\solid} E^{\solid}(R_2^{\delta}, \Gamma_2)) \simeq L_{K(h)}^{\solid} (E^{\solid} (R_1,\Gamma_1) \otimes_{\SW(k)^{\delta}} ^{\solid} E(R_2, \Gamma_2)^{\delta}).\]
  Concretely, the right hand side may be computed as the profinite hypersheafification of the presheaf
  \[S \mapsto L_{K(h)} (E(R_1(S),\Gamma_1) \otimes_{\SW(k)} E(R_2, \Gamma_2)).\]
  By \Cref{thm:alg-tensor}, this identifies with the presheaf
  \[S \mapsto E((R_1 (S),\Gamma_1) \notimes_k (R_2,\Gamma_2)),\]
  which is $E^{\solid} ((R_1, \Gamma_1) \notimes^{\solid} _{k^{\delta}} (R_2,\Gamma_2))$ by \Cref{prop:solid-LT}.

  Upgrading this to an equivalence of $\Ani^{\cond}$-enriched functors may be done similarly to the proof of \Cref{prop:Kh-alg-top-solid}.
\end{proof}

\subsection{Examples and applications} \label{subsec:exapp}
We now discuss some examples and applications of the above theory that will be used in the sequel.

\subsubsection{Examples}
\label{ssec:examples}
We now introduce some important examples of solid $\LKN$ pairs.  In these examples, $k$ refers to an algebraic extension of $\F_p$ which contains $\F_{p^t}$ for all $t$ which appear as a height in the context of the example discussed, and $\Gamma_h$ refers to a 1-dimensional connected $p$-divisible group over $k$ all of whose geometric endomorphisms are $k$-rational: for example, the $p$-divisible group associated to the Honda formal group of height $h$ when $\mathbb{F}_{p^h} \subseteq k$.

\begin{example}
\label{exm:usualLT}
    The most fundamental example starts with the $\LKN$ pair $(k, \Gamma_h)$.
    Then $E_h = E(k, \Gamma_h)$ is the usual Lubin--Tate theory of height $h$, and $E_h ^{\solid} = E(k^\delta,\Gamma_h)$ is its solid refinement as constructed by Clausen \cite[Section 12]{Clausen2025arxiv}.
    Explicitly, $E(k^\delta,\Gamma_h)$ is the $L_{K(h)}^{\solid}$-localization of $E_h ^{\delta}$ by \Cref{lem:ht-h-discrete}.

    The Morava stabilizer group, when viewed as a constant pro-finite \'etale group, acts on the pair $(k,\Gamma_h)$. This implies that there are compatible actions of $\G_h (S)$ on $(k^\delta (S), \Gamma_h)$ for each profinite $S$, so that there is a condensed action of the profinite group $\G_h$ on the solid $\LKN$ pair $(k^{\delta}, \Gamma_h)$, where we make sense of this using the $\Ani^{\cond}$-enrichment of $\LKN^{\solid}$.

    Thus, by $\Ani^{\cond}$-enriched functoriality of $E^{\solid}$, the Morava stabilizer group $\G_h$ acts on the solid $\Einf$-ring $E_h ^\solid$.
    Since $E_h^\solid$ is discrete as a solid object in $K(h)$-local $\Einf$-rings, this arises from a unique smooth action of $\G_h$ on $E_h^{\solid}$ by \Cref{prop:smoothdiscrete}.
\end{example}

\begin{example} \label{exm:kh-1}
    Of particular interest to the chromatic splitting conjecture in coheight $1$ is the solid spectrum $L_{K(h-1)} ^{\solid} E_h ^{\solid}$, which is endowed with an induced action of the profinite group $\G_h$.
    By \Cref{prop:Kh-alg-top-solid}, we have a $\G_h$-equivariant identification
    \[L_{K(h-1)}^{\solid} E_h ^\solid \simeq E(K, \Gamma_{h,K}^{\mathrm{un}})\]
    where $(K,\Gamma_{h,K}^{\mathrm{un}}) = L_{K(h-1)}^{\solid, \alg} (k,\Gamma_h)$.
    Examining the definition of $L_{K(h-1)}^{\solid, \alg}$, we see that $K \cong k((u_{h-1}))$, with the solid structure arising from the $u_{h-1}$-adic topology, and that the action of $\G_h$ on $(K,\Gamma_{h,K}^{\mathrm{un}})$ is obtained from its action on $(R^{\cl}_{(k, \Gamma_h)},\Gamma_h^{\un})$ via quotienting out by $I_{h-1}$ and inverting $u_{h-1}$.

    We note that the action of $\G_h$ is on the $K(h-1)$-local solid $\Einf$-ring $L_{K(h-1)}^{\solid} E_h ^\solid$ is no longer smooth: indeed, the action of $\G_h$ on the solid ring $\pi_0 (L_{K(h-1)}^{\solid} E_h ^\solid) / I_{h-1} \cong K \cong k((u_{h-1}))$ is not smooth.
\end{example}

\begin{example} \label{exm:B}
    The $\Einf$-ring $\mathbb{B} = L_{K(h-1)} (L_{K(h)}E_h \otimes_{\SW(k)} E_{h-1})$, which naturally enriches to the solid $\Einf$-ring $\mathbb{B}^{\solid} = L_{K(h-1)} ^{\solid} (L_{K(h-1)} ^{\solid} E_h ^{\solid} \otimes_{\SW(k)} ^{\solid} E_{h-1} ^{\solid})$, may be viewed as a solid Lubin--Tate theory in an equivariant way.

    First, we note that
    $L_{K(h-1)} ^{\solid} (L_{K(h-1)} ^{\solid} E_h ^{\solid} \otimes_{\SW(k)} ^{\solid} E_{h-1} ^{\solid})$
    admits an action of the profinite group $\G_{h,h-1} = \G_h \times_{\Gal(k/\F_p)} \G_{h-1}$
    Indeed, this is a consequence of the fact that the cospan of solid $\Einf$-rings
    $L_{K(h-1)} ^{\solid} E_h ^{\solid} \leftarrow \SW(k) \to E_{h-1} ^{\solid}$
    is acted on by the span of profinite groups
    $\G_h \to \Gal(k/\F_p) \leftarrow \G_{h-1}$.

    Let $(L,\Gamma_{h,L}^{\un}) = (K,\Gamma_{h,K}^{\un}) \notimes_{k^\delta}^\solid (k, \Gamma_{h-1})$.
    Then $(L,\Gamma_{h,L}^{\un})$ is equipped with an action of $\G_{h,h-1}$.
    By construction, we see that $K \otimes_{K(\ast)^\delta} L(\ast)^\delta \to L$ is an isomorphism and that $K(\ast) \to L(\ast)$ is the pro-finite \'etale cover classifying isomorphisms $(\Gamma_{h,K}^{\un})^{\circ} \cong \Gamma_{h-1}$, i.e., it is the Gross--Torii extension discussed in the introduction.
    By \cite[Theorem 3.5(1)(a)]{GrossRamification}, $K(\ast) \to L(\ast)$ is a totally ramified $\overline{\G}_{h-1}$-Galois extension of fields.
    As a consequence, $K \to L$ is a $\overline{\G}_{h-1}^{\mathrm{sm}}$-torsor in the sense of \cite[Definition 6.2.1]{anDR}; in particular, the action of $\overline{\G}_{h-1}$ on $L$ is smooth.
\end{example}

\begin{example}
\label{ExBtoBhat}
    We continue to view $L$ as a solid ring; let $\widehat{L} = (\mathcal{O}_L)^{\wedge} _{u_{h-1}} [u_{h-1}^{-1}]$ be its uniform completion, also considered as a solid ring.
    Here $\mathcal{O}_L$ is the solid ring of integers of $L$, which may be recovered functorially from $L$ as its solid ring of power bounded elements. \cite[Definition 2.6.1(3)]{anDR}
    The solid ring $\widehat{L}$ arises from $\widehat{L}$ viewed as a topological ring with the $u_{h-1}$-adic topology; this is a consequence of the fact that $\mathcal{O}_L / u_{h-1}$ is discrete as a solid ring.
    
    Since $\widehat{L}(\ast)$ is a perfect field by \Cref{PropL}, and both the Frobenius and its inverse are continuous for the $u_{h-1}$-adic topology, we find that $\widehat{L}(T)$ is perfect for all profinite sets $T$.
    It follows that $\widehat{L}(\ast) \to \widehat{L}(T)$ is flat and relatively perfect for all $T$, so that the pair $(\widehat{L},\Gamma_{h-1,\widehat{L}})$ is a solid LKN pair of height $h-1$. We may therefore define the solid $\Einf$-ring 
    \begin{equation}
        \label{eq:Bhatsolid}
\widehat{\mathbb{B}}^\solid = E^\solid(\widehat{L}, \Gamma_{h-1,\widehat{L}}).
    \end{equation} 

    We now construct a $\G_{h,h-1}$-equivariant map of solid $\Einf$-rings
    \[\mathbb{B}^{\solid} \to \widehat{\mathbb{B}}^\solid\]
    which realizes the natural map $L \to \widehat{L}$ on the level of $\pi_0$ modulo $I_{h-1}$.
    We do so by composing two maps
    \[E^\solid (L, \Gamma^{\un} _{h,L}) \to E^\solid (\widehat{L}, \Gamma^{\un} _{h,\widehat{L}})  \to E^\solid (\widehat{L}, (\Gamma^{\un} _{h,\widehat{L}})^\circ) \simeq E^\solid (\widehat{L}, \Gamma_{h-1,\widehat{L}}).\]

    The first map arises from the $\G_{h,h-1}$-equivariant map of solid $\LKN$ pairs $(L, \Gamma^{\un} _{h,L}) \to E^\solid (\widehat{L}, \Gamma^{\un} _{h,\widehat{L}})$ obtained from the $\G_{h,h-1}$-equivariant map of solid rings $L \to \widehat{L}$ and base changing the $p$-divisible group.

    To construct the second map, we use the fact there is a canonical splitting $\Gamma^{\un} _{h,\widehat{L}} \cong (\Gamma^{\un} _{h,\widehat{L}}) ^{\circ} \oplus (\Gamma^{\un} _{h,\widehat{L}}) ^\et$ over the perfect field $\widehat{L}$ (and similarly over $\widehat{L}(T)$).
    It follows that, given a deformation of $(\Gamma^{\un} _{h,\widehat{L}}) ^{\circ}$, one may take its direct sum with (the unique deformation of) $(\Gamma^{\un} _{h,\widehat{L}}) ^\et$ to obtain a deformation of $\Gamma^{\un} _{h,\widehat{L}}$.
    By the universal property of Lubin--Tate theory, this determine natural $\G_{h,h-1}(T)$-equivariant maps of $\Einf$-rings
    $E(\widehat{L}(T),\Gamma^{\un}_{h,\widehat{L}(T)}) \to E(\widehat{L}(T),(\Gamma^{\un}_{h,\widehat{L}(T)})^{\circ})$
    that fit together into a $\G_{h,h-1}$-equivariant map of solid $\Einf$-rings
    $E^{\solid} (\widehat{L},\Gamma^{\un}_{h,\widehat{L}}) \to E^\solid (\widehat{L}, (\Gamma^{\un} _{h,\widehat{L}})^\circ)$.
\end{example}

\subsubsection{Applications}
In this section, we study in more depth the above examples and the fixed points spectra for their profinite group actions.

\begin{rmk}
    In the following, we use the notion of descendable morphism of algebras introduced and studied by Mathew \cite{Mathew2016,AkhilNilpotence}: A morphism of commutative algebra objects $A \to B$ in a stable presentably symmetric monoidal category $\CC$ is said to be \emph{descendable} if the thick tensor ideal generated by $B$ contains $A$. Importantly, descendability is preserved under exact lax symmetric monoidal functors (\cite[Corollary 4.13]{MNN2017}).
\end{rmk}

We begin by recalling the following result, which is a consequence of the proof of the smash product theorem.

\begin{thm}[{{Hopkins--Ravenel \cite{RavenelNilpotence},\cite[Theorem 4.18]{AkhilNilpotence}}}] \label{thm:HRdesc}
    The map of $\Einf$-rings $L_{h} \SS \to E_h$ is descendable.
\end{thm}

\begin{cor} \label{cor:solidHRdesc}
  The map of solid $\Einf$-rings $L_{K(h)} ^{\solid} \SS \to E_h ^{\solid}$ is descendable.
\end{cor}

\begin{proof}
  This map is obtained from $L_{h} \SS \to E_h$ by applying $L_{K(h)}^{\solid} ((-)^{\delta}) \colon \mathrm{Sp} \to \mathrm{Sp} ^{\solid}$. Since this is an exact and lax symmetric monoidal functor, it preserves descendability.
\end{proof}

\begin{proposition} \label{prop:solid-Galois}
  The map $L^{\solid} _{K(h)} \SS \to E_h ^{\solid}$ is a faithful $K(h)$-local solid $\G_h$-Galois extension.
  In other words, the natural maps
  \[L^{\solid} _{K(h)} \SS \to (E_h ^{\solid})^{h \G_h}\]
  and
  \[L_{K(h)}^{\solid} (E_h ^{\solid} \otimes_{L^{\solid} _{K(h)} \SS} ^\solid E_h ^{\solid}) \to C (\G_h, E_h ^{\solid})\]
  are equivalences of solid $\Einf$-rings.
\end{proposition}

\begin{proof}
  We first prove that the second map is an equivalence.
  By \Cref{prop:solid-tensor-thm} and \Cref{rmk:prof-cotensor-preserved}, the second map is obtained by applying $E^{\solid} (-)$ to the canonical morphism of solid $\LKN$ pairs
  \[(k, \Gamma_h) \notimes_{\F_p} (k, \Gamma_h) \to C(\G_h,(k, \Gamma_h)).\]
  This map is an isomorphism by the definition of $\notimes_{\F_p}$: $C(\G_h,k)$ classifies isomorphisms from $\iota_1 ^* \Gamma_h$ to $\iota_2 ^* \Gamma_h$ over $k \otimes_{\F_p} k \cong C(\Gal(k/\F_p),k)$, where $\iota_1 , \iota_2 \colon k \hookrightarrow k \otimes_{\F_p} k$ are the two inclusions.

  Given this, the first map identifies with the natural map from $L^{\solid} _{K(h)} \SS$ to the \v{C}ech nerve of $L^{\solid} _{K(h)} \SS \to E_h ^{\solid}$, which is an isomorphism by descendability.
\end{proof}

\begin{cor} \label{cor:kh-1-still-descent}
    The map of $K(h-1)$-local solid $\Einf$-rings $L_{K(h-1)} ^\solid L_{K(h)} ^\solid \SS \to (L_{K(h-1)} ^\solid E_h ^{\solid})^{h\G_h}$ is an equivalence.
\end{cor}

\begin{proof}
  By \Cref{prop:solid-Galois}, there is an equivalence
  \[L_{K(h)}^{\solid} \SS \simeq \Tot C(\G_h^{\bullet}, E_h ^\solid).\]
  Apply $L_{K(h-1)}^{\solid}$ to this equivalence.
  It follows from descendability that the totalization commutes past the $L_{K(h-1)}^{\solid}$.
  Then the proof of \Cref{prop:Kh-alg-top-solid} implies that there are natural equivalences
  \[L_{K(h-1)}^\solid C(\G_h^{\bullet}, E_h ^\solid) \simeq C(\G_h^{\bullet}, L_{K(h-1)}^\solid E_h ^\solid),\]
  from which the statement follows.
  More precisely, we use the statement that, for a solid $\LKN$ pair $(R,\Gamma)$, $L_{K(h-1)}^{\solid} E^{\solid} (R,\Gamma)$ is computed by
  \[S \mapsto L_{K(h-1)} E(R(S),\Gamma).\]
  In particular, no additional sheafification is necessary.
\end{proof}

\begin{rmk} \label{rmk:kh-1-general-G}
    Let $G \subseteq \G_h$ denote a closed subgroup. By \Cref{cor:solidHRdesc}, the factorization
    \[L_{K(h)} ^{\solid} \SS \to (E_h ^\solid) ^{h G} \to E_h^{\solid}\]
    implies that $(E_h ^\solid) ^{h G} \to E_h^{\solid}$ is descendable as a map of solid $K(h)$-local $\Einf$-rings. Moreover, it is a faithful $G$-Galois extension by \Cref{prop:solid-Galois}.
    It therefore follows as in \Cref{cor:kh-1-still-descent} that
    \[L_{K(h-1)}^\solid ((E_h ^\solid) ^{h G}) \to (L_{K(h-1)}^\solid E_h ^\solid)^{hG}\]
    is an equivalence.
\end{rmk}

The following result will be useful in what comes below.

\begin{proposition} \label{prop:hG-E-commute}
    Let $(S,\Gamma)$ denote a solid $\LKN$ pair with an action of a profinite group $G$, and $(R,\Gamma)$ denote another solid $\LKN$ pair, which we equip with the trivial action of $G$.
    Then if $(R,\Gamma) \to (S,\Gamma)$ is a $G$-equivariant map such that the map $R \to S^{hG}$ is an equivalence and the map $R_0 \to R(\ast) \to S(\ast)$ is flat and relatively perfect for some choice of witness $(R_0, \Gamma_0) \to (R(\ast), \Gamma)$, the natural map
    $E^{\solid} (R,\Gamma) \to (E^{\solid} (S,\Gamma))^{h G}$
    is an equivalence of solid $\Einf$-rings.
\end{proposition}

\begin{proof}
    Associated to the action of $G$ on $(S,\Gamma)$ and the map $(R,\Gamma) \to (S,\Gamma)$ is an augmented cosimplicial solid $\LKN$ pair
    \begin{center}
        \begin{tikzcd}
            (R,\Gamma) \ar[r] & (S,\Gamma) \ar[r, shift left = 3] \ar[r, shift right = 3] & C(G,(S,\Gamma)) \ar[l] \ar[r] \ar[r, shift left = 6] \ar[r, shift right = 6] & C(G \times G, (S,\Gamma)) \ar[l, shift left = 3] \ar[l, shift right = 3] & \cdots
        \end{tikzcd}
    \end{center}
    which is a limit diagram by assumption.
    Since the solid Lubin--Tate theory functor preserves cotensors by profinite sets \Cref{rmk:prof-cotensor-preserved}, 
    we find that there is a natural equivalence
    \[E^{\solid} (S,\Gamma)^{hG} \simeq \Tot C(G^\bullet, E^{\solid} (S,\Gamma)) \simeq \Tot E^{\solid} (C(G^{\bullet},(S,\Gamma))),\]
    so that it suffices to prove that $E^{\solid} (-)$ takes the above diagram to a limit diagram.
    Evaluating at each profinite set $T$, this is a consequence of \Cref{prop:tot-preserved}.
\end{proof}

\begin{proposition} \label{prop:not-just-smooth}
  The natural map $L_{K(h-1)}^{\solid} E_h^{\solid} \to (\mathbb{B}^\solid)^{h \overline{\G}_{h-1}}$ is an equivalence.
\end{proposition}

\begin{proof}
  By \Cref{prop:hG-E-commute}, it suffices to prove that $K \to L^{h\overline{\G}_{h-1}}$ is an equivalence.
  Since $K \to L$ is a $\overline{\G}_{h-1} ^{\sm}$-torsor, the map $K \to L^{h\overline{\G}_{h-1}^{\sm}}$ is an equivalence.
  It follows from \Cref{rslt:smooth-F-p-reps-embed-fully-faithfully} that the map $L^{h\overline{\G}_{h-1}^{\sm}} \to L^{h\overline{\G}_{h-1}}$ is an equivalence, so we are done.
\end{proof}

\begin{rmk}
    Note, however, that the map $L_{K(h-1)}^{\solid} E_h^{\solid} \to \mathbb{B}^{\solid}$ is \emph{not} a $K(h-1)$-local solid $\overline{\G}_{h-1}$-Galois extension.
    Instead, it is instead of a $K(h-1)$-local solid $\overline{\G}_{h-1}^{\mathrm{sm}}$-Galois extension.
    This is consequence of the fact, pointed out in \Cref{exm:B}, that $K \to L$ is a $\overline{\G}_{h-1}^{\mathrm{sm}}$-Galois extension, not a $\overline{\G}_{h-1}$-Galois extension.
\end{rmk}

Combining \Cref{cor:kh-1-still-descent} and \Cref{prop:not-just-smooth}, we obtain the following corollary.

\begin{cor}
    The natural map $L_{K(h-1)}^{\solid} L_{K(h)}^\solid \SS \to (\mathbb{B}^{\solid})^{h\G_{h,h-1}}$ is an equivalence of solid $\Einf$-rings.
\end{cor}

\newpage
\section{The chromatic splitting conjecture and its variants}\label{sec:CSC}

The main purpose of this section is to present the proof of \cref{thmx:main}. Assembling the results from the previous sections, the key ingredient is \cref{thm:bhatfactor1} which identifies the $\G_h^1$ homotopy fixed points of $\widehat{\B}$. As we will explain at the end of the section, similar techniques allow us to reduce the proof of \cref{thmx:height2} on the strong form of the chromatic splitting conjecture to a cohomological criterion (\cref{PropDeperfection}), which will then be studied in the next section. 

We begin with a quick review of the key chromatic players and the chromatic splitting conjecture(s). 

\subsection{Review of the chromatic splitting conjecture}\label{ssec:csc}

In this expository section, we briefly review the chromatic splitting conjecture and its variants in order to put our results in context. 

\begin{conj}[Hopkins' chromatic splitting conjecture]\label{conj:CSC}
    For every $h \geq 1$ and sufficiently generic
    prime $p$, there is an equivalence of spectra
    \begin{equation}\label{eq:CSC}
        L_{h-1}L_{K(h)}\SS^0 \cong 
            \bigoplus_{0\leq j\leq h \atop 1 \leq i_1 < \ldots < i_j \leq h} \left (\bigotimes_{k=1}^j L_{h-i_k}\SS_p^{1-2i_k} \right ),
    \end{equation}
    where the tensor product is taken in the category of $L_{h-1}\SS_p^0$-modules. 
\end{conj}

\begin{rmk}
    The original form of Hopkins' chromatic splitting conjecture as stated above is known to fail for $p=2=h$, as already discussed in the introduction, see also \cref{rmk:CSCknowncases}. Currently, we do not have a formulation of the conjecture that works for all heights and primes uniformly, i.e., predicts the necessary amendments at small primes for heights $h \geq 3$. However, for all $h \leq 2$ and all primes, the canonical map $L_{h-1}\SS_p^0 \to L_{h-1}L_{K(h)}\SS_p^0$ is known to be split; the \emph{weak chromatic splitting conjecture} asks whether this is always the case.
\end{rmk}

There are two kinds of variants of the splitting conjecture that one may naturally take into consideration: On the one hand, one can isolate specific intermediate chromatic heights by applying $L_{K(t)}$ for $0 \leq t < n$ and analyze the resulting conjecture there. On the other hand, one might hope to at least construct maps from certain parts of the right hand side of \eqref{eq:CSC} to the left side in order to realize them as summands of $L_{h-1}L_{K(h)}\SS^0$, without establishing a complete decomposition. Below, both directions will be studied in tandem. 

In the first direction, by construction, for every pair of non-negative integers $t,m$ we have canonical equivalences of functors
\[
    L_{K(t)}L_m \simeq 
    \begin{cases}
        L_{K(t)} & 0 \leq t \leq m \\
        0 & m<t.
    \end{cases}
\]
Therefore, fixing $0 \leq t \leq h$ and applying $L_{K(t)}$ to \eqref{eq:CSC} results in a coheight $h-t$ variant of the chromatic splitting conjecture. The case most relevant to the present paper is when $t = h-1$. The coheight 1 chromatic splitting conjecture thus reads:

\begin{conj}[Hopkins' chromatic splitting conjecture in coheight 1]\label{conj:CSCcoht1}
    For every $h \geq 1$ and sufficiently generic prime $p$, there is an equivalence of spectra
    \begin{equation}\label{eq:CSCcoht1}
        L_{K(h-1)}L_{K(h)}\SS^0 \cong L_{K(h-1)}\SS_p^0 \oplus L_{K(h-1)}\SS_p^{-1}.
    \end{equation}
\end{conj}
Note that, when $h>1$, there is a canonical equivalence $L_{K(h-1)}\SS \simeq L_{K(h-1)}\SS_p$, so we may remove the $p$-completion on the spheres on the right hand side of \eqref{eq:CSCcoht1}. 

\begin{rmk}\label{rmk:CSCknowncases}
    Let us summarize the state of the art on the splitting conjecture. More detailed discussions along with comprehensive lists of references can be found in \cite{BarthelBeaudry,Barthel2020,beaudry2022chromatic,BSSW}.
        \begin{enumerate}
            \item At height $h=1$, the chromatic splitting conjecture is explained in \cite[\S4.3]{BarthelBeaudry}.
            \item For $h=2$ and odd primes, the chromatic splitting conjecture was proven through work of many authors, including Shimomura--Yabe (\cite{ShimomuraYabe1995}), Hopkins (unpublished), and revisited by Behrens (\cite{Behrens2012}) and Lader (\cite{Lader2013}) for $p>3$ and by Goerss, Henn, Karamanov, Mahowald, Rezk, and Shimomura \cite{Shimomura2000,GHMR2005,HKM2013,GHM2014} for $p=3$.
            \item At $h=2$ and $p=2$, the chromatic splitting conjecture does not hold, and a modified version containing two additional (shifted) $\SS^0/p$-summands is proven in \cite{beaudry2022chromatic}. This implies that the coheight 1 chromatic splitting conjecture also fails for $(h,p) = (2,2)$, while the coheight 2 chromatic splitting conjecture holds here.
            \item The coheight $h$ chromatic splitting conjecture has been established for all heights $h\geq 1$ in \cite{BSSW}.
            \item There have been important advances regarding the chromatic splitting conjecture since the writing of this paper. Let us only comment here on how they relate to the chromatic splitting conjecture in coheight $1$; we leave a discussion of the bigger picture to another day.
            First, Agn\`es Beaudry and TS (\cite{hexagon:2610.00138v1}) have reported on an AI-generated disproof of the weak chromatic splitting conjecture in coheight $1$ for primes which are small relative to $h$.
            Second, another AI-generated document (\cite{OAI:filtration}) claims a proof of the strong chromatic splitting conjecture in coheight one when $p > h+1$, by proving that the hypothesis of \Cref{PropDeperfection} holds in this case.
        \end{enumerate}
\end{rmk}

In the second direction, Devinatz and Hopkins \cite{DevinatzHopkins2004} have constructed an element $\zeta_h \in \pi_{-1}L_{K(h)}\SS$ for every $h\geq 1$ which conjecturally supports the right summand of \eqref{eq:CSCcoht1} upon $K(h-1)$-localization. Geometrically, this element arises through the \emph{Iwasawa extension} $\hS{h}$ of $L_{K(h)}\SS$, also known as the \emph{(local) half-sphere}. In order to recall the construction, consider the \emph{determinant} 
\[ 
\det \colon \G_h \to \Z_p^{\times},    
\]
of $\G_h$, constructed carefully in \cite[Definition 3.16]{BarthelBeaudry}. Composing the determinant with the quotient map to $\Z_p^{\times}/\mu \cong \Z_p$ gives a homomorphism (sometimes called the \emph{reduced determinant}),
    \begin{equation}\label{eq:zetan}
        \zeta_h \colon \G_h \longrightarrow \Z_p,
    \end{equation}
where $\mu=C_2$ if $p=2$ and $\mu = C_{p-1}$ if $p$ is odd. The map $\zeta_h$ corresponds to a class $\zeta_h \in H_c^1(\G_h, \Z_p)$. Denote the corresponding kernels by
    \begin{equation}\label{eq:detsubgroups}
        S\G_h\coloneqq \ker(\det) \quad \text{and} \quad \G_h^1 \coloneqq \ker(\zeta_h).
    \end{equation}
Note that the homomorphism $\zeta_h$ is split surjective since $\Z_p$ is topologically free; therefore we obtain a decomposition $\G_h \cong \G_h^1 \rtimes \Z_p$ which trivializes to a product when $p$ does not divide $h$. 

By the Devinatz--Hopkins theorem \cite[Theorem 1]{DevinatzHopkins2004}, $E_h$ forms a $K(h)$-local pro-Galois extension of $L_{K(h)}\SS$ with group $\G_h$, in the sense of Rognes \cite{RognesGalois}.  Taking $\G_h^1$-homotopy fixed points, we get:
    
\begin{proposition}\label{Prop:ChromaticIwasawa}
    The $\EE_{\infty}$-ring spectrum $E_h^{h\G_h^1}$ is a $K(h)$-local pro-Galois extension of $L_{K(h)}\SS$ with group $\Z_p$. In particular, there is a canonical equivalence $L_{K(h)}\SS \simeq \hS{h}^{h\Z_p}$ of $\EE_{\infty}$-ring spectra.
\end{proposition}

This result provides the motivation for: 

    \begin{defn}\label{def:iwasawasphere}
        The $\EE_{\infty}$-ring spectrum $\hS{h}\coloneqq E_h^{h\G_h^1}$ is called the Iwasawa extension of $L_{K(h)}\SS$, also known as the \emph{chromatic half sphere (at height $h$)}.
    \end{defn}

\begin{rmk}\label{rmk:halfspherevariants}
    \cite[Appendix A]{Behrens2006} gives a ``cellular'' description of $E_h^{h\G_h^1}$, where the cells are attached along $\zeta_h$. Other models of the chromatic half-sphere and the related $\Z_p^{\times}$-extension of $\SS_{K(h)}$ have been constructed in \cite{Westerland2017} and \cite{CSY2024}.
\end{rmk}

Choosing an element $\psi \in \G_h$ such that $\zeta_h(\psi)$ is a topological generator of $\Z_p$, we thus obtain a fiber sequence
\begin{equation}\label{eq:halfspherefiberseq} 
    L_{K(h)}\SS \xrightarrow{\iota} \hS{h} \xrightarrow{\psi-1}  \hS{h} \xrightarrow{\delta} L_{K(h)}\SS^1,
\end{equation}
justifying the name ``chromatic half-sphere'' for $\hS{h}$. The composite $\SS \to \hS{h} \xrightarrow{\delta} L_{K(h)}\SS^1$ yields the desired class $\zeta_h \in \pi_{-1}L_{K(h)}\SS$, lifting the cocycle $\zeta_h \in H_c^1(\G_h, \Z_p)$ via the homotopy fixed point spectral sequence. We propose:

\begin{conj}[The half sphere chromatic splitting conjecture in coheight 1]\label{conj:strongCSC} 
Assume that $p > h \ge 2$. The unit map $L_{K(h-1)}\SS \to L_{K(h-1)}\hS{h}$ of $\EE_\infty$-ring spectra is an isomorphism. 
\end{conj}

\begin{rmk}\label{rem:strongCSCh1}
    The half sphere chromatic splitting conjecture in coheight 1 is not true when $h=1$. Indeed, consider for simplicity the case that $p>2$. In this case $\hS{1} \simeq L_p$ is an Adams summand of $p$-complete $K$-theory, so $L_{K(0)}\hS{1} \simeq \Q \otimes L_p \simeq \Q_p[u^{\pm 1}]$ for $u$ in degree $2(p-1)$, while $L_{K(0)}\SS_p \cong \Q_p$.
\end{rmk}

\begin{rmk}\label{rem:strongCSCh2generic}
    For $h=2$ and $p>2$, the half sphere chromatic splitting conjecture in coheight 1 is known to be true. For $p>3$, this is proven en route to the computation of $\pi_* \SS_{K(2)}$ and is thus due to Shimomura--Yabe \cite{ShimomuraYabe1995}, Hopkins (unpublished), and Behrens \cite{Behrens2012}. Since the identification of $L_{K(1)}\hS{2}$ does not seem to be stated explicitly in these sources, we include a brief explanation, using Lader's thesis \cite{Lader2013}. Indeed, from his Corollaire 4.4 we deduce 
        \[
            v_1^{-1}H_{\cts}^*(\G_2^1,E_*/p) \cong \F_p[v_1^{\pm}]\{1,(h_0)_0\}
        \]
    and from his Lemme 4.36 one can show that the element $(h_0)_0$ detects the image of $\alpha_1$. So, the unit $\SS \to \hS{2}=E_2^{h\G_2^1}$ is an equivalence after applying $L_{K(1)}(\SS/p \otimes -)$. This gives the claim.
    
    For $p=3$, the result is a special case of \cite[Theorem 5.7]{GHM2014} by Goerss--Henn--Mahowald, based on previous computations of Shimomura \cite{Shimomura2000}, Goerss--Henn--Maho-wald--Rezk \cite{GHMR2005}, and Henn--Karamanov--Mahowald \cite{HKM2013}.
\end{rmk}

\begin{rmk}\label{rem:strongCSCh2p2}
    At height $h=2$ and $p=2$, the $K(1)$-localization of the Iwasasa sphere has been computed in work of Beaudry, Goerss, and Henn \cite{beaudry2022chromatic}; it is given by
        \[
            L_{K(1)}\hS{2} \simeq L_{K(1)}(\SS \oplus \SS^{-2}/2).
        \]
    Since this is not explicitly stated in this form in their paper, we explain how to extract it. Indeed, there is a map
        \[
            1 + x \colon L_{K(1)}(\SS \oplus \SS^{-2}/2) \to L_{K(1)}\SS_{K(2)} \to L_{K(1)}\hS{2},
        \]
    where $x$ is defined in their \cite[Definition 8.4.1]{beaudry2022chromatic}. The proof of \cite[Theorem 8.4.9]{beaudry2022chromatic} shows that this map induces an isomorphism on homotopy groups after tensoring with  the finite type $1$ spectrum $Y =  \SS/2 \otimes \SS/\eta$, from which the claim follows.
\end{rmk}

The next result explains the connection between our strong version of the coheight~1 splitting conjecture to the one by Hopkins:

\begin{proposition}
    When $h\ge 2$, \cref{conj:strongCSC} implies \cref{conj:CSCcoht1}.
\end{proposition}
\begin{proof}
    Indeed, the fiber sequence \eqref{eq:halfspherefiberseq} 
    \begin{enumerate}
        \item exhibits an equivalence of $\EE_{\infty}$-ring spectra $L_{K(h)}\SS \simeq \hS{h}^{h\Z_p}$, and
        \item shows that taking $\Z_p$-fixed points on $\hS{h}$ is a finite (homotopy) limit and thus commutes with localization.
    \end{enumerate}
It follows that the canonical map $L_{K(h-1)}\SS \to L_{K(h-1)}L_{K(h)}\SS \xrightarrow{L_{K(h-1)}\iota} L_{K(h-1)}\hS{h}$ induces, upon passing to $\Z_p$-fixed points a map of $\EE_{\infty}$-ring spectra
    \[
        (L_{K(h-1)}\SS)^{h\Z_p} \to L_{K(h-1)}L_{K(h)}\SS.
    \]
Here, the $\Z_p$-action on the domain is trivial, thus $(L_{K(h-1)}\SS)^{h\Z_p} \simeq L_{K(h-1)}\SS \oplus L_{K(h-1)}\SS^{-1}$. \Cref{conj:strongCSC} implies that the displayed map is an equivalence, hence so is the composite $L_{K(h-1)}\SS \oplus L_{K(h-1)}\SS^{-1} \to L_{K(h-1)}L_{K(h)}\SS$.
\end{proof}

\subsection{The proof of \cref{thmx:main} and complements}

Our next goal is to assemble the proof of \cref{thmx:main} from the material established above. We begin by recalling the key players in the proof, referring to \eqref{eq:allgroups} and \eqref{eq:allgroups2} for the definition of the groups in use. Let $k/\F_p$ be an algebraic extension field containing $\F_{p^{(h-1)h}}$. Write $K\isom k\laurentseries{u_{h-1}}$ for the Laurent series field over $k$ and consider the field extensions viewed as solid rings throughout:
    \begin{equation}\label{eq:toriiextfields}
        K \xrightarrow{\G_{h-1}} L \to \widehat{L}.
    \end{equation}
Here, $L/K$ is the totally ramified $\overline{\G}_{h-1}$-extension constructed by Gross \cite{GrossRamification} and $\widehat{L}/L$ is the $u_{h-1}$-adic completion of $L$. As such, $\widehat{L}$ is a perfectoid field of characteristic $p$, isomorphic to $k\laurentseries{u_{h-1}^{1/p^\infty}}$, see \cref{PropL}. The field $K$ comes equipped with a natural $\G_h$-action which extends to an action of $\G_{h,h-1}$ on $L$ and $\widehat{L}$. 

As explained in Section \ref{sec:solidLT}, the diagram \eqref{eq:toriiextfields} admits a spectral deformation to a diagram of solid $\EE_{\infty}$-ring spectra equipped with $\G_{h,h-1}$-action:
    \begin{equation}\label{eq:toriiextspectra}
        \xymatrix{
        L_{K(h-1)}^{\solid}E_h^{\solid} \ar[r] & \B^{\solid} \ar[r] & \widehat{\B}^{\solid},
        }
    \end{equation}
where the action of $\overline{\G}_{h-1}$ on $L_{K(h-1)}^{\solid}E_h^{\solid}$ is taken to be trivial. 
The sequence \eqref{eq:toriiextfields} is recovered from $\eqref{eq:toriiextspectra}$ by taking $\pi_0$ and reducing modulo the Landweber ideal $I$. The solid $\EE_{\infty}$-ring spectrum $\B^{\solid}$ is equivalent to $L_{K(h-1)}^{\solid}(E_h^{\solid} \otimes_{\SW(k)}^{\solid} E_{h-1}^{\solid})$, hence there is also a diagram of solid $\EE_{\infty}$-ring spectra with $\G_{h,h-1}$-action:
    \begin{equation}\label{eq:toriiextspectra2}
        \xymatrix{
        E_{h-1}^{\solid} \ar[r] & \B^{\solid} \ar[r] & \widehat{\B}^{\solid},
        }
    \end{equation}
where the action of $\overline{\G}_h$ on $E_{h-1}$ is taken to be trivial. We now apply continuous homotopy fixed points to \eqref{eq:toriiextspectra2} for the group $S\overline{\G}_h$ to obtain a diagram of solid $\EE_{\infty}$-ring spectra admitting an action of the quotient $\G_{h,h-1}/S\overline{\G}_h\isom \Z_p^\times \times\G_{h-1}$:
    \begin{equation}\label{eq:phimap}
        \xymatrix{
        \varphi_h\colon E_{h-1}^{\solid} \ar[r] & (\B^{\solid})^{hS\overline{\G}_h} \ar[r] & (\widehat{\B}^{\solid})^{hS\overline{\G}_h}.
        }
    \end{equation}
The next step is to take homotopy fixed points under the subgroup $(\Z/p)^\times\subset\Z_p^\times$. 

We recall from \cref{ThmCohomologyOfLBody} the fiber sequence of derived $k$-modules with $\G_{h-1}\times \Z_p^\times$-action:
\begin{equation*}
    k \to \widehat{L}^{hS\G_h} \to k[-h] 
\end{equation*}
Here $\overline{\G}_{h-1}\times\Z_p^\times$ acts on the $k[-h]$ term via the character $(g,a)\mapsto \det g^{-h}a^{h-1}$ modulo $p$.  

\begin{lemma}\label{lem:heyerjacquetgeneric}
    Let $\mu= (\Z/p)^{\times}$ be the torsion subgroup of $\Z_p^\times$. Then $k[-h]^{h\mu}=0$.
\end{lemma}
\begin{proof}  The action of the cyclic group $\mu$ on $k[-h]$ is nontrivial if and only if $p-1\nmid h-1$.  Assume this is the case: then $k[-h]^{h\mu}$ is computed by a two-term complex with two copies of $k$ and boundary map $g-1$ for a generator $g\in\mu$;  since $g-1$ acts invertibly, we have $k[-h]^{h\mu}=0$.
\end{proof}

\begin{thm}\label{thm:bhatfactor1}
    Assume $h\geq 2$ and $p$ odd, and consider the map $\varphi_h\colon E_{h-1}^{\solid} \to (\widehat{\B}^{\solid})^{hS\overline{\G}_h}$ of solid $\EE_{\infty}$-ring spectra equipped with $\Z_p^{\times} \times \G_{h-1}$-action constructed above. If $p-1\nmid h-1$, then taking homotopy fixed points with respect to the factor $(\Z/p)^{\times} \subseteq \Z_{p}^{\times}$ induces an equivalence
        \[
            \xymatrix{\varphi_h^{h(\Z/p)^\times}\colon E_{h-1}^{\solid} \ar[r]^-{\simeq} & (\widehat{\B}^{\solid})^{h\overline{\G}_h^1}}
        \]
    of solid $\EE_{\infty}$-ring spectra with residual $\Z_p\times \G_{h-1}$-action.
\end{thm}
\begin{proof}
    The map
    \[E_{h-1} ^{\solid} \to \widehat{\B}^{\solid}\]
    arises from the $\G_{h,h-1}$-equivariant map of solid $\LKN$ pairs
    \[(k, \Gamma_{h-1}) \to (\widehat{L}, \Gamma_{h-1,\widehat{L}}).\]
    It therefore follows from \Cref{prop:hG-E-commute} that $\varphi_h^{h(\Z/p)^\times}$ is an equivalence if the map
    $k \to \widehat{L}^{h\overline{\G}_{h}^1}$
    is an equivalence.
    To see that this holds when $p-1 \nmid h-1$, recall the fiber sequence of $\Z_p^{\times} \times \G_{h-1}$-modules established in \cref{ThmCohomologyOfLBody}, see \eqref{EqFiberSeqk}:
        \[
            k \to \widehat{L}^{hS\overline{\G}_h} \to k[-h],
        \]
    and apply \Cref{lem:heyerjacquetgeneric}.
\end{proof}

In the following arguments, we will repeatedly use that homotopy fixed points with respect to the group appearing as the pullback
    \[
        \xymatrix{{\G}_h^1\times_{\Gal(k/\F_p)} \G_{h-1} \ar[r] \ar[d] & \G_{h-1} \ar[d] \\
        {\G}_h^1 \ar[r] & \Gal(k/\F_p)}
    \]
can be computed iteratively by first taking homotopy fixed points with respect to the fiber in the horizontal (respectively in the vertical direction) and then the quotient. 

We are ready for the proof of the first part of \cref{thmx:main}, which we restate here (in an equivalent formulation) for convenience:

\begin{thm}\label{thm:transchromaticIwasawa}
    Let $h\geq 2$, and assume $p-1\nmid h-1$. The $K(h-1)$-localized Iwasawa extension admits a $\Z_p$-equivariant augmentation $L_{K(h-1)}\hS{h} \to \SS_{K(h-1)}$. 
\end{thm}
\begin{proof}
    Since $p$ is odd, we have a decomposition $\Z_p^{\times} \cong (\Z/p)^{\times}\times \Z_p$. Consider the maps of solid $\EE_{\infty}$-ring spectra with $\Z_p\times\G_{h-1}$-action 
        \begin{equation}\label{eq:reducedphimap}
            \varphi_h^{h(\Z/p)^{\times}}\colon E_{h-1}^{\solid} \to (\B^{\solid})^{h\overline{\G}_h^1} \to (\widehat{\B}^{\solid})^{h\overline{\G}_h^1}
        \end{equation}
    derived from \eqref{eq:phimap} by taking homotopy fixed points with respect to the $(\Z/p)^{\times}$-factor in $\Z_p^{\times}$. By \cref{thm:bhatfactor1}, the composite is an equivariant equivalence of solid $\EE_{\infty}$-ring spectra. 
    
    We will study the effect of taking homotopy $\G_{h-1}$-fixed points on \eqref{eq:reducedphimap}. On the one hand, the middle term identifies as
        \begin{align*}
            ((\B^{\solid})^{h\overline{\G}_h^1})^{h\G_{h-1}} & \simeq ((\B^{\solid})^{h\overline{\G}_{h-1}})^{h\G_{h}^1} && \text{commutativity of homotopy fixed points} \\
            & \simeq (L_{K(h-1)}^{\solid}E_h^{\solid})^{h\G_{h}^1} && \text{by  \cref{prop:not-just-smooth}} \\
            & \simeq L_{K(h-1)}^{\solid}((E_h^{\solid})^{h\G_{h}^1}) && \text{by \cref{rmk:kh-1-general-G}} \\ 
            & = L_{K(h-1)}\hS{h},
        \end{align*}
    compatibly with the residual $\Z_p$-action. On the other hand, the Devinatz--Hopkins theorem supplies a $\Z_p$-equivariant equivalence $L_{K(h-1)}\SS \simeq E_{h-1}^{h\G_{h-1}}$. Substituting these identifications into \eqref{eq:reducedphimap}, we obtain maps of $\Z_p$-equivariant $\EE_{\infty}$-ring spectra
        \[
            L_{K(h-1)}\SS \to L_{K(h-1)}\hS{h} \to ((\widehat{\B}^{\solid})^{h\overline{\G}_h^1})^{h\G_{h-1}}(*),
        \]
    whose composite is an equivalence. This completes the proof. 
\end{proof}

Corollary \ref{corx:main} is then an immediate consequence:

\begin{cor}\label{cor:transchromaticIwasawa}
    Let $h\geq 2$ and assume $p-1\nmid h-1$. Then the map of $\EE_{\infty}$-ring spectra
        \[
            \SS_{K(h-1)}^{B\Z_p} \to L_{K(h-1)}\SS_{K(h)}
        \]
    admits a retract.
\end{cor}
\begin{proof}
    This follows from \cref{thm:transchromaticIwasawa} upon applying homotopy fixed points with respect to the remaining $\Z_p$-action. Indeed, the action on the domain is trivial, while for the codomain we use \cref{Prop:ChromaticIwasawa} and invoke the presentation \eqref{eq:halfspherefiberseq} to commute homotopy fixed points and localization. 
\end{proof}

Let us complement the above results at generic primes ($p-1 \nmid h-1$) by some partial results in the case of odd primes $p$ with $p-1 \mid h-1$ and $h\geq 2$. Recall the map \eqref{eq:phimap} of solid $\EE_{\infty}$-ring spectra with residual $\Z_p^{\times} \times \G_{h-1}$-action
    \[
        \xymatrix{
        \varphi_h\colon E_{h-1}^{\solid} \ar[r] & (\B^{\solid})^{hS\overline{\G}_h} \ar[r] & (\widehat{\B}^{\solid})^{hS\overline{\G}_h}.
        }
    \]
We are now interested in the effect of taking homotopy fixed points on the composite with respect to the first factor. 

\begin{thm}\label{thm:bhatfactor2}
    Let $h \geq 2$, $p>2$, and assume $p-1\mid h-1$. The cofiber $\mathcal{C}$ of $\varphi_h^{h(\Z/p)^{\times}} \colon E_{h-1}^{\solid}\to (\widehat{\B}^{\solid})^{h \overline{\G}_h^1}$ is an invertible $E_{h-1}^{\solid}$-module.
\end{thm}
\begin{proof}
    Since the $(\Z/p)^{\times}$-action under consideration is trivial on the $p$-complete spectrum $E_{h-1}^{\solid}$, we have $(E_{h-1}^{\solid})^{h(\Z/p)^{\times}} \simeq E_{h-1}^{\solid}$, while we also have $(\widehat{\B}^{\solid})^{h \overline{\G}_h^1} \simeq ((\widehat{\B}^{\solid})^{hS\overline{\G}_h})^{h(\Z/p)^{\times}}$. With these identifications, it therefore remains to understand the cofiber of the induced map 
        \[
            \xymatrix{
            \varphi_h^{h(\Z/p)^{\times}}\colon E_{h-1}^{\solid}\to (\widehat{\B}^{\solid})^{h\overline{\G}_h^1}.    
            }
        \]
    In order to prove that $\mathcal{C}$ is invertible as an $E_{h-1}^{\solid}$-module, it suffices to prove that $(\pi_*{\mathcal{C}}) \otimes_{\pi_* E_{h-1} ^\solid}^L \pi_* E_{h-1} ^{\solid}/I$ is invertible as a graded $(\pi_* E_{h-1}^\solid)/I$-module, where $I$ is the Landweber ideal of $E_{h-1}$. Consider the fiber sequence of $\Z_p^{\times} \times \G_{h-1}$-modules established in \cref{ThmCohomologyOfLBody}, see \eqref{EqFiberSeqk}:
        \[
            k \to \widehat{L}^{hS\overline{\G}_h} \to k[-h].
        \]
    In the case that $p-1\mid h-1$, the $\Z_p^\times$-action on the Heyer--Jacquet module is trivial, and taking derived $(\Z/p)^{\times}$-invariants results in a fiber sequence
        \[
            k\to \widehat{L}^{h\overline{\G}_h^1} \to k[-h]. 
        \]
    Letting $u \in (\pi_{-2} E_{h-1})/I$ denote a periodicity class, its image in $(\pi_{-2} \widehat{\B}^{\solid})/I$ is invariant under the action of $\overline{\G}_h^1$.
    We thus find that there is a fiber sequence
    \[k[u^{\pm 1}] \to (\widehat{L}[u^{\pm 1}])^{h\overline{\G}_h^1} \to k[u^{\pm 1}][-h],\]
    from which one may deduce that $(\pi_*{\mathcal{C}}) \otimes_{\pi_* E_{h-1} ^\solid}^L \pi_* E_{h-1} ^{\solid}/I \cong k[u^{\pm 1}][-h]$ is an invertible graded $(\pi_* E_{h-1}^{\solid})/I$-module, and hence that $\mathcal{C}$ is itself an invertible module over $E_{h-1}^{\solid}$.
\end{proof}

\begin{cor}
    Under the assumptions of \cref{thm:bhatfactor2}, the cofiber of the unit map $L_{K(h-1)}^\solid \SS \to (\widehat{\B}^\solid)^{h (\G_h^1\times_{\Gal(k/\F_p)} \G_{h-1})}$ is an invertible $K(h-1)$-local solid spectrum.
\end{cor}
\begin{proof}
    The cofiber in question is $\mathcal{C}^{h\G_{h-1}}$ with $\mathcal{C}$ as in \Cref{thm:bhatfactor2}.
    It follows from the computation of the derived mod $I$ homotopy groups of $\mathcal{C}$ that it is discrete as a solid $K(h-1)$-local spectrum.
    By \Cref{prop:smoothdiscrete}, it follows that the $\G_{h-1}$-action arises from a unique $\G_{h-1}^{\sm}$-action, and that the natural map $\mathcal{C}^{h\G_{h-1}^{\sm}} \to \mathcal{C}^{h\G_{h-1}}$ is an equivalence.

    Now, it follows from \Cref{cor:solidHRdesc}, \Cref{prop:solid-Galois}, and \Cref{prop:smoothdiscrete} that
    \[L_{K(h-1)} ^\solid \SS \to E_{h-1}^\solid\]
    is a descendable $K(h-1)$-local $\G_{h-1}^{\sm}$-Galois extension.
    Since $C^{\sm} (\G_{h-1}, -)$ commutes with the solid tensor product, the $E_{h-1}^{\solid}$-semilinear $\G_{h-1}^{\sm}$-action on $\mathcal{C}$ may be viewed as a descent datum for this map.
    The result now follows from the fact that module categories, and hence invertible modules, satisfy descent for descendable maps.
\end{proof}

\begin{rmk}
    We finish this subsection with a few comments on possible extensions of \cref{thm:bhatfactor2}.
        \begin{enumerate}
            \item Related to \cref{rem:HJMquestions}, it may be possible to identify the $\G_{h-1}$-equivariant invertible $E_{h-1}^{\solid}$-module $\mathcal{C}$ from \cref{thm:bhatfactor2}. On underlying spectra, such modules are specified by pairs $(\mathcal{L},\phi)$, where $\mathcal{L}$ is a $\G_{h-1}$-equivariant invertible $E_{h-1}$-module, and $\phi$ is a $\G_{h-1}$-equivariant  automorphism of $\mathcal{L}$. 
            The object $\mathcal{L}$ represents a class in the $K(h-1)$-local Picard group $\Pic_{h-1,p}$.  Its homotopy $\pi_*\mathcal{L}$ is a $\G_{h-1}$-equivariant invertible $\pi_*E_{h-1}$-module; i.e., it is a Morava module in the sense of \cite{HopkinsPicard}.  The main theorem of \cite{HopkinsPicard} is a complete description of the group of Morava modules, and it would be interesting to determine which class is represented by $\pi_*\mathcal{L}$.
            \item In the case that $p$ is odd and $p-1\mid h-1$, the map $\varphi_h^{h\Z_p^{\times}}$ is not an isomorphism, but one could still ask whether it admits a retract, which is to say, whether $\mathcal{C}\to \Sigma E_{h-1}^{\solid}$ is null.  If so, then $\SS_{K(h-1)}^{B\Z} \to L_{K(h-1)}\SS_{K(h)}$ has a retract as well, albeit a priori only as spectra. 
            \item The case $p=2$ merits further study. In that case, \cite[Theorem 1.1.9]{beaudry2022chromatic} finds an isomorphism $L_{K(1)}\SS_{K(2)}\isom \SS_{K(1)}^{B\Z} \oplus \SS^{-2}/2^{B\Z}$, where $\SS/2$ is the mod 2 Moore spectrum.  We suspect that the $\SS^{-2}/2^{B\Z}$ term must be related to the Heyer--Jacquet module.  
        \end{enumerate}
\end{rmk}

\subsection{Equivalent statements to the strong chromatic splitting conjecture}\label{ssec:eqstrongCSC} 

Having settled \cref{thmx:main}, we return in this final subsection to the strong form of the conjecture (\Cref{conj:strongCSC}), which states that for $h\geq 2$, the unit map $L_{K(h-1)}\SS \to L_{K(h-1)}\hS{h}$ is not just split but in fact an equivalence for $p>h$. We remind the reader that this extra hypothesis is expected to be necessary, see \Cref{rmk:astra-compute-kernel} and \cref{rem:strongCSCh2p2}.

We first give several reformulations of the strong chromatic splitting conjecture in coheight 1. Recall the maps $E_{h-1}^{\solid} \to \B^{\solid}$ and $\B^{\solid} \to \widehat{\B}^{\solid}$ from \eqref{eq:toriiextspectra2}.  Also let $\SS_{K(h),\infty}^\solid$ be the solid version of the Iwasawa sphere:  this is the homotopy fixed points of $E_h^{\solid}$ under $\G_h^1$.  

\begin{lemma}\label{lem:TFAEsCSC} 
    Assume that $h\geq 2$ and $p-1\nmid h-1$.  The following are equivalent:
        \begin{enumerate} 
            \item The unit map $\iota\colon L_{K(h-1)}^\solid \SS \to L_{K(h-1)}^\solid\SS_{K(h),\infty}^\solid$ is an equivalence;
            \item the induced map $E_{h-1}^{\solid}\to (\B^{\solid})^{h\overline{\G}_h^1}$ is an equivalence;
            \item the induced map $(\B^{\solid})^{h\overline{\G}_h^1} \to (\widehat{\B}^{\solid})^{h\overline{\G}_h^1}$ is an equivalence; 
            \item the induced map $(\B^{\solid})^{h({\G}_h^1\times_{\Gal(k/\F_p)} \G_{h-1})}\to (\widehat{\B}^{\solid})^{h({\G}_h^1\times_{\Gal(k/\F_p)} \G_{h-1})}$ is an equivalence.
        \end{enumerate}
\end{lemma}
\begin{proof} As in \cref{thm:bhatfactor1}, we observe that the composition
        \[
            \varphi^{h(\Z/p)^{\times}}\colon E_{h-1}^{\solid} \to (\B^{\solid})^{h\overline{\G}_{h}^1} \to (\widehat{\B}^{\solid})^{h\overline{\G}_{h}^1}
        \]
    of solid $\EE_{\infty}$-ring spectra with residual $\Z_p\times \G_{h-1}$-action is an equivalence.  This immediately gives $(2) \iff (3)$.  Likewise, applying homotopy $\G_{h-1}$-fixed points supplies an equivalence
        \[
            \varphi^{h((\Z/p)^{\times}\times \G_{h-1})}\colon L_{K(h-1)}\SS\simeq (E_{h-1}^{\solid})^{h\G_{h-1}} \to (\B^{\solid})^{h({\G}_h^1\times_{\Gal(k/\F_p)} \G_{h-1})} \to (\widehat{\B}^{\solid})^{h({\G}_h^1\times_{\Gal(k/\F_p)} \G_{h-1})},
        \]
    using the exact sequence $1 \to \overline{\G}_h^1 \to {\G}_h^1\times_{\Gal(k/\F_p)} \G_{h-1} \to \G_{h-1} \to 1$.
    The proof of \cref{thm:transchromaticIwasawa} establishes an equivalence $L_{K(h-1)}\hS{h} \simeq (\B^{\solid})^{h({\G}_h^1\times_{\Gal(k/\F_p)} \G_{h-1})}(*)$, so we deduce that $(1) \iff (4)$. Since $(3) \implies (4)$ follows from applying  homotopy $\G_{h-1}$-fixed points, it remains to prove $(1) \implies (2)$.

    Consider the effect of tensoring the map in $(1)$ with $E_{h-1}^\solid$ and then applying $L_{K(h-1)}^\solid$:
        \begin{eqnarray*}  
            L_{K(h-1)}(E_{h-1}^\solid\otimes^\solid \iota)\colon E_{h-1}^{\solid} 
            &\to& 
            L_{K(h-1)}^\solid( E_{h-1}^{\solid} {\otimes}^{\solid}L_{K(h-1)}^\solid \hS{h}^\solid)\\
            &=& L_{K(h-1)}(E_{h-1}^{\solid} {\otimes}^{\solid}_{\SW(k)} L_{K(h-1)}(E_{h}^{\solid})^{h\overline{\G}_{h}^1}).
        \end{eqnarray*}
      Now we claim that the map 
        \begin{equation} 
    \label{eq:projection_formula}
         L_{K(h-1)}(E_{h-1}^{\solid} {\otimes}^{\solid}_{\SW(k)} L_{K(h-1)}(E_{h}^{\solid})^{h\overline{\G}_{h}^1}) 
        \to L_{K(h-1)}^\solid
        (E_{h-1}\otimes^{\solid}_{\SW(k)} L_{K(h-1)}^{\solid}E_h^\solid)^{h\overline{\G}_h^1}=(\B^{\solid})^{h\overline{\G}_{h}^1} 
        \end{equation}
        is an equivalence.  This would establish $(1) \implies (2)$, as if  $\iota$ is an  equivalence, then so is $L_{K(h-1)}(E_{h-1}^\solid\otimes^\solid \iota)$.  
        
    The argument goes as follows:  the Devinatz--Hopkins extension $L_{K(h)}^\solid \SS\to E_h^{\solid}$ is  $\G_h$-Galois and descendable (\cref{cor:solidHRdesc}), and therefore the intermediate extension $(E_h^{\solid})^{h\overline{\G}_h^1} \to E_h^\solid$ is $\overline{\G}_h^1$-Galois and descendable.  This implies that the Amitsur complex with terms $C((\overline{\G}_h^1)^\bullet, E_h^\solid)$ converges quickly to $(E_h^\solid)^{h\overline{\G}_h^1}$ \cite[Proposition 2.26]{AkhilNilpotence}.  Quick convergence is preserved by exact functors \cite[Proposition 2.20]{AkhilNilpotence}.  In particular, upon applying the exact functor $X\mapsto L_{K(h-1)}^\solid(E_{h-1}^\solid \otimes^\solid_{\SW(k)} L_{K(h-1)}^\solid X)$, we find that the cosimplicial complex with terms
    \begin{equation} 
    \label{eq:amitsur_terms}
    L_{K(h-1)}^\solid (E_{h-1}^\solid \otimes^\solid_{\SW(k)} L_{K(h-1)}^{\solid} C((\overline{\G}_h^1)^\bullet, E_h^\solid)) 
    \end{equation}
    has limit $L_{K(h-1)}(E_{h-1}^{\solid} {\otimes}^{\solid}_{\SW(k)} L_{K(h-1)}(E_{h}^{\solid})^{h\overline{\G}_{h}^1})$.
   The terms in \eqref{eq:amitsur_terms} are built out of the $S$-points of $E_h^\solid$ for $S$ profinite by applying the operations $L_{K(h-1)}^\solid$ and tensoring with $E_{h-1}^\solid=E^\solid(k^\delta,\Gamma_{h-1})$. By   \cref{prop:Kh-alg-top-solid} and \cref{prop:solid-tensor-thm}, these operations can be performed on the level of LKN pairs, where they commute with taking $S$-points by design (see Definition \ref{defn:notimes-solid}).  Therefore the complex in \eqref{eq:amitsur_terms} is equivalent to $C((\overline{\G}_h^1)^\bullet, \B^\solid)$, whose limit is $(\B^\solid)^{h\overline{\G}_h^1}$.  
\end{proof}

In order to deduce a cohomological criterion for the strong chromatic splitting conjecture in coheight 1, we require an analogue of the identification 
    \begin{equation}\label{eq:BGh-1}
        (\B^{\solid})^{h\overline{\G}_{h-1}}\isom E^{\solid}(K,\fg_h^{\un}\otimes_{\pi_0E_h} K) \isom L_{K(h-1)}^{\solid} E_h^{\solid}
    \end{equation}
    from Proposition \ref{prop:not-just-smooth} for $\widehat{\B}^{\solid}$. We begin with an auxiliary lemma.

\begin{lemma}\label{LemPerfectoidFields}
    Let $F$ be a perfectoid field, and let $E/F$ be a (possibly infinite) Galois extension with profinite group $G$. Let $\widehat{E}$ be the completion of $E$ with respect to the unique metric on the field $E$ which extends the metric on $F$.  Then the natural map $F\to H^*_{\cts}(G,\widehat{E})$ is an isomorphism of graded $F$-algebras.
\end{lemma}
\begin{proof} 
    The morphism $\Spa \widehat{E} \to \Spa F$ is a perfectoid pro-\'etale $G$-torsor, and so by \cite[Proposition 3.6.3]{BSSW} we have an isomorphism in the derived category of condensed $\OO_F$-modules:
        \begin{equation}\label{EqGTorsorInFields}
            R\Gamma((\Spa F)_{\proet}, \widehat{\OO}^+) \isom R\Gamma((\Spa \widehat{E})_{\proet},\widehat{\OO}^+)^{hG}. 
        \end{equation}
    By \cite[Theorem 1.2]{ScholzeDiamond}, whenever $X=\Spa(R,R^+)$ is a perfectoid affinoid space, the map $R^+\to H^*(X_{\proet},\OO^+)$ is an almost isomorphism, where $\OO^+=\widehat{\OO}^+$ in our conventions. (The cited result is about the finer v-topology, which suffices for our purposes.)  Therefore, base changing \eqref{EqGTorsorInFields} from $\OO_F$ to $F$ yields an isomorphism
        \[ 
            F\to H^*_{\cts}(G,\widehat{E}),
        \] 
    as claimed. 
\end{proof}

To state the desired identification of $(\widehat{\B}^{\solid})^{h\overline{\G}_{h-1}}$, we introduce the notation
    \[ 
        K_p \coloneqq \bigcup_{n\geq 1} K^{1/p^n} 
    \]
for the colimit perfection of $K=k\laurentseries{u_{h-1}}$.  Let 
    \[
        \val_K\from K\to \Z \cup\set{\infty} 
    \] 
denote the $u_{h-1}$-adic valuation on $K$, and also its $\Z[1/p]$-valued extension to $K_p$.  Let $\widehat{K_p}$ be the completion of $K_p$ with respect to $\val_K$.   Then $\widehat{K}_p$ is a perfectoid subfield of $\widehat{L}$ admitting a continuous action of $\G_h$.  A theorem of Ax \cite{Ax1970} shows that the subfield of $\overline{\G}_{h-1}$-invariants in $\widehat{L}$ is exactly $\widehat{K}_p$. 

\begin{lemma}\label{LemmaFixedSpectra} 
    We have equivalences of solid spectra with $\G_h$-action: 
    \begin{enumerate}
        \item $(\widehat{\B}^{\solid})^{h\overline{\G}_{h-1}} \isom E^{\solid}(\widehat{K}_p,(\Gamma_{h,\widehat{K}_p}^{\un})^{\circ})$
        \item  $(\hB^{\solid})^{h({\G}_h^1\times_{\Gal(k/\F_p)} \G_{h-1})}\isom E^{\solid}(\widehat{K}_p,(\Gamma_{h,\widehat{K}_p}^{\un})^{\circ})^{h\G_h^1}$. 
 \end{enumerate}
\end{lemma}
\begin{proof} 
    For (1):  Recall from \cref{ExBtoBhat} that $E^{\solid}(\widehat{L},(\fg_{h,\widehat{L}}^{\un})^{\circ}) \isom E^{\solid}(\widehat{L},\fg_{h-1}\otimes_k \widehat{L}) =\hB^{\solid}$.  Since $\widehat{K}_p$ is the $\overline{\G}_{h-1}$-fixed subfield of $\widehat{L}$, we obtain a natural map 
        \[
            E^{\solid}(\widehat{K}_p,(\Gamma_{h,\widehat{K}_p}^{\un})^{\circ})\to E^{\solid}(\widehat{L},(\fg_{h,\widehat{L}}^{\un})^{\circ})^{h\overline{\G}_{h-1}} \isom (\widehat{\B}^{\solid})^{h\overline{\G}_{h-1}}
        \]
    from the equivariant map of solid $\LKN$ pairs $(\widehat{K}_p,(\Gamma_{h,\widehat{K}_p}^{\un})^{\circ}) \to (\widehat{L},(\fg_{h,\widehat{L}}^{\un})^{\circ})$.
    This map is an equivalence by \Cref{LemPerfectoidFields} and \Cref{prop:hG-E-commute}.

    Part (2) follows from (1) via the following manipulation:
        \[ 
            (\widehat{\B}^{\solid})^{h({\G}_h^1\times_{\Gal(k/\F_p)} \G_{h-1})}
                \isom ((\widehat{\B}^{\solid})^{h\overline{\G}_{h-1}})^{h\G_h^1}
                \isom E^{\solid}(\widehat{K}_p,(\Gamma_{h,\widehat{K}_p}^{\un})^{\circ})^{h\G_h^1},
        \] 
    using the exact sequence $1 \to \overline{\G}_{h-1} \to {\G}_h^1\times_{\Gal(k/\F_p)} \G_{h-1} \to \G_h^1 \to 1$.
\end{proof}

With this identification in hand, we obtain an algebraic criterion for the strong chromatic splitting conjecture in coheight 1:

\begin{proposition}\label{PropDeperfection} 
    Assume that the natural map of bigraded solid rings
        \[
            H^*_{\cts}(\G_h^1,\pi_*E_h^\solid \otimes_{\pi_0E_h^\solid} K)\to H^*_{\cts}(\G_h^1,\pi_*E_h^\solid \otimes_{\pi_0E_h^\solid} \widehat{K}_p) 
        \] 
    is an isomorphism. 
    Then the map $(\B^{\solid})^{h ({\G}_h^1\times_{\Gal(k/\F_p)} \G_{h-1})}\to (\widehat{\B}^{\solid})^{h({\G}_h^1\times_{\Gal(k/\F_p)} \G_{h-1})}$ is an equivalence.  If furthermore $p-1\nmid h-1$, then the unit map $\iota\colon \SS_{K(h-1)} \to L_{K(h-1)}\hS{h}$ is an equivalence.
\end{proposition}
\begin{proof}
    By \eqref{eq:BGh-1} and \cref{LemmaFixedSpectra}, this reduces to showing that 
        \[ 
            \pi_*(E^\solid(K,\Gamma_{h,K}^{\un})^{h\G_h^1}) \to \pi_* (E^\solid(\widehat{K}_p,(\Gamma_{h,\widehat{K}_p}^{\un})^{\circ})^{h\G_h^1}) 
        \]
    is an isomorphism. There is a morphism of $\G_h^1$-fixed point spectral sequences converging to this map, with $E_2$-pages
        \[ 
            H^*_{\cts}(\G_h^1,\pi_*E^\solid (K,\Gamma_{h,K}^{\un})) \to H^*_{\cts}(\G_h^1,\pi_*E^\solid (\widehat{K}_p,(\Gamma_{h,\widehat{K}_p}^{\un})^{\circ})). 
        \]
    It suffices to prove that this induces an isomorphism.
    Since the coefficients are $I$-complete, it suffices to check that this is an isomorphism after cofibering out by $(p,v_1, \dots, v_{h-2})$. The resulting map
        \[ 
            H^*_{\cts}(\G_h^1,\pi_*E^\solid (K,\Gamma_{h,K}^{\un})/I) \to H^*_{\cts}(\G_h^1,\pi_*E^\solid (\widehat{K}_p,(\Gamma_{h,\widehat{K}_p}^{\un})^{\circ})/I). 
        \]
    identifies with the map of the proposition, so the claim follows.
    
    The second part of the proposition follows from \cref{lem:TFAEsCSC}.
\end{proof}

The next section develops techniques that allow us to verify the algebraic criterion of \cref{PropDeperfection} at height 2, thereby proving the strong chromatic splitting conjecture at height 2 and odd primes.

\newpage
\section{An acyclicity criterion and the case of height 2}\label{sec:height2}

In this final section, we will develop methods to study the cohomological criterion (\cref{PropDeperfection}) for the strong form of the chromatic splitting conjecture established in the previous section. It relies on an acyclicity criterion (\cref{ThmAcyclicityOne}) based on an elaboration of Lazard's work, which we verify at height 2 and odd primes. The essential computational input for the latter are formulas for the action of the Morava stabilizer group due to Salch.  

\subsection{The Lazard complex}
Assume $p-1\nmid h-1$.  \Cref{PropDeperfection} reduces the strong splitting conjecture to the statement that a certain $\G_h^1$-module has trivial continuous cohomology, namely the graded $K$-vector space $(\widehat{K}_p/K)\omega^{\ast}$.  It suffices to replace $\G_h^1$ with any closed normal subgroup (such as $\overline{\G}_h^1$), as one sees by applying the Lyndon--Hochschild--Serre spectral sequence.  

Lazard \cite{Lazard} invented powerful techniques to control the continuous cohomology of a certain class of profinite groups $G$ (``$p$-valued groups of finite rank'').  Crucially, $\overline{\G}_h^1$ contains a basis of open subgroups lying in this class.  

Suppose $G$ is a $p$-valued group of finite rank, and let $V$ be a $\Z_p\powerseries{G}$-module. Lazard shows that the continuous cohomology of $H^i_{\cts}(G,V)$ is computed by a finite complex, each of whose terms is a direct sum of finitely many copies of $V$.  In particular $H^i_{\cts}(G,V)$ vanishes for $i$ large enough. The observation underlying this construction is that $\Z_p\powerseries{G}$ has a natural filtration, whose associated graded ring is a polynomial algebra $P$ over $\F_p$.  There is a standard finite free resolution of $\F_p$ as a $P$-module, namely the Koszul resolution.  

We state here the result (Lazard calls it ``Serre's lemma'') which lifts free resolutions from a graded ring to any filtered ring which gives rise to it.  First, some preparations concerning resolutions.  Suppose $\Omega$ is any ring, and $A$ is an $\Omega$-algebra equipped with an augmentation homomorphism $\varepsilon\from A\to \Omega$.  A {\em split resolution of $\Omega$ by free $A$-modules} is a pair $(X\to \Omega,s)$, where $X\to \Omega$ represents a complex of $\Omega$-modules of the form
\[ \dots \to X_2\to X_1\to X_0=A\stackrel{\varepsilon}{\to} \Omega \] 
with each $X_i$ free over $A$, and $s$ is a nullhomotopy of the complex $X\to\Omega$.  For such a pair, $\Omega$ is equivalent to $X$ in the derived category of $A$-modules.   

Now suppose $\Omega$ is a complete filtered ring (such as $\Z_p$ or any discrete ring), and that $A$ is a complete filtered $\Omega$-algebra.
We assume that the augmentation ideal $I$ is free-complete over $\Omega$ \cite[I.2.1.17]{Lazard}. The discussion of split resolutions above carries over in the categories of filtered rings and graded rings.  In particular, if $(X\to \Omega,s)$ is a split resolution of $\Omega$ by free $A$-modules, then its associated graded $(\Gr X\to \Gr \Omega,\Gr s)$ is a split resolution of $\Gr \Omega$ by free $\Gr A$-modules.  

\begin{lemma}[Serre's lemma, {\cite[V.2.1.1, Corollaire V.2.1.1.3, Corollaire V.2.1.6]{Lazard}}] \label{LemmaSerre} Let $(Y\to \Gr \Omega,\overline{s})$ be a resolution of $\Gr \Omega$ by graded free $\Gr A$-modules with discrete grading.  
\begin{enumerate}
    \item There exists a split resolution $(X\to \Omega,s)$ of $\Omega$ by filtered free $A$-modules, whose associated graded is $(Y\to \Gr \Omega, \overline{s})$.
    \item The resolution $X\to\Omega$ is unique up to isomorphism (but the isomorphism need not commute with the nullhomotopies). 
    \item If each term of $Y$ is a free $\Gr A$-module of finite rank, then each term of $X$ is a free $A$-module of the same rank.
\end{enumerate}
\end{lemma}

The pair $(X\to \Omega,s)$ is constructed recursively through a limit procedure that uses the completeness of $A$ in an essential way.

\subsection{Continuous cohomology of uniform pro-$p$ groups with coefficients in characteristic $p$}\label{ssec:lazardmodp}

Lazard's results on continuous cohomology apply to the class of complete $p$-valued groups.  At the cost of a little generality, we focus on a convenient subclass, the uniform pro-$p$ groups (also called uniformly powerful pro-$p$ groups).  

\begin{defn} [{\cite[Definitions 3.1 and 4.1]{DDMS}}] Let $d\geq 0$.  A {\em uniform pro-$p$ group of dimension $d$} is a pro-$p$ group $G$ which satisfies these criteria:
\begin{itemize}
    \item $G$ is topologically finitely generated.
    \item If $p$ is odd, then $G/\overline{G^p}$ is abelian.  If $p=2$, then $G/\overline{G^4}$ is abelian.  (Here $G^n$ is the subgroup generated by $n$th powers of elements of $G$, and $\overline{G^n}$ refers to the topological closure.)
    \item Writing $G=P_1(G)\supset P_2(G)\supset P_3(G)\supset \cdots$ for the lower $p$-series of $G$, we have $[P_i(G):P_{i+1}(G)]=p^d$ for all $i\geq 1$.  
\end{itemize}
\end{defn}

We refer the reader to \cite{DDMS} for a comprehensive account of uniform pro-$p$ groups and their properties, some of which we list here.  Suppose $G$ is a uniform pro-$p$ group of dimension $d$.  
Then the $p$th power map $x\mapsto x^p$ induces an isomorphism $P_i(G)/P_{i+1}(G)\to P_{i+1}(G)/P_{i+2}(G)$ for all $i\geq 1$ \cite[Definition 4.1(iii)']{DDMS}. Furthermore, $P_i(G)$ is exactly the subset of $p^{i-1}$st powers in $G$ \cite[Lemma 4.10]{DDMS}.  Finally, 
if $g_1,\dots,g_d$ are topological generators of $G$, then the map (not necessarily a homomorphism) $\Z_p^{\oplus d}\to G$ carrying $(a_1,\dots,a_d)$ onto $g_1^{a_1}\cdots g_d^{a_d}$ is a homeomorphism \cite[Theorem 4.17]{DDMS}.   For each $i\geq 1$, this homeomorphism factors through a bijection $(\Z/p^{i-1}\Z)^{\oplus d} 
\to G/P_{i}(G)$. 

\begin{lemma} 
\label{LemmaCommensurable}
Let $G$ be a uniform pro-$p$ group of dimension $d$, and let $\F_p[G]$ be its group algebra with $\F_p$-coefficients.  Let $I=\ker(\F_p[G]\to\F_p)$ be the augmentation ideal, and for all $i\geq 1$, let $J_i=\ker(\F_p[G]\to\F_p[G/P_i(G)])$.  Then:
\begin{enumerate}
    \item $J_i\subset I^{p^{i-1}}$.
    \item $I^{p^{d(i-1)}}\subset J_i$.
\end{enumerate}
Therefore the collections of ideals $I^i$ and $J_i$ induce the same topology on $\F_p[G]$.
\end{lemma}

\begin{proof} For (1), we proceed by induction on $i\geq 2$, the case $i=1$ being trivial.  It suffices to show that $g-1\in I^i$ for $g\in P_i(G)$.   As noted above, $P_i(G)$ is exactly the subset of $p^{i-1}$th powers in $G$, so $g=h^p$ for some $h\in P_{i-1}(G)$, and then $g-1=(h-1)^p\in J_{i-1}^p\subset I^{p^{i-1}}$ by the inductive hypothesis. 

For (2), consider the action of the pro-$p$ group $G$ on the $p^{d(i-1)}$-dimensional $\F_p$-vector space $V=\F_p[G/P_i(G)]$.  The image of $G$ in $GL(V)$ is contained in a Sylow $p$-subgroup, which up to conjugacy is the subgroup of upper-triangular unipotent matrices of dimension $p^{d(i-1)}$.  Given $p^{d(i-1)}$ such matrices $I+M_j$, the product of the $M_j$ is 0.  Therefore for any $g_1,\dots,g_{p^{d(i-1)}}\in G$, the element $(g_1-1)\cdots (g_{p^{d(i-1)}}-1)$ acts trivially on $V$, which is to say it lies in $J_i$.
\end{proof}  

\begin{lemma} \label{LemmaUniformProps}
Let $G$ be a uniform pro-$p$ group with topological generators $g_1,\dots,g_d$.  Let $I\subset \F_p[G]$ be the augmentation ideal. 

\begin{enumerate}
    \item For all $n\geq 0$ and all $g\in G$ we have $[I^n,g]\subset I^{n+2}$, where $[x,y]=xy-yx$ is the Lie bracket on $\F_p[G]$.  
    \item For $n\geq 0$, let $W_n\subset\F_p[G]$ denote the $\F_p$-span of all elements of the form 
    \[ (g_1-1)^{a_1}\cdots (g_d-1)^{a_d},\]
    where the $a_j\geq 0$ satisfy $\sum_j a_j = n$.  Then for all $n\geq 0$ we have
    \[ I^n=W_n+I^{n+1}. \] 
\end{enumerate}
\end{lemma}

\begin{proof}
For (1):  We proceed by induction on $n$.  For the base case, $I^0=\F_p[G]$ is the $\F_p$-span of elements $h\in G$, and then $[h,g]=hg-gh=(hgh^{-1}g^{-1}-1)gh\in J_2\subset I^p\subset I^2$ by \cref{LemmaCommensurable}(1).  

When $n\geq 1$, suppose $x\in I^{n-1}$ and $y\in I$.  Then 
\[ [xy,g]=x[y,g]+[x,g]y \]
lies in $I^{n+2}$ by the inductive hypothesis.

For (2):  Again we proceed by induction on $n$.  The case $n=0$ being trivial, we assume $n\geq 1$.  By the inductive hypothesis,
\[  I^n=I^{n-1}I=(W_{n-1}+I^n)(W_1+I^2)\subset W_{n-1}W_1+I^{n+1},  \] 
and so we are reduced to the claim that $W_{n-1}W_1\subset W_n+I^{n+1}$.  The product $W_{n-1}W_1$ is the $\F_p$-span of elements of the form 
\[ (g_1-1)^{a_1}\cdots (g_d-1)^{a_d}(g_i-1), \]
with $\sum_j a_j=n-1$.
The difference between the above element and 
\[ (g_1-1)^{a_1}\cdots (g_i-1)^{a_i+1}\cdots (g_d-1)^{a_d}\in W_n \] equals
\[ (g_1-1)^{a_1}\cdots (g_i-1)^{a_i}[x,g_i] \]
for $x=(g_{i+1}-1)^{a_{i+1}}\cdots (g_d-1)^{a_d}\in I^{\sum_{j\geq i+1} a_j}$.  By (1), this difference lies in $I^{n+1}$.
\end{proof}

For a general profinite group $G$, the completed group algebra $\F_p\powerseries{G}$ is defined as the limit of $\F_p[G/H]$ for $H$ running through a cofinal system of open normal subgroups of $G$.  If $G$ is a uniform pro-$p$ group, then one can take the lower $p$-series to be the system of open normal subgroups, and then \cref{LemmaCommensurable} shows that $\F_p\powerseries{G}$ is the completion of $\F_p[G]$ with respect to the augmentation ideal.

The following theorem identifies the associated graded $\F_p$-algebra of the filtration on $\F_p\powerseries{G}$ induced by the augmentation ideal.  

\begin{thm} \label{ThmSymmetric} Let $G$ be a uniform pro-$p$ group of dimension $d$.   Let $W=G/P_2(G)=G/G^p$, so that $W$ is an $\F_p$-vector space of dimension $d$.  With respect to the filtration on $\F_p\powerseries{G}$ induced by powers of the augmentation ideal, the graded $\F_p$-algebra $\Gr \F_p\powerseries{G}$ is isomorphic to $\Symm W$, the symmetric $\F_p$-algebra on $W$.  Given $g\in G$ with image $w\in W$, the isomorphism sends $g-1$ to $w$.  
\end{thm}

\begin{proof} Since $\F_p\powerseries{G}$ is the $I$-adic completion of $\F_p[G]$, there is an isomorphism $\Gr\F_p[G]\to \Gr \F_p\powerseries{G}$, where both $\F_p[G]$ and $\F_p\powerseries{G}$ are filtered by their augmentation ideals.  

The map $G\to I/I^2$ carrying $g$ to $g-1$ is a group homomorphism onto a $p$-torsion abelian group, so it factors through an $\F_p$-linear map $W\to I/I^2$. This map extends uniquely to a map of graded $\F_p$-algebras $\Symm W\to \Gr \F_p[G]$, which carries a polynomial $f(w_1,\dots,w_d)$ onto $f(g_1-1,\dots,g_d-1)$.  We claim this map is an isomorphism.  Its $n$th graded piece is the surjective map $\Symm^n W\to I^n/I^{n+1}$, which factors through the space $W_n\subset \F_p[G]$ from \cref{LemmaUniformProps}.  

We must show that $\Symm^n W\to I^n/I^{n+1}$ is also injective.  Suppose $f_n(w_1,\dots,w_d)\in \Symm^n W$ lies in the kernel.  Applying \cref{LemmaUniformProps} repeatedly, we find for all $N\geq n+1$ a relation
\[ \sum_{j=n}^{N-1} f_j(g_1-1,\dots,g_d-1)\in I^N, \] 
where $f_j$ is a homogeneous polynomial of degree $j$. Choosing $N=p^{d(i-1)}$ for a large enough $i$, we have $I^{N}\subset \ker(\F_p[G]\to \F_p[G/P_i(G)])$ by \cref{LemmaCommensurable}(2).

The point now is that $G/P_i(G)$ is exactly the set of classes of elements $g_1^{a_1}\cdots g_d^{a_d}$, where $0\leq a_j\leq p^{i-1}-1$ for each $j$.  Therefore such elements form a $\F_p$-basis for $\F_p[G/P_i(G)]$.  A routine calculation shows that the $(g_1-1)^{a_1}\cdots (g_d-1)^{a_d}$ also form a basis.  In particular these elements are linearly independent modulo $I^{N}$.  This implies that $f_j=0$ for all $j$, and in particular $f_n=0$.  
\end{proof}

\begin{rmk} We learned \cref{ThmSymmetric} from notes of Pstr\k{a}gowski\footnote{Available at \url{https://www.kurims.kyoto-u.ac.jp/~piotr/291y_notes.pdf}, see Theorem 1.2.}.  The proof that appears here is nothing but an adaptation of those notes to the characteristic $p$ setting.  The theorem is a special case of a result of Lazard \cite[A.3, Th\'eor\`eme 3.5]{Lazard}, which for a general pro-$p$ group $G$ identifies $\Gr \F_p\powerseries{G}$ with the universal enveloping algebra of a certain restricted Lie algebra $\Gr G$ associated with $G$.  When $G$ is a uniform pro-$p$ group, $\Gr G$ is abelian and freely generated by its degree 1 part.
\end{rmk}

We intend to apply Serre's Lemma (\cref{LemmaSerre}) to the situation in \cref{ThmSymmetric}.  The graded algebra $\Symm W$ is a polynomial algebra over $\F_p$.  There is a canonical resolution of $\F_p$ in free modules over a polynomial algebra, namely the Koszul complex.  We review it now.

Let $k$ be any field, and let $W$ be a finite-dimensional $k$-vector space, with dimension $r$.  Then $\Symm W$ is a polynomial $k$-algebra with $r$ generators.  The Koszul complex 
\[ \Kos(W)\colon   \Kos^r(W)\to \Kos^{r-1}(W)\to \cdots \to \Kos^1(W)\to \Kos^0(W) \] 
is a complex of finite-dimensional free $\Symm W$-modules with $n$th term
\[ \Kos^n(W)=\Symm W\otimes \bigwedge^n W \] 
for $0\leq n\leq r$, and differential determined by the formula
\[ d(1\otimes x_1 \cdots x_n) = \sum_{i=1}^n (-1)^{i+1}x_i \otimes x_1\cdots \hat{x}_i \cdots x_n \]
for $x_1,\dots,x_n\in W$.  

Importantly, the augmented complex $\Kos(W)\to k$ admits a nullhomotopy.  The key observation is that if $W=W'\oplus W''$ is a direct sum of $k$-vector spaces, then $\Kos(W)\isom \Kos(W')\otimes \Kos(W'')$ in the category of complexes of $k$-vector spaces.  Therefore a basis for $W$ presents $\Kos(W)$ as a tensor product of complexes of the form $\Kos(k)$.  In that case, the augmented complex $\Kos(k)^+=[\Kos(k)\to k]$ admits a nullhomotopy $s$: 
\[
\xymatrix{
k[x] \ar[r]^{\times x} & k[x] \ar[r] \ar@/^1.0pc/[l]^{s_1} & k \ar@/^1.0pc/[l]^{s_0} }
\]
Here $s_0$ is the inclusion of the constants into $k[x]$, and $s_1(f(x))=(f(x)-f(0))/x$.  Tensoring together a copy of $(\Kos(k)^+,s)$ for each vector in a basis for $W$, we obtain the pair $(\Kos(W)^+,s_W)$, where $s_W$ is a null-homotopy of $\Kos(W)^+$.  (Strictly speaking, $s_W$ depends on the choice of an ordered basis for $W$.)  This shows that $\Kos(W)$ is indeed a split resolution of $k$ by free $\Symm W$-modules.

Consequently, if $V$ is a $\Symm W$-module, the Ext-groups  $\Ext^n_{\Symm W}(
k,V)$ are computed by the complex
\[ C_{\Kos}(\Symm W,V)\colon C^0_{\Kos}(\Symm W,V)\to C^1_{\Kos}(\Symm W,V)\to \cdots \to C^r_{\Kos}(\Symm W,V), \] 
where the $n$th term is 
\[ C^n_{\Kos}(\Symm W,V) = \Hom(\bigwedge^n W, V). \] 

\begin{defn}  Let $G$ be a uniform pro-$p$ group of dimension $d$, and let $W=G/G^p$, an $\F_p$-vector space of dimension $d$.  The {\em Lazard complex} $\Laz(G)$ is a complex of finite free $\F_p\powerseries{G}$-modules 
\[ \Laz^d(G)\to \Laz^{d-1}(G)\to \cdots\to \Laz^1(G)\to \Laz^0(G)\]
with $n$th term
\[ \Laz^n(G)=\F_p\powerseries{G} \otimes \bigwedge^n W\] 
for $0\leq n\leq d$, and the {\em Lazard resolution} is a split resolution of $\F_p$ by finite free $\F_p\powerseries{G}$-modules, consisting of a pair $(\Laz(G)\to \F_p,s_G)$.  The Lazard resolution arises by applying Serre's Lemma: it is the unique resolution whose associated graded is the Koszul resolution $(\Kos(W)\to \F_p,s_W)$, which in turn comes from choosing an ordered basis of $W$. 

Given a representation of $\F_p\powerseries{G}$ on a topological $\F_p$-vector space $V$, let $C_{\Laz}(G,V)$ denote the cochain complex obtained by applying $\Hom(-,V)$ to the complex $\Laz(G)$.  Thus 
\[ C^n_{\Laz}(G,V)=\Hom(\bigwedge^n W,V). \] 
\end{defn}

Since the Lazard complex really is a resolution of the trivial module $\F_p$, the complex $C_{\Laz}(G,V)$ computes the continuous cohomology $H^*_{\cts}(G,V)$.  We have arrived at the characteristic $p$ version of Lazard's theorem.

\begin{thm}\label{ThmLazard} Let $G$ be a uniform pro-$p$ group of finite dimension, and let $V$ be a representation of $\F_p\powerseries{G}$ on a topological $\F_p$-vector space $V$.  The continuous cohomology $H^*_{\cts}(G,V)$ is computed by a complex $C_{\Laz}(G,V)$ whose terms are direct sums of finitely many copies of $V$:
\[ C_{\Laz}^n(G,V)=\Hom(\bigwedge^n G/G^p,V).\] 
\end{thm}

In this context, cohomology of trivial $G$-modules is particularly simple.

\begin{cor} \label{CorCohoTrivial} 
Let $G$ be a uniform pro-$p$ group of dimension $d$, and let $V$ be the trivial representation of $G$ on a topological $\F_p$-vector space.  Then for each $n\in\Z$,
\[ H^n_{\cts}(G,V)\isom \Hom(\bigwedge^n G/G^p,V)\isom V^{\oplus  {d\choose{n}}}. \] 
\end{cor}

\begin{proof} The differentials $\Laz^n(G)\to \Laz^{n-1}(G)$ have coefficients in the augmentation ideal $I$, which means that the differentials in the complex $C_{\Laz}(G,V)$ are zero.   
\end{proof}

It is difficult to write down the differential $\Laz^n(G)\to \Laz^{n-1}(G)$ directly, other than to say that its coefficients lie in $I$.  But the reduction of the differential modulo $I^2$ is known:
\begin{equation}
    \label{eq:differentialmodI2}
 d(1\otimes w_1 \cdots w_n) \equiv  \sum_{i=1}^n (-1)^{i+1} (g_i-1) \otimes w_1\cdots \hat{w}_i \cdots w_n \pmod{I^2} 
 \end{equation} 
where $g_i\in G$ is a lift of $w_i\in W$. 
This calculation comes from the fact that 
the associated graded of $\Laz(G)$ is $\Kos(W)$.

\subsection{An acyclicity criterion using one basis element}

Let $V$ be a complete ultrametric $\F_p$-vector space:  this means $V$ is separated and complete for an ultrametric valuation $\val_V\from V\to \R\cup\set{\infty}$.  

For a linear operator $A\in \End V$, we write
\[ \val_V(A)=\inf_{v\in V} \val_V(Av)-\val_V(v) \in \R\cup\set{\pm \infty}\] 
for its operator valuation.  The operator is bounded (respectively, a strict contraction) if $\val_V(A)>-\infty$ (respectively, $\val_V(A)>0$).  A bounded operator is continuous.  (The converse is false;  the Frobenius map on $V=\F_p\laurentseries{t}$ with its $t$-adic valuation is continuous but unbounded.) 

Now suppose $G$ is a uniform pro-$p$ group of finite dimension which acts linearly on $V$.  We say that the action of $G$ on $V$ is strictly unipotent if for all $g\in G$, the operator $(g-1)$ is a strict contraction on $V$.  In this situation $V$ becomes an $\F_p\powerseries{G}$-module.
 
We state here a condition which assures the acyclicity of the $\F_p\powerseries{G}$-module $V$.

\begin{thm}
\label{ThmAcyclicityOne} Let $G$ be a uniform pro-$p$ group of dimension $d$, with topological generators $g_1,\dots,g_d$, and let $V$ be a complete ultrametric $\F_p$-vector space admitting a strictly unipotent action of $G$. Assume that the operator $(g_1-1)$ admits a bounded inverse $(g_1-1)^{-1}$ on $V$. 
Then $H^*_{\cts}(G,V)=0$.
\end{thm}
\begin{proof}  First we observe that if $G'\subset G$ is an open subgroup, it suffices to show that $H^*_{\cts}(G',V)=0$.  Indeed, since $G$ is pro-$p$, it is pro-nilpotent, and so $G'$ is subnormal in $G$.  Repeated application of the Lyndon--Hochschild--Serre spectral sequence gives $H^*_{\cts}(G,V)=0$.

Let $L$ be the Lie algebra of $G$, with $\Z_p$-basis $x_1,\dots,x_d$ corresponding to $g_1,\dots,g_d$. By \cite[Theorem 9.10]{DDMS}, uniform open subgroups of $G$ correspond to finite-index powerful Lie subalgebras of $L$ of finite index.  For $N\geq 0$, it is easily checked that $x_1,p^Nx_2,\dots,p^Nx_d$ span such a Lie subalgebra.  Let $G'\subset G$ be the corresponding uniform open subgroup corresponding to $L'$;  then $G'$ has topological generators $g_1,g_2^{p^N},\dots,g_d^{p^N}$.  By our hypotheses, $\val_V(g_1-1)^{-1}>-\infty$ and $\val_V(g_i-1)>0$, and so we can choose $N$ large enough so that for $i=2,\dots,d$,
\[ p^N\val_V(g_i-1) > -\val_V(g_1-1)^{-1}. \]
Note that $\val_V(g_i^{p^N}-1)\geq p^N \val_V(g_i-1)$.  
In light of the reduction of the previous paragraph, we may relabel $G'$ as $G$, and assume henceforth that $(g_1-1)^{-1}(g_i-1)$ is a strict contraction on $V$ for $i=2,\dots,d$.

Let $I$ be the augmentation ideal of $\F_p\powerseries{G}$. The continuous cohomology $H^*_{\cts}(G,V)$ is computed by the Lazard complex $C_{\Laz}(G,V)$.  We will construct a nullhomotopy of $C_{\Laz}(G,V)$.  Recall the notation $W=G/G^p$; this is an $\F_p$-vector space with distinguished basis $w_1,\dots,w_d$ in correspondence with $g_1,\dots,g_d$.  We separate out the first basis vector from the rest:
\[ W = W_1\oplus W_c  \]  
where $W_1$ is spanned by the single vector $w_1$, and the complement $W_c$ is spanned by the others.  Accordingly we may decompose 
\[ \bigwedge^n W = \left[\bigwedge^{n-1} W_c \otimes W_1 \right] \oplus \bigwedge^n W_c \] 
and also
\[ C^n_{\Laz}(G,V)=
\left[ V\otimes \bigwedge^{n-1} W_c^* \otimes W_1^* \right]\oplus \left[ V\otimes \bigwedge^n W_c^*\right] , \] 
where the $\ast$ indicates $\F_p$-linear dual.  

With respect to this decomposition, the differential $d^n\from C^n_{\Laz}(G,V)\to C^{n+1}_{\Laz}(G,V)$ takes the block form 
\[ d^n = \begin{pmatrix} 
A_n & B_n \\ 
C_n & D_n
\end{pmatrix}
\]
We now apply \eqref{eq:differentialmodI2}, which identifies $d^n$ modulo $I^2$.  Consider in particular the operator \[ B_n \from V^{\oplus{d-1\choose n}} \isom V\otimes \bigwedge^n W_c^\ast \to V\otimes \left[\bigwedge^n W_c^\ast \otimes W_1^\ast \right]\isom V^{\oplus{d-1\choose n}}. \] 
This can be written 
\[ B_n = (g_1-1)\operatorname{Id}+\gamma, \]
where $\operatorname{Id}$ is the identity operator on $V^{\oplus{d-1 \choose n}}$, and $\gamma$ is a matrix with entries in $I^2$.  

We claim that $B_n$ is invertible on $V^{\oplus{d-1\choose n}}$.  Since $(g_1-1)$ is invertible by hypothesis, this reduces to the invertibility of the operator $\operatorname{Id}+\gamma(g_1-1)^{-1}$.  Formally, the inverse of this operator is 
\begin{equation}
\label{EqSumContraction}
\sum_{j=0}^\infty (-1)^j (\gamma(g_1-1)^{-1})^j, 
\end{equation}
which converges precisely because each entry of the matrix $\gamma(g_1-1)^{-1}$ is a strict contraction by hypothesis.  

Now we claim that 
\[ s^n = \begin{pmatrix} 0 & 0 \\ B_{n-1}^{-1} & 0 \end{pmatrix} \from C^{n}_{\Laz}(G,V) \to C^{n-1}_{\Laz}(G,V) \] 
is a nullhomotopy of $C_{\Laz}(G,V)$.  Indeed,
\begin{eqnarray*}
    d^{n-1}s^n+s^{n+1}d^n&=& 
\begin{pmatrix}
    A_{n-1} & B_{n-1} \\
    C_{n-1} & D_{n-1} 
\end{pmatrix}
\begin{pmatrix}
0 & 0 \\
B_{n-1}^{-1} & 0 
\end{pmatrix}
+
\begin{pmatrix}
0 & 0 \\
B_n^{-1} & 0
\end{pmatrix}
\begin{pmatrix}
    A_n & B_n \\
    C_n & D_n
\end{pmatrix}
\\
&=& 
\begin{pmatrix}
    1 & 0 \\
D_{n-1}B_{n-1}^{-1}+B_n^{-1}A_n & 1 
\end{pmatrix}
\end{eqnarray*}
Examining the relation $d^nd^{n-1}=0$, we see $A_nB_{n-1}+B_nD_{n-1}=0$, which implies $D_{n-1}B_{n-1}^{-1}+B_n^{-1}A_n=0$.  
\end{proof}

\subsection{The Nottingham group} 
\Cref{PropDeperfection} reduces the strong splitting conjecture to the statement $H^*_{\cts}(G,(\widehat{K}_p/K)\omega^{\otimes n})=0$ for all $n\in\Z$, where $G$ is a uniform pro-$p$ normal subgroup of the Morava stabilizer group $\G_h$.  In this subsection we gather some abstract results relevant to this situation.

Suppose $k$ is a perfect field of characteristic $p$, and let  
\[ K = k\laurentseries{t} \]
be the field of Laurent series over $k$.  It is
separated and complete for the $t$-adic valuation $\val_K$. 
Let $\omega=\OO_Ku$ be a free $\OO_K$-module of rank 1, and let $\omega_K=\omega\otimes_{\OO_K}K=Ku$.  From these data we construct the graded $K$-algebra 
\[ 
    \bigoplus \omega_K^{\otimes n}. 
\] 
Each graded piece $\omega_K^{\otimes n}=Ku^n$ comes equipped with a valuation $\val_{Ku^n}$ which assigns the value of 0 to the basis vector $u^n$.

The {\em Nottingham group} $\operatorname{Nott}(k)$ is the group of strictly unipotent automorphisms of the $k$-algebra $K$ which act trivially on the associated graded of $K$.  In other words, an element $g\in \operatorname{Nott}(k)$ is a {\em substitution}
\[ g(t)= t(1+a_1t+a_2t^2+\dots), \] 
with each $a_i\in k$. 
The group of all strictly unipotent automorphisms of the graded $k$-algebra $\bigoplus_{n\in\Z} \omega_K^{\otimes n}$ whose restriction to $K$ lies in $\operatorname{Nott}(k)$ is a split extension of $\operatorname{Nott}(k)$ by the topologically unipotent subgroup $\OO_K^1$ of $\OO_K^\times$, which we might call the extended Nottingham group $\operatorname{Nott}'(k)$:
\[
\operatorname{Nott}'(k)=\operatorname{Nott}(k) \ltimes \OO_K^1. 
\] 
An element $g\in \operatorname{Nott}'(k)$ is specified by two power series in $\OO_K^1$, namely $g(t)/t$ and $g(u)/u$.  

Consider the following $k$-vector spaces $V$ derived from $Ku^n$:
\begin{itemize}
    \item The quotient $Ku^{pn}/(Ku^n)^p=(K/K^p)u^{pn}$.
    \item The base changes $K^{1/p^j}\otimes_K K\omega^{\otimes n}=K^{1/p^j}u^n$ for each $j\geq 0$,
    \item Their union $K_pu^n$,
    \item The completion $\widehat{K}_pu^n$,
    \item The quotient $(\widehat{K}_p/K)u^n$.
\end{itemize}
Each of the above $V$ carries an action of $\operatorname{Nott}'(k)$ and an invariant valuation $\val_V$, which assigns the value $a/p^j$ to $t^{a/p^j}u^n$.
The action is strictly unipotent on each $K^{1/p^j}u^n$, but not on $K_pu^n$ or the spaces derived from it. 

Now suppose $G$ is a uniform pro-$p$ group of finite dimension admitting a homomorphism $G\to\operatorname{Nott}'(k)$.  
We state a variation on \cref{ThmAcyclicityOne} which can be used to establish acyclicity for the action of $G$ on $(\widehat{K}_p/K)u^n$.

\begin{proposition} 
\label{prop:acyclicity_one}
Let $G$ be a uniform pro-$p$ group of dimension $d$, with topological generators $g_1,\dots,g_d$.  Let $G\to\operatorname{Nott}'(k)$ be a homomorphism.  Assume that the operator $(g_1-1)$ admits an inverse $(g-1)^{-1}$ on $(K/K^p)u^{pn}$ for all $n\in\Z$, whose operator valuation is bounded below independently of $n$.  Then for all $n\in\Z$ we have $H^*_{\cts}(G,(\widehat{K}_p/K)u^n)=0$.
\end{proposition}

\begin{proof}
 Since the elements $g_2,\dots,g_d$ act strictly unipotently on $K$ and $Ku$, there exists $B>0$ for which $\val_K(g_i-1)\geq B$ and $\val_{Ku}(g_i-1)\geq B$ for $i=2,\dots,d$.  It is straightforward to check that for all $n\in\Z$, the operator valuation $\val_{Ku^n}(g_i-1)$ is also bounded below by $B$.  

By hypothesis, there exists $C\geq 0$ for which $\val_{(K/K^p)u^{pn}}(g-1)^{-1}\geq -C$ for all $n\in\Z$.  Choose $N$ large enough so that for $i=2,\dots,d$,
we have $p^NB>(p/(p-1))C$.   Replacing $G$ with the subgroup generated by $g_1,g_2^{p^N},\dots,g_d^{p^N}$, we may assume the operator $(g_1-1)^{-1}$ exists on $(K/K^p)u^{pn}$ for all $n\in\Z$, with operator valuation bounded below by $-C$, and that for $i=2,\dots,d$ and all $n\in\Z$ we have $\val_{Ku^n}(g_i-1)\geq B$, with $B>(p/(p-1))C$.  
 
For each $m\geq 1$, the ultrametric space $(K^{1/p^m}/K)u^n$ admits a finite $\operatorname{Nott}'(k)$-stable filtration by closed subspaces $(K^{1/p^i}/K)u^n$ for $1\leq i\leq m$, with subquotients
\[(K^{1/p}/K)u^n ,\dots, (K^{1/p^m}/K^{1/p^{m-1}})u^n;\] these are isomorphic to 
\[ (K/K^p)u^{pn},\dots,(K/K^p)u^{p^mn},\] 
respectively, via $f^{1/p^i} u^n\mapsto f u^{p^in}$. 

With respect to a splitting of this filtration of $(K^{1/p^m}/K)u^n$, the operator $(g_1-1)$ may be written as an upper-triangular matrix with bounded entries, whose $i$th diagonal entry admits an inverse with operator valuation bounded below by $-C/p^i$.  Therefore $(g_1-1)^{-1}$ exists on $(K^{1/p^m}/K)u^n$;  its matrix is upper-triangular, and for $j\geq i$, the $ij$th entry of this matrix has operator valuation bounded below by
\[ -C/p^i-C/p^{i+1}-\cdots-C/p^j \geq -C/p^{i-1}(p-1); \]
ie, entries in the $i$th row are bounded below by $-C/p^{i-1}(p-1)$.  In particular, the operator norm of $(g_1-1)^{-1}$ is bounded below independently of $m$.  

Whereas for each $i=2,\dots,d$, the matrix of $(g_i-1)$ is upper-triangular, and entries in the $i$th column are bounded below by $B/p^i$.
Therefore the entries of $(g_i-1)(g_1-1)^{-1}$ have operator norm bounded by $B/p^i-C/p^{i-1}(p-1)>0$, implying that $(g_i-1)(g_1-1)^{-1}$ is a strict contraction on $(K^{1/p^m}/K)u^n$. 

Applying the results of the above two paragraphs for all $m\geq 1$, we find that $(g_1-1)^{-1}$ exists as a bounded operator on $(K_p/K)u^n$, which extends to the completion $(\widehat{K}_p/K)u^n$.  We also find that for any $\gamma\in I^2$ (with $I$ the augmentation ideal of $\F_p\powerseries{G}$), the terms of the series \eqref{EqSumContraction} have operator norms tending to 0 on $(K^{1/p^m}/K)u^n$ for each $m\geq 1$, and so the operator $\operatorname{Id}+(g_1-1)^{-1}\gamma$ has bounded inverse on $(K_p/K)u^n$, which extends to $(\widehat{K}_p/K)u^n$.  Therefore $(g_1-1)\operatorname{Id}+\gamma$ is invertible on the complete ultrametric space $(\widehat{K}_p/K)u^n$.  This is exactly the input needed to prove the acyclicity of $(\widehat{K}_p/K)u^n$, as in the proof of \cref{ThmAcyclicityOne}.
\end{proof}

Let us call an element $g\in \operatorname{Nott}'(k)$ ``good'' if it satisfies the hypothesis imposed on $g_1$ in \cref{prop:acyclicity_one}:  $(g-1)$ admits an inverse on $(K/K^p)u^{pn}$ for all $n\in\Z$, whose operator valuation is bounded below independently of $n$.  
Certainly there are many ``bad'' elements, such as $g(t)=t+t^p$, for then $(g-1)$ annihilates the image of $t^{1/p}$ in $\widehat{K}_p/K$.  As a further example, the Nottingham group contains many $p$-torsion elements, for instance the geometric series $g(t)=t(1-t)^{-1}$ or more generally $g(t)=t(1-\alpha t^Q)^{-1/Q}$ for $Q$ an integer prime to $p$ and $\alpha\in k^\times$.  For such elements, $(g-1)^p=g^p-1=0$ in the characteristic $p$ group algebra, so $(g-1)$ cannot admit an inverse. 

We highlight a simple condition guaranteeing ``goodness'', requiring little accuracy.  

\begin{proposition}
\label{PropGoodGuy}
Let $g$ be an element of $\operatorname{Nott}'(k)$, subject to the following assumptions.

In the case $p\neq 2$:  Let $Q\equiv 1\pmod{p}$ be a positive integer, and assume $g$ takes the following form:
\begin{eqnarray*}
 g(t)&=&t(1+\alpha t^Q + \beta t^{2Q}+\cdots)\\
 g(u)&=&(1+\gamma t^r+\cdots)u
 \end{eqnarray*}
where $\alpha,\beta,\gamma\in k$, and where $\cdots$ indicates terms of higher valuation.  Assume:
\begin{enumerate}
    \item $\alpha\neq 0$,
    \item $\beta\neq \alpha^2$,
    \item $r>2Q/p$.
\end{enumerate}

In the case $p=2$:  let $0<P<2Q<2P$, with $P$ odd, and assume $g$ takes the following form:
\begin{eqnarray*}
    g(t)&=&t(1+\alpha t^P+\beta t^{2Q}+\cdots) \\
    g(u)&=&(1+\gamma t^r+\cdots)u
\end{eqnarray*}
where $\alpha,\beta,\gamma\in k$.  
Assume:
\begin{enumerate}
    \item $\beta\neq 0$,
    \item $r>Q$.
\end{enumerate}
Then for all $n\in\Z$, the operator $(g-1)$ is invertible on $(K/K^p)u^{pn}$, and $\val_{(K/K^p)u^{pn}} (g-1)^{-1} =-2Q$.  In particular, $g$ is ``good''.
\end{proposition}

\begin{proof}  Start with the case that $p$ is odd.
Consider the action of $(g-1)$ on a monomial $t^\ell$ for $p\nmid \ell$.  We have
\[ 
    g(t^\ell) = t^\ell(1+\alpha t^Q+\dots)^\ell = t^\ell + \ell\alpha t^{\ell+Q}+\dots,
\] 
As long as $\ell\not\equiv -1\pmod{p}$, the dominant term of $(g-1)t^\ell$ is $\ell\alpha t^{\ell+Q}$.  On the other hand if $\ell=-1+pk$, then:
\begin{eqnarray*}
g(t^\ell)&=&t^\ell(1+\alpha t^Q+\beta t^{2Q}+\dots)^{-1}(1+\alpha^p t^{pQ}+\beta^p t^{2pQ}+\dots)^k \\ 
&=& t^\ell(1+\alpha t^Q+\beta t^{2Q})^{-1} +O(t^{\ell+2Q+1}) \\
&=& t^\ell(1-\alpha t^Q + (\alpha^2-\beta)t^{2Q})+O(t^{\ell+2Q+1}) \\
&\equiv& t^\ell + (\alpha^2-\beta)t^{\ell+2Q} +O(t^{\ell+2Q+1}).
\end{eqnarray*}
In the second equality:  we have used $p>2$ to ignore the factor with the $k$th power.  The final congruence is modulo $K^p$; note that $\ell+Q\equiv 0\pmod{p}$.  In summary, the dominant term of $(g-1)t^\ell$ as a class in $K/K^p$ is 
\[ (g-1)t^\ell=\begin{cases}
    \ell\alpha t^{\ell+Q}+\cdots,& \ell\not\equiv -1\pmod{p} \\
    (\alpha^2-\beta) t^{\ell+2Q}+\cdots,& \ell\equiv -1\pmod{p} 
\end{cases}
\] 
Now let $n\in\Z$, and consider the action of $(g-1)$ on $(K/K^p)u^{pn}$ and in particular on the elements $t^\ell u^{pn}$ for $p\nmid \ell$.  We have
\[ (g-1)(t^\ell u^{pn})=((g-1)t^\ell)u^{pn}+t^\ell(g-1)u^{pn}
+((g-1)t^\ell)((g-1)u^{pn})
. \] 
Since $(g-1)u^{pn}\in t^{pr}\OO_Ku^{pn}$, and $pr>2Q$, the first term in the displayed sum dominates, which is to say
\[ (g-1)t^\ell u^{pn}=\begin{cases}
    \ell \alpha t^{\ell+Q}u^{pn}+\cdots,& \ell\not\equiv -1\pmod{p} \\
    (\alpha^2-\beta) t^{\ell+2Q}u^{pn}+\cdots,& \ell\equiv -1\pmod{p}. 
\end{cases}
\] 

In the case $p=2$:  Under the listed hypothesis, we consider the action of $(g-1)$ on a monomial $t^\ell$ for $\ell=1+2k$ odd.  We have
\begin{eqnarray*}
 g(t^\ell)&=&t^{\ell}(1+\alpha t^P+\beta^{2Q}+\dots)(1+\alpha^2 t^{2P}+\cdots)^k  \\
 &=& t^{\ell}(1+\alpha t^P+\beta^{2Q}+\cdots)\\
 &\equiv& t^{\ell} + \beta t^{\ell+2Q} +\cdots,
\end{eqnarray*}
where the second equality comes from $2Q<2P$, and the final congruence is modulo $K^2$, noting that $\ell+P$ is even.  Again let $n\in\Z$, and consider the action of $(g-1)$ on $(K/K^2)u^{2n}$ and in particular on the element $t^\ell u^{2n}$.  By a similar calculation as in the $p$ odd case, we have
\[ (g-1)t^{\ell} u^{2n} = \beta t^{\ell+2Q}u^{2n} +\cdots. \] 

In either case, the key observation is that for every integer $m$ prime to $p$, there exists a unique $\ell\geq m-2Q$ such that the dominant term of $(g-1)t^\ell u^{pn}$ is a unit times $t^m u^{pn}$.  This suggests the following recursive procedure for determining $(g-1)^{-1}v$ for $v\in (K/K^p)u^{pn}$:  If $\val_{(K/K^p)u^{pn}} v=m$, there exists some $a\in k$ and $\ell\geq m -2Q$ for which $\val_{(K/K^p)u^{pn}} (v-a(g-1)t^\ell u^{pn})> m$; repeat the process with $v-a(g-1)t^\ell u^{pn}$ in the place of $v$.  One obtains $(g-1)^{-1}v=at^{\ell}u^{pn} + \cdots$.  The result is an operator $(g-1)^{-1}$ on $(K/K^p)u^{pn}$ with valuation $-2Q$.  
\end{proof}

\subsection{The representation of the Morava stabilizer group in the Nottingham group}
Assume $h\geq 2$, and recall our conventions regarding the Morava stabilizer group:  $D_h$ is the division algebra over $\Q_p$ with invariant $1/h$, with ring of integers $\OO_{D_h}$ and unit group $\overline{\G}_h=\OO_{D_h}^\times$.  Let $\zeta_h\from \overline{\G}_h\to \Z_p$ be the $\Z_p$-component of the determinant map, and let $\overline{\G}_h^1$ be its kernel. 

Let
\[ G = 1+q\OO_{D_h}\subset \G_h, \]
where as usual $q=p$ if $p$ is odd and $q=4$ if $p=2$.  

\begin{lemma}
$G$ is a uniform pro-$p$ group of dimension $h^2$.  Denoting by $P_i(G)$ its lower $p$-series, we have $P_i(G)=1+qp^{i-1}\OO_{D_h}$ for $i\geq 1$.  

Let $G^1=G\cap \overline{\G}_h^1$.  Then $G^1$ is a uniform pro-$p$ group of dimension $h^2-1$. 
\end{lemma}

\begin{proof} The claim about $G$ appears as \cite[Example 2.15]{BGHS2022}, but can also be deduced from earlier work of Henn \cite{Henn1998}.  

The claim about $G^1$ follows from the following statement:  if $G$ is a uniform pro-$p$ group of dimension $d$, and $\zeta\from G\to\Z_p$ is a surjective homomorphism, then $H\coloneqq \ker \zeta$ is a uniform pro-$p$ group of dimension $d-1$.  Indeed, in this situation we have $[H,H]\subset [G,G]\cap H\subset G^q\cap H\subset H^q$, so that $H$ is powerful.  It is also torsion-free since $G$ is;  therefore it is uniform by \cite[Theorem 4.8]{DDMS}.  Finally, note that $H/H^p$ is the kernel of the surjective map $G/G^p\to \Z/p$ induced by $\zeta$;  the dimension of $H$ is the $\F_p$-dimension of $H/H^p$, which must be $d-1$.
\end{proof}

Let us identify $(\pi_0E_h)/(p,u_1,\dots,u_{h-2})$ with the power series ring $\OO_K=k\powerseries{t}$ via $u_{h-1} \mapsto t$, and $(\pi_2E_h)/(p,u_1,\dots,u_{h-2})$ with the invertible $\OO_K$-module $\omega=\OO_Ku$.  Then $\G_h$ acts continuously on \[ (\pi_*E_h)/(p,u_1,\dots,u_{h-2}) \isom \bigoplus_{n\in\Z} \omega^{\otimes n},\]
and the action of the unipotent subgroup $G\subset \G_h$ induces a homomorphism into the extended Nottingham group: 
\[ G\to \operatorname{Nott}'(k). \] 

\begin{thm} \label{ThmActionEstimate} For all $g\in 1+p^i\OO_D$ and $n\in \Z$ we have:
\[ \val_{K\omega^{\otimes n}}(g-1)\geq p^{(h-1)i}.\] 
\end{thm}
\begin{proof}
Let us write $\fg_0$ for the formal group of height $h$ over $k$, and $\fg$ for its deformation to $\OO_K$.  The defining quality of $\OO_K$ as a quotient of the Lubin--Tate deformation ring is that the $p$-series $[p]_{\fg}$ factors through $F^{h-1}$, where $F$ is the relative Frobenius.  Similarly, the power $p^i$ factors through $F^{(h-1)i}$.  Observe that the $p^{(h-1)i}$ power Frobenius map on $\OO_K/t^{p^{(h-1)i}}$ factors through $\OO_K/t\isom k$, so that $F^{(h-1)i}$ factors through $\fg_0\otimes_k \OO_K/t^{p(h-1)i}$.  Therefore if $A\in \OO_D$ is an endomorphism of $\fg_0$, then $p^iA$ lifts to an endomorphism of $\fg\otimes_{\OO_K} \OO_K/t^{p^{(h-1)i}}$ which factors through $\fg_0\otimes_k \OO_K/t^{p^{(h-1)i}}$.  Consequently $g\coloneqq 1+p^iA$ lifts to an automorphism $\gamma$ of $\fg\otimes \OO_K/t^{p^{(h-1)i}}$, which means exactly that $g$ acts trivially on $\OO_K/t^{p^{(h-1)i}}$.  From this it follows that $\val_K(g-1)\geq p^{(h-1)i}$.   

The derivative of $\gamma$ is an endomorphism of $\Lie \fg\otimes \OO_K/t^{p^{(h-1)i}}$ which factors through the derivative of $g$ on $\fg_0\otimes \OO_K/t^{p^{(h-1)i}}$, which is the identity;  this gives the desired result on differentials. 
\end{proof}

\begin{rmk} This argument owes a debt to \cite{GrossCanonical}.
\end{rmk}

\begin{rmk}
\label{rmk:goodness}
Since $(g-1)^p=g^p-1$ in the group algebra over $k$, we have
\[ \val_K (g^p-1)\geq p\val_K(g-1) \] 
for all $g\in \operatorname{Nott}(k)$.  If $g$ is a generic element, one might expect that as $n\to \infty$, the quantity $\val_K(g^{p^n}-1)$ is $O(p^n)$.  
For elements $g$ of the uniform pro-$p$ group $G$, we have $g^{p^i}\in 1+qp^i\OO_D$, and then \cref{ThmActionEstimate} shows that 
$\val_K(g^{p^i}-1)\geq p^{(h-1)i}$.  Thus in height $h\geq 3$, elements of $\operatorname{Nott}(k)$ arising from the Morava stabilizer group are quite exceptional.

Suppose $g\in \operatorname{Nott}(k)$ is an element for which an inverse operator $(g-1)^{-1}$ exists on $V=K/K^p$ and has operator norm bounded below by $-B$.  Then for all nonzero $v\in V$ we have $\val_V(g-1)v-\val_V v\leq B$, which is to say, $\val_V(g-1)\leq B$.  Inductively we see $\val_V(g^{p^n}-1)=\val_V(g-1)^{p^n} \leq p^nB$.  If $g$ arose from the action of an element of the Morava stabilizer group at height $h$, then by \cref{ThmActionEstimate} we have $\val_K(g^{p^n}-1)\geq p^{(h-1)n}$, which implies $\val_V(g^{p^n}-1)\geq p^{(h-1)n}$.  We conclude that if $h\geq 3$, then no element of the Morava stabilizer group is ``good''. 
\end{rmk}

\subsection{The case of height 2 and the proof of \cref{thmx:height2}}
In the case of height $h=2$, Salch \cite{Salch} has discovered formulas for the action of the Morava stabilizer group $\G_2$ on $\pi_*E_2$.  For our purposes we only need an approximate formula for the action of one particular element of $\G_2$, which lies in the subgroup $W(\F_{p^2})^\times$ of unit Witt vectors.

Throughout we use the ``Hazewinkel coordinates'' $\pi_*E_2\isom W(k)\powerseries{u_1}[u^{-\pm 1}]$;  here $u\in \pi_{-2}E_2$ is a generator for that module.   
We let $t$ denote the image of $u_1$ in $\pi_0E_2/p$, so that $\pi_0E_2/p\isom k\powerseries{t}\coloneqq \OO_K$, and $\pi_*E_2/p\isom k\powerseries{t}[u^{\pm 1}]$. 
 
\begin{proposition}[{\cite[Corollary 6.4 and Corollary 7.3]{Salch}, \cite[Theorem 6.7]{Salch} for $p = 3$, personal communication from Salch for $p=2$}] 
\label{PropSalch} Suppose $p$ is odd. Let 
$g\in W(\F_{p^2})^\times\subset \G_2$ be an element of the form $g\equiv 1+p[a] \pmod{p^2}$, where $a\in \F_{p^2}$.  Let $\alpha=a^p-a$.  Then 
\begin{eqnarray*} g(t)&=&t(1+\alpha t^{p+1}+(\alpha^2+\alpha)t^{2(p+1)} +\cdots) \\
g(u)&=&(1-\alpha t^{p+1}+\cdots)u. 
\end{eqnarray*}

Now suppose $p=2$.  Let $g\in W(\F_4)^\times\subset \G_2$ be the element $1+4[a]$, where $a\in\F_4^\times$ is a generator of that field.  Then
\begin{eqnarray*}
g(t)&=&t(1+t^9+t^{12}+t^{15}+t^{18}+a^2t^{21}+at^{24}+\cdots) \\
g(u)&=&(1+t^9+t^{12}+t^{15}+t^{18}+a^2t^{21}+at^{24}+\cdots)u
\end{eqnarray*}
\end{proposition}

\begin{thm} \label{ThmDeperf2} For all primes $p$ and all $n\in\Z$, the map
\[H^*_{\cts}(\G_2^1,\pi_nE_2\otimes_{\pi_0E_2} K)\to H^*_{\cts}(\G_2^1,\pi_nE_2\otimes_{\pi_0E_2} \widehat{K}_p) \] 
is an isomorphism.  
\end{thm}

\begin{proof} As observed before, it suffices to show that for all $n\in\Z$, the module $\pi_{2n} E_2\otimes_{\pi_0E_2} (\widehat{K}_p/K)=(\widehat{K}_p/K)u^n$ is acyclic for the action of the closed normal uniform pro-$p$ subgroup $G^1\subset \G_2^1$.

First assume $p$ is odd. Let $a\in\F_{p^2}$ be any nonzero element of trace 0, so that     
$\alpha=a^p-a=a^p+a-2a=-2a\neq 0$.  Let $g_1\in G^1$ be any element with $g\equiv 1+p[a] \pmod{p^2\OO_{D_2}}$.  Using \cref{PropSalch}, we find that $(g_1-1)$ satisfies the hypotheses of \cref{PropGoodGuy}, with $Q=p+1$ and $r=p+1$.  

Now assume $p=2$.  Let $a\in\F_4$ be a generator.  Then $\OO_{D_2}$ is generated over $W(\F_4)$ by an element $\pi$ with $\pi^2=2$ satisfying $\pi[a]=[a^2]\pi$.  The uniform pro-2 group $G^1$ has generators $g_1,g_2,g_3$, with 
\begin{eqnarray*}
    g_1&=&(1+4[a])^2/13\\
    g_2&=&-(1+2\pi)^2/7 \\
    g_3&=&-(1+2[a]\pi)^2/7.
\end{eqnarray*}
(It is easy to check these lie in $G^1$.  To see that they generate $G^1$, it is enough to show that they generate $G^1/(G^1)^2$.  
The group $G^1/(G^1)^2$ is identified as an $\F_2$-vector space via $g\mapsto (g-1)/4$ with traceless elements of $\OO/2$, and now one easily checks that the $g_i$ map to a basis of this vector space.)

 Using \cref{PropSalch}, and noting that scalars in $\Z_2^\times$ act trivially on $t$ and $u$, we find 
 \begin{eqnarray*}
g_1(t)&=&t(1+t^{21}+t^{24}+t^{27}+\cdots)\\g_1(u)&=&(1+t^{21}+\cdots)u.
 \end{eqnarray*}
 We see that $(g_1-1)$ satisfies the hypotheses of \cref{PropGoodGuy}, with $P=21$, $Q=12$ and $r=21$.  

In both cases, we apply \cref{PropGoodGuy} to find that the operator $(g_1-1)$ admits an inverse on $(K/K^p)u^{pn}$ with operator valuation bounded below independently of $n$.
The $G^1$-acyclicity of $(\widehat{K_p}/K)u^n$ now follows from \cref{prop:acyclicity_one}.
\end{proof}

\begin{cor}
    For all primes $p$ and $h=2$, the map $(\B^{\solid})^{h ({\G}_h^1\times_{\Gal(k/\F_p)} \G_{h-1})}\to (\widehat{\B}^{\solid})^{h({\G}_h^1\times_{\Gal(k/\F_p)} \G_{h-1})}$ is an equivalence.  If $p$ is odd, then the strong splitting conjecture is true at height 2:  the unit map $\SS_{K(1)}\to L_{K(1)}\SS_{K(2),\infty}$ is an equivalence.  
\end{cor}

\begin{proof} The first claim follows from \cref{PropDeperfection} and \cref{ThmDeperf2}.  The second claim follows from \cref{lem:TFAEsCSC}. 
\end{proof}

We conclude with some comments about the strong splitting conjecture in height 3 for primes $p\geq 5$.  
Approximate formulas for the action of the height 3 Morava stabilizer group $\G_3$ on $\pi_0E_3$ have been obtained by Davis \cite{DavisAction}.  It would be interesting to see if these formulas are accurate enough to prove the strong splitting conjecture.  The same strategy which worked in height 2 cannot work in height 3, because as we have noted, no elements of $\G_3$ are ``good''.  We only remark that computer calculations using Davis' formulas suggest that the quotient $K/K^p$ admits no fixed vectors under the action of $\G_3$.  

After much of this article was written, we asked an AI model to investigate the hypothesis in \cref{PropDeperfection} for heights $h\geq 3$.  To our surprise, the model claims to have disproven the hypothesis for primes $p\leq h-1$. We therefore have reason to believe that the statement of \cref{ConjStrongSplitting} is false for such primes, even if $p-1$ does not divide $h-1$.   

\newpage
\appendix

\section{Actions of profinite groups}\label{app:actions}

Throughout this article, actions of profinite groups play a central role.
It is therefore important to understand precisely what is meant by an action of a profinite group.
On the one hand, when working with condensed mathematics, there is a native notion of action of a profinite group, viewed as a condensed group.
On the other hand, there is the notion of a smooth profinite group action, and which makes sense without any condensed structure.
Since the interplay between these notions plays an important role in our work, we will carefully compare and contrast them here.
Our main results are:
\begin{enumerate}
    \item A proof that, for discrete condensed objects, smooth actions are the same as condensed actions (\cref{prop:smoothdiscrete}).
    \item A comparison between different notions of smooth action that appear in the literature (\cref{prop:smooth-same}).
    \item A proof that, in the $p$-adically complete setting and with respect to a compact $p$-adic Lie group, solid objects with a smooth action embed fully faithfully into solid objects with a condensed action (\cref{prop:compactpadicsmooth}).
\end{enumerate}
This appendix was significantly influenced by \cite{solidRep} and \cite{solidRep2}.

\begin{conv}
    In this section, we let $\CC$ denote a stable presentably symmetric monoidal category. For example, $\CC$ could be the category of spectra, or the category of $K(n)$-local spectra.
\end{conv}

\begin{defn}
    We define the category of condensed objects in $\CC$ to be
    \[\CC^{\cond} \coloneqq \Sp^{\cond} \otimes \CC.\]
    Equivalently, this is the category of accessible $\CC$-valued hypersheaves on $\Prof$.
\end{defn}

There is a fully faithful discrete object functor
\[(-)^\delta \colon \CC \hookrightarrow \CC^{\cond}\]
which is left adjoint to the forgetful functor evaluating at the point
\[(-)(\ast) \colon \CC^{\cond} \to \CC.\]
On the level of sheaves, $(-)^\delta$ is described as follows: given $X \in \CC$ and $S = \varprojlim_{\alpha} S_\alpha \in \Prof$, we have
$X^{\delta}(S) = \varinjlim_{\alpha} X(S_{\alpha})$.

\subsection{Condensed actions}

\begin{conv}
    From now on, we let $\DD$ denote a commutative $\Sp^{\cond}$-algebra in $\Pr^L$: in other words, we let $\iota \colon \Sp^{\cond} \to \DD$ denote a symmetric monoidal left adjoint.
    For example, $\DD$ could be $\CC^{\cond}$ or $\Sp^{\solid}$.
\end{conv}

\begin{rmk}
    Strictly speaking, $\Sp^\cond$ is not a presentable category, but instead is a big presentable category in the sense of \cite[Definition 1.2.1]{brink}. We choose to simply ignore the subtleties arising from this. The reader who is uncomfortable with this is invited to replace all instances of condensed spectra in this work with light condensed spectra, and note that no meaningful change is effected.
\end{rmk}

\begin{construction}
    There are symmetric monoidal functors
    \[(\Prof, \times) \to (\Sp^{\cond}, \otimes) \to (\DD, \otimes).\]
    Since a profinite group $G$ is a cocommutative bialgebra in $(\Prof, \times)$, we obtain a cocommutative bialgebra $\mathbf{1}[G]$ in $(\DD, \otimes).$
\end{construction}

By \cite[Variant 2.2.3]{Arpon}, if $A \in \DD$ is a cocommutative bialgebra, then $\LMod_{A} (\DD)$ and the forgetful functor $\LMod_{A} (\DD) \to \DD$ are naturally equipped with symmetric monoidal structures.

\begin{defn}
    We define the category of objects in $\DD$ with condensed $G$-action to be $\DD^{BG} = \LMod_{\mathbf{1}[G]} (\DD)$. By the above, this is a symmetric monoidal category equipped with a symmetric monoidal forgetful functor
    \[\DD^{BG} \to \DD.\]
    This determines a functor $\Grp(\Prof)^{op} \to \CAlg(\Pr^L)_{/\DD}$.
\end{defn}

\begin{rmk} \label{rmk:G-finite-agree}
    When $G$ is a finite group, it follows from \cite[Example 2.2.9]{Arpon} that $\DD^{BG}$ is the usual category of objects with $G$-action.
\end{rmk}

\begin{rmk}
    Let $\DD \to \DD'$ denote a morphism in $\CAlg(\Pr^L)_{\Sp^{\cond}/}$. Then the natural functor
    \[\DD' \otimes_{\DD} \DD^{BG} \to (\DD')^{BG}\]
    is an equivalence by \cite[Theorem 4.8.4.6]{HA}.
    In particular, we have
    \[\DD \otimes_{\Sp^{\cond}} (\Sp^{\cond})^{BG} \simeq \DD^{BG}.\]
\end{rmk}

For the purpose of comparison with smooth actions, it will be convenient to describe $\DD^{BG}$ comonadically. To incorporate its symmetric monoidal structure into this description, we require the following definition.

\begin{defn}
    Let $\Fun^{\lax}(\DD,\DD)$ denote the category of lax symmetric monoidal functors from $\DD$ to $\DD$. Composition endows $\Fun^{\lax}(\DD,\DD)$ with the structure of of a monoidal category. We define the category $\coMon^{\lax} (\DD)$ of lax symmetric monoidal comonads to be the category of coalgebras in $\Fun^{\lax} (\DD,\DD)$ with respect to this monoidal structure.
\end{defn}

\begin{rmk}
    We note that the forgetful functors
    \[\coMon^{\lax} (\DD) \to \Fun^{\lax}(\DD,\DD) \to \Fun(\DD,\DD)\]
    preserve filtered colimits.
\end{rmk}

\begin{thm}
    Suppose that $C$ is a lax symmetric monoidal comonad on a category $\mathcal{E}$. Then the category $\mathrm{coAlg}_{C} (\mathcal{E})$ is canonically equipped with a symmetric monoidal structure, and the forgetful functor $\mathrm{coAlg}_{C} (\mathcal{E}) \to \mathcal{E}$ with the structure of a symmetric monoidal functor.

    Conversely, if $\mathcal{F} \to \mathcal{E}$ is a symmetric monoidal comonadic functor, then the corresponding comonad on $C$ on $\mathcal{E}$ is canonically equipped with the structure of a lax symmetric monoidal comonad, and there is an identification $\mathcal{F} \simeq \mathrm{coAlg}_{C} (\mathcal{E})$ in $\CAlg(\Pr^L)_{/\EE}$.
\end{thm}

\begin{proof}
    The first statement is \cite[Theorem 1.11]{Heine} applied to $\mathcal{E}^{op}$, and the second statement follows from \cite[Theorem 5.2]{spectral-lie}.
\end{proof}

\begin{construction}
    Let $G$ denote a profinite group.
    Since the forgetful functor
    \[\DD^{BG} \to \DD\]
    preserves limits and colimits and is conservative, it is comonadic.
    Since it is symmetric monoidal the comonad is lax symmetric monoidal.
    We let $C(G,-)$ denote the corresponding lax symmetric monoidal comonad.

    This determines a functor
    \[C \colon \Grp(\Prof)^{op} \to \coMon^{\lax} (\DD).\]
\end{construction}

\begin{rmk} \label{rmk:colim-preserve}
    When $\DD = \CC^{\cond}$ or $\Sp^{\solid}$, then the underlying functor of $C(G,-)$ is given on the level of sheaves on $\Prof$ by the formula $C(G,X) (S) = X(G \times S)$.
    In the case $\DD = \Sp^{\solid}$, this implies that the functor commutes with $C(G,-)$ with filtered colimits, hence all colimits.
    Indeed, the above formula reduces this to statement that $X \mapsto X(S)$ preserves filtered colimits for a fixed profinite set $S$, which in turn reduces to the statement that $\SS[S]^{\solid}$ is compact as a solid spectrum.
    This holds because $\SS[S]^{\solid}$ is isomorphic to a product of copies of the unit for any $S$: in particular, it is a summand of $\SS[T]^{\solid}$ for an extremally disconnected profinite set $T$.
    For such a $T$, the condensed spectrum $\SS[T]$ is compact, and solidification preserves compactness because its right adjoint commutes with small colimits.
\end{rmk}

\begin{rmk}
    When $\DD = \Sp^{\solid}$, this
\end{rmk}

\subsection{Smooth actions}

We now describe the category of objects with a smooth $G$-action.

\begin{construction}
    The assignment taking a finite group $G$ to the comonadic symmetric monoidal forgetful functor
    \[\CC^{BG} \to \CC\]
    determines a functor
    \[C \colon \Grp(\Fin)^{op} \to \coMon^{\lax} (\CC).\]
    In fact, since the functor $C(G,-)$ preserves colimits, the functor $C$ factors through the full subcategory $\coMon^{L, \lax} (\CC)$ consisting of colimit-preserving lax symmetric monoidal comonads.
    This category is presentable, hence admits all colimits.
    By left Kan extending along
    \[\Grp(\Fin)^{op} \hookrightarrow \Ind(\Grp(\Fin)^{op}) \simeq \Grp(\Prof)^{op},\]
    we obtain a functor
    \[C^{\sm} \colon \Grp(\Prof)^{op} \to \coMon^{L, \lax} (\CC).\]
    For $G = \varprojlim_\alpha G / U_{\alpha} \in \Grp(\Prof)$ and $X \in \C$, this is given by the formula
    \[C^{\sm} (G,X) \simeq \varinjlim_\alpha C(G/U_{\alpha}, X).\]
\end{construction}

\begin{defn}
    Given a profinite group $G$, we define the symmetric monoidal category of objects with smooth $G$-action in $\CC$ to be
    \[\CC^{BG^{\sm}} \coloneqq \Coalg_{C^{\sm} (G,-)} (\CC).\]
    This comes with a symmetric monoidal forgetful functor $\CC^{BG^{\sm}} \to \CC.$
\end{defn}

\begin{rmk}
    Let $G$ denote a finite group. Since $\mathbf{1}[G]$ is dualizable, the endofunctor $C(G,-)$ is naturally identified with the functor $C(G, \mathbf{1}) \otimes -$, and $C(G, \mathbf{1})$ is the dual of $\mathbf{1}[G]$. This induces a commutative bialgebra structure on $C(G,\mathbf{1})$ and identifies $\CC^{BG}$ with comodules over $C(G,\mathbf{1})$.

    Taking filtered colimits, we are able to extend this to profinite $G = \varprojlim_\alpha G/U_{\alpha}$: the symmetric lax-monoidal comonad $C^{\sm} (G,-)$ is isomorphic to that associated to the commutative bialgebra $C^{\sm} (G, \mathbf{1}) \simeq \varinjlim_{\alpha} C(G/U_{\alpha},\mathbf{1}),$, and $C^{BG^{\sm}}$ is the category of comodules over this bialgebra.
\end{rmk}

\begin{rmk}
    More generally, if one wants to mix smooth and condensed actions, one may combine the above two constructions in the following way.
    Suppose that $C(G,-)$ preserves colimits in $\DD$ for all profinite groups $G$: for example, we could have $\DD = \Sp^{\solid}$ as in \Cref{rmk:colim-preserve}.
    Starting with $C \colon \Grp(\Prof)^{\op} \to \coMon^{L, \lax}(\DD)$, one may left Kan extend to obtain a functor
    $C^{\mix} \colon (\Pro(\Grp(\Prof)))^{\op} \to \coMon^{L, \lax}(\DD)$.
    Thus, given a pro-profinite group $G$, we may consider the symmetric monoidal category of coalgebras over the lax symmetric monoidal comonad $C^{\mix} (G, -)$.
    The ``outer'' pro-direction acts smoothly, while the ``inner'' pro-direction acts via a condensed action.

    For example, if we have two profinite groups $G = \varprojlim_\alpha G/U_{\alpha}$ and $H$, an element of $((\DD)^{BH})^{BG^{\sm}}$ is the same thing as a coalgebra for the comonad associated to the pro-profinite group $\varprojlim_{\alpha} (G/U_{\alpha} \times H)$.
\end{rmk}

To conclude, we prove a base-change property for categories of objects with smooth $G$-action. The key step is the following lemma, which is due to an AI system.
In the following lemma, we relax our conventions on $\CC$ and $\DD$.

\begin{lemma} \label{lem:chatgpt}
    Let $\CC$ denote a stable presentable category, and let $C \in \coMon^L (\CC)$ denote a colimit-preserving comonad on $\CC$.
    Given a second presentable stable category $\mathcal{D}$, $C \otimes \id_{\mathcal{D}} \in \coMon^L (\CC \otimes \mathcal{D})$, and there is a canonical comparison functor
    \[\Coalg_{C} (\CC) \otimes \mathcal{D} \to \Coalg_{C \otimes \id_{\mathcal{D}}} (\CC \otimes \mathcal{D})\]
    If $\mathcal{D}$ is dualizable as a presentable stable category (for example, if $\DD$ is compactly generated), then this functor is an equivalence of categories.

    If $\CC$ and $\DD$ are symmetric monoidal and the comonad $C$ is lax symmetric monoidal, then so is $C \otimes \id_{\DD}$, and the comparison functor naturally upgrades to a symmetric monoidal functor.
\end{lemma}

\begin{proof}
    The left adjoint
    \[\Coalg_{C} (\CC) \otimes \mathcal{D} \to \CC \otimes \mathcal{D}\]
    has associated comonad $C \otimes \id_{\mathcal{D}}$, so it suffices to prove that it's comonadic. If $\CC$ and $\DD$ are symmetric monoidal and $C$ is lax symmetric monoidal, then this functor is symmetric monoidal and its associated comonad $C \otimes \id_{\mathcal{D}}$ lifts to a lax symmetric monoidal comonad.

    To do so, we use dualizability to identify it with
    \[G\colon \Fun^{L} (\mathcal{D}^\vee, \Coalg_{C} (\CC)) \to \Fun^{L} (\mathcal{D}^\vee, \CC).\]
    We will verify comonadicity using the co-Barr--Beck--Lurie theorem.
    This requires us to prove that the functor is conservative, preserves colimits, and preserves $G$-split totalizations.
    It is clear from the outset that it preserves colimits.
    It is conservative because
    \[ H\colon \Coalg_{C} (\CC) \to \CC\]
    is conservative and isomorphisms are checked levelwise in functor categories.

    Let $X^\bullet \in \Fun(\Delta, \Fun^{L} (\mathcal{D}^\vee, \Coalg_{C} (\CC)))$ be $G$-split.
    Since $G(X^{\bullet})$ is split, its totalization is preserved by all functors.
    In particular, it is preserved by evaluation at any object of $\mathcal{D}^\vee$.
    It follows that the totalization of $G(X^{\bullet})$ is equivalent to its levelwise totalization.
    Now, since $X^{\bullet}$ is $G$-split, it levelwise $H$-split.
    Moreover, $G$ preserves levelwise $H$-split totalizations.
    It therefore suffices to show that the totalization of $X^{\bullet}$ is computed levelwise, i.e., that its levelwise totalization is a colimit-preserving functor.
    This follows from the fact that $H \colon \Coalg_{C} (\CC) \to \CC$ creates colimits and the totalization of $G(X^{\bullet})$ is computed levelwise.
\end{proof}

\begin{cor}
    Suppose that $\CC$ is a dualizable stable presentably symmetric monoidal $\infty$-category.
    Then the natural comparison functor
    \[\Sp^{BG^{\sm}} \otimes_{\Sp} \CC \to \CC^{BG^{\sm}}\]
    is an equivalence of symmetric monoidal categories.
\end{cor}

\subsection{Comparison between smooth and condensed actions I}

\begin{construction}
It follows from \Cref{rmk:G-finite-agree} that the two functors
\[C, C^{\sm} \colon \Grp(\Prof)^{op} \to \coMon^{\lax} (\DD)\]
agree after restriction to $\Grp(\Fin)^{op}$.
It therefore follows from the universal property of left Kan extension that there is a natural transformation of functors $C^{\sm} \to C$. This induces a diagram of symmetric monoidal left adjoints
\begin{center}
    \begin{tikzcd}
        & \DD & \\
        (\DD)^{BG^\sm} \ar[rr] \ar[ur, "\for"] & & (\DD)^{BG} \ar[ul,"\for",swap] \\
        & \DD \ar[lu, "\triv^{\sm}"] \ar[ru, "\triv", swap] & 
    \end{tikzcd}
\end{center}
the bottom part of which produces a diagram of lax symmetric monoidal right adjoints
\begin{center}
    \begin{tikzcd}
        \DD^{BG^\sm} \ar[rd, "(-)^{hG^{\sm}}", swap] & & \DD^{BG} \ar[ll, "(-)^{G-\sm}", swap] \ar[ld, "(-)^{hG}"] \\
        & \DD.  &
    \end{tikzcd}
\end{center}
The functor $(-)^{G_\sm}$ is the functor of (derived) $G$-smooth vectors.
\end{construction}

\begin{rmk}
    The forgetful functor
    \[\DD^{BG^\sm} \to \DD^{BG}\]
    is \emph{not} fully faithful in general, and the natural transformation of functors on $\DD^{BG^\sm}$
    \[(-)^{hG^\sm} \to (-)^{hG}\]
    often fails to be an isomorphism. However, see \Cref{prop:smoothdiscrete} and \Cref{prop:compactpadicsmooth}.
\end{rmk}

We are now ready to state our first result.

\begin{proposition}\label{prop:smoothdiscrete}
    The lax symmetric monoidal comonads $C^{\sm} (G,-)$  and $C(G,-)$ on $\CC^{\cond}$ restrict to the essential image of the embedding $(-)^{\delta} \colon \CC \hookrightarrow \CC^{\cond}$,
    where they are both naturally isomorphic to the lax symmetric monoidal comonad $C^{\sm} (G,-)$ on $\CC$.

    In particular, we find that:
    \begin{enumerate}
        \item If $X \in (\CC^{\cond})^{BG}$ and the underlying object of $X$ is discrete, then $X$ canonically lifts to an object of $(\CC^{\cond})^{BG^{\sm}}$.
        \item If $X \in (\CC^{\cond})^{BG^{\sm}}$ and the underlying object of $X$ is discrete, then the unit map $X \to X^{G-\sm}$ is an isomorphism, and hence the comparison map
        $X^{hG^{\sm}} \to X^{hG}$ is an isomorphism.
    \end{enumerate}
\end{proposition}

\begin{proof}
    Since $(-)^\delta$ is a symmetric monoidal left adjoint, it is elementary to check that $C^{\sm} (G,-)$ on $\CC^{\cond}$ restricts along $(-)^{\delta}$ to $C^{\sm}(G,-)$ on $\CC$.
    Given $X \in \CC$, it remains to check that  the natural map
    \[C^{\sm}(G,X^{\delta}) \to C(G,X^\delta)\]
    is an equivalence, which is immediate from evaluating at any given profinite set $S$ and expanding out the definitions.

    Consequence (1) is then immediate, and (2) follows from writing out for $X \in \CC^{BG^{\sm}}$
    \[X^{\delta} \simeq \Tot C^{\sm}(G^{\bullet+1}, X^{\delta}) \to \Tot C^{\sm}(G,C(G^{\bullet},X^{\delta})) \simeq \Tot C(G^{\bullet+1}, X^{\delta})^{G-\sm} \simeq (X^{\delta})^{G-\sm}\]
    and using that $C^{\sm}(G^{\bullet},X^{\delta}) \to C(G^{\bullet},X^\delta)$ is an equivalence.
\end{proof}

Before we proceed further, we prove two results giving conditions for when condensed and smooth homotopy fixed points commute with filtered colimits.

\begin{lemma} \label{lem:smooth-fp-filt-colims}
    Suppose $\CC$ is equipped with a $t$-structure that is compatible with filtered colimits and for which $\CC_{\leq 0}$ is compactly generated by cotruncated objects.
    Then $C^{\sm} (G,-)$ restricts to a comonad on $\CC_{\leq 0}$, and the functor
    \[(-)^{hG^{\sm}}\colon (\CC_{\leq 0})^{BG^{\sm}} \to \CC_{\leq 0}\]
    commutes with filtered colimits.
\end{lemma}

\begin{proof}
    Given a profinite set $S \simeq \varprojlim_i S_i$, the functor
    \[C^{\sm} (S, -) \simeq \varinjlim_i C(S_i, -)\]
    preserves $\CC_{\leq 0}$ because the $t$-structure is assumed compatible with filtered colimits.
    Given $X \in (\CC_{\leq 0})^{BG^{\sm}}$, write
    \[X^{hG^{\sm}} \simeq \Tot C^{\sm} (G^{\bullet}, X).\]
    To check commutation with filtered colimits of coconnective objects, we just need to commute the filtered colimit past the totalization.
    This holds by \cite[Lemma 3.7(1)]{arc}.
\end{proof}

\begin{lemma} \label{lem:fp-filt-colims}
    Suppose that $\DD$ is equipped with a $t$-structure, and that:
    \begin{enumerate}
        \item the $t$-structure is compatible with filtered colimits.
        \item the category $\DD_{\leq 0}$ is compactly generated by cotruncated objects.
        \item for each $G \in \Grp(\Prof)$, the functor $C(G,-)$ on $\DD$ restricts to $\DD_{\leq 0}$ and commutes with filtered colimits in $\DD_{\leq 0}$.
    \end{enumerate}
    Then the functor
    \[(-)^{hG} \colon (\DD_{\leq 0})^{BG} \to \DD_{\leq 0}\]
    commutes with filtered colimits.
\end{lemma}

\begin{proof}
    Given $X \in (\DD_{\leq 0})^{BG}$, write
    \[X^{hG} \simeq \Tot C (G^{\bullet}, X).\]
    The functor $C(G^\bullet, -)$ commutes with filtered colimits in $\DD_{\leq 0}$ by parts (1) and (3).
    The totalization commutes with filtered colimits of coconnective objects by part (2) and \cite[Lemma 3.7(1)]{arc}, so we're done.
\end{proof}

\begin{example}
    The assumptions of the above lemma are satisfied by $\Sp^{\solid}$ and its standard $t$-structure.
\end{example}

\subsection{Models for smooth actions}

In the case that $\CC = \Sp$, we record several equivalent definitions of the category of spectra with smooth $G$-action. In particular, we compare our definition with various categories of sheaves on the site of finite continuous $G$-sets, which is the language that Devinatz--Hopkins used to formulate the action of the Morava stabilizer group on Lubin--Tate theory. \cite{DevinatzHopkins2004}

\begin{defn}
    Let $G$ denote a profinite group. Then we let $\mathcal{T}_G$ denote the category of finite continuous $G$-sets. We view $\mathcal{T}_G$ as a site by declaring a family $\{X_i \to X\}_{i \in I}$ to be a cover if it is jointly surjective.
\end{defn}

We will compare the following symmetric monoidal categories:
\begin{enumerate}
    \item $\Sp^{BG^{\sm}}$, the category of spectra with smooth $G$-action. Equivalently, the category of $C^{\sm} (G, \SS)$-comodules in spectra.
    \item $\Shv^{\hyp} (\mathcal{T}_G; \Sp)$, the category of hypercomplete sheaves of spectra on $\mathcal{T}_G$, with the pointwise symmetric monoidal structure.
    \item $\widehat{\Shv}(\mathcal{T}_G; \Sp)$, the left completion of the category of sheaves of spectra on $\mathcal{T}_G$ with respect to the usual $t$-structure, with the pointwise symmetric monoidal structure.
    \item $D(BG)$, the value of Heyer--Mann's six-functor formalism on the condensed anima $BG$. \cite{HeyerMann6FF}
\end{enumerate}

\begin{rmk}
    Note that the left completion of the category of sheaves of spectra on a topos $\mathcal{T}$ is not necessarily a localization of $\Shv(\mathcal{T}; \Sp)$.
    Indeed, by definition $\widehat{\Shv}(\mathcal{T}; \Sp) = \varprojlim_n \Shv(\mathcal{T}; \Sp_{\leq n})$, and if $\{X_n\}$ is an object of this category then it may very well be the case that the natural map
    \[\tau_{\leq n} \varprojlim_k X_k \to X_n\]
    is not an equivalence.

    Moreover, since a sheaf of $n$-truncated spectra is automatically a hypersheaf, the left completions of $\Shv(\mathcal{T}; \Sp)$ and $\Shv(\mathcal{T}^{\hyp}; \Sp)$ agree.
\end{rmk}

Three of these categories agree for all profinite groups $G$. We start with a lemma.

\begin{proposition} \label{prop:smooth-same}
    For any profinite group $G$, there are natural symmetric monoidal equivalences of categories:
    \[\Sp^{BG^{\sm}} \simeq \widehat{\Shv}(\mathcal{T}_G; \Sp) \simeq D(BG).\]
\end{proposition}

\begin{proof}
    We start by proving the equivalence
    \[\Sp^{BG^{\sm}} \simeq \widehat{\Shv}(\mathcal{T}_G; \Sp).\]
    Let $G = \varprojlim_{\alpha} G / U_{\alpha}$.
    First, we note that by \cite[Construction 4.5]{hyperdescent} there is an equivalence
    \[\Shv(\mathcal{T}_G; \Sp) \simeq \varinjlim_{\alpha} \Sp^{BG/U_\alpha}\]
    in $\CAlg(\Pr^L)$, where the colimit is taken along the inflation functors.
    Inflation along $G \to G/U_\alpha$ therefore induces a symmetric left adjoint
    \[\Shv(\mathcal{T}_G; \Sp) \simeq \varinjlim_{\alpha} \Sp^{BG/U_\alpha} \to \Sp^{BG^{\sm}}.\]
    Concretely, this left adjoint is computed by $\{X_{\alpha}\} \mapsto \varinjlim_{\alpha} X_{\alpha}$, and its right adjoint is computed by $X \mapsto \{X^{hU_{\alpha}^{\sm}}\}$.
    These functors are right $t$-exact, hence induce an adjoint pair on categories of $n$-truncated objects
    \[\Shv(\mathcal{T}_G; \tau_{\leq n} \Sp) \leftrightarrows (\tau_{\leq n} \Sp)^{BG^{\sm}}.\]
    By taking the limit over $n$, the desired result reduces to proving that this is an adjoint equivalence.
    We start by proving that the counit map
    \[\varinjlim_{\alpha} X^{hU_{\alpha}^{\sm}} \to X\]
    is an equivalence for $X \in (\tau_{\leq n} \Sp)^{BG^{\sm}}.$
    This is equivalent to the colimit-comparison map
    \[\varinjlim_{\alpha} ((X \otimes C^{\sm} (G/U_{\alpha},\SS))^{hG^{\sm}}) \to (X \otimes C^{\sm} (G,\SS))^{hG^{\sm}},\]
    which is an equivalence because $(-)^{hG^{\sm}}$ commutes with filtered colimits of $n$-truncated objects by \Cref{lem:smooth-fp-filt-colims}.
    This implies that
    \[(\tau_{\leq n} \Sp)^{BG^{\sm}} \to \Shv(\mathcal{T}_G; \tau_{\leq n} \Sp)\]
    is fully faithful.
    To conclude, it suffices to show that the functor
    \[\Shv(\mathcal{T}_G; \tau_{\leq n} \Sp) \to (\tau_{\leq n} \Sp)^{BG^{\sm}}\]
    is conservative---equivalently, its composite with $(\tau_{\leq n} \Sp)^{BG^{\sm}} \to \tau_{\leq n} \Sp$ is.
    Since this functor is colimit-preserving, it suffices to prove that the induced functor on hearts
    \[\Shv(\mathcal{T}_G; \mathrm{Ab}) \to \mathrm{Ab}\]
    is conservative. This follows from \cite[Section III.9]{MacLaneMoerdijk}, which identifies $\Shv(\mathcal{T}_G; \mathrm{Ab})$ with the category of discrete abelian groups with continuous $G$-action.
    
    Now we deal with $D(BG)$.
    For this, we write
    \[D(BG) \simeq \Tot(D(G^{\bullet})) \simeq \Tot (\varinjlim D((G/U_\alpha)^\bullet)),\]
    as symmetric monoidal categories.
    On the other hand, we have
    $\Sp^{BG/U_{\alpha}} \simeq \Tot (\Sp^{(G/U_{\alpha})^{\bullet}}) \simeq \Tot (D((G/U_{\alpha})^{\bullet})),$ whence we get a natural symmetric monoidal limit-colimit interchange functor 
    \[\Shv(\mathcal{T}_G; \Sp) \simeq \varinjlim \Sp^{BG/U_{\alpha}} \to D(BG).\]
    Since $D(BG)$ is left complete, and this functor is right $t$-exact, this factors through a symmetric monoidal functor
    \[\Comod_{C^{\sm}(G,\SS)} (\Sp) \simeq (\Sp)^{BG^{\sm}} \simeq \widehat{\Shv}(\mathcal{T}_G; \Sp) \to D(BG).\]
    Using that the composite functor 
    \[\Comod_{C^{\sm}(G,\SS)} (\Sp) \to D(BG) \to \Sp\]
    is the forgetful one, one checks using \cite[Corollary 4.7.5.3]{HA} that
    $\Comod_{C^{\sm}(G,\SS)} (\Sp) \to D(BG)$
    is an equivalence, as desired.
    The key observation is that $D(G^\bullet) = \Mod_{C^{\sm}(G^\bullet,\SS)} (\SS).$
\end{proof}

When $G$ is of virtual finite cohomological dimension, it was proven by Clausen and Mathew that the left-completed category of sheaves agrees with the category of hypercomplete sheaves.

\begin{proposition}[{\cite[Corollary 4.28]{hyperdescent}}] \label{prop:hyp-post-same}
    Suppose that $G$ is a profinite group of virtual finite cohomological dimension.
    Then $\Shv^{\hyp} (\mathcal{T}_G; \Sp)$ is left complete, so that 
    \[\Shv^{\hyp} (\mathcal{T}_G; \Sp) \simeq \widehat{\Shv}(\mathcal{T}_G; \Sp).\]
\end{proposition}

\begin{example}
    Uniform pro-$p$ groups are of finite cohomological dimension, so compact $p$-adic Lie groups are of virtual finite cohomological dimension.
\end{example}

Finally, let us record some useful corollaries of the above.

\begin{cor} \label{cor:hyper-smooth-equiv}
    Suppose that $\CC$ is a dualizable stable presentably symmetric monoidal category, and that $G = \varprojlim_{\alpha} G / U_{\alpha}$ is of virtual finite cohomological dimension.
    Then the natural symmetric monoidal left adjoint
    \[\Shv (\mathcal{T}_G; \CC) \simeq \varinjlim_{\alpha} \CC^{BG/U_\alpha} \to \CC^{BG^{\sm}}\]
    factors through an equivalence
    \[\Shv^{\hyp} (\mathcal{T}_G; \CC) \simeq \CC^{BG^{\sm}}\]
    of presentably symmetric monoidal categories.
\end{cor}

\begin{proof}
  This is a consequence of \Cref{prop:smooth-same} and \Cref{prop:hyp-post-same} when $\CC = \Sp$. The case of general $\CC$ follows from \Cref{lem:chatgpt} and the fact that
  \[\Shv^{\hyp} (\mathcal{T}_G; \Sp) \otimes_{\Sp} \CC \simeq \Shv^{\hyp} (\mathcal{T}_G; \CC).\]
\end{proof}

\begin{example}
    For example, if $G$ is a compact $p$-adic Lie group, this implies that
    \[(\Sp_{K(n)})^{BG^{\sm}} \simeq \Sp^{BG^{\sm}} \otimes \Sp_{K(n)} \simeq \Shv^{\hyp} (\mathcal{T}_G; \Sp_{K(n)})\]
    since $\Sp_{K(n)}$ is compactly generated and hence dualizable.
\end{example}

\begin{cor} \label{cor:smooth-colim-fixed}
    Suppose that $\CC$ is a dualizable stable presentably symmetric monoidal category, and that $G = \varprojlim_{\alpha} G / U_{\alpha}$ is of virtual finite cohomological dimension. Then, for all $X \in \CC^{BG^{\sm}}$,
    \[\varinjlim_{\alpha} X^{hU_{\alpha}^{\sm}} \to X\]
    is an equivalence.
\end{cor}

\begin{proof}
  This map is question is the counit of the left adjoint
    \[\Shv (\mathcal{T}_G; \CC) \simeq \varinjlim_{\alpha} \CC^{BG/U_\alpha} \to \CC^{BG^{\sm}}.\]
  This is a Bousfield localization by \Cref{cor:hyper-smooth-equiv}, hence its counit is an isomorphism.
\end{proof}

\begin{rmk} \label{rmk:G-smooth-colimit}
    Under the assumptions of \Cref{cor:smooth-colim-fixed}, for $X \in \DD^{BG}$ there is a natural equivalence
    \[\varinjlim_{\alpha} X^{hU_{\alpha}} \simeq X^{G-\sm}\]
    by applying \Cref{cor:smooth-colim-fixed} to $X^{G-\sm}$ and using the equivalences $X^{hU_{\alpha}} \simeq (X^{G-\sm})^{hU_{\alpha}^{\sm}}.$
\end{rmk}

\subsection{Comparison between smooth and condensed actions II}

In the case of compact $p$-adic Lie groups and $p$-adically complete solid coefficients, objects with smooth action embed fully faithfully into objects with condensed action.

\begin{defn}
  Let $\Sp^{\solid} _p$ denote the category of $p$-adically complete solid spectra, i.e., the category of solid spectra $X$ such that $X \simeq \varprojlim_i X/p^i$.
\end{defn}

\begin{rmk}
  Note that, under this definition, $\mathbb{Q}_{p}$ with its usual topology is \emph{not} a $p$-adically complete solid spectrum.
\end{rmk}

\begin{proposition}\label{prop:compactpadicsmooth}
    Suppose that $G$ is a compact $p$-adic Lie group. The forgetful functor
    \[(\Sp_p^{\solid})^{BG^\sm} \to (\Sp_{p}^{\solid})^{BG}\]
    is fully faithful.
    In other words, given $X \in (\Sp_p^{\solid})^{BG^\sm}$, the unit map
    \[X \to X^{G-\sm}\]
    is an equivalence.
\end{proposition}

\begin{proof}
    By \Cref{rmk:G-smooth-colimit}, we have
    \[X^{G-\sm} \simeq \varinjlim_{\alpha} X^{hU_{\alpha}}.\]
    Under this identification, the unit map is identified with
    \[X \simeq \varinjlim_{\alpha} X^{hU_{\alpha}^{\sm}} \to \varinjlim_{\alpha} X^{hU_{\alpha}}\]
    is an equivalence.
    This may be checked after replacing $G$ by an open subgroup.
    We may therefore assume that $G$ is uniform pro-$p$ group.
    Cofibering out by $p$ and using $p$-completeness, we may assume that the endomorphism $p^2$ of $X$ is null.
    Now, write $X \simeq \varprojlim_n \tau_{\leq n} X.$
    Since $G$ is of finite cohomological dimension for abelian groups, it is of finite cohomological dimension for solid abelian groups by \cite[Theorem 3.6]{solidProf}.
    We may therefore commute the limit and colimit to obtain the following equivalence:
    \[\varinjlim_{\alpha} X^{hU_{\alpha}} \simeq \varprojlim_n \varinjlim_{\alpha} (\tau_{\leq n} X)^{hU_{\alpha}},\]
    so that we may reduce to the case that $X$ is $n$-truncated.

    Now, writing $X \simeq \varinjlim_{\alpha} X^{hU_{\alpha}^{\sm}}$, and using that $(-)^{hU_{\alpha}}$ commutes with filtered colimits of $n$-truncated solid spectra by \Cref{lem:fp-filt-colims}, we may assume, replacing $G$ by an open subgroup, that the action of $G$ on $X$ is trivial.
    
    Next, writing $X \simeq \varinjlim_m \tau_{\geq m}X$ and again using that $(-)^{hU_{\alpha}}$ commutes with filtered colimits of $n$-truncated solid spectra by \Cref{lem:fp-filt-colims},
    we are able to reduce to the case that $X$ is a solid abelian group in which $p^2=0$.
    Considering the $p$-adic filtration, we reduce to the case that $X$ is a solid $\F_p$-vector space.

    Using the solid version of Lazard's complex \cite[Theorem 5.7]{solidRep}, we find that the natural map
    \[\F_p ^{hU_{\alpha}} \otimes_{\F_p} X \to X^{hU_{\alpha}}\]
    is an isomorphism.
    It therefore suffices to assume that $X = \F_p$. Since $\F_p$ is discrete,
    \[\F_p \to (\F_p)^{G-\sm}\]
    is an isomorphism by \Cref{prop:smoothdiscrete}, as desired.
\end{proof}

The following are immediate corollaries of the above proposition.

\begin{cor} \label{rslt:smooth-F-p-reps-embed-fully-faithfully}
    Let $G$ denote a compact $p$-adic Lie group and let $X \in (\Sp^{\solid} _p)^{BG^{\sm}}$.
    Then the canonical map
    \[X^{hG^{\sm}} \to X^{hG}\]
    is an equivalence.
\end{cor}

\begin{cor}
    Let $G$ denote a compact $p$-adic Lie group and let $X \in (\Sp^{\solid}_p)^{BG}$.
    Then the $G$-action on $X$ is smooth if and only if the canonical map
    \[\varinjlim_{\alpha} X^{hU_{\alpha}} \to X\]
    is an equivalence.
\end{cor}

\begin{rmk}
    \cref{rslt:smooth-F-p-reps-embed-fully-faithfully} was communicated to the second author by Rodríguez Camargo at a conference in 2023.
\end{rmk}

\begin{rmk}
    In the above proposition, it is important that we only work with $p$-adically complete coefficients. For $\mathbb{Q}_{p,\solid}$-modules, it is false.
    Indeed, one has
    $\mathbb{Q}_p \simeq (\mathbb{Q}_p)^{h \Z_p^{\sm}} \not\simeq (\mathbb{Q}_p)^{h\Z_p} \simeq \mathbb{Q}_p \oplus \mathbb{Q}_p[-1].$
\end{rmk}

\bibliographystyle{amsalpha}
\bibliography{bibfile}

@article{BhattMorrowScholze,
  title={Integral p-adic Hodge theory},
  author={Bhatt, Bhargav and Morrow, Matthew and Scholze, Peter},
  journal={Publications math{\'e}matiques de l'IH{\'E}S},
  volume={128},
  number={1},
  pages={219--397},
  year={2018},
  publisher={Springer}
}

@article {CDNIntegral,
    AUTHOR = {Colmez, Pierre and Dospinescu, Gabriel and Nizio{\l} , Wies{\l}awa},
     TITLE = {Integral {$p$}-adic \'{e}tale cohomology of {D}rinfeld symmetric
              spaces},
   JOURNAL = {Duke Math. J.},
  FJOURNAL = {Duke Mathematical Journal},
    VOLUME = {170},
      YEAR = {2021},
    NUMBER = {3},
     PAGES = {575--613},
      ISSN = {0012-7094},
   MRCLASS = {11G18 (14F30)},
  MRNUMBER = {4255044},
MRREVIEWER = {Lance Edward Miller},
       DOI = {10.1215/00127094-2020-0084},
       URL = {https://doi-org.ezproxy.bu.edu/10.1215/00127094-2020-0084},
}

@incollection {BarthelBeaudry,
    AUTHOR = {Barthel, Tobias and Beaudry, Agn\`es},
     TITLE = {Chromatic structures in stable homotopy theory},
 BOOKTITLE = {Handbook of homotopy theory},
    SERIES = {CRC Press/Chapman Hall Handb. Math. Ser.},
     PAGES = {163--220},
 PUBLISHER = {CRC Press, Boca Raton, FL},
      YEAR = {[2020] \copyright 2020},
      ISBN = {978-0-815-36970-7},
   MRCLASS = {55P42 (55N20 55N22 55P60 55Q45 55Q51)},
  MRNUMBER = {4197985},
MRREVIEWER = {Constanze\ Roitzheim},
}

@incollection {Barthel2020,
    AUTHOR = {Barthel, Tobias},
     TITLE = {A short introduction to the telescope and chromatic splitting
              conjectures},
 BOOKTITLE = {Bousfield classes and {O}hkawa's theorem},
    SERIES = {Springer Proc. Math. Stat.},
    VOLUME = {309},
     PAGES = {261--273},
 PUBLISHER = {Springer, Singapore},
      YEAR = {[2020] \copyright 2020},
      ISBN = {978-981-15-1588-0; 978-981-15-1587-3},
   MRCLASS = {55P42},
  MRNUMBER = {4100659},
MRREVIEWER = {Hadi\ Zare},
       DOI = {10.1007/978-981-15-1588-0\{_}9}

@article{beaudry2022chromatic,
  title={Chromatic splitting for the {K}(2)--local sphere at p= 2},
  author={Beaudry, Agnes and Goerss, Paul G and Henn, Hans-Werner},
  journal={Geometry \& Topology},
  volume={26},
  number={1},
  pages={377--476},
  year={2022},
  publisher={Mathematical Sciences Publishers}
}

@misc{FarguesScholze,
      title={Geometrization of the local Langlands correspondence}, 
      author={Laurent Fargues and Peter Scholze},
      year={2024},
      eprint={2102.13459},
      archivePrefix={arXiv},
      primaryClass={math.RT},
      url={https://arxiv.org/abs/2102.13459}, 
}

@ARTICLE{BSSW,
       author = {{Barthel}, Tobias and {Schlank}, Tomer M. and {Stapleton}, Nathaniel and {Weinstein}, Jared},
        title = "{On the rationalization of the $K(n)$-local sphere}",
      journal = {arXiv e-prints, accepted for publication in {Annals of Mathematics}},
         year = 2024,
        month = feb,
          eid = {arXiv:2402.00960},
        pages = {arXiv:2402.00960},
          doi = {10.48550/arXiv.2402.00960},
archivePrefix = {arXiv},
       eprint = {2402.00960},
 primaryClass = {math.AT},
       adsurl = {https://ui.adsabs.harvard.edu/abs/2024arXiv240200960B}
}

@ARTICLE{HopkinsPicard,
       author = {{Barthel}, Tobias and {Schlank}, Tomer M. and {Stapleton}, Nathaniel and {Weinstein}, Jared},
        title = "{On Hopkins' Picard group}",
      journal = {arXiv e-prints, accepted for publication in {Acta Mathematica}},
         year = 2024,
        month = jul,
          eid = {arXiv:2407.20958},
        pages = {arXiv:2407.20958},
          doi = {10.48550/arXiv.2407.20958},
archivePrefix = {arXiv},
       eprint = {2407.20958},
 primaryClass = {math.AT},
       adsurl = {https://ui.adsabs.harvard.edu/abs/2024arXiv240720958B}
}

@article {LauSmoothness,
    AUTHOR = {Lau, Eike},
     TITLE = {Smoothness of the truncated display functor},
   JOURNAL = {J. Amer. Math. Soc.},
  FJOURNAL = {Journal of the American Mathematical Society},
    VOLUME = {26},
      YEAR = {2013},
    NUMBER = {1},
     PAGES = {129--165},
      ISSN = {0894-0347,1088-6834},
   MRCLASS = {14F30 (14L05)},
  MRNUMBER = {2983008},
MRREVIEWER = {Eva\ Viehmann},
       DOI = {10.1090/S0894-0347-2012-00744-9},
       URL = {https://doi.org/10.1090/S0894-0347-2012-00744-9},
}

@article {ScholzeWeinstein,
    AUTHOR = {Scholze, Peter and Weinstein, Jared},
     TITLE = {Moduli of {$p$}-divisible groups},
   JOURNAL = {Camb. J. Math.},
  FJOURNAL = {Cambridge Journal of Mathematics},
    VOLUME = {1},
      YEAR = {2013},
    NUMBER = {2},
     PAGES = {145--237},
      ISSN = {2168-0930,2168-0949},
   MRCLASS = {14L05 (11G25 14C30 14D20)},
  MRNUMBER = {3272049},
MRREVIEWER = {Rui\ Miguel\ Saramago},
       DOI = {10.4310/CJM.2013.v1.n2.a1},
       URL = {https://doi.org/10.4310/CJM.2013.v1.n2.a1},
}

@article {SchneiderStuhler,
    AUTHOR = {Schneider, P. and Stuhler, U.},
     TITLE = {The cohomology of {$p$}-adic symmetric spaces},
   JOURNAL = {Invent. Math.},
  FJOURNAL = {Inventiones Mathematicae},
    VOLUME = {105},
      YEAR = {1991},
    NUMBER = {1},
     PAGES = {47--122},
      ISSN = {0020-9910,1432-1297},
   MRCLASS = {11F85 (14G20)},
  MRNUMBER = {1109620},
MRREVIEWER = {Gerd\ Faltings},
       DOI = {10.1007/BF01232257},
       URL = {https://doi.org/10.1007/BF01232257},
}

@book {BerkeleyLectures,
    AUTHOR = {Scholze, Peter and Weinstein, Jared},
     TITLE = {Berkeley lectures on {$p$}-adic geometry},
    SERIES = {Annals of Mathematics Studies},
    VOLUME = {207},
 PUBLISHER = {Princeton University Press, Princeton, NJ},
      YEAR = {2020},
     PAGES = {x+250},
      ISBN = {978-0-691-20209-9; 978-0-691-20208-2; 978-0-691-20215-0},
   MRCLASS = {14G45 (14A15 14F30 14G22 14G35 14M15)},
  MRNUMBER = {4446467},
MRREVIEWER = {Lance\ Edward\ Miller},
}

@article {MingjiaZhang,
    AUTHOR = {Zhang, Mingjia},
     TITLE = {A {PEL}-type {I}gusa stack and the {$p$}-adic geometry of
              {S}himura varieties},
   JOURNAL = {Camb. J. Math.},
  FJOURNAL = {Cambridge Journal of Mathematics},
    VOLUME = {14},
      YEAR = {2026},
    NUMBER = {2},
     PAGES = {377--486},
      ISSN = {2168-0930,2168-0949},
   MRCLASS = {14G35 (11G18 14G45)},
  MRNUMBER = {5073031},
       DOI = {10.4310/cjm.260514231850},
       URL = {https://doi.org/10.4310/cjm.260514231850},
}

@ARTICLE{ScholzeDiamond,
       author = {{Scholze}, Peter},
        title = "{Etale cohomology of diamonds}",
      journal = {arXiv e-prints},
         year = 2017,
        month = sep,
          eid = {arXiv:1709.07343},
        pages = {arXiv:1709.07343},
          doi = {10.48550/arXiv.1709.07343},
archivePrefix = {arXiv},
       eprint = {1709.07343},
 primaryClass = {math.AG},
       adsurl = {https://ui.adsabs.harvard.edu/abs/2017arXiv170907343S}
}

@book{MessingCrystals,
	address = {Berlin, Heidelberg},
	series = {Lecture {Notes} in {Mathematics}},
	title = {The {Crystals} {Associated} to {Barsotti}-{Tate} {Groups}: with {Applications} to {Abelian} {Schemes}},
	volume = {264},
	copyright = {http://www.springer.com/tdm},
	isbn = {978-3-540-05840-3 978-3-540-37475-6},
	shorttitle = {The {Crystals} {Associated} to {Barsotti}-{Tate} {Groups}},
	url = {http://link.springer.com/10.1007/BFb0058301},
	urldate = {2024-12-09},
	publisher = {Springer Berlin Heidelberg},
	author = {Messing, William},
	year = {1972},
	doi = {10.1007/BFb0058301},
}

@misc{GleasonIvanovZillinger,
      title={Meromorphic vector bundles on the {F}argues--{F}ontaine curve}, 
      author={Ian Gleason and Alexander B. Ivanov and Felix Zillinger},
      year={2025},
      eprint={2307.00887},
      archivePrefix={arXiv},
      primaryClass={math.AG},
      url={https://arxiv.org/abs/2307.00887}, 
}

@misc{stacks-project,
  author       = {The {Stacks project authors}},
  title        = {The Stacks project},
  howpublished = {\url{https://stacks.math.columbia.edu}},
  year         = {2025},
}

@ARTICLE{AM-Descent-on-vStacks,
       author = {{Ansch{\"u}tz}, Johannes and {Mann}, Lucas},
        title = "{Descent for solid quasi-coherent sheaves on perfectoid spaces}",
      journal = {arXiv e-prints},
         year = 2024,
        month = mar,
          eid = {arXiv:2403.01951},
        pages = {arXiv:2403.01951},
          doi = {10.48550/arXiv.2403.01951},
archivePrefix = {arXiv},
       eprint = {2403.01951},
 primaryClass = {math.AG},
       adsurl = {https://ui.adsabs.harvard.edu/abs/2024arXiv240301951A}
}

@misc{ALBM-6ff-on-FF-curve,
      title={A 6-functor formalism for solid quasi-coherent sheaves on the {F}argues-{F}ontaine curve}, 
      author={Johannes Anschütz and Arthur-César Le Bras and Lucas Mann},
      year={2024},
      eprint={2412.20968},
      archivePrefix={arXiv},
      primaryClass={math.AG},
      url={https://arxiv.org/abs/2412.20968}, 
}

@article {Ivanov-arc-descent,
    AUTHOR = {Ivanov, Alexander B.},
     TITLE = {{Arc-descent} for the perfect loop functor and {$p$}-adic
              {D}eligne-{L}usztig spaces},
   JOURNAL = {{J. Reine Angew. Math.}},
  FJOURNAL = {Journal f\"ur die Reine und Angewandte Mathematik. [Crelle's
              Journal]},
    VOLUME = {794},
      YEAR = {2023},
     PAGES = {1--54},
      ISSN = {0075-4102,1435-5345},
   MRCLASS = {14G45 (14D24 20G05 22E50)},
  MRNUMBER = {4529410},
MRREVIEWER = {Alessandra\ Bertapelle},
       DOI = {10.1515/crelle-2022-0060},
       URL = {https://doi.org/10.1515/crelle-2022-0060},
}

@book {RapoportZink,
    AUTHOR = {Rapoport, M. and Zink, Th.},
     TITLE = {Period spaces for {$p$}-divisible groups},
    SERIES = {Annals of Mathematics Studies},
    VOLUME = {141},
 PUBLISHER = {Princeton University Press, Princeton, NJ},
      YEAR = {1996},
     PAGES = {xxii+324},
      ISBN = {0-691-02782-X; 0-691-02781-1},
   MRCLASS = {14G20 (11G18 14F30 14L05 14M15 20G05 20G25)},
  MRNUMBER = {1393439},
MRREVIEWER = {Robert\ E.\ Kottwitz},
       DOI = {10.1515/9781400882601},
       URL = {https://doi.org/10.1515/9781400882601},
}

@incollection {TatePDiv,
    AUTHOR = {Tate, J. T.},
     TITLE = {{$p$}-divisible groups},
 BOOKTITLE = {Proc. {C}onf. {L}ocal {F}ields ({D}riebergen, 1966)},
     PAGES = {158--183},
 PUBLISHER = {Springer, Berlin-New York},
      YEAR = {1967},
   MRCLASS = {14.50},
  MRNUMBER = {231827},
MRREVIEWER = {M.\ Miyanishi},
}

@ARTICLE{DavisAction,
       author = {{Davis}, Andr{\'e}},
        title = "{A Computation of the Action of the Morava Stabilizer Group on the Lubin-Tate Deformation Ring}",
      journal = {arXiv e-prints},
         year = 2022,
        month = jan,
          eid = {arXiv:2201.08913},
        pages = {arXiv:2201.08913},
          doi = {10.48550/arXiv.2201.08913},
archivePrefix = {arXiv},
       eprint = {2201.08913},
 primaryClass = {math.AT},
       adsurl = {https://ui.adsabs.harvard.edu/abs/2022arXiv220108913D}
}

@article {Lazard,
    AUTHOR = {Lazard, Michel},
     TITLE = {Groupes analytiques {$p$}-adiques},
   JOURNAL = {Inst. Hautes \'Etudes Sci. Publ. Math.},
  FJOURNAL = {Institut des Hautes \'Etudes Scientifiques. Publications
              Math\'ematiques},
    NUMBER = {26},
      YEAR = {1965},
     PAGES = {389--603},
      ISSN = {0073-8301,1618-1913},
   MRCLASS = {14.50},
  MRNUMBER = {209286},
       URL = {http://www.numdam.org/item?id=PMIHES_1965__26__389_0},
}

@book{DDMS,
    AUTHOR = {Dixon, J. D. and du Sautoy, M. P. F. and Mann, A. and Segal,
              D.},
     TITLE = {Analytic pro-{$p$}-groups},
    SERIES = {London Mathematical Society Lecture Note Series},
    VOLUME = {157},
 PUBLISHER = {Cambridge University Press, Cambridge},
      YEAR = {1991},
     PAGES = {x+251},
      ISBN = {0-521-39580-1},
   MRCLASS = {20E18 (20G30)},
  MRNUMBER = {1152800},
MRREVIEWER = {Andy\ R.\ Magid},
}

@ARTICLE{Salch,
       author = {{Salch}, Andrew},
        title = "{The action of the Morava stabilizer group on the coefficients of Morava E-theory at height 2}",
      journal = {arXiv e-prints},
         year = 2025,
        month = mar,
          eid = {arXiv:2503.04686},
        pages = {arXiv:2503.04686},
          doi = {10.48550/arXiv.2503.04686},
archivePrefix = {arXiv},
       eprint = {2503.04686},
 primaryClass = {math.AT},
       adsurl = {https://ui.adsabs.harvard.edu/abs/2025arXiv250304686S}
}

@article {GrossCanonical,
    AUTHOR = {Gross, Benedict H.},
     TITLE = {On canonical and quasicanonical liftings},
   JOURNAL = {Invent. Math.},
  FJOURNAL = {Inventiones Mathematicae},
    VOLUME = {84},
      YEAR = {1986},
    NUMBER = {2},
     PAGES = {321--326},
      ISSN = {0020-9910,1432-1297},
   MRCLASS = {14L05 (11S31 14H25)},
  MRNUMBER = {833193},
MRREVIEWER = {K.\ Shiratani},
       DOI = {10.1007/BF01388810},
       URL = {https://doi.org/10.1007/BF01388810},
}

@article {ShimomuraYabe1995,
    AUTHOR = {Shimomura, Katsumi and Yabe, Atsuko},
     TITLE = {The homotopy groups {$\pi_*(L_2S^0)$}},
   JOURNAL = {Topology},
  FJOURNAL = {Topology. An International Journal of Mathematics},
    VOLUME = {34},
      YEAR = {1995},
    NUMBER = {2},
     PAGES = {261--289},
      ISSN = {0040-9383},
   MRCLASS = {55Q45},
  MRNUMBER = {1318877},
MRREVIEWER = {Jin\ Kun\ Lin},
       DOI = {10.1016/0040-9383(94)00032-G},
       URL = {https://doi.org/10.1016/0040-9383(94)00032-G},
}

@article {LubinTate,
    AUTHOR = {Lubin, Jonathan and Tate, John},
     TITLE = {Formal moduli for one-parameter formal {L}ie groups},
   JOURNAL = {Bull. Soc. Math. France},
  FJOURNAL = {Bulletin de la Soci\'et\'e{} Math\'ematique de France},
    VOLUME = {94},
      YEAR = {1966},
     PAGES = {49--59},
      ISSN = {0037-9484},
   MRCLASS = {14.50},
  MRNUMBER = {238854},
MRREVIEWER = {J.\ Dieudonn\'e},
       URL = {http://www.numdam.org/item?id=BSMF_1966__94__49_0},
}

@article {Behrens2012,
    AUTHOR = {Behrens, Mark},
     TITLE = {The homotopy groups of {$S_{E(2)}$} at {$p\geq 5$} revisited},
   JOURNAL = {Adv. Math.},
  FJOURNAL = {Advances in Mathematics},
    VOLUME = {230},
      YEAR = {2012},
    NUMBER = {2},
     PAGES = {458--492},
      ISSN = {0001-8708,1090-2082},
   MRCLASS = {55Q45 (55T15)},
  MRNUMBER = {2914955},
MRREVIEWER = {Paul\ G.\ Goerss},
       DOI = {10.1016/j.aim.2012.02.023},
       URL = {https://doi.org/10.1016/j.aim.2012.02.023},
}

@article {MR2794609,
    AUTHOR = {Torii, Takeshi},
     TITLE = {{$K(n)$}-localization of the {$K(n+1)$}-local
              {$E_{n+1}$}-{A}dams spectral sequences},
   JOURNAL = {Pacific J. Math.},
  FJOURNAL = {Pacific Journal of Mathematics},
    VOLUME = {250},
      YEAR = {2011},
    NUMBER = {2},
     PAGES = {439--471},
      ISSN = {0030-8730,1945-5844},
   MRCLASS = {55P42 (55Q51 55T15 55T25)},
  MRNUMBER = {2794609},
MRREVIEWER = {Vigleik\ Angeltveit},
       DOI = {10.2140/pjm.2011.250.439},
       URL = {https://doi.org/10.2140/pjm.2011.250.439},
}

@article {Torii2010,
    AUTHOR = {Torii, Takeshi},
     TITLE = {On {$E_\infty$}-structure of the generalized {C}hern
              character},
   JOURNAL = {Bull. Lond. Math. Soc.},
  FJOURNAL = {Bulletin of the London Mathematical Society},
    VOLUME = {42},
      YEAR = {2010},
    NUMBER = {4},
     PAGES = {680--690},
      ISSN = {0024-6093},
   MRCLASS = {55P43 (55N22)},
  MRNUMBER = {2669689},
MRREVIEWER = {Andrey Yu. Lazarev},
       DOI = {10.1112/blms/bdq026},
       URL = {https://doi.org/10.1112/blms/bdq026},
}

@article {ToriiComparison,
    AUTHOR = {Torii, Takeshi},
     TITLE = {Comparison of {M}orava {$E$}-theories},
   JOURNAL = {Math. Z.},
  FJOURNAL = {Mathematische Zeitschrift},
    VOLUME = {266},
      YEAR = {2010},
    NUMBER = {4},
     PAGES = {933--951},
      ISSN = {0025-5874,1432-1823},
   MRCLASS = {55N22 (19L10 55N20 55S05)},
  MRNUMBER = {2729298},
MRREVIEWER = {John\ Hunton},
       DOI = {10.1007/s00209-009-0605-9},
       URL = {https://doi-org.ezproxy.bu.edu/10.1007/s00209-009-0605-9},
}

@ARTICLE{Clausen2025arxiv,
       author = {{Clausen}, Dustin},
        title = "{Duality and linearization for p-adic lie groups}",
      journal = {arXiv e-prints},
         year = 2025,
        month = jun,
          eid = {arXiv:2506.18174},
        pages = {arXiv:2506.18174},
          doi = {10.48550/arXiv.2506.18174},
archivePrefix = {arXiv},
       eprint = {2506.18174},
 primaryClass = {math.AT},
       adsurl = {https://ui.adsabs.harvard.edu/abs/2025arXiv250618174C}
}

@incollection {MR1320994,
    AUTHOR = {Hovey, Mark},
     TITLE = {Bousfield localization functors and {H}opkins' chromatic
              splitting conjecture},
 BOOKTITLE = {The \v Cech centennial ({B}oston, {MA}, 1993)},
    SERIES = {Contemp. Math.},
    VOLUME = {181},
     PAGES = {225--250},
 PUBLISHER = {Amer. Math. Soc., Providence, RI},
      YEAR = {1995},
      ISBN = {0-8218-0296-8},
   MRCLASS = {55P42 (55N20 55N22 55P60)},
  MRNUMBER = {1320994},
       DOI = {10.1090/conm/181/02036},
       URL = {https://doi.org/10.1090/conm/181/02036},
}

@book {RavenelNilpotence,
    AUTHOR = {Ravenel, Douglas},
    TITLE = {Nilpotence and periodicity in stable homotopy theory},
    SERIES = {Ann. of Math. Stud.},
    NUMBER = {128},
    PUBLISHER = {Princeton University Press, Princeton, NJ},
    YEAR = {1992},
    PAGES = {xiv+209 pp.},
    ISBN = {0-691-02572-X},
}

@unpublished{HeyerMann6FF,
    title = {Six Functor Formalisms and Representations},
    author = {Claudius Heyer and Lucas Mann},
    year = {2024},
}

@unpublished {MannThesis,
    title = {A $p$-Adic 6-Functor Formalism in Rigid-Analytic Geometry},
    author = {Lucas Mann},
    year = 2022,
    note = {\url{https://arxiv.org/abs/2206.02022}}
}

@article {HansenMannJL,
	title = {$p$-adic sheaves on classifying stacks, and the $p$-adic {J}acquet-{L}anglands correspondence},
	author = {David Hansen and Lucas Mann},
	year = 2022,
	note = {Preprint. Available at \url{https://arxiv.org/pdf/2207.04073.pdf}},
}

@ARTICLE{ZavyalovAbstractPoincare,
       author = {{Zavyalov}, Bogdan},
        title = "{Poincar{\'e} Duality in abstract 6-functor formalisms}",
      journal = {arXiv e-prints},
         year = 2023,
        month = jan,
          eid = {arXiv:2301.03821},
        pages = {arXiv:2301.03821},
          doi = {10.48550/arXiv.2301.03821},
archivePrefix = {arXiv},
       eprint = {2301.03821},
 primaryClass = {math.AG},
       adsurl = {https://ui.adsabs.harvard.edu/abs/2023arXiv230103821Z}
}

@misc{ClausenScholzeCondensed,
      title={Lectures on Condensed Mathematics}, 
      author={Peter Scholze},
      year={2026},
      eprint={2605.03658},
      archivePrefix={arXiv},
      primaryClass={math.NT},
      url={https://arxiv.org/abs/2605.03658}, 
}

@article {FarguesFontaineCourbes,
    AUTHOR = {Fargues, Laurent and Fontaine, Jean-Marc},
     TITLE = {Courbes et fibr\'es vectoriels en th\'eorie de {H}odge
              {$p$}-adique},
      NOTE = {With a preface by Pierre Colmez},
   JOURNAL = {Ast\'erisque},
  FJOURNAL = {Ast\'erisque},
    NUMBER = {406},
      YEAR = {2018},
     PAGES = {xiii+382},
      ISSN = {0303-1179,2492-5926},
      ISBN = {978-2-85629-896-1},
   MRCLASS = {14F30 (11F80 11G25)},
  MRNUMBER = {3917141},
MRREVIEWER = {Giovanni\ Rosso},
}

@article{ECII,
    author = {Lurie, Jacob},
    title = {Elliptic {C}ohomology {II}: {O}rientations},
    year = {2018},
    note = {\newline {A}vailable at \href{https://www.math.ias.edu/~lurie/papers/Elliptic-II.pdf}{https://www.math.ias.edu/~lurie/papers/Elliptic-II.pdf}},
}

@article{SAG,
    author = {Lurie, Jacob},
    title = {Spectral Algebraic Geometry},
    year = {2018},
    note = {\newline {A}vailable at \href{https://www.math.ias.edu/~lurie/papers/SAG-rootfile.pdf}{https://www.math.ias.edu/~lurie/papers/SAG-rootfile.pdf}},
}

@article{HA,
    author = {Lurie, Jacob},
    title = {Higher Algebra},
    year = {2017},
    note = {\newline {A}vailable at \href{https://www.math.ias.edu/~lurie/papers/HA.pdf}{https://www.math.ias.edu/~lurie/papers/HA.pdf}},
}

@article {DrinfeldCoverings,
    AUTHOR = {Drinfeld, V. G.},
     TITLE = {Coverings of {$p$}-adic symmetric domains},
   JOURNAL = {Funkcional. Anal. i Prilo\v zen.},
  FJOURNAL = {Akademija Nauk SSSR. Funkcional\cprime nyi Analiz i ego
              Prilo\v zenija},
    VOLUME = {10},
      YEAR = {1976},
    NUMBER = {2},
     PAGES = {29--40},
      ISSN = {0374-1990},
   MRCLASS = {14L05 (10D05 12B25 14G20)},
  MRNUMBER = {422290},
MRREVIEWER = {Michiel\ Hazewinkel},
}

@article{hansen2017degenerating,
  title={Degenerating vector bundles in p-adic Hodge theory},
  author={Hansen, David},
  journal={J. Inst. Math. Jussieu, to appear},
  year={2017}
}

@article{Wintenberger1979,
  author = {Wintenberger, Jean-Pierre},
  title = {{Extensions de Lie et groupes d'automorphismes de corps locaux de caract{\'{e}}ristique p}},
  journal = {Comptes rendus hebdomadaires des s{\'{e}}ances de l'Acad{\'{e}}mie des sciences. S{\'{e}}ries A et B, Sciences math{\'{e}}matiques et Sciences physiques},
  volume = {288},
  series = {A},
  pages = {367--370},
  year = {1979},
  mrnumber = {526137},
  zbl = {0475.12020}
}

@incollection{FarguesFontaine2014,
  author = {Fargues, Laurent and Fontaine, Jean-Marc},
  title = {Vector bundles on curves and p-adic Hodge theory},
  booktitle = {Automorphic Forms and Galois Representations},
  volume = {2},
  pages = {17--104},
  series = {London Mathematical Society Lecture Note Series},
  publisher = {Cambridge University Press},
  address = {Cambridge},
  year = {2014},
  editor = {Diamond, Fred and Kassaei, Payman L. and Kim, Minhyong},
}

@article {MR2004428,
    AUTHOR = {Torii, Takeshi},
     TITLE = {On degeneration of one-dimensional formal group laws and
              applications to stable homotopy theory},
   JOURNAL = {Amer. J. Math.},
  FJOURNAL = {American Journal of Mathematics},
    VOLUME = {125},
      YEAR = {2003},
    NUMBER = {5},
     PAGES = {1037--1077},
      ISSN = {0002-9327,1080-6377},
   MRCLASS = {55N22 (55P42 55P60)},
  MRNUMBER = {2004428},
MRREVIEWER = {Mark\ Hovey},
       URL =
              {http://muse.jhu.edu.ezproxy.bu.edu/journals/american_journal_of_mathematics/v125/125.5torii.pdf},
}

@ARTICLE{Heyer-GeometricLemma,
       author = {{Heyer}, Claudius},
        title = "{The Geometrical Lemma for Smooth Representations in Natural Characteristic}",
      journal = {arXiv e-prints},
         year = 2023,
        month = mar,
          eid = {arXiv:2303.14721},
        pages = {arXiv:2303.14721},
          doi = {10.48550/arXiv.2303.14721},
archivePrefix = {arXiv},
       eprint = {2303.14721},
 primaryClass = {math.RT},
       adsurl = {https://ui.adsabs.harvard.edu/abs/2023arXiv230314721H}
}

@article{Heyer2023Left,
  author    = {Heyer, Claudius},
  title     = {The Left Adjoint of Derived Parabolic Induction},
  journal   = {Mathematische Zeitschrift},
  volume    = {305},
  number    = {3},
  pages     = {60},
  year      = {2023},
  doi       = {10.1007/s00209-023-03385-5},
  eprint    = {2204.11581},
  archivePrefix = {arXiv},
  primaryClass = {math.RT},
  url       = {https://arxiv.org/abs/2204.11581}
}

@article {Mathew2016,
    AUTHOR = {Mathew, Akhil},
     TITLE = {The {G}alois group of a stable homotopy theory},
   JOURNAL = {Adv. Math.},
  FJOURNAL = {Advances in Mathematics},
    VOLUME = {291},
      YEAR = {2016},
     PAGES = {403--541},
      ISSN = {0001-8708,1090-2082},
   MRCLASS = {14F35 (11F11 11S20 14F20 18D10 18G55 55P43 55U35)},
  MRNUMBER = {3459022},
MRREVIEWER = {Rui\ Miguel\ Saramago},
       DOI = {10.1016/j.aim.2015.12.017},
       URL = {https://doi.org/10.1016/j.aim.2015.12.017},
}

@article {MNN2017,
    AUTHOR = {Mathew, Akhil and Naumann, Niko and Noel, Justin},
     TITLE = {Nilpotence and descent in equivariant stable homotopy theory},
   JOURNAL = {Adv. Math.},
  FJOURNAL = {Advances in Mathematics},
    VOLUME = {305},
      YEAR = {2017},
     PAGES = {994--1084},
      ISSN = {0001-8708,1090-2082},
   MRCLASS = {55P91 (55P42)},
  MRNUMBER = {3570153},
MRREVIEWER = {Gregory\ Z.\ Arone},
       DOI = {10.1016/j.aim.2016.09.027},
       URL = {https://doi.org/10.1016/j.aim.2016.09.027},
}

@incollection {AkhilNilpotence,
    AUTHOR = {Mathew, Akhil},
     TITLE = {Examples of descent up to nilpotence},
 BOOKTITLE = {Geometric and topological aspects of the representation theory
              of finite groups},
    SERIES = {Springer Proc. Math. Stat.},
    VOLUME = {242},
     PAGES = {269--311},
 PUBLISHER = {Springer, Cham},
      YEAR = {2018},
      ISBN = {978-3-319-94033-5; 978-3-319-94032-8},
   MRCLASS = {18G80},
  MRNUMBER = {3901164},
MRREVIEWER = {Paul\ Balmer},
       DOI = {10.1007/978-3-319-94033-5\_11},
       URL = {https://doi.org/10.1007/978-3-319-94033-5_11},
}

@article {DevinatzHopkins2004,
    AUTHOR = {Devinatz, Ethan S. and Hopkins, Michael J.},
     TITLE = {Homotopy fixed point spectra for closed subgroups of the
              {M}orava stabilizer groups},
   JOURNAL = {Topology},
  FJOURNAL = {Topology. An International Journal of Mathematics},
    VOLUME = {43},
      YEAR = {2004},
    NUMBER = {1},
     PAGES = {1--47},
      ISSN = {0040-9383},
   MRCLASS = {55P43 (55N22 55T15)},
  MRNUMBER = {2030586},
MRREVIEWER = {R.\ E.\ Stong},
       DOI = {10.1016/S0040-9383(03)00029-6},
       URL = {https://doi.org/10.1016/S0040-9383(03)00029-6},
}

@article {GHM2014,
    AUTHOR = {Goerss, Paul G. and Henn, Hans-Werner and Mahowald, Mark},
     TITLE = {The rational homotopy of the {$K(2)$}-local sphere and the
              chromatic splitting conjecture for the prime 3 and level 2},
   JOURNAL = {Doc. Math.},
  FJOURNAL = {Documenta Mathematica},
    VOLUME = {19},
      YEAR = {2014},
     PAGES = {1271--1290},
      ISSN = {1431-0635,1431-0643},
   MRCLASS = {55Q45 (55P15 55P42)},
  MRNUMBER = {3312144},
MRREVIEWER = {Sa\"{\i}d\ Zarati},
}

@article {HKM2013,
    AUTHOR = {Henn, Hans-Werner and Karamanov, Nasko and Mahowald, Mark},
     TITLE = {The homotopy of the {$K(2)$}-local {M}oore spectrum at the
              prime 3 revisited},
   JOURNAL = {Math. Z.},
  FJOURNAL = {Mathematische Zeitschrift},
    VOLUME = {275},
      YEAR = {2013},
    NUMBER = {3-4},
     PAGES = {953--1004},
      ISSN = {0025-5874,1432-1823},
   MRCLASS = {55Q52 (55P15)},
  MRNUMBER = {3127044},
MRREVIEWER = {Jie\ Wu},
       DOI = {10.1007/s00209-013-1167-4},
       URL = {https://doi.org/10.1007/s00209-013-1167-4},
}

@article {Shimomura2000,
    AUTHOR = {Shimomura, Katsumi},
     TITLE = {The homotopy groups of the {$L_2$}-localized mod {$3$} {M}oore
              spectrum},
   JOURNAL = {J. Math. Soc. Japan},
  FJOURNAL = {Journal of the Mathematical Society of Japan},
    VOLUME = {52},
      YEAR = {2000},
    NUMBER = {1},
     PAGES = {65--90},
      ISSN = {0025-5645,1881-1167},
   MRCLASS = {55Q45 (55P42 55Q52 55T15)},
  MRNUMBER = {1727130},
MRREVIEWER = {R.\ E.\ Stong},
       DOI = {10.2969/jmsj/05210065},
       URL = {https://doi.org/10.2969/jmsj/05210065},
}

@article {GHMR2005,
    AUTHOR = {Goerss, P. and Henn, H.-W. and Mahowald, M. and Rezk, C.},
     TITLE = {A resolution of the {$K(2)$}-local sphere at the prime 3},
   JOURNAL = {Ann. of Math. (2)},
  FJOURNAL = {Annals of Mathematics. Second Series},
    VOLUME = {162},
      YEAR = {2005},
    NUMBER = {2},
     PAGES = {777--822},
      ISSN = {0003-486X,1939-8980},
   MRCLASS = {55Q45 (55P60)},
  MRNUMBER = {2183282},
MRREVIEWER = {Mark\ J.\ Behrens},
       DOI = {10.4007/annals.2005.162.777},
       URL = {https://doi.org/10.4007/annals.2005.162.777},
}

@article {Behrens2006,
    AUTHOR = {Behrens, Mark},
     TITLE = {A modular description of the {$K(2)$}-local sphere at the
              prime 3},
   JOURNAL = {Topology},
  FJOURNAL = {Topology. An International Journal of Mathematics},
    VOLUME = {45},
      YEAR = {2006},
    NUMBER = {2},
     PAGES = {343--402},
      ISSN = {0040-9383},
   MRCLASS = {55Q40 (11F23 14H52 55N34 55Q51 55S05)},
  MRNUMBER = {2193339},
MRREVIEWER = {Mark\ Hovey},
       DOI = {10.1016/j.top.2005.08.005},
       URL = {https://doi.org/10.1016/j.top.2005.08.005},
}

@article {Westerland2017,
    AUTHOR = {Westerland, Craig},
     TITLE = {A higher chromatic analogue of the image of {$J$}},
   JOURNAL = {Geom. Topol.},
  FJOURNAL = {Geometry \& Topology},
    VOLUME = {21},
      YEAR = {2017},
    NUMBER = {2},
     PAGES = {1033--1093},
      ISSN = {1465-3060,1364-0380},
   MRCLASS = {55P42 (19L20 55N15 55P20 55Q51)},
  MRNUMBER = {3626597},
MRREVIEWER = {Markus\ Szymik},
       DOI = {10.2140/gt.2017.21.1033},
       URL = {https://doi.org/10.2140/gt.2017.21.1033},
}

@article {CSY2024,
    AUTHOR = {Carmeli, Shachar and Schlank, Tomer M. and Yanovski, Lior},
     TITLE = {Chromatic cyclotomic extensions},
   JOURNAL = {Geom. Topol.},
  FJOURNAL = {Geometry \& Topology},
    VOLUME = {28},
      YEAR = {2024},
    NUMBER = {8},
     PAGES = {3511--3564},
      ISSN = {1465-3060,1364-0380},
   MRCLASS = {55P42 (18N60)},
  MRNUMBER = {4843744},
MRREVIEWER = {Geoffrey\ M. L. Powell},
       DOI = {10.2140/gt.2024.28.3511},
       URL = {https://doi.org/10.2140/gt.2024.28.3511},
}

@article {BGHS2022,
    AUTHOR = {Beaudry, Agn\`es and Goerss, Paul G. and Hopkins, Michael J.
              and Stojanoska, Vesna},
     TITLE = {Dualizing spheres for compact {$p$}-adic analytic groups and
              duality in chromatic homotopy},
   JOURNAL = {Invent. Math.},
  FJOURNAL = {Inventiones Mathematicae},
    VOLUME = {229},
      YEAR = {2022},
    NUMBER = {3},
     PAGES = {1301--1434},
      ISSN = {0020-9910,1432-1297},
   MRCLASS = {55P92 (20E18 22E41 55R25 55R50 55U30)},
  MRNUMBER = {4462627},
MRREVIEWER = {Hans-Werner\ Henn},
       DOI = {10.1007/s00222-022-01120-1},
       URL = {https://doi.org/10.1007/s00222-022-01120-1},
}

@article {Henn1998,
    AUTHOR = {Henn, Hans-Werner},
     TITLE = {Centralizers of elementary abelian {$p$}-subgroups and
              mod-{$p$} cohomology of profinite groups},
   JOURNAL = {Duke Math. J.},
  FJOURNAL = {Duke Mathematical Journal},
    VOLUME = {91},
      YEAR = {1998},
    NUMBER = {3},
     PAGES = {561--585},
      ISSN = {0012-7094,1547-7398},
   MRCLASS = {20J05 (20E18)},
  MRNUMBER = {1604171},
MRREVIEWER = {Norberto\ Gavioli},
       DOI = {10.1215/S0012-7094-98-09121-9},
       URL = {https://doi.org/10.1215/S0012-7094-98-09121-9},
}

@article {RognesGalois,
    AUTHOR = {Rognes, John},
     TITLE = {Galois extensions of structured ring spectra. {S}tably
              dualizable groups},
   JOURNAL = {Mem. Amer. Math. Soc.},
  FJOURNAL = {Memoirs of the American Mathematical Society},
    VOLUME = {192},
      YEAR = {2008},
    NUMBER = {898},
     PAGES = {viii+137},
      ISSN = {0065-9266,1947-6221},
   MRCLASS = {55P43 (55M05 55P35 57T05)},
  MRNUMBER = {2387923},
MRREVIEWER = {Alberto\ Cavicchioli},
       DOI = {10.1090/memo/0898},
       URL = {https://doi.org/10.1090/memo/0898},
}

@phdthesis{Lader2013,
  TITLE = {{Une r{\'e}solution projective pour le second groupe de Morava pour p $\ge$ 5 et applications}},
  AUTHOR = {Lader, Olivier},
  URL = {https://theses.hal.science/tel-00875761},
  NUMBER = {2013STRAD028},
  SCHOOL = {{Universit{\'e} de Strasbourg}},
  YEAR = {2013},
  MONTH = Oct,
  TYPE = {Theses},
  HAL_ID = {tel-00875761},
  HAL_VERSION = {v1},
}

@book {Araki,
    AUTHOR = {Araki, Sh\^or\^o},
     TITLE = {Typical formal groups in complex cobordism and {$K$}-theory},
    SERIES = {Lectures in Mathematics, Department of Mathematics, Kyoto
              University},
    VOLUME = {No. 6},
 PUBLISHER = {Kinokuniya Book Store Co., Ltd., Tokyo},
      YEAR = {1973},
     PAGES = {v+101 pp. (loose errata)},
   MRCLASS = {57D90 (14L05 55B15)},
  MRNUMBER = {375354},
MRREVIEWER = {A.\ Dold},
}

@article {Ax1970,
    AUTHOR = {Ax, James},
     TITLE = {Zeros of polynomials over local fields---{T}he {G}alois
              action},
   JOURNAL = {J. Algebra},
  FJOURNAL = {Journal of Algebra},
    VOLUME = {15},
      YEAR = {1970},
     PAGES = {417--428},
      ISSN = {0021-8693},
   MRCLASS = {12.30 (10.14)},
  MRNUMBER = {263786},
MRREVIEWER = {M.\ Miyanishi},
       DOI = {10.1016/0021-8693(70)90069-4},
       URL = {https://doi.org/10.1016/0021-8693(70)90069-4},
}

@article {HoveyStrickland,
    AUTHOR = {Hovey, Mark and Strickland, Neil P.},
     TITLE = {Morava {$K$}-theories and localisation},
   JOURNAL = {Mem. Amer. Math. Soc.},
  FJOURNAL = {Memoirs of the American Mathematical Society},
    VOLUME = {139},
      YEAR = {1999},
    NUMBER = {666},
     PAGES = {viii+100},
      ISSN = {0065-9266,1947-6221},
   MRCLASS = {55P60 (55N22 55P42 55T15)},
  MRNUMBER = {1601906},
MRREVIEWER = {J.\ P. C. Greenlees},
       DOI = {10.1090/memo/0666},
       URL = {https://doi.org/10.1090/memo/0666},
}

@incollection {GrossRamification,
    AUTHOR = {Gross, Benedict H.},
     TITLE = {Ramification in {$p$}-adic {L}ie extensions},
 BOOKTITLE = {Journ\'ees de {G}\'eom\'etrie {A}lg\'ebrique de {R}ennes
              ({R}ennes, 1978), {V}ol. {III}},
    SERIES = {Ast\'erisque},
    VOLUME = {65},
     PAGES = {81--102},
 PUBLISHER = {Soc. Math. France, Paris},
      YEAR = {1979},
   MRCLASS = {12B10 (14L05)},
  MRNUMBER = {563473},
MRREVIEWER = {Helmut\ Koch},
}

@article{anDR,
    author = {Esteban Rodríguez Camargo, Juan},
    title = {The analytic {D}e {R}ham stack in rigid geometry},
    year = {2024},
    note = {\href{https://arxiv.org/abs/2401.07738}{arXiv:2401.07738}},
}

@article {solidRep,
    AUTHOR = {Rodrigues Jacinto, Joaqu\'{\i}n and Rodr\'{\i}guez Camargo,
              Juan Esteban},
     TITLE = {Solid locally analytic representations of {$p$}-adic {L}ie
              groups},
   JOURNAL = {Represent. Theory},
  FJOURNAL = {Representation Theory. An Electronic Journal of the American
              Mathematical Society},
    VOLUME = {26},
      YEAR = {2022},
     PAGES = {962--1024},
      ISSN = {1088-4165},
   MRCLASS = {11F85 (22E41 22E50)},
  MRNUMBER = {4475468},
MRREVIEWER = {Anton\ Deitmar},
       DOI = {10.1090/ert/615},
       URL = {https://doi.org/10.1090/ert/615},
}

@article{solidRep2,
    author = {Rodrigues Jacinto, Joaquín and Esteban Rodríguez Camargo, Juan},
    title = {Solid locally analytic representations},
    year = {2026},
    note = {\href{https://arxiv.org/abs/2305.03162}{arXiv:2305.03162}},
}

@incollection {Faltings2002,
    AUTHOR = {Faltings, Gerd},
     TITLE = {A relation between two moduli spaces studied by {V}. {G}.
              {D}rinfeld},
 BOOKTITLE = {Algebraic number theory and algebraic geometry},
    SERIES = {Contemp. Math.},
    VOLUME = {300},
     PAGES = {115--129},
 PUBLISHER = {Amer. Math. Soc., Providence, RI},
      YEAR = {2002},
      ISBN = {0-8218-3267-0},
   MRCLASS = {14L05 (14G22)},
  MRNUMBER = {1936369},
MRREVIEWER = {Noriko\ Yui},
       DOI = {10.1090/conm/300/05145},
       URL = {https://doi.org/10.1090/conm/300/05145},
}

@article {CaraianiScholze,
    AUTHOR = {Caraiani, Ana and Scholze, Peter},
     TITLE = {On the generic part of the cohomology of non-compact unitary
              {S}himura varieties},
   JOURNAL = {Ann. of Math. (2)},
  FJOURNAL = {Annals of Mathematics. Second Series},
    VOLUME = {199},
      YEAR = {2024},
    NUMBER = {2},
     PAGES = {483--590},
      ISSN = {0003-486X,1939-8980},
   MRCLASS = {11R39 (14G35 14G45)},
  MRNUMBER = {4713019},
MRREVIEWER = {Nguy\cftil en Qu\^oc Th\'ang},
       DOI = {10.4007/annals.2024.199.2.1},
       URL = {https://doi-org.ezproxy.bu.edu/10.4007/annals.2024.199.2.1},
}

@article {BirkbeckExtensions,
    AUTHOR = {Birkbeck, Christopher and Feng, Tony and Hansen, David and
              Hong, Serin and Li, Qirui and Wang, Anthony and Ye, Lynnelle},
     TITLE = {Extensions of vector bundles on the {F}argues-{F}ontaine
              curve},
   JOURNAL = {J. Inst. Math. Jussieu},
  FJOURNAL = {Journal of the Institute of Mathematics of Jussieu. JIMJ.
              Journal de l'Institut de Math\'ematiques de Jussieu},
    VOLUME = {21},
      YEAR = {2022},
    NUMBER = {2},
     PAGES = {487--532},
      ISSN = {1474-7480,1475-3030},
   MRCLASS = {14D24 (14G45)},
  MRNUMBER = {4386821},
MRREVIEWER = {Alexander\ Boris\ Ivanov},
       DOI = {10.1017/S1474748020000183},
       URL = {https://doi-org.ezproxy.bu.edu/10.1017/S1474748020000183},
}

@book{LocallyConvexSpacesBook,
  title = {Locally Convex Spaces over Non-{{Archimedean}} Valued Fields},
  author = {{Perez-Garcia}, C.},
  year = 2010,
  series = {Cambridge Studies in Advanced Mathematics},
  number = {119},
  publisher = {Cambridge University Press},
  address = {New York},
  collaborator = {Schikhof, Wilhelmus Hendricus},
  isbn = {978-0-521-19243-9 978-0-511-72571-5},
  langid = {english}
}

@misc{Bosco,
      title={On the $p$-adic pro-\'etale cohomology of Drinfeld symmetric spaces}, 
      author={Guido Bosco},
      year={2023},
      eprint={2110.10683},
      archivePrefix={arXiv},
      primaryClass={math.NT},
      url={https://arxiv.org/abs/2110.10683}, 
}

@book{TopologicalVectorSpaces-Schaefer-Wolff,
  title = {Topological {{Vector Spaces}}},
  author = {Schaefer, H. H. and Wolff, M. P.},
  year = 1999,
  series = {Graduate {{Texts}} in {{Mathematics}}},
  volume = {3},
  publisher = {Springer New York},
  address = {New York, NY},
  doi = {10.1007/978-1-4612-1468-7},
  urldate = {2026-05-29},
  copyright = {http://www.springer.com/tdm},
  isbn = {978-1-4612-7155-0 978-1-4612-1468-7}
}

@article {hyperdescent,
    AUTHOR = {Clausen, Dustin and Mathew, Akhil},
     TITLE = {Hyperdescent and \'etale {$K$}-theory},
   JOURNAL = {Invent. Math.},
  FJOURNAL = {Inventiones Mathematicae},
    VOLUME = {225},
      YEAR = {2021},
    NUMBER = {3},
     PAGES = {981--1076},
      ISSN = {0020-9910,1432-1297},
   MRCLASS = {18F25 (14F20 55N15)},
  MRNUMBER = {4296353},
MRREVIEWER = {Barry\ H.\ Dayton},
       DOI = {10.1007/s00222-021-01043-3},
       URL = {https://doi.org/10.1007/s00222-021-01043-3},
}

@article {Arpon,
    AUTHOR = {Raksit, Arpon},
     TITLE = {Hochschild homology and the derived de {R}ham complex
              revisited},
   JOURNAL = {Selecta Math. (N.S.)},
  FJOURNAL = {Selecta Mathematica. New Series},
    VOLUME = {32},
      YEAR = {2026},
    NUMBER = {3},
     PAGES = {Paper No. 61, 96},
      ISSN = {1022-1824,1420-9020},
   MRCLASS = {14A30 (13D03 14F40)},
  MRNUMBER = {5090707},
       DOI = {10.1007/s00029-026-01154-w},
       URL = {https://doi.org/10.1007/s00029-026-01154-w},
}

@article {arc,
    AUTHOR = {Bhatt, Bhargav and Mathew, Akhil},
     TITLE = {The arc-topology},
   JOURNAL = {Duke Math. J.},
  FJOURNAL = {Duke Mathematical Journal},
    VOLUME = {170},
      YEAR = {2021},
    NUMBER = {9},
     PAGES = {1899--1988},
      ISSN = {0012-7094,1547-7398},
   MRCLASS = {14F20 (14F06 14G22)},
  MRNUMBER = {4278670},
MRREVIEWER = {Lance\ Edward\ Miller},
       DOI = {10.1215/00127094-2020-0088},
       URL = {https://doi.org/10.1215/00127094-2020-0088},
}

@book {MacLaneMoerdijk,
    AUTHOR = {Mac Lane, Saunders and Moerdijk, Ieke},
     TITLE = {Sheaves in geometry and logic},
    SERIES = {Universitext},
      NOTE = {A first introduction to topos theory,
              Corrected reprint of the 1992 edition},
 PUBLISHER = {Springer-Verlag, New York},
      YEAR = {1994},
     PAGES = {xii+629},
      ISBN = {0-387-97710-4},
   MRCLASS = {03G30 (18B25 54B40)},
  MRNUMBER = {1300636},
MRREVIEWER = {M.\ Makkai},
}

@article{solidProf,
    author = {Ged Corob Cook and Max Gheorghiu and Sofía Marlasca Aparicio and Thomas A. Wasserman},
    title = {Solid duality for profinite groups},
    year = {2026},
    note = {\href{https://arxiv.org/abs/2607.16446}{arXiv:2607.16446}},
}

@misc{KedlayaLiuImperfect,
  author        = {Kedlaya, Kiran S. and Liu, Ruochuan},
  title         = {Relative $p$-adic {H}odge theory, {II}: {I}mperfect period rings},
  elocation-id  = {arXiv:1602.06899},
  eprint        = {1602.06899},
  archivePrefix = {arXiv},
  primaryClass  = {math.NT},
  year          = {2019},
  note          = {v3; major changes to the treatment of pseudocoherent modules}
}

@article{SerreCohDim,
  title = {Sur La Dimension Cohomologique Des Groupes Profinis},
  author = {Serre, Jean-Pierre},
  year = 1965,
  month = jun,
  journal = {Topology},
  volume = {3},
  number = {4},
  pages = {413--420},
  issn = {00409383},
  doi = {10.1016/0040-9383(65)90006-6},
  urldate = {2026-09-18}
}

@article{Heine,
    author = {Heine, Hadrian},
    title = {A duality between monads and monadic morphisms},
    year = {2026},
    note = {\href{https://arxiv.org/abs/1712.00555v6}{arXiv:1712.00555v6}},
}

@article{spectral-lie,
    author = {Blans, Max and Heuts, Gijs},
    title = {A characterization of the spectral {L}ie operad},
    year = {2026},
    note = {\href{https://arxiv.org/abs/2606.17313}{arXiv:2606.17313}},
}

@article{brink,
    author = {Brink, Emma},
    title = {Condensed {G}roup {C}ohomology},
    year = {2025},
    note = {\href{https://arxiv.org/abs/2512.03648}{arXiv:2512.03648}},
}

@misc{hexagon:2610.00138v1,
  title={Weak chromatic splitting fails at n = p for every odd prime p},
  key={hexagon:2610.00138v1},
  howpublished={Hexagon preprint hexagon:2610.00138v1},
  year={2026},
  note={Contributors: Agnes Beaudry; Tomer Schlank},
  url={https://hexagonmath.org/2610.00138v1}
}

@misc{OAI:filtration,
  author = {{OpenAI}},
  title = {{Filtered chromatic splitting at generic primes}},
  howpublished = {OpenAI Math Release preprint
                  \href{https://github.com/openai/math/blob/main/preprints/Filtered-chromatic-splitting-at-generic-primes-September-25-2026/paper.pdf}{OAI:Filtered-chromatic-splitting-at-generic-primes-September-25-2026}},
  year = {2026}
}

\end{document}